\documentclass[10pt,reqno]{article}
\usepackage{amsfonts,amsrefs,latexsym,amsmath, amssymb, mathrsfs, verbatim}
\usepackage{url,color}
\usepackage{pifont}
\usepackage{upgreek}
\usepackage{subfigure}
\usepackage{fancyhdr}
\usepackage[linktocpage]{hyperref}
\usepackage{calligra}
\usepackage{marvosym}
\usepackage[percent]{overpic}
\usepackage{pict2e}
\usepackage{tensor}
\usepackage{amsthm}
\usepackage[hypcap=false]{caption}
\usepackage{xcolor}
\usepackage{wasysym}
\usepackage{accents}
\usepackage{enumitem}
\usepackage{tikz}

\usetikzlibrary{shapes.geometric, arrows}
\usetikzlibrary{mindmap}
\setlist[enumerate]{label*=\arabic*.}
\numberwithin{equation}{section}

\topmargin - .23 in

\newtheorem{theorem}{Theorem}[section]
\newtheorem{proposition}{Proposition}[section]
\newtheorem{lemma}[proposition]{Lemma}
\newtheorem{corollary}[proposition]{Corollary}
\newtheorem{notation}[proposition]{Notation}

\theoremstyle{definition}
\newtheorem{definition}{Definition}[section]
\newtheorem{remark}[proposition]{Remark}

\DeclareMathAlphabet{\mathcalligra}{T1}{calligra}{m}{n}
\DeclareFontShape{T1}{calligra}{m}{n}{<->s*[2.2]callig15}{}

\newcommand{\Ent}{s}

\newcommand{\gradEnt}{S}
\newcommand{\vgradientEnt}{\vec{\gradEnt}}
\newcommand{\pressure}{p}
\newcommand{\density}{\varrho}
\newcommand{\backgroundDensity}{\overline{\density}}
\newcommand{\logDensity}{\uprho}
\newcommand{\logDensitybar}{\underline{\uprho}}

\newcommand{\vort}{\upomega}
\newcommand{\vvort}{\vec{\vortrenormalized}}
\newcommand{\vortrenormalized}{\Omega}
\newcommand{\vVort}{\vec{\vortrenormalized}}
\newcommand{\angvort}{{\vortrenormalized \mkern-12mu / } \, }
\newcommand{\anggradEnt}{{\gradEnt \mkern-12mu / } \, }
\newcommand{\speed}{c}

\newcommand{\velbar}{\underline{v}}
\newcommand{\angvbar}{\underline{\angv}}

\newcommand{\entbar}{\underline{\Ent}}

\newcommand{\pressurebar}{\underline{\pressure}}

\newcommand{\backgroundDensitybar}{\underline{\backgroundDensity}}

\newcommand{\speedbar}{\underline{\speed}}

\newcommand{\inivel}{\mathring{v}}
\newcommand{\iniangv}{\mathring{\angv}}

\newcommand{\inient}{\mathring{\Ent}}
\newcommand{\inigradEnt}{\mathring{\gradEnt}}

\newcommand{\inipressure}{\mathring{\pressure}}

\newcommand{\inilogDensity}{\mathring{\logDensity}}

\newcommand{\inivortrenormalized}{\mathring{\vortrenormalized}}

\newcommand{\inispeed}{\mathring{\speed}}

\newcommand{\inivelbar}{\underline{\mathring{v}}}
\newcommand{\iniangvbar}{\underline{\iniangv}}

\newcommand{\inientbar}{\underline{\mathring{\Ent}}}

\newcommand{\inilogDensitybar}{\underline{\mathring{\logDensity}}}

\newcommand{\inispeedbar}{\underline{\mathring{\speed}}}

\newcommand{\DivGradEnt}{\mathcal{D}}
\newcommand{\vectvort}{\vec{\modivort}}
\newcommand{\modivort}{\mathcal{C}}

\newcommand{\iniDivGradEnt}{\mathring{\DivGradEnt}}

\newcommand{\inimodivort}{\mathring{\modivort}}

\newcommand{\variables}{\Psi}
\newcommand{\vvariables}{\vec{\variables}}

\newcommand{\inivvariables}{\mathring{\vvariables}}
\newcommand{\inivvariablesbar}{\underline{\mathring{\vvariables}}}

\newcommand{\norm}[1]{\left\lVert#1\right\rVert}
\newcommand{\onenorm}[4]{\norm{#1}_{L_{#2}^{#3}(#4)}}
\newcommand{\twonorms}[6]{\norm{#1}_{L_{#2}^{#3}L_{#4}^{#5}(#6)}}

\newcommand{\holdernorm}[5]{\norm{#1}_{C_{#2}^{#3,#4}(#5)}}
\newcommand{\semiholdernorm}[5]{\norm{#1}_{\dot{C}_{#2}^{#3,#4}(#5)}}
\newcommand{\holdertwonorms}[7]{\norm{#1}_{L_{#2}^{#3}C_{#4}^{#5,#6}(#7)}}

\newcommand{\sobolevnorm}[3]{\norm{#1}_{H^{#2}(#3)}}
\newcommand{\abs}[1]{\left\vert#1\right\vert}

\newcommand{\gtauabs}[1]{\left\vert#1\right\vert_{\tilde{g}}}
\newcommand{\inigfabs}[1]{\left\vert#1\right\vert_{\mathring{\gfour}}}

\newcommand{\sgabs}[1]{\left\vert#1\right\vert_{\gsphere}}
\newcommand{\seabs}[1]{\left\vert#1\right\vert_{\esphere}}

\newcommand{\dive}{\mbox{\upshape div}}
\newcommand{\curl}{\mbox{\upshape curl}}
\newcommand{\antisymmetic}{\upepsilon}

\newcommand{\cartiasphere}{\theta^{A}}
\newcommand{\sphere}{\theta}

\newcommand{\sangle}{\theta}

\newcommand{\ExistenceRegionthreearg}[3]{\mathcal{M}_{#1}^{#2,#3}}
\newcommand{\region}{\mathcal{M}}
\newcommand{\Rendallregionwithgeometry}{\hat{\mathcal{M}}}
\newcommand{\Rendallexistencetime}{\uptau_{\mbox{\upshape Rendall}}}
\newcommand{\timecutoffRendallregionwithgeometry}{\uptau_{\mbox{\upshape local}}}
\newcommand{\extendedregion}{\mathcal{O}}
\newcommand{\hypersur}{\mathcal{H}}
\newcommand{\hypersurbar}{\underline{\mathcal{H}}}

\newcommand{\inicartiasphere}{\mathring{\theta}^A}

\newcommand{\inihypersur}{\mathring{\hypersur}}
\newcommand{\inihypersurbar}{\mathring{\hypersurbar}}
\newcommand{\iniLunit}{\mathring{\Lunit}}
\newcommand{\iniuLunit}{\mathring{\uLunit}}
\newcommand{\uiniLunit}{\tilde{\Lunit}}
\newcommand{\uiniuLunit}{\tilde{\uLunit}}
\newcommand{\iniLgeo}{\mathring{\mathrm{F}}}
\newcommand{\iniLbargeo}{\underline{\mathring{\mathrm{F}}}}

\newcommand{\St}{\Sigma_t}
\newcommand{\Stau}{\Sigma_{\uptau}}

\newcommand{\suubarprime}{S_{u,\ubarprime}}
\newcommand{\suubar}{S_{u,\ubar}}

\newcommand{\suprimeubar}{S_{\uprime,\ubar}}
\newcommand{\Stwo}{\mathbb{S}^2}
\newcommand{\stu}{S_{t,u}}
\newcommand{\sTauu}{S_{\uptau,u}}

\newcommand{\iniepsilon}{\Lambda_0}
\newcommand{\RBAinfty}{\widetilde{\BAinfty}}
\newcommand{\BAinfty}{\Delta_0}
\newcommand{\iniBAinfty}{\mathring{\Delta}_0}
\newcommand{\iniLtwo}{\mathcal{O}_0}
\newcommand{\iniLtwoS}{\mathring{\Lambda}_{0}}
\newcommand{\LOTsmoothfunction}{\mathrm{f}_{(l.o.t.)}}
\newcommand{\gensmoothfunction}{\mathrm{f}}
\newcommand{\quadsmoothfunction}{\mathscr{Q}}
\newcommand{\linsmoothfunction}{\mathscr{L}}

\newcommand{\Ricfour}[2]{\mathbf{Ric}_{#1 #2}}
\newcommand{\Riemfour}[4]{\mathbf{Riem}_{#1 #2 #3 #4}}

\newcommand{\spherenormal}{\mathbf{N}}

\newcommand{\inilapse}{\mathring{\lapse}}
\newcommand{\boundinilapse}{K}
\newcommand{\invfolia}{\upmu}
\newcommand{\iniinvfolia}{\mathring{\invfolia}}
\newcommand{\uinvfolia}{\underline{\upmu}}
\newcommand{\iniuinvfolia}{\mathring{\uinvfolia}}
\newcommand{\nullapse}{\uplambda}
\newcommand{\lapse}{\upnu}

\newcommand{\ueta}{\underline{\upeta}}

\newcommand{\hchi}{\hat{\upchi}}
\newcommand{\uchi}{\underline{\upchi}}
\newcommand{\huchi}{\hat{\uchi}}

\newcommand{\ACC}{\mathcal{A}}
\newcommand{\iniACC}{\mathring{\ACC}}
\newcommand{\GaussCurvature}{\mathbf{K}}
\newcommand{\Jen}[1]{^{(#1)} \mkern-3mu \mathbf{J}}
\newcommand{\Jenarg}[2]{{^{(#1)} \mkern-3mu \mathbf{J}^{#2}}}

\newcommand{\deform}[1]{{^{(#1)} \mkern-1mu \pmb{\pi}}}

\newcommand{\gfour}{\mathbf{g}}
\newcommand{\inigfour}{\mathring{\gfour}}
\newcommand{\inigfourinverse}{\mathring{\gfour}^{-1}}

\newcommand{\gt}{g}
\newcommand{\dfour}[2]{G_{#1#2}}

\newcommand{\gtau}{\tilde{g}}
\newcommand{\gsphere}{g \mkern-8.5mu / }
\newcommand{\inigsphere}{\ring{g} \mkern-8.5mu / }
\newcommand{\esphere}{e \mkern-8.5mu / }

\newcommand{\Chfour}{\pmb{\Gamma}}

\newcommand{\volume}{\upsilon}

\newcommand{\gvol}{\varpi_{\gfour}}

\newcommand{\tvol}{\varpi_\gt}
\newcommand{\flattvol}{\varpi_{e}}
\newcommand{\tauvol}{\varpi_{\gtau}}

\newcommand{\spherevol}{\varpi_{\gsphere}}
\newcommand{\flatspherevol}{\varpi_{\esphere}}

\newcommand{\sphereproject}{{\Pi \mkern-12mu / } \, }
\newcommand{\Sigmatproject}{\underline{\Pi}}
\newcommand{\Sigmatauproject}{\tilde{\Pi}}

\newcommand{\mytr}{{\mbox{\upshape tr}}}

\newcommand{\gtr}{\mytr_{\gsphere}}
\newcommand{\inigtr}{\mytr_{\inigsphere}}

\newcommand{\Timelike}{\mathbf{T}}
\newcommand{\Trenomal}{\breve{\mathbf{T}}}
\newcommand{\materialderivative}{\mathbf{B}}
\newcommand{\Lunit}{L}
\newcommand{\Lrenormal}{\breve{\Lunit}}
\newcommand{\uLunit}{\underline{L}}
\newcommand{\uLrenormal}{\breve{\uLunit}}
\newcommand{\Lgeo}{\Lunit_{(Geo)}}
\newcommand{\uLgeo}{\uLunit_{(Geo)}}
\newcommand{\uomega}{\underline{\omega}}

\newcommand{\p}{\partial}
\newcommand{\pfour}{\pmb\partial}

\newcommand{\ubar}{\underline{u}}
\newcommand{\ubarprime}{\underline{u}^\prime}
\newcommand{\uprime}{u^\prime}

\newcommand{\ubarfinite}{\hat{\ubar}}
\newcommand{\ufinite}{\hat{u}}

\newcommand{\inimaterialderivative}{\mathbf{\mathring{B}}}

\newcommand{\inisphere}{\mathring{\sphere}}
\newcommand{\iniu}{\mathring{u}}
\newcommand{\iniubar}{\mathring{\ubar}}
\newcommand{\inip}{\mathring{\pmb{\p}}}

\newcommand{\embed}{\mathcal{E}}
\newcommand{\embedding}{\mathcal{I}}
\newcommand{\uembedding}{\mathcal{\underline{I}}}

\newcommand{\embedkappa}{\upkappa}
\newcommand{\embeddomain}{\mathcal{D}}

\newcommand{\embedgeo}{\Phi}

\newcommand{\Dfour}{\mathbf{D}}
\newcommand{\angv}{{v \mkern-10mu / }}
\newcommand{\angD}{ {\Dfour \mkern-14mu / \,}}
\newcommand{\angB}{ {\materialderivative \mkern-12mu / \,}}
\newcommand{\iniangB}{ {\mathring{\materialderivative} \mkern-12mu / \,}}
\newcommand{\angp}{ {\p \mkern-10mu / }}
\newcommand{\angV}{ {V \mkern-14mu / }}
\newcommand{\taup}{\tilde{\p}}
\newcommand{\iniangp}{ {\mathring{\p} \mkern-10mu / }}
\newcommand{\angnabla}{ {\nabla \mkern-14mu / \,} }

\newcommand{\angdfour}{{G \mkern-14mu / \,}}

\newcommand{\angdiv}{\mbox{\upshape{div} $\mkern-17mu /$\,}}
\newcommand{\angcurl}{\mbox{\upshape{curl} $\mkern-17mu /$\,}}
\newcommand{\anglap}{ {\Delta \mkern-12mu / \, } }

\newcommand{\rgeo}{\widetilde{r}}

\newcommand{\timefunction}{\uptau}

\newcommand{\Lie}{\mathcal{L}}

\newcommand{\angLie}{ { \mathcal{L} \mkern-10mu / } }

\newcommand{\boxg}{\square_{\gfour}}

\newcommand{\diff}{\mathrm d}

\newcommand{\average}[1]{\overline{#1}}

\newcommand{\Energy}{\mathbb{F}}
\newcommand{\KEnergy}{\mathbb{K}}

\newcommand{\uEnergy}{\underline{\mathbb{F}}}

\newcommand{\bEnergy}{\mathbb{E}}

\newcommand{\FR}{\mathfrak{R}}
\newcommand{\FRBar}{\underline{\FR}}
\newcommand{\FV}{\mathfrak{V}}
\newcommand{\FVBar}{\underline{\FV}}
\newcommand{\FS}{\mathfrak{S}}
\newcommand{\FSBar}{\underline{\FS}}
\newcommand{\FW}{\mathfrak{W}}
\newcommand{\FWBar}{\underline{\FW}}

\newcommand{\mypmH}{z}

\title{The Characteristic Initial Value Problem for the 3D Compressible Euler Equations}

\date{\today}      
\begin{document}	
	\author{
		Jared Speck\thanks{Vanderbilt University, Nashville, TN, USA. jared.speck@vanderbilt.edu.
		J.\, Speck gratefully acknowledges support from NSF grants DMS-2349575 and DMS-2054184.} and
		Sifan Yu\thanks{National University of Singapore, Singapore. sifanyu@nus.edu.sg}}
	\maketitle
	
	\begin{abstract}
		We prove local existence and uniqueness, propagation of Sobolev regularity, and continuation criteria
		results for the characteristic initial value problem for the $3D$ compressible Euler equations.
		This is the Cauchy problem with initial data given on a pair of transversally intersecting (acoustically) null
		hypersurfaces. The region of existence in our main theorem 
		is a future-neighborhood of the union of the data null hypersurfaces,
		where the size of the neighborhood depends on various norms of the data as well as the structure of the
		data null hypersurfaces.
		We allow for any equation of state with positive sound speed, and we do not make any symmetry,
		irrotationality, or isentropicity assumptions.
		Our analysis relies on the geometric formulation of the flow
		derived in \cite{ANewFormulationofthe3DCompressibleEulerEquationswithDynamicEntropy},
		which decomposes the equations into a ``wave-part'' 
		coupled to a ``transport-div-curl'' part.
		Compared to Rendall's celebrated solution of the  
		characteristic initial value problem for wave equations
		\cite{ReductionoftheCharacteristicInitialValueProblem},
		substantial new ideas are needed to handle the fact that the
		compressible Euler equations are more complicated than pure wave equations.
		\\		
		\\
		\textbf{Keywords:} characteristic initial value problem; double-null foliation; energy estimates; 
		eikonal function; local well-posedness; continuation criteria; vectorfield method\\
		\\
		\textbf{Mathematics Subject Classification:} 35Q31; Secondary: 35L67, 76L05. \\
		\\
	\end{abstract}
	\tableofcontents
	\newpage

	\section{Introduction}
	\label{S:INTRO}
\begin{figure} \centering    
		\subfigure[In Rendall's work \cite{ReductionoftheCharacteristicInitialValueProblem} on wave equations and Einstein's equations. ] {
			\label{fig3.1}     
			\includegraphics[width=7cm,height=7cm]{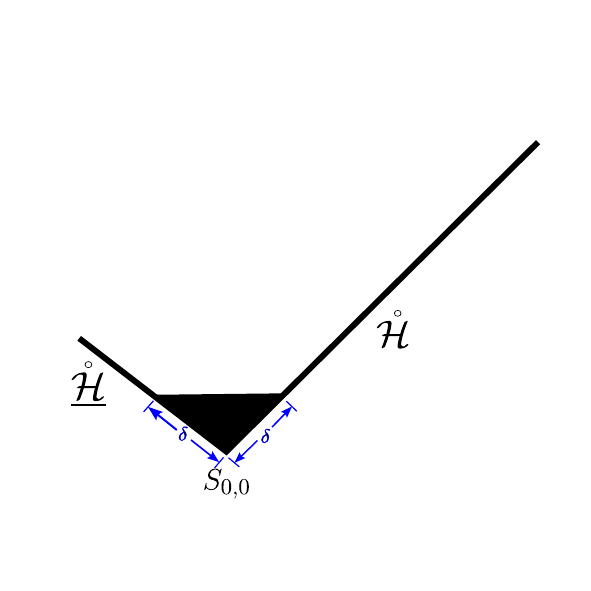}  
		}     
		\subfigure[In this article for compressible Euler equations.] {
			\label{fig3.2}     
			\includegraphics[width=7cm,height=7cm]{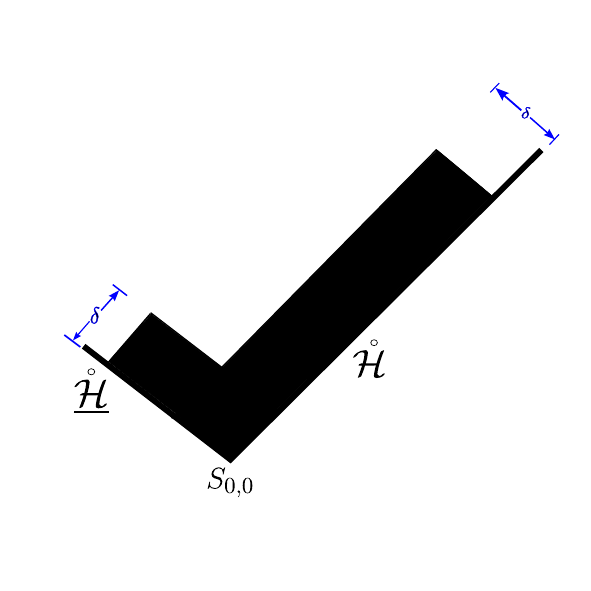}  
		}  
		\caption{The regions of local existence and uniqueness in double-null coordinates.}  
		\label{fig3}   
	\end{figure}

\begin{figure}\centering    
		\subfigure[In Rendall's work \cite{ReductionoftheCharacteristicInitialValueProblem}. ] {
			\label{fig1.1}     
			\includegraphics[width=6cm,height=8.5cm]{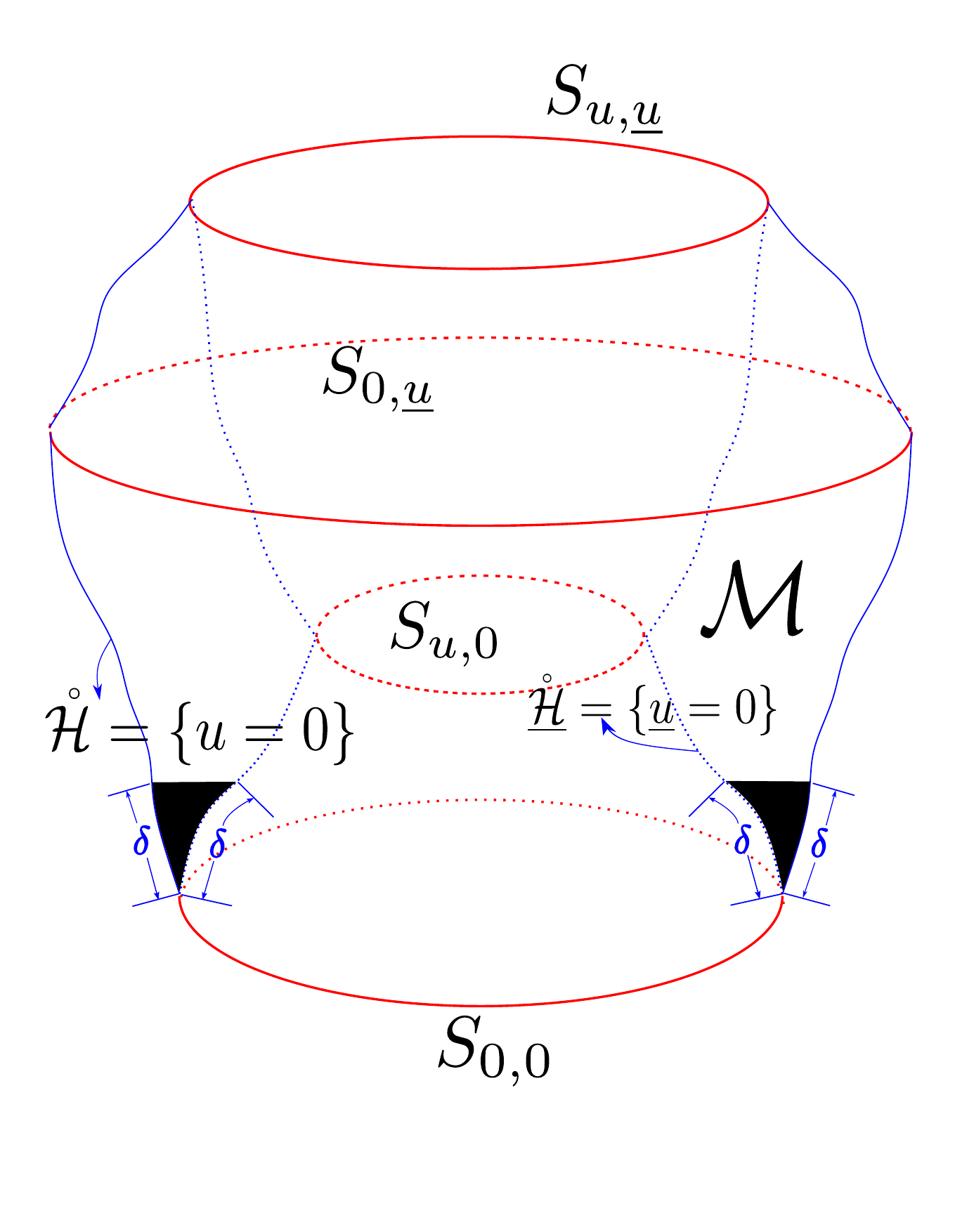}  
		}     
		\subfigure[In this article.] {
			\label{fig1.2}     
			\includegraphics[width=6cm,height=8.5cm]{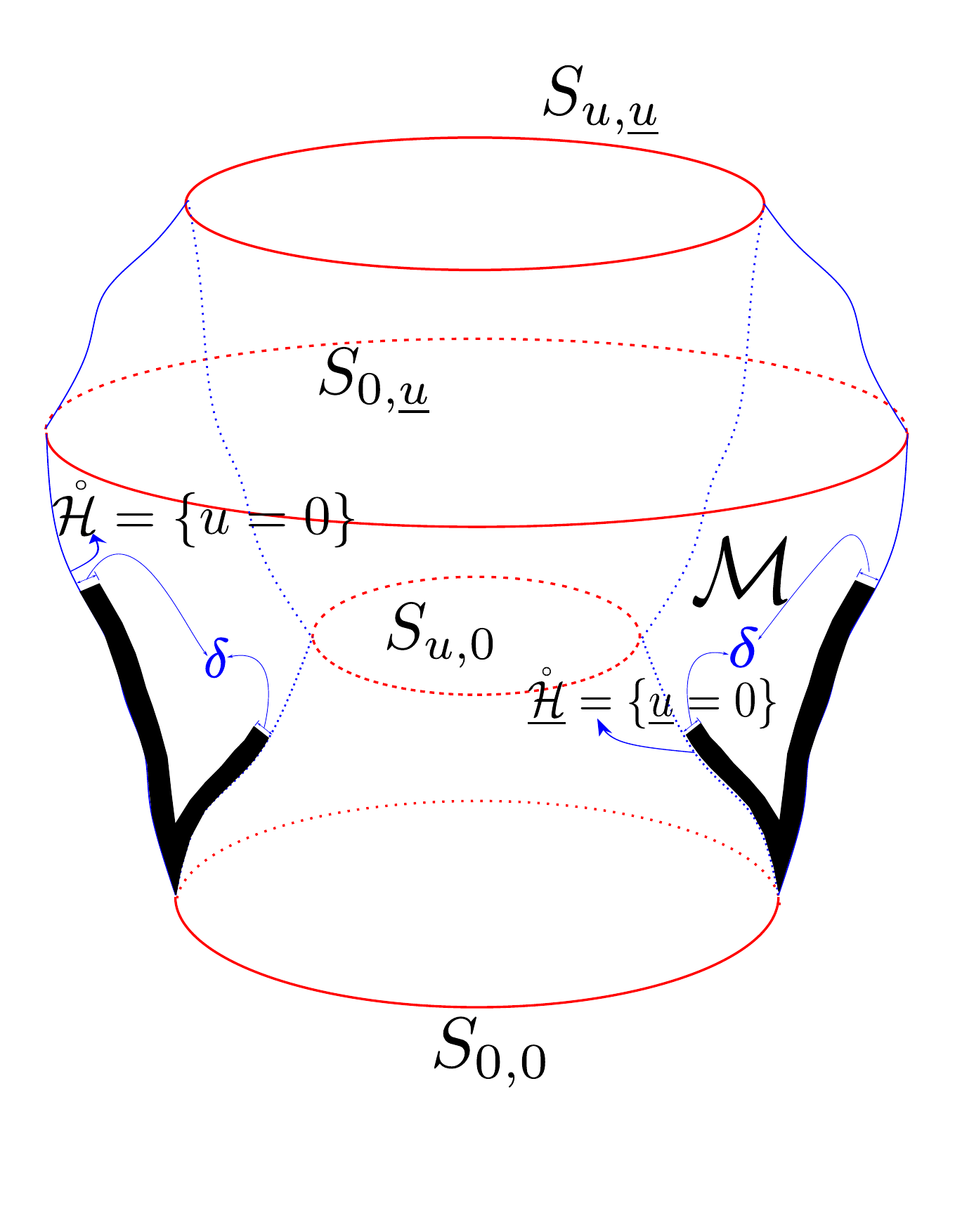}  
		}  
		\caption{The regions of local existence and uniqueness illustrated in Cartesian spacetime.}  
		\label{fig1}   
	\end{figure}

	In this article, we study the Cauchy problem for the $3D$ compressible Euler equations with vorticity and entropy
	given suitable \emph{characteristic initial data}, i.e., data on a pair of transversally intersecting
	cones that are null with respect to the data-dependent intrinsic geometry of the equations
	(also known as ``sound cones'' or ``characteristics'' - see Fig.\,\ref{fig1}) and whose intersection is diffeomorphic to $\mathbb{S}^2$.
	Our results can roughly be summarized as follows:
	\begin{theorem}[Main theorem, heuristic version]
	\label{TH:MAINTHEMVERYQUICK}
		The characteristic initial value problem (CharIVP)
		for the $3D$ compressible Euler equations \eqref{EQ:introEuler}, 
		with suitable $C^{\infty}$ data prescribed on a pair of transversally intersecting
		sound cones that are characteristic with respect to the data, 
		has a local, unique, $C^{\infty}$ solution in a \underline{future-neighborhood of the cones}. See Fig.\,\ref{fig3},
		where the region of existence is depicted relative to acoustical null coordinates $(u,\ubar,\theta^1,\theta^2)$,
		and Fig.\,\ref{fig1}, where the region of existence is depicted relative to standard Cartesian coordinates.
		Here, $u$ and $\ubar$ are acoustical eikonal functions, i.e., solutions to the acoustic eikonal equation \eqref{eq: eikonal equations},
		their level sets are null hypersurfaces (i.e., characteristics), and $(\theta^1,\theta^2)$ are local ``angular'' coordinates on the 
		intersection of the level sets of $u$ and $\ubar$, which are diffeomorphic to $\mathbb{S}^2$.
		Moreover, various kinds of Sobolev regularity for the data
		are propagated to the solution by the flow, and the solution can be continued as long as it
		satisfies appropriate a priori estimates and the cones don't degenerate.
		In particular, even though we have assumed the data are $C^{\infty}$ to infer existence, 
		our a priori energy estimates close at finite Sobolev regularity.
	\end{theorem}
	We refer readers to Theorem~\ref{TH:intromain2} for a more precise statement of the main results
	and to Theorems~\ref{TH:localwellposedness} and \ref{TH:continuationcriteria} for fully detailed statements.
	
	The main novelties of this paper can be summarized as follows:
	\begin{itemize}
		\item Compared to Rendall's celebrated solution \cite{ReductionoftheCharacteristicInitialValueProblem} 
			of the CharIVP for quasilinear wave equations of type 
			$(\gfour^{-1})^{\alpha \beta}(\Psi) \partial_{\alpha} \partial_{\beta} \Psi = \mathcal{N}(\Psi,\partial \Psi)$, 	
			to prove our results, we must overcome interesting new challenges, tied in particular to the difficulty of solving for the fluid's transversal 
			derivatives on the initial hypersurfaces. In a nutshell, 
			solving for the transversal derivatives is a much more subtle task for fluids compared to wave equations;
			see Prop.\,\ref{P:DETERMININGTRANSVERALDERIVATIVESANDSATISFYINGEQUATIONSONINITIALNULLHYPERSURFACES}.
		\item Compressible Euler flow is a \emph{multi-speed} hyperbolic system, i.e., in the physical space picture,
			its characteristics consist of sound cones as well as
			the flow lines of the material derivative vectorfield 
			$\materialderivative :=\p_t + v^a\p_a$, where $v$ is the velocity.
			Handling the multi-speed aspect of the system requires
			substantial new ideas compared to the pure wave case. To this end, we rely on
			the geometric formulation of compressible Euler flow
			derived in \cite{ANewFormulationofthe3DCompressibleEulerEquationswithDynamicEntropy},
			where the barotropic case was first handled in 
			\cite{ThehiddennullstructureofthecompressibleEulerequationsandapreludetoapplications};
			see Sect.\,\ref{SSS:NEWFORMULATIONSCHEMATIC} for an overview of this formulation.
		\item Our existence-region is a future-neighborhood of the subset of the initial null hypersurfaces 
			on which the solution jet 
			(i.e., the array comprising the solution and all of its partial derivatives -- including transversal ones -- to all orders) is well-defined,
			i.e., the same kind of region treated by Luk in his study \cite{GRCharIVP} of CharIVP for
			Einstein's equations. 
			In contrast, Rendall's approach for wave equations yielded existence only in a future-neighborhood
			of the intersection of the two null hypersurfaces; see Figs.\,\ref{fig3} and \ref{fig1}.
		\item Our existence-region is covered by a double-null foliation, i.e., solutions $u$ and $\ubar$ to the acoustic eikonal equation:
			\begin{align}  \label{E:INTROEIK}
					(\gfour^{-1})^{\alpha\beta}\p_\alpha u\p_\beta u=&0,
					&& (\gfour^{-1})^{\alpha\beta}\p_\alpha \ubar\p_\beta \ubar=0,
				\end{align}
			where $\gfour$ is the acoustical metric (see Def.\,\ref{DE:introacousticalmetric}), 
			a fluid-dependent Lorentzian metric tied to the intrinsic causal structure of the flow. 
			\eqref{E:INTROEIK} implies that the level sets of $u$ and $\ubar$ are characteristic for the compressible Euler equations.
			In our main results, the level sets of $u$ and $\ubar$ intersect transversally in surfaces diffeomorphic to $\mathbb{S}^2$.
			This is the first paper that rigorously constructs double-null foliations in the context of well-posedness theory for
			multi-dimensional compressible fluids.
		\item 
			There is a caveat compared to Einstein's equations: for the compressible Euler equations, 
			even if the fluid data along the initial null hypersurfaces are $C^{\infty}$
			in tangential directions, 
			the derivatives of the solution in directions transversal to the data-null hypersurfaces can still blow up,
			thanks to Riccati-type nonlinearities; see equation \eqref{E:WAVEEQUATIONFORDENSITYINADAPATEDCOORDIANTESALONGINITIALNULL}. 
		\item Admissible fluid data must satisfy \emph{constraint} equations (see \eqref{E:constraint1}--\eqref{E:constraint2}) 
			along each null hypersurface,
			a feature that does not occur in the case of Rendall's quasilinear wave equations.
			 Moreover, we view the initial null hypersurfaces as part of the data that also must satisfy constraint equations 
			(see, e.g., \eqref{E:FIRSTINITIALEMBEDDINGISNULL}--\eqref{E:FIRSTINITIALEMBEDDINGNULLDIRECTIONORTHOGONALTOANGULARDIRECTIONS})
			coupled to the fluid data, 
			because the question of whether a given hypersurface in $\mathbb{R}^{1+3}$ is null depends on the fluid data along the hypersurface itself.
		\item The fluid data must satisfy various ``corner compatibility conditions'' at the intersection of the two initial null hypersurfaces
			(see Sect.\,\ref{S:FLUIDCHARDATAFURTHERASPECTS}),
			which is an issue that does not arise in the case of Rendall's equations.
		\item We are able to close our a priori estimates by commuting the equations only with Cartesian coordinate partial derivatives $\pmb{\partial}$, 
			as opposed to geometric derivatives tied to $u$ and $\ubar$. This leads to a \emph{simplified approach} to a priori estimates.
			We use the eikonal functions $u$ and $\ubar$ only to control the ``shape'' of the region of existence. For this reason,
			to close our proof, we only need to control $u$ and $\ubar$ at the low derivative levels, i.e., 
			$L^{\infty}$-control of $\pmb{\partial}^{\leq 3} u$ and $\pmb{\partial}^{\leq 3} \ubar$ is more than enough 
			(we have not attempted to optimize the number of derivatives). Such $L^{\infty}$ estimates are easy to derive by treating \eqref{E:INTROEIK}
			as (fully nonlinear) transport equations with $C^{\infty}$ coefficients;
			see the proof of Prop.\,\ref{PR:BoundforGeometricQuantities}.
		\item Compared to standard estimates, we propagate one extra degree worth of Sobolev regularity for the vorticity and entropy;
			see Theorem~\ref{TH:UPTOTOPORDERENERGYESTIMATES}, where the extra regularity is afforded by the spacetime integrals ``$\KEnergy$.''
			This extra derivative is important for applications to the formation of shocks and the structure of maximal globally hyperbolic developments
			\cite{TheEmergenceoftheSingularBoundaryfromtheCrease,TheEmergenceoftheCauchyHorizonfromtheCrease}.
			To derive this improved regularity, 
			we use the aforementioned geometric formulation of compressible Euler flow
			from \cite{ANewFormulationofthe3DCompressibleEulerEquationswithDynamicEntropy}
			as well as the remarkable $L^2$-type integral identities derived in \cite{IntegralIdentities}.
		\item Unlike in the case of Einstein's equations, for the compressible Euler equations, the sound cones are necessarily embedded into Cartesian spacetime.
			The sound cones can therefore degenerate (e.g., pinch off) with respect to the Cartesian differential structure as the flow progresses. 
			Hence, the size of the region of existence in our main results is
			limited not just by the behavior of the fluid (which could blow up due to the nonlinearities), 
			but also by the behavior of the null hypersurfaces (which could collapse).
			This is an issue for even the characteristic initial value problem for the standard linear wave equation if the 
			data are posed on highly distorted null cones (as opposed to flat light cones).
			In the context of Figs.\,\ref{fig3} and \ref{fig1}, the existence width
			$\updelta > 0$ might have to be chosen to be small to avoid the formation of degeneracies in the characteristics.
	\end{itemize}

	\begin{remark}[We do not construct the data or the initial null hypersurfaces here]
		\label{R:DONOTCONSTRUCTDATA}
		In this paper, we do not ``construct'' characteristic data, but rather, we assume that it exists and satisfies various 
		constraint and compatibility equations (which necessarily are satisfied by any $C^{\infty}$ solution).
		It is an interesting open problem to construct such data\footnote{See Section \ref{SS:WANGYUYUDATACONSTRUCTION} for a discussion of the work \cite{WangYuYu25} on constructing initial data on one characteristic hypersurface in $(t,u,\sphere^1,\sphere^2)$ coordinates.}, which have to be constructed
		in conjunction with the corresponding null hypersurfaces on which they are posed.
		Moreover, the data and null hypersurfaces cannot be independently constructed, for a given hypersurface
		can be null only for ``some'' initial data (and some hypersurfaces, such as $\lbrace t = 0 \rbrace$, can never be null).
	\end{remark}

	Before providing further context for our results, we first introduce the equations.

	\subsection{The Equations - First Versions}
	\label{SS:EQUATIONSFIRSTVERSION}	
	
	\subsubsection{A standard first-order formulation}
	\label{SSS:STANDARDFIRSTORDER}
	We first provide a standard first-order formulation of the $3D$ compressible Euler equations.
	We will use this formulation in conjunction with the geometric wave-transport formulation
	from \cite{ANewFormulationofthe3DCompressibleEulerEquationswithDynamicEntropy}, which we
	further discuss in Sect.\,\ref{SSS:NEWFORMULATIONSCHEMATIC}.
	Specifically, for $(t,x) \in \mathbb{R} \times \mathbb{R}^3$, the $3D$ compressible Euler equations 
	can be formulated as a first-order PDE system in the 
	velocity $v: \mathbb{R} \times \mathbb{R}^3 \rightarrow \mathbb{R}^3$, 
	the density $\varrho: \mathbb{R} \times \mathbb{R}^3 \rightarrow [0,\infty)$, 
	and the entropy $\Ent: \mathbb{R} \times \mathbb{R}^3 \rightarrow \mathbb{R}$,
	which accounts for thermodynamics. 
	Instead of working with the density, we prefer to work with the \emph{logarithmic density}:
	\begin{align}
		\logDensity
		& := 
		\ln \frac{\varrho}{\bar{\varrho}},
	\end{align}
	where $\bar{\varrho} > 0$ is a constant. We will avoid regions with vanishing density and thus
	$\logDensity$ will be well-defined.
	To close the equations, we assume an \emph{equation of state}, that is, that the pressure $p$ is a given function
	of $(\varrho,\Ent)$. We allow for any equation of state with positive sound speed $\speed$,
	defined by:
	\begin{align} \label{E:SPEEDDEF}
		\speed
		& :=\sqrt{\frac{\p \pressure}{\p\density}},
	\end{align}
	where $\frac{\p \pressure}{\p\density}$ is the partial derivative of the equation of state with respect to density at fixed entropy.
	We will typically regard $p$ and $\speed$ as functions of $(\logDensity,\Ent)$.
	
	Given the above considerations, relative to the standard 
	``Cartesian coordinates'' $ x^0 :=t, x^1,x^2,x^3$ on ``Cartesian space'' $\mathbb{R}^{1+3}$, the equations can be written as follows,\footnote{We use the Einstein summation convention in that repeated adjacent indices are summed over their respective ranges, which are $1-3$ for Latin indices and $0-3$ for Greek indices.}
	where $\partial_{\alpha} := \frac{\partial}{\partial x^{\alpha}}$:
	\begin{subequations}\label{EQ:introEuler}
		\begin{align}
			\label{EQ:introBrho}\materialderivative\logDensity&=-\dive v,\\
			\label{EQ:introBv}\materialderivative v^i&= - \speed^2\delta^{ia}\p_a\logDensity-e^{-\logDensity}\frac{\pressure_{;\Ent}}{\backgroundDensity}\delta^{ia}\p_a\Ent,\\
			\label{EQ:introBs}\materialderivative s&=0,
		\end{align}
	\end{subequations}
	where $\delta^{ab}$ is the Kronecker delta, 
	\begin{align} \label{E:MATDERIV}
		\materialderivative
		& :=\p_t + v^a\p_a 
	\end{align}
	is the material derivative vectorfield,
	$\dive v := \partial_a v^a$ is the Euclidean divergence of the velocity,
	and on the RHS of \eqref{EQ:introBv}, $\pressure_{;\Ent}$ denotes 
	the partial derivative of the equation of state with respect to $\Ent$ at fixed $\logDensity$.

	\subsubsection{A new formulation of the flow - schematic version}
	\label{SSS:NEWFORMULATIONSCHEMATIC}
	Our work heavily relies on a new formulation of the equations derived in \cite{ANewFormulationofthe3DCompressibleEulerEquationswithDynamicEntropy,ThehiddennullstructureofthecompressibleEulerequationsandapreludetoapplications}, where the authors found that the flow can be decomposed into a ``sound-wave-part'' (``wave-part'' for short) for $(\logDensity,\Ent,v)$ and a ``transport-div-curl-part'' (``transport-part'' for short) for the specific vorticity $\vortrenormalized$ and the entropy gradient $\gradEnt$ (see Defs.\,\ref{D:LOGDENSITYANDSPECIFICFORTICITY} and \ref{D:entgradient} for their precise definitions). 
	This geometric formulation can be schematically expressed as follows:\footnote{Throughout the paper, schematic spatial partial derivatives and space-time partial derivatives are denoted by $\p$ and $\pfour$, respectively. Additionally, the following schematic notations are consistently used, where $A,B,C$ represent arrays of variables:
		
		\begin{itemize}
			\item $\linsmoothfunction[A](B)$ denotes any scalar-valued function that is linear in the components of $B$ with coefficients that are a function of the components of $A$.
			\item $\quadsmoothfunction[A](B,C)$ denotes any scalar-valued function that is quadratic in the components of $B$ and $C$ with coefficients that are a function of the components of $A$.
	\end{itemize}}\\

	\textbf{Geometric wave equations}	
	\begin{align}
		\label{EQ:intro1}\boxg\variables=\linsmoothfunction(\vvariables)[\vectvort,\DivGradEnt]+\quadsmoothfunction(\vvariables)[\pfour\vvariables,\pfour\vvariables].
	\end{align}
	
	\textbf{Transport equations}
	\begin{subequations}\label{EQ:intro2}
		\begin{align}
			\label{EQ:intro2.1}	&\materialderivative\vortrenormalized^i=\linsmoothfunction(\vvariables,\vvort,\vgradientEnt)[\pfour\vvariables],\\ \label{EQ:intro2.2}&\materialderivative\gradEnt^i=\linsmoothfunction(\vvariables,\vgradientEnt)[\pfour\vvariables].
		\end{align}
	\end{subequations}
	
	\textbf{Transport-Div-Curl system}
	\begin{subequations}\label{EQ:intro3}
		\begin{align}
			\label{EQ:intro3.1}&\materialderivative \modivort^i=\mathfrak{F}_{(\modivort^i)}:=\quadsmoothfunction(\vvariables)[\pfour\vvariables,(\pfour\vortrenormalized,\pfour\vgradientEnt,\pfour\vvariables)]+\quadsmoothfunction(\vgradientEnt)[\pfour\vvariables,\pfour\vvariables]+\linsmoothfunction(\vvariables,\vvort,\vgradientEnt)[\pfour\vvariables,\pfour\vgradientEnt],\\
			\label{EQ:intro3.2}&\materialderivative \DivGradEnt=\mathfrak{F}_{(\DivGradEnt)}:=\linsmoothfunction(\vvariables,\vgradientEnt)[\pfour\vortrenormalized]+\quadsmoothfunction(\vvariables)[\pfour\vvariables,(\pfour\vgradientEnt,\pfour\vvariables)]+\quadsmoothfunction(\vgradientEnt)[\pfour\vvariables,\pfour\vvariables]+\linsmoothfunction(\vvariables,\vvort,\vgradientEnt)[\pfour\vvariables],\\
			\label{EQ:intro3.3}& \curl S=0,\\ 
			\label{EQ:intro3.4}& \dive \vortrenormalized =\linsmoothfunction(\vortrenormalized)[\pfour\vvariables],
		\end{align}
	\end{subequations}
	where $\vvariables:=(v^1,v^2,v^3,\logDensity,\Ent)$, $\gfour=\gfour(\vvariables)$ is a solution-dependent Lorentzian metric which governs the geometry of sound waves, $\modivort\simeq  \curl \vortrenormalized$ and $\DivGradEnt\simeq\text{div} \gradEnt$ are special ``modified fluid variables.'' See (\ref{DE:vortrenomal}) for the definition of $\vortrenormalized$, (\ref{DE:ent}) for the definition of $\gradEnt$, (\ref{DE:CD}) for the definition of $\modivort,\DivGradEnt$, (\ref{D:acousticmetric}) for the definition of $\gfour$, and Theorem~\ref{T:geometricformulation} for the precise geometric formulation of the equations. This formulation reveals miraculous regularity and geometric properties of the flow, which play a fundamental role in our current study. 
These geometric properties are not visible in the first-order formulation \eqref{EQ:introEuler}.
	

	\subsection{Further Context}
	\label{SS:FURTHERCONTEXT}
	
	\subsubsection{Previous work on CharIVP}
	\label{SSS:PREVIOUSCHARIVPWORK}
	Besides Rendall's \cite{ReductionoftheCharacteristicInitialValueProblem} and Luk's \cite{GRCharIVP} aforementioned papers 
	on the CharIVP for quasilinear wave equations and Einstein's equations, we also mention
	\cite{aCpCrT2016}, which in addition treats special classes of symmetric hyperbolic systems that exhibit a special ``splitting'' property relative
	to a characteristic coordinate system; the compressible Euler equations do not fall under the scope of \cite{aCpCrT2016}.
	The only other paper on CharIVP for the compressible Euler equations is
	Lisibach's work \cite{CharivpLisibach} on the spherically symmetric problem for barotropic equations of state, i.e., 
	when the pressure is a function of density alone (as opposed to the general case we treat here, in which the pressure can depend on density and entropy). The spherically symmetric problem, while interesting, is substantially easier because it can be treated using transport equation techniques, without the need for energy estimates; see Sect.\,\ref{SS:intro1D} for further discussion of the symmetric problem.
	
	\subsubsection{Connections to related problems}
	\label{SSS:CONNECTIONS}
	Although CharIVP is of interest in itself, 
	another key motivation for solving it is that the setup can be
	of immense value for studying the long-time behavior of solutions. Inspiring examples for the Einstein-vacuum equations
	include Christodoulou's proof of the formation of trapped surfaces \cite{TheFormationofBlackHolesinGR} 
	and the Dafermos--Luk proof of the $C^0$ nonlinear stability of the Cauchy horizon inside Kerr black holes 
	\cite{StabilityofKerrBlackHole}.
	
	In the context of compressible fluid flow, CharIVP has deep connections to the study of shock formation, maximal globally hyperbolic developments,
	and weak solutions, as we now explain. We first note that even if the data are $C^{\infty}$, 
	one cannot generally hope to avoid singularities for all time.
	That is, it is well-known that even in the irrotational and isentropic case, the compression of sound waves can cause shocks (gradient-blowup) to develop in flows launched by smooth initial data. Moreover, in more than one space dimension and away from symmetry, with or without dynamic vorticity and entropy, these singularities are known to be stable; see \cite{Theformationofshocksin3-dimensionalfluids,CompressibleflowandEuler'sequations,ShockFormationandVorticityCreationfor3DEuler,Smalldatashocks,Shockformationinsolutionstothe2DcompressibleEulerequationsinthepresenceofnon-zerovorticity,Speck2016,ShkollerCPAM,LukSpeck24,MiaoYu17}. 
We also highlight that shocks are not the only kinds of singularities that can develop; in \cite{Rodnianski20191,Rodnianski20192}, the authors proved that implosion singularities (infinite density blowup) can develop starting from $C^{\infty}$ data, which is a more severe singularity compared to shocks (gradient-blowup).
	
In view of finite speed of propagation, one can continue to study the classical solution 
in spacetime regions that contain points to the future of the constant-time hypersurface of first blowup.
The holy grail object in this regard is a \emph{maximal globally hyperbolic development} (MGHD) of the data.
Roughly speaking, an MGHD is a largest possible globally hyperbolic solution + region on which the solution is classical and determined by the initial data; 
we refer to the survey article \cite{abbresciaspecknotices} and the appendices of \cite{AbbresciaSpeckYu2026} for a detailed discussion of MGHDs.
Here, ``globally hyperbolic'' means that the manifold on which the solution exists contains a Cauchy hypersurface $\Sigma$, i.e., a surface such that every inextendible curve that is timelike with respect to the acoustical metric $\gfour$ (see Def.\,\ref{DE:introacousticalmetric}) intersects $\Sigma$ precisely once; clearly the notion of an MGHD relies on the fact that the causal structure of the equations is governed by the Lorentzian metric $\gfour$.
For $3D$ irrotational and isentropic relativistic fluids, in the
breakthrough work \cite{Theformationofshocksin3-dimensionalfluids}, Christodoulou constructed a large, implicit portion of the MGHD up to the boundary
for open sets of shock-forming solutions arising from data near that of a trivial solution with vanishing velocity and constant density.
\cite{Theformationofshocksin3-dimensionalfluids} was the first
work to study the structure of an MGHD boundary for a shock-forming solution to a multi-dimensional quasilinear hyperbolic PDE without symmetry assumptions.
In the $2D$ isentropic case,
\cite{ShkollerVicol2024} constructed an $\mathcal{O}(\epsilon)$-size portion of the boundary of the MGHD for open sets of shock-forming solutions
whose data has gradients of size $\mathcal{O}(1/\epsilon)$ and satisfy a convexity condition.
In the $3D$ case with vorticity and entropy, 
\cite{TheEmergenceoftheSingularBoundaryfromtheCrease,TheEmergenceoftheCauchyHorizonfromtheCrease}
constructed an $\mathcal{O}(1)$-size portion of the boundary of the MGHD up to the boundary for open sets of shock-forming solutions that satisfy a weaker
``transversal convexity'' condition, which is essentially the minimal condition needed for the stability of the structure of the MGHD boundary. In particular, in 
\cite{ShkollerVicol2024,TheEmergenceoftheCauchyHorizonfromtheCrease}, part of the MGHD boundary is a Cauchy horizon, 
which is a null hypersurface that emanates from a co-dimension $2$ manifold where the fluid's gradient blows up.
To follow the solution up the Cauchy horizon, \cite{TheEmergenceoftheCauchyHorizonfromtheCrease} used arguments based on a priori estimates 
and continuation criteria on double-null foliations,
i.e., it used some of the results of this article.
See also the work \cite{AbbresciaSpeckYu2026} in spherical symmetry for global solutions as well as unique maximal global hyperbolic developments that terminate in shocks
for open sets of small initial data.

One can also study a class of weak solutions that develop starting from the co-dimension $2$ manifolds of first gradient blowup mentioned in the previous paragraph. 
Describing the transition of the solution from classical to weak as well as constructing the shock hypersurface, across which the weak solution exhibits a jump discontinuity
dictated by Rankine--Hugoniot jump conditions,
is called the shock development problem. In his celebrated monograph
\cite{TheShockDevelopmentProblem}, Christodoulou solved the ``restricted'' shock development problem without symmetry assumptions.
Here, the word ``restricted'' means that Christodoulou ignored the jump in entropy and vorticity that generally occurs in fluids across shock hypersurfaces. That is, he solved an idealized problem by constructing a weak solution to a quasilinear wave equation that agrees with the compressible Euler equations in regions where the flow is classical, irrotational and isentropic. In \cite{TheShockDevelopmentProblem}, it is shown that the ``initial data'' for the shock development problem
can be taken to be a classical MGHD for a solution that exhibits gradient blowup, and whose boundary consists of a characteristic surface portion along which the solution is smooth (the aforementioned Cauchy horizon), as well as another ``singular'' portion along which the fluid gradient blows up. In particular, characteristic hypersurfaces play a fundamental role in setting up the shock development problem.
Relatedly, characteristic surfaces play a fundamental role in the study of
shock reflection and the shock interaction problems; see \cite{lisibach2021shock,lisibach2022shock} for such problems in plane symmetry and \cite{wang2023shock} for the shock interaction problem in spherical symmetry.
We also mention that \cite{SimultaneousDevelopmentofShocksandCusps} solved the shock development problem in azimuthal symmetry with initial data given near the initial shock

\subsection{Remarks on the Standard Cauchy Problem for the Compressible Euler Equations}
The solution to the standard Cauchy problem for the compressible Euler equations, where initial data is given on a spacelike hypersurface, 
is a classic result that can be directly established \cite{kato,Kato1976} for the first-order formulation \eqref{EQ:introEuler} of the flow
via an iteration scheme, a priori energy estimates, Sobolev embedding, and the fractional Sobolev–Moser calculus. 
In \cite{3DCompressibleEuler,Alterproof,Roughsolutionsofthe3-DcompressibleEulerequations}, based on the aforementioned geometric formulation from 
Sect.\,\ref{SSS:NEWFORMULATIONSCHEMATIC},
the authors proved low-regularity well-posedness, with the regularity assumption on the fluid variables $(v^i,\logDensity,\Ent)$ being $\frac{1}{2}$-less than that in the classical results \cite{kato,Kato1976}. The proofs were based on the Strichartz estimates and energy estimates; see \cite{AGeoMetricApproach,Equationsd'ondesquasilineairesetestimationdeStrichartz,Equationsd'ondesquasilinearireseteffetdispersif,ImprovedLocalwellPosedness,RoughsolutionstotheEinsteinvacuumequations,Strichartzestimatesforoperatorswithnonsmooth,Strichartzestimatesforsecondorderhyperbolicoperators,counterexamplestolocalexistenceforquasilinearwaveequations,Sharplocalwell-posednessresultsforthenonlinearwaveequation,ACommutingVectorfieldsApproachtoStrichartz} for further insights into the development of the low-regularity problem in quasilinear wave equations. A parallel development exists for relativistic fluids. For example, Disconzi--Speck  \cite{TheRelativistivEulerEquations} derived a geometric formulation of the relativistic Euler equations with vorticity and dynamic entropy, and \cite{yu2022rough} proved low-regularity local well-posedness.

	\subsection{``CharIVP'' in Spherical Symmetry}\label{SS:intro1D}
	In \cite{CharivpLisibach}, the author considers the ``CharIVP'' in spherical symmetry for barotropic (i.e., the pressure is a function of density alone) 
	compressible Euler equations, where the characteristic system can be written as the following coupled system of transport equations:\\
	
	\textbf{On $\hypersur:=\{u=\text{const}\}$:}
	\begin{subequations}
		\begin{align}
			\label{EQ:Riemann1}
			\frac{\partial}{\partial \ubar} \mathcal{R}_{+}
				& =-\frac{\speed(\mathcal{R}_{+},\mathcal{R}_{-})}{r}(\mathcal{R}_{+}-\mathcal{R}_{-}),\\
		\frac{\partial}{\partial \ubar} r & = 
		\frac{1}{2}(\mathcal{R}_{+}-\mathcal{R}_{-})+\speed(\mathcal{R}_{+},\mathcal{R}_{-}).
		\end{align}
	\end{subequations}

	\textbf{On $\hypersurbar:=\{\ubar=\text{const}\}$:}
	\begin{subequations}
		\begin{align}
			\label{EQ:Riemann2}
			\frac{\partial}{\partial u} \mathcal{R}_{-} & = -\frac{\speed(\mathcal{R}_{+},\mathcal{R}_{-})}{r}(\mathcal{R}_{+}-\mathcal{R}_{-}),\\
			\frac{\partial}{\partial u} r
			& = \frac{1}{2}(\mathcal{R}_{+}-\mathcal{R}_{-})-\speed(\mathcal{R}_{+},\mathcal{R}_{-}),
		\end{align}
	\end{subequations}
	where $\mathcal{R}_{+},\mathcal{R}_{-}$ are a pair of Riemann invariants (introduced in \cite{Riemann1860}), $r$ is the Euclidean radial variable, 
	$u,\ubar$ are the geometric coordinates, and $\speed$ is the speed of sound. In the above system, one has the freedom to prescribe $\mathcal{R}_{-}$ on $\inihypersur$ ($\{u=0\}$) and $\mathcal{R}_{+}$ on $\inihypersurbar$ ($\{\ubar=0\}$). 
	
	The above system yields important insights about the initial data 
	for the compressible Euler equations:
	\begin{itemize}
		\item The initial hypersurfaces and the initial fluid data are coupled, so one has to construct the initial hypersurfaces and initial fluid data simultaneously. In other words, the characteristic initial hypersurfaces can/should be viewed as a part of the initial data.
		\item A crucial difference compared to the case of the standard Cauchy problem (in which the data are
		prescribed on a space-like hypersurface) is that in the characteristic case, \textbf{some components of the
			initial data are constrained}, that is, determined by the free components by virtue of equations obtained by restricting the PDE to the initial characteristic hypersurfaces.
	\end{itemize}
	
	\subsection{Further Discussion of the Initial Data and the Initial Embeddings} \label{SS:DISCUSSIONDATAANDEMBEDDINGS}
	As we alluded to in Section \ref{SS:intro1D}, for CharIVP, we can only freely prescribe
	``part'' of the initial fluid data; the remaining part is constrained by the compressible Euler equations and should be determined by 
	``constraint PDEs'' restricted to the initial hypersurfaces. In addition, since the notion of a characteristic hypersurface is coupled to 
	the fluid variables, any attempt to construct initial characteristic hypersurfaces is necessarily coupled to the construction of fluid data along them.
	Furthermore, away from symmetry, 
	in the process of constructing initial fluid data and initial hypersurfaces, there is a natural phenomenon of loss of derivatives for the transversal derivatives along null hypersurfaces. This issue plays only a minor role in the present article since all of our data are $C^{\infty}$, 
	and our a priori estimates show that there is no further loss of derivatives from the data to the solution.
	
	\subsubsection{A model problem}\label{SSS:AModel}
	As a model problem, consider the standard linear wave equation $\p_t^2\variables-\Delta \variables=0$, which can be expressed as:
	\begin{align}\label{EQ:Lwave}
		\Lunit(\uLunit\variables) & = \anglap\variables + l.o.t,
	\end{align} 
	where $\Lunit:=\p_t+\p_r, \uLunit:=\p_t-\p_r$, and $\anglap$ is the angular Laplacian (on the Euclidean round spheres). 
	For the wave equation, we can freely prescribe data for $\variables$ on $\inihypersur=\{t-r=0\}$ and $\inihypersurbar=\{t+r=0\}$ (in Figure \ref{fig1}, we have displayed analogous ``curved cones'' for the ``CharIVP'' for $3D$ compressible Euler flow). In particular, unlike the fluid case, there are no constraint equations for the wave equation data.
	Note that tangential derivatives at all orders are also prescribed. 
	In particular, because $\inihypersur$ and $\inihypersurbar$ intersect transversally, 
	on the initial sphere $S_{0,0}=\inihypersur\bigcap\inihypersurbar$, all the first derivatives of $\variables$ can be computed. In characteristic initial value problems, the transversal derivatives along the data hypersurfaces have to be computed by solving an evolution equation determined by the PDE. For example, for the linear wave equation, 
	(\ref{EQ:Lwave}) implies that along $\inihypersur$, $\uLunit\variables$ satisfies a transport equation with $\anglap\variables$ as a source term, where the data of $\uLunit\variables$ on $S_{0,0}$ is known. This equation loses derivatives in the sense that along a data null hypersurface ruled by integral curves of $\Lunit$, we can control
the data of $\uLunit\variables$ only by integrating \eqref{EQ:Lwave} along the integral curves of $\Lunit$, and this equation
features the derivative-losing source term $\anglap\variables$. 
This is an example of a phenomenon that typically occurs for multi-dimensional ``CharIVPs'' for hyperbolic PDEs: some transversal derivatives (specifically, $\uLunit\variables$ along $\inihypersur$ in the case of the linear wave equation) have to be solved for via transport equations that lose derivatives. In particular, for 3D compressible Euler, we find that a similar but more complicated phenomenon occurs, and moreover, some of the tangential derivatives must obey constraint equations (again, such constraint equations do not play a role in the above analysis for the linear wave equation); see (\ref{E:constraint1}) for the precise constraint equation on $\inihypersur$. 
In this paper, the derivative loss occurs only in a ``soft way,'' because our data are $C^{\infty}$. That is, our evolution estimates do not lose any Sobolev regularity compared to the initial data. Moreover, since our data are $C^{\infty}$, the loss of a finite amount of Sobolev regularity for the transversal derivatives of the data is not important in the sense that along the initial null hypersurfaces, all transversal derivatives to all orders are $C^{\infty}$.
	
	\subsubsection{Characteristic initial data and initial embeddings for compressible fluids }\label{SSS:introInitialEmbeddingsFORCOMPRESSIBLEFLUIDS}
	To further explain our approach to proving our main results, we first discuss the acoustic geometry for the compressible Euler flow, that is, the geometry of the fluid characteristics. The geometry of the flow is intimately connected to the acoustical metric $\gfour$ and its inverse $\gfour^{-1}$, which can be expressed 
	relative to the Cartesian coordinates as follows, where $\materialderivative =\p_t + v^a\p_a$ is the material derivative:
	\begin{subequations}\label{DE:introacousticalmetric}
		\begin{align}
			\gfour&:=-\diff t \otimes \diff t + \speed^{-2}\sum\limits^3_{a=1}(\diff x^a-v^a\diff t) \otimes (\diff x^a-v^a\diff t),\\
			\gfour^{-1}&:=-\materialderivative \otimes \materialderivative+\speed^2\sum\limits^3_{a=1}\p_a \otimes\p_a.
		\end{align}					
	\end{subequations}
	The initial data for ``CharIVP'' is a pair of transversally intersecting null hypersurfaces (null with respect to the acoustical metric $\gfour$ defined in (\ref{DE:introacousticalmetric})) embedded in Cartesian spacetime and denoted by $\inihypersur$ and $\inihypersurbar$, as well as 
	fluid data $(v^i,\density,\Ent)$ on $\inihypersur$ and $\inihypersurbar$; see Figure \ref{fig1}. We stress that 
	\textbf{the embeddings are part of the data.} Given this setup, we solve for appropriate transversal derivatives of the data in terms of tangential derivatives by solving transport equations that lose derivatives, much like in the case of the linear wave equation. However, compared to the linear wave equation, the analysis is more complicated due to the tensorial structure of the 3D compressible Euler equations and the presence of vorticity and entropy. 
	To overcome the above difficulties, we use both the first-order formulation (\ref{EQ:introEuler}) and the geometric formulation (\ref{EQ:intro1})-(\ref{EQ:intro3}) and proceed as follows:
	\begin{enumerate}[label=\arabic*.]
		\item We assume we are given 
		initial $C^{\infty}$ embeddings $\embedding$ and $\uembedding$ of the initial hypersurfaces $\inihypersur$ and  $\inihypersurbar$ respectively
		into Cartesian spacetime and $C^{\infty}$ fluid data along each hypersurface 
		such that: \textbf{1)}  $\inihypersur$ and  $\inihypersurbar$ are null with respect to the acoustical metric $\gfour$ 
		and \textbf{2)}: The fluid data on $\inihypersur$ and $\inihypersurbar$ satisfy the aforementioned constraint equations (see \eqref{E:constraint}), which are imposed on them by the Euler equations. 
		These assumptions are coupled in the sense that the constraint equations depend on the embeddings of $\inihypersur$ and $\inihypersurbar$ into the Cartesian spacetime, and the initial embeddings depend on the fluid variables (see Section \ref{SS:initialembedding2}). Moreover, 
		we assume that the fluid data satisfy various corner compatibility conditions\footnote{See Section \ref{SS:CORNERcompatibilityCONDITIONS} for the corner compatibility conditions in geometric coordinates.} at $S_{0,0}$; these conditions are needed to ensure that
		the compressible Euler equations \eqref{EQ:introEuler} are satisfied to all orders at $S_{0,0}$ and that the solution's derivatives to all orders
		are well-defined at $S_{0,0}$.
		\item Let $(\inilogDensity,\inivel^i,\inient)$ denote the fluid data on $\inihypersur$.
		We solve for the transversal derivative of $\inient$ in terms of tangential derivatives of $\inient$ using (\ref{EQ:introBs}). Then, we express the transversal derivative of $\mathring{v}^i$ and $\iniDivGradEnt$ in terms of tangential derivatives of $\mathring{v}^i$ and all derivatives of $\inilogDensity$ as follows (see Section \ref{SSS:DETERMININGTRANSVERALDERIVATIVESANDSATISFYINGEQUATIONSONINITIALNULLHYPERSURFACES} for $\iniDivGradEnt$ and (\ref{E:e3v}) for $\mathring{e}_3\mathring{v}^i$):
		\begin{align}\label{EQ:introe3v}
			\mathring{e}_3\mathring{v}^i,\iniDivGradEnt=\pfour(\inilogDensity,\inient)+(\mathring{e}_4,\iniangp)(\mathring{v}^i,\inigradEnt),
		\end{align}
		where $\mathring{e}_3$ is the transversal derivative and $\mathring{e}_4,\iniangp$ are the tangential derivatives on $\inihypersur$.
		\item We solve for $\mathring{e}_3\inilogDensity$ by using the wave equation (\ref{EQ:intro1}) for $\logDensity$ (see (\ref{E:boxgdensity})) with source terms $(\mathring{e}_4,\iniangp)(\mathring{v}^i,\inilogDensity,\inigradEnt)$ (see (\ref{EQ:e4e3rho})). Then, we solve for $\mathring{e}_3\mathring{v}^i,\iniDivGradEnt$ by using (\ref{EQ:introe3v}). Finally, we algebraically solve for $\inimodivort$ by the wave equation (\ref{E:boxgv}) for $v^i$; see Section \ref{SSS:inimodivort}.
		\item We carry out analogous analysis on the other data null hypersurface $\inihypersurbar$.
		\item We solve for the higher-order transversal derivatives using a similar procedure, based on commuting the equations; 
		see Prop.\,\ref{P:DETERMININGTRANSVERALDERIVATIVESANDSATISFYINGEQUATIONSONINITIALNULLHYPERSURFACES}.
		Moreover, Props.\,\ref{P:DETERMININGTRANSVERALDERIVATIVESANDSATISFYINGEQUATIONSONINITIALNULLHYPERSURFACES} and \ref{P:CORNERCOMPATIBILITYCONDITIONS}
		show that this procedure is well-defined and yields $C^{\infty}$ data on $\inihypersur \cup \inihypersurbar$
		such that the compressible Euler equations \eqref{EQ:introEuler}
		are satisfied to all orders on $\inihypersur \cup \inihypersurbar$.
	\end{enumerate}
	
	\subsubsection{The initial data construction in \cite{WangYuYu25}}
	\label{SS:WANGYUYUDATACONSTRUCTION}
	Recall that in the present paper, we have \emph{assumed} that we are given suitable $C^{\infty}$ data and initial null hypersurfaces; 
	we have not made any attempt to construct them.
	In contrast, a construction of initial data on a \emph{single} null hypersurface using $(t,u,\sphere^1,\sphere^2)$ coordinates was achieved by the second author with Wang and Yu in \cite{WangYuYu25}, where $t$ is the standard Cartesian time function and $u$ is an acoustic eikonal function.
	Because of the availability of the time function $t$, in \cite{WangYuYu25}, the sphere-normal vectorfield $\spherenormal$ and spherical-tangent vectorfields $e_A$ are not fluid-dependent because the induced metric on $\St$ is flat 
(just renormalize \eqref{DE:introacousticalmetric} by $\speed^2$). Consequently, equations guaranteeing that the hypersurface is null can be solved in conjunction with the aforementioned fluid constraint equation, i.e., one can construct large families of fluid data and null hypersurfaces.
We emphasize that in this article, we use the coordinates 
$(u,\ubar,\sphere^1,\sphere^2)$, where $u$ and $\ubar$ are a pair of eikonal functions,
and the spheres $\suubar$ given by the intersection of the level sets of $u$ and $\ubar$ are not generally 
subsets of the constant-Cartesian-time slice $\St$. This significantly increases the difficulty of any attempt to carry out a construction of data/null hypersurfaces since, without the availability of the time function $t$, 
all the vectorfields in the null frame are solution-dependent 
(in  \cite{WangYuYu25}, only one out of the four frame vectorfields was solution-dependent). Moreover, in the present article, we have 
initial data on both incoming and outgoing characteristic hypersurfaces. One cannot simply construct data separately on each null hypersurface and then union them together because satisfying the compressible Euler equations \eqref{EQ:introEuler} at all orders on the initial sphere $S_{0,0}:=\hypersur_{0}\cap\hypersurbar_0$ requires the aforementioned corner compatibility conditions 
(see Sect.\,\ref{SS:CORNERcompatibilityCONDITIONS}). In total, it remains an interesting open problem to carry out the construction of admissible $3D$ compressible fluid initial data on a pair of transversally intersecting null hypersurfaces.

	\subsection{Main Theorem and Ideas of the Proof} \label{SS:introAproof}
	Given the setup of the initial fluid data and the initial null hypersurfaces as discussed in Section \ref{SS:DISCUSSIONDATAANDEMBEDDINGS},
	we can now state a more precise version of our main results.
	
	\begin{theorem}[Local existence and uniqueness in a neighborhood of the data null hypersurfaces]\label{TH:intromain2}
		Consider $C^{\infty}$ fluid data for the compressible Euler equations \eqref{EQ:introEuler}
		that are prescribed on the union 
		$\inihypersur \cup \inihypersurbar$ as described in Section \ref{SSS:introInitialEmbeddingsFORCOMPRESSIBLEFLUIDS},
		where $\inihypersur$ and $\inihypersurbar$ are $C^{\infty}$ characteristic hypersurface portions 
		that intersect transversally in a set $S_{0,0}$ that is diffeomorphic to $\mathbb{S}^2$.
		Assume the fluid data satisfy the Euler constraints \eqref{E:constraint} and the 
		corner compatibility conditions
		from Def.\,\ref{D:CORNERCOMPATIBILITYCONDITIONSONDATA}.
		Let $\hypersur_{0}([0,\ubarfinite])$ and $\hypersurbar_{0}([0,\ufinite])$
		denote portions of $\inihypersur$ and $\inihypersurbar$ respectively such that all derivatives to all orders
		of the solution (including in transversal directions) are well-defined along them (see Remark~\ref{R:BLOWUPONDATAHYPERSURFACES}).
		Assume that for some integer $N \geq 17$,
		on $\hypersur_{0}([0,\ubarfinite]) \cup \hypersurbar_{0}([0,\ufinite])$,
		the initial data satisfy the norm bounds stated in Sect.\,\ref{SS:iniregularity2}.
	
		Then there exists a number $\updelta > 0$ and a corresponding
		future-neighborhood $\mathcal{M}^{\ufinite,\ubarfinite}_\updelta$
		of $\hypersur_{0}([0,\ubarfinite]) \cup \hypersurbar_{0}([0,\ufinite]) \subset \mathbb{R}^{1+3}$
		such that the following hold on $\mathcal{M}^{\ufinite,\ubarfinite}_\updelta$:
		\begin{itemize}
			\item There exists a unique $C^{\infty}$ solution $(\rho,v^1,v^2,v^3,s)$
			to the compressible Euler equations \eqref{EQ:introEuler}.
			\item The initial embeddings yield initial data for the geometric coordinates $(u,\ubar,\theta^1,\theta^2)$,
			such that there exist unique $C^{\infty}$ solutions to the following initial value problems:
			\begin{subequations}\label{eq: eikonal equations}
				\begin{align} 
					(\gfour^{-1})^{\alpha\beta}\p_\alpha u\p_\beta u=&0,&
					(\gfour^{-1})^{\alpha\beta}\p_\alpha \ubar\p_\beta \ubar=&0,\\
					(\gfour^{-1})^{\alpha\beta}\p_\alpha u\p_\beta \theta^1=&0,&
					(\gfour^{-1})^{\alpha\beta}\p_\alpha u\p_\beta \theta^2=&0.
				\end{align}
			\end{subequations}
			Moreover, the level sets of $u$ and $\ubar$ intersect transversally 
			in two-dimensional $g$-spacelike submanifolds 
			$S_{u',\ubar'} := \lbrace u = u' \rbrace \cap \lbrace \ubar = \ubar' \rbrace$,
			all of which are diffeomorphic to $\mathbb{S}^2$.
			\item The change of variables map $\embedgeo$ from geometric coordinates to Cartesian coordinates, i.e.,
			\begin{align} \label{DE:ChangeofVariableMap}
			\embedgeo(u,\ubar,\theta^1,\theta^2)
				& := (t,x^1,x^2,x^3),
			\end{align}
			is a $C^{\infty}$ diffeomorphism from the domain
			$[0,\updelta] \times [0,\ubarfinite] \times \mathbb{S}^2 \bigcup [0,\ufinite] \times [0,\updelta] \times \mathbb{S}^2$
			onto the image set $\mathcal{M}^{\ufinite,\ubarfinite}_\updelta$ in Cartesian coordinate space.
			\item For $k \approx \frac{N}{2}$, the Cartesian coordinate partial derivatives of the solution satisfy $L^2$-type energy estimates in terms of the data, 
			which can be schematically captured as:
			\begin{align} \label{E:SCHEMATICENERGYESTIMATES}
				\norm{\pfour^{\leq k} (\logDensity,v^i,s,\vortrenormalized,\gradEnt)}_{L^2(\mathcal{M}^{\ufinite,\ubarfinite}_\updelta)} 
				& \lesssim \norm{\pfour^{\leq k} (\logDensity,v^i,s,\vortrenormalized,\gradEnt)}_{L^2(\hypersur_{0}([0,\ubarfinite]) \cup \hypersurbar_{0}([0,\ufinite]))}.
			\end{align}
	\end{itemize}
	\end{theorem}
		We stress that the actual energy estimates are rather involved and are not perfectly captured by the schematic estimates in
	\eqref{E:SCHEMATICENERGYESTIMATES}; see Sect.\,\ref{S:APRIORIESTIMATESUPTOTOPORDER} for the details.
	The one issue that is well-captured by \eqref{E:SCHEMATICENERGYESTIMATES} is that the same number of derivatives appear on the left and right-hand sides,
	i.e., the a priori energy estimates for the solution do not lose any derivatives relative to the right-hand side.
	We refer readers to Section \ref{S:Lastslice} for detailed statements and proofs of our local existence and uniqueness result as well as continuation criteria.

	\subsubsection{Local existence and uniqueness in a neighborhood of the initial sphere}\label{SSS:neighborhood}
	We recall Rendall's result \cite{ReductionoftheCharacteristicInitialValueProblem}, in which he proved local existence and uniqueness for the ``CharIVP'' for quasilinear wave equations and Einstein vacuum equations in a neighborhood of the initial sphere $S_{0,0}:=\inihypersur\bigcap\inihypersurbar$ (see Figure \ref{fig3}). To prove his result, Rendall cleverly transforms CharIVP into a standard Cauchy problem (where initial data are given on a spacelike hypersurface), such that the solution to the new problem coincides with the ``CharIVP'' in a neighborhood of $S_{0,0}$. Below, we will explain how we adapt his approach to the compressible Euler equations, ensuring a (small) region where a solution uniquely exists; see Section \ref{S:localexistence} for detailed proofs.
	
	First, rewriting \eqref{EQ:introEuler} in the form of $A^{\alpha}[\vvariables]\p_{\alpha}\vvariables=0$, where, for each $\alpha=0,1,2,3$, $A^\alpha$ is a solution-dependent matrix with $A^0$ positive definite, we consider the following standard Cauchy problem for the unknown $\tilde{\vvariables}$:
	\begin{subequations}
		\begin{align} \label{E:introSPLITFIRSTORDEREULERINTOEXTENDEDDATAPLUSREMAINDER}
			A^{\alpha}[\vvariables_{(\text{ext})}+\tilde{\vvariables}]
			\p_{\alpha}
			\left(\vvariables_{(\text{ext})}
			+
			\tilde{\vvariables}
			\right)
			& = F,
			\\
			\tilde{\vvariables} \restriction_{\widetilde{\Sigma}_0^{\mathcal{N}}}
			& = 0,
			\label{E:intro0SPACELIKEDATAFORSPLITFIRSTORDEREULERINTOEXTENDEDDATAPLUSREMAINDER}
		\end{align}
	\end{subequations}
	where $\vvariables_{(\text{ext})}$ denotes appropriately constructed Whitney-extended fluid functions\footnote{The Whitney-extended fluid variables solve compressible Euler equations \eqref{EQ:introEuler} at all orders on, and only on, the initial characteristic hypersurfaces $\hypersur_{0}\cup\hypersurbar_{0}$.} of the prescribed initial data, $\mathcal{N}$ is a neighborhood of $S_{0,0}$, $\widetilde{\Sigma}_0^{\mathcal{N}}:=\embed(\iniu+\iniubar=0)\cap\mathcal{N}$ and 
	\begin{align} \label{E:introINHOMOGENEOUSTERMFORSPLITFIRSTORDEREULERINTOEXTENDEDDATAPLUSREMAINDER}
		F
		& := 
		\begin{cases}
			0, & \mbox{in } \extendedregion_{(\text{in})},
			\\
			A^{\alpha}[\vvariables_{(\text{ext})}]
			\p_{\alpha}
			\vvariables_{(\text{ext})},
			&  \mbox{in } \extendedregion_{(\text{out})}.
		\end{cases}
	\end{align}
	In the above,
	\begin{subequations}
		\begin{align}
			\extendedregion_{(\text{in})} 
			& := \embed\left(\lbrace \iniu+\iniubar \geq 0 \rbrace \cap \lbrace  \iniu \geq 0 \rbrace \cap 
			\lbrace \iniubar \geq 0 \rbrace \right) \cap \mathcal{N},
			\label{E:introINTERIORRENDALLREGION} \\
			\extendedregion_{(\text{out})} 
			& := \embed\left( \lbrace \iniu+\iniubar \geq 0 \rbrace 
			\cap 
			(\lbrace  \iniu < 0 \rbrace 
			\cup 
			\lbrace \iniubar < 0 \rbrace) 
			\right) \cap \mathcal{N},
			\label{E:introEXTERIORRENDALLREGION} \\
			\extendedregion 
			& := 
			\extendedregion_{(\text{in})}  \cup \extendedregion_{(\text{out})}  
			=
			\embed\left( \lbrace \mathring{\uptau} \geq 0 \rbrace \right) \cap \mathcal{N}.
			\label{E:introTOTALRENDALLREGION}
		\end{align}
	\end{subequations}
	By construction, the source-term $F$ in \eqref{E:introINHOMOGENEOUSTERMFORSPLITFIRSTORDEREULERINTOEXTENDEDDATAPLUSREMAINDER} is $C^{\infty}$.
	Then, by local existence and uniqueness theory for standard Cauchy problems, we have $\tilde{\vvariables}=0$ in $\extendedregion_{(\text{out})}$. Moreover, $\vvariables := \vvariables_{(\text{ext})}+\tilde{\vvariables}$ solves the compressible
	Euler equations \eqref{EQ:introEuler} in a future neighborhood of $S_{0,0}$, $\region_{(\text{in})}\subset\extendedregion_{(\text{in})}$, and takes on the given initial data
	on $(\inihypersur \bigcup \inihypersurbar) \cap \region_{(\text{in})}$.
	
	The above argument relies on the Whitney extension theorem (see Lemma \ref{L:Whitneyextension}), which is not applicable with $H^k$ Sobolev-class data. 
	Compared to Rendall's wave equations, for compressible fluids, it is much more difficult to set up the machinery needed to apply Whitney's extension theorem. The reason is that $3D$ compressible fluid flow does not ``cleanly'' separate into ingoing and outgoing sound waves, so it is not as straightforward to solve for transversal derivatives, i.e., it is much harder compared to the model problem from Sect.\ref{SSS:AModel}.  
	To achieve this, we combine the new formulation of the flow from Sect.\,\ref{SSS:NEWFORMULATIONSCHEMATIC} with the standard first-order formulation
	and make careful use of the characteristic geometry;
	see Sect.\,\ref{S:FLUIDCHARDATAFURTHERASPECTS} and Sect.\,\ref{SS:WHITNEY} for the details.
	We stress that for the argument to work, the function
	$F$ in \eqref{E:introINHOMOGENEOUSTERMFORSPLITFIRSTORDEREULERINTOEXTENDEDDATAPLUSREMAINDER} needs to be $C^{\infty}$, which requires that the prescribed initial fluid data are smooth and satisfy Euler equations at all orders. Furthermore, significant difficulties arise in smoothing the solution due to the change of characteristics, and compatibility issues occur at $S_{0,0}$. Therefore, in this article, we only consider $C^{\infty}$ data and solutions. Nevertheless, our a priori estimates for Sobolev norms close at a finite Sobolev regularity level.
	
\subsubsection{Causal geometry}
	To extend the region of existence and uniqueness from a neighborhood of the initial sphere to a neighborhood of the initial characteristic hypersurfaces for ``CharIVP,'' we use the new formulation of the flow discussed above in Section \ref{SSS:NEWFORMULATIONSCHEMATIC}. 
	In Section \ref{SS:CONSTRUCTIONOFACOUSTICEIKONALFUNCTIONSANDGEOMETRICANGULARCOORDINATESINRENDALLREGION}, we construct a pair of acoustical eikonal functions with transversal level sets, denoted by $u$ and $\ubar$, i.e., solutions to:
	\begin{align}
		(\gfour^{-1})^{\alpha\beta}\p_\alpha u\p_\beta u=&0,&
		(\gfour^{-1})^{\alpha\beta}\p_\alpha \ubar\p_\beta \ubar=&0,
	\end{align} 
	where $\gfour$ and $\gfour^{-1}$ are the acoustical metric and its inverse defined in (\ref{DE:introacousticalmetric}). The levels of $u$ and $\ubar$ are characteristic hypersurfaces for the fluid evolution equations. In particular, we construct $u$ and $\ubar$ such that $\{u=0\}=\inihypersur$, and $\{\ubar=0\}=\inihypersurbar$. 
	We prove existence and uniqueness in the shaded region $\{0\leq u\leq\updelta,0\leq\ubar\leq \ubarfinite\}\bigcup\{0\leq \ubar\leq\updelta,0\leq u\leq \ufinite\}$ in Figure \ref{fig3}, compared to $\{0\leq u\leq\updelta,0\leq\ubar\leq \updelta\}$ in Section \ref{SSS:neighborhood}, based on a priori estimates for the fluid variables and the acoustic geometry
	(i.e., for the eikonal functions). 
	We note that the region of existence is similar to the case in \cite{GRCharIVP} for Einstein vacuum equations. We also derive continuation criteria, i.e., conditions on the solution and geometry guaranteeing that one can extend the fluid solution and acoustic geometry for a short time. 
	
	Our analysis involves three types of regions, as illustrated in Figure \ref{Figure4}. $\region_1$ denotes the ``triangle'' region bounded below by transversally intersecting initial null hypersurfaces. $\region_2$ denotes the ``parallelogram'' region bounded below by the top (spacelike) boundary of $\region_{1}$, with null lateral boundaries. $\region_3$ denotes the ``trapezoid'' region bounded below by the top (spacelike) boundary of $\region_2$, with null lateral boundaries. We derive a priori estimates on all three of these kinds of regions in order to conclude our main existence result.

	\subsubsection{Space-time transport-elliptic estimates}\label{SSS:elliptic}
	The transport-div-curl part of the compressible Euler flow from Sect.\,\ref{SSS:NEWFORMULATIONSCHEMATIC}
	lies outside of the scope of the traditional vectorfield method (developed for wave equations) on spatially compact regions due to the difficult boundary integrals. Despite this difficulty, in \cite{IntegralIdentities}, the authors provided remarkable integral identities, which 
effectively allow one to extend the vectorfield method to the compressible Euler equations by yielding	
transport-elliptic estimates for the ``transport-part'' of the solution on arbitrary globally hyperbolic domains.
 In particular, in \cite{IntegralIdentities}, $\hypersur,\hypersurbar$ are allowed to be acoustically spacelike or acoustically null. 
 Specifically, the methods of \cite{IntegralIdentities} yield the following
integral identity\footnote{We present the schematic version of the identity in $\region_1$; we refer readers to Lemma \ref{L:integralidentitiesforS} for the full details of the identities.} for $V\in\{\vortrenormalized,\gradEnt\}$ (see Section \ref{SS:morefluid} for $\vortrenormalized,\gradEnt$) in $\region_1$:
	\begin{align}\label{ES:introIntegral}
		\int_{\region}\abs{\pfour V}^2\diff\gvol+\int_{S_{u,0}}\abs{\angV}^2+\int_{S_{0,\ubar}}\abs{\angV}^2=&~\int_{S_{0,0}}\abs{\angV}^2+\int_\region \text{good terms}\\
		\notag& +\int_{\hypersur}\text{tangential terms}+\int_{\hypersurbar}\text{tangential terms},
	\end{align}
	where $\angV$ is the $\gfour$-orthogonal projection of $V$ onto the 2-sphere $\suubar$. The ``good terms'' in (\ref{ES:introIntegral}) consist of $\modivort,\DivGradEnt$, which can be controlled by transport equations (\ref{EQ:intro3.1})-(\ref{EQ:intro3.2}). We emphasize that the boundary terms on $\hypersur$ and $\hypersurbar$ consist of $\vortrenormalized,\gradEnt$, and tangential derivatives of fluid variables and null frames. This tangential structure is crucial since ``wave-part'' energies on the null hypersurfaces are semi-coercive, that is, one can only control tangential derivatives on null hypersurfaces by using the wave equation (\ref{EQ:intro1}). The estimates on $\region_2$ and $\region_3$ follow a similar approach.
	
 Combining the integral identity (\ref{ES:introIntegral}) with relatively standard transport estimates  
	based on (\ref{EQ:intro3.1})-(\ref{EQ:intro3.2}), 
	we gain one derivative for the vorticity and entropy gradient compared to standard estimates:
	\begin{align}
		\int_{\region}\abs{\pfour\vortrenormalized,\pfour\gradEnt}^2\diff\gvol\lesssim\int_{\region} \abs{\pfour(v^i,\logDensity,\Ent)}^2\diff\gvol.
	\end{align}
	This is a gain of one derivative compared to viewing $\pfour\vortrenormalized,\pfour\gradEnt \sim \pfour\p(v^i,\Ent)$,
	which follows by definition.
	
	\subsubsection{Energy estimates}
	\label{SSS:introEnergyEstimates}
	The local construction near $S_{0,0}$ described in Sect.\,\ref{SSS:neighborhood}
	provides the starting point for our
	analysis. The next step is to construct the eikonal functions and geometric coordinates 
	near $S_{0,0}$. This step is relatively standard and we will not discuss it here; see Sect.\,\ref{SS:CONSTRUCTIONOFACOUSTICEIKONALFUNCTIONSANDGEOMETRICANGULARCOORDINATESINRENDALLREGION}
	for the details.
	To extend the solution to a larger region,
	we establish a series of a priori estimates that simultaneously control the fluid variables and the acoustic geometry (the eikonal functions). 
	To derive these estimates, we combine the
	vectorfield method for the wave equations \eqref{EQ:intro1}
	with standard transport estimates for the transport equations \eqref{EQ:intro2} and
	spacetime transport-elliptic estimates for the equations \eqref{EQ:intro3}, as described in
	Section~\ref{SSS:elliptic}.
	
	To derive energy estimates for the wave equations, we use a future-directed timelike multiplier
	$\Timelike$, which is adapted to the acoustic geometry; see Sect.\,\ref{SS:energies}. 
	The associated energies control all first
	derivatives of the wave variables on the spacelike hypersurfaces
	$\Stau$, where $\uptau=u+\ubar$, while the null fluxes control only
	derivatives tangent to the corresponding characteristic hypersurfaces.
	The principal difficulty is that we have to close these wave equation energy estimates together with those for
	the transport variables. The wave equations contain the modified variables
	$\modivort$ and $\DivGradEnt$, whose transport equations enjoy a gain\footnote{Estimating these derivatives directly from their
		definitions would require an additional derivative of the fluid variables.} in regularity. The 
	transport-elliptic estimates described in
	Section~\ref{SSS:elliptic} recover spacetime $L^2$ control
	of vorticity $\vortrenormalized$ and entropy gradient $\gradEnt$.
	Crucially, the derivative terms on the null boundaries have the tangential
	structure needed to be estimated by the available fluxes. This closes the
	coupled estimates without an additional loss of derivatives in the
	evolution.
	
	We obtain higher-order estimates by commuting the equations of Sect.\,\ref{SSS:NEWFORMULATIONSCHEMATIC}
	and the eikonal equations
	with the Cartesian spacetime
	derivatives $\pfour$, and carrying out an induction argument on the number of derivatives. At order
	$k$, the energies control $k$ derivatives of
	$\vvariables=(\logDensity,v^1,v^2,v^3,\Ent)$ on spacelike slices and
	$k-1$ derivatives of $\modivort$ and $\DivGradEnt$. The accompanying
	transport-elliptic estimates control $k$ derivatives of
	$(\vortrenormalized,\gradEnt)$ in spacetime $L^2$. Thus, the gain of one
	derivative for the latter variables is realized in the spacetime estimates.
	Applying these arguments successively in the regions
	$\region_1$, $\region_2$, and $\region_3$, and using Young's and
	Gr\"onwall's inequalities, yields schematically
	\begin{align}\label{ES:introgeneral}
		\tilde{\bEnergy}_{(k)}
		+\tilde{\Energy}_{(k)}
		+\tilde{\uEnergy}_{(k)}
		+\sum_{r=1}^{3}\tilde{\KEnergy}_{(k)}
		\lesssim
		\iniLtwo\exp\!\left(\RBAinfty(\updelta_u+\updelta_{\ubar})\right),
		\qquad 1\leq k\leq\frac{N-1}{2},
	\end{align}
	where $\tilde{\bEnergy}_{(k)}$ is the energy below-k-th-order for the fluid variables on $\Stau$, $\tilde{\Energy}_{(k)},\tilde{\uEnergy}_{(k)}$ are the fluid energy fluxes on $\hypersur$ and $\hypersurbar$, respectively, and $\tilde{\KEnergy}_{(k)}$ is the spacetime energy for transport variables; see Section \ref{SS:energies} for their precise definitions.
	In \eqref{ES:introgeneral}, $\iniLtwo$ measures the initial
	energy size, and $\updelta_u=u_1-u_0$ and
	$\updelta_{\ubar}=\ubar_1-\ubar_0$ are the null coordinate lengths of the
	region. The precise estimates are stated in
	Theorem~\ref{TH:UPTOTOPORDERENERGYESTIMATES}.
	
	An essential feature of this argument is that only one of these two lengths must be
	sufficiently small. The estimates therefore apply in a thin neighborhood
	of either initial characteristic hypersurface, with its prescribed
	longitudinal extent kept fixed. Combining the energy estimates with
	the Sobolev embedding on $S_{u,\ubar}$, we can integrate in the  ``short direction''
	to improve the fluid bootstrap assumptions. The eikonal and geometric
	transport equations then control the null lapses and Ricci coefficients
	and preserve the ellipticity of the induced sphere metrics. These
	arguments close the bootstrap for the fluid and acoustic geometry; see
	Section~\ref{S:APRIORIESTIMATESUPTOTOPORDER} for the details.
	
	\subsubsection{The main result: Existence and uniqueness in a neighborhood of the initial hypersurfaces}\label{SSS:neighborhood2}
	Given the a priori estimates discussed in Section \ref{SSS:introEnergyEstimates}, 
	it is not too difficult to extend the region of existence and uniqueness from a neighborhood of the initial sphere (established in Section \ref{SSS:neighborhood}) to a neighborhood of the initial hypersurfaces. We refer the reader to Figs.\ref{fig3} and \ref{fig1} for a comparison of these two regions. In this section, we give an overview of the extension argument; see Section \ref{S:Lastslice} for the details.
	
	We discuss only how to extend the solution to a neighborhood of $\inihypersur$; the solution can be extended to a neighborhood of $\inihypersurbar$ using nearly identical arguments. More precisely, we will explain why the solution can be extended to the $\inihypersur$-neighborhood 
	$[0,\updelta]\times[0,\ubarfinite]\times\Stwo$ on which we have already shown that a priori estimates hold. 
	We refer readers to Fig.\,\ref{Figure3} for an illustration of the regions that play a role in the argument.
	Suppose for the sake of deriving a contradiction
	that the maximal region of existence and uniqueness is a smaller region $[0,\updelta') \times[0,\ubar')\times\Stwo$, 
	where $0<\updelta'<\updelta$, $0<\ubar'<\ubarfinite$. Our a priori estimates can be extended to show the following ``soft fact:'' 
	the solution extends as a $C^{\infty}$ solution to the closure $\widetilde{\mathcal{M}}_0:=[0,\updelta']\times[0,\ubar']\times\Stwo$;
	see Sect.\,\ref{SS:SmoothnessoftheSolutions} and the proof of Theorem~\ref{TH:localwellposedness}.
	Hence, using ``Rendall-type'' arguments from Section \ref{SSS:neighborhood}, 
	we can deduce that there exists an ``extra'' region of existence and uniqueness $\widetilde{\mathcal{M}}_1:=[\updelta',\updelta'+\varepsilon]\times[0,\varepsilon]\times\Stwo$. Notice that we obtain a spacelike slice $\Sigma_{\timefunction'}$, where $\timefunction'=\updelta'+\varepsilon$, on which we can apply the well-posedness theory for the standard Cauchy problem with initial data on the spacelike slice $\Sigma_{\timefunction'}$. 
This standard theory, when combined with our a priori estimates, 
allows us to extend the region of existence and uniqueness to $[0,\updelta'+\varepsilon]\times[0,\ubar']\times\Stwo$,
which in Fig.\,\ref{Figure3} is $\widetilde{\mathcal{M}}_0 \cup \widetilde{\mathcal{M}}_1 \cup \widetilde{\mathcal{M}}_2$.
Then, using the same strategy, we can further extend the region of existence to $[0,\updelta'+\varepsilon]\times[0,\ubar'+\varepsilon]\times\Stwo$,
which in Fig.\,\ref{Figure3} is $\widetilde{\mathcal{M}}_0 \cup \widetilde{\mathcal{M}}_1 \cup \widetilde{\mathcal{M}}_2 \cup \widetilde{\mathcal{M}}_3$.
This contradicts the maximality of $\widetilde{\mathcal{M}}_0$
and shows the region of existence can be extended to the full region where the priori estimates hold. We refer the reader to 
	Theorem~\ref{TH:localwellposedness} for the detailed statement and proof.
	Our continuation result, namely Theorem~\ref{TH:continuationcriteria}, can be proved using similar ideas.
	
	\subsection{Paper Outline}
	The remainder of the paper is organized as follows:
	\begin{itemize}
		\item In Section \ref{S:Notations}, we introduce notation that we use throughout the paper. We then present a standard first-order formulation of the compressible Euler equations.
		\item In Section \ref{S:Geometricformulation}, we state the geometric formulation of the compressible Euler equations introduced in \cite{ANewFormulationofthe3DCompressibleEulerEquationswithDynamicEntropy}. The equations are provided as Theorem~\ref{T:geometricformulation}.
		\item In Section \ref{S:ACOUSTICGEOMETRY}, we introduce various
		geometric constructions tied to characteristic surfaces.
		\item In Section \ref{S:CID}, we state our assumptions on the characteristic initial data and the embeddings of the initial characteristic hypersurfaces into Cartesian spacetime. 
		\item In Section \ref{S:FLUIDCHARDATAFURTHERASPECTS}, we solve for transversal derivatives of the initial data along the initial null hypersurfaces such that the initial data that we constructed satisfy the first-order Euler equations at all orders. We also identify corner compatibility conditions that the data must satisfy along
			the intersection of the initial null hypersurfaces.
		\item In Section \ref{S:localexistence}, we establish the local existence and uniqueness in a neighborhood of the initial sphere $S_{0,0}$. 
		We then construct acoustical eikonal functions in a perhaps smaller neighborhood.
		\item In Section \ref{S:geometry}, we introduce acoustic geometry based on the acoustical eikonal functions constructed 
		in Section \ref{S:localexistence}. We define additional geometric quantities, including various regions, the connection coefficients of the null frame, projections, and volume forms.
		\item In Section \ref{S:REGULARITYOFDATA}, 
			we derive estimates for various norms of the fluid and acoustic geometry along the initial null hypersurfaces.
		\item In Section \ref{S:SectionAssumptions}, 
		we state the Sobolev norms on the initial data that we will use in the proof of our main a priori estimates.
			We also state various bootstrap assumptions that we will use in deriving a priori estimates.
		\item  In Section \ref{S:ENERGYANDENERGYIDENTITIES}, we define various energies on constant-geometric-time hypersurfaces $\Stau$, null hypersurfaces $\hypersur,\hypersurbar$, and the space-time region $\region$. Then, in Section \ref{SS:energyinequalities}, 
		we present preliminary energy identities and estimates
		for the equations of Theorem~\ref{T:geometricformulation}.
		Finally, in Section \ref{SS:elliptic}, using the approach of \cite{IntegralIdentities}, we derive 
		spacetime elliptic estimates for the specific vorticity and entropy gradient.
		These allow us to gain one derivative worth of Sobolev regularity for the specific vorticity and entropy gradient compared to standard estimates.
		\item In Section \ref{S:BASELEVELENERGYESTIMATES}, we derive base-level a priori energy estimates for solutions to the equations of 
		Theorem~\ref{T:geometricformulation},
			i.e., energy estimates for the un-commuted equations.
		\item In Section \ref{S:APRIORIESTIMATESUPTOTOPORDER}, in a neighborhood of the initial hypersurfaces,
			we derive higher-order a priori estimates for solutions to the equations of Theorem~\ref{T:geometricformulation}.
		Moreover, under a smallness assumption on the width of the existence-region, 
		we improve the bootstrap assumptions (stated in Section \ref{S:SectionAssumptions}). 
		\item  In Section \ref{S:Lastslice}, 
		we prove our main theorem on local existence and uniqueness in a neighborhood of the initial hypersurfaces. In addition, 
		in Section \ref{SS:continuationcriterion}, we provide continuation criteria.
	\end{itemize}

	\section{Notation and First-order Formulation of the Compressible Euler Equations}\label{S:Notations}
	\subsection{Notation}\label{Notation}

	We use the following notation throughout the paper:
	\begin{itemize}
		\item Greek ``space-time" indices take on the values $0,1,2,3$, while the Latin ``spatial" indices take on the values $1,2,3$. 
		\item We use the Einstein summation convention.
		\item $\antisymmetic_{\alpha\beta\gamma\updelta}$ and $\antisymmetic_{ijk}$ denote the fully antisymmetric symbols normalized by $\antisymmetic_{0123}=1$ and $\antisymmetic_{123}=1$.
		\item We raise and lower indices for tensorfields using the acoustical metric\footnote{See Definition \ref{D:acousticmetric} for the definition of $\gfour$ and its inverse.} $\gfour$ and its inverse. For example, $V^\alpha=(\gfour^{-1})^{\alpha\beta}V_\beta$, $V_\beta=\gfour_{\alpha\beta}V^\alpha$.
		\item We denote $\St:=\{(t^\prime,x^1,x^2,x^3)\in\mathbb{R}^{1+3}|t\equiv t^\prime\}$ as the standard constant-time slice.
		\item $\linsmoothfunction[A](B)$ denotes any scalar-valued function that is linear in the components of $B$ with coefficients that are a function of the components of $A$.
		\item $\quadsmoothfunction[A](B,C)$ denotes any scalar-valued function that is quadratic in the components of $B$ and $C$ with coefficients that are a function of the components of $A$.
		\item We denote the schematic spatial partial derivatives by $\p$ and the schematic space-time partial derivatives by $\pfour$.
		\item We denote $\pfour^{\leq k}\phi$ any linear combination of any components of $\pfour^j\phi$ with $j\leq k$.
		\item We denote $A\lesssim B$ as $A\leq C\cdot B$ for some universal constant $C$.
		\item We denote $\gensmoothfunction$ the smooth function of its argument. We allow $\gensmoothfunction$ to vary from line to line.
	\end{itemize}
	We refer readers to Appendix \ref{AP:notations} for more notations that we use throughout the article.
	\subsection{Fluid Variables, Equation of State and Speed of Sound}\label{D:speedofsound}
	We denote $v^i$ with $i=1,2,3$ the fluid velocity, $\density$ the fluid density, $\Ent$ the fluid entropy per particle.
	
	\begin{definition}[Partial derivatives with respect to $\density$ and $\Ent$]\label{DE:SoundSpeed}
		If $Q$ is a quantity that can be expressed as a function of $(\density,\Ent)$, then
		\begin{align}
			Q_{;\density}&=\left.Q_{;\density}(\density,\Ent):=\frac{\p Q}{\p\density}\right\vert_\Ent,\\
			Q_{;\Ent}&=\left.Q_{;\Ent}(\density,\Ent):=\frac{\p Q}{\p\Ent}\right\vert_\density.
		\end{align}
	\end{definition}
	
	We assume an equation of state of the form $p=p(\density,\Ent)$. The speed of sound $\speed$ is defined as follows:
	\begin{equation}
		\speed:=\sqrt{\left.\frac{\p p}{\p\density}\right\vert_\Ent}.
	\end{equation}
	In this article, we assume that $\speed\geq\speed_0>0$ for some constant $\speed_0$. 
	
	\subsection{More Fluid Variables}\label{SS:morefluid}
	\begin{definition}[Divergence and Curl of Vectorfields] For a $\St$-tangent vectorfield $V$, we define its divergence $\dive V$ and curl $\curl V$ as follows:
		\begin{align}
			\dive V&:=\p_aV^a,&
			(\curl V)^i&:=\antisymmetic_{iab}\p_aV^b.
		\end{align}
	\end{definition}
	\begin{definition}[Vorticity] We define the vorticity $\vort:\mathbb{R}^{1+3}\rightarrow\mathbb{R}^3$ as follows:
		\begin{align}
			\vort^i:=(\curl v)^i.
		\end{align}		
	\end{definition}
	
	We find it convenient to work with the logarithmic density $\logDensity$ and the specific vorticity $\vortrenormalized$. 
	To define these, we fix a constant background density $\backgroundDensity$ such that $\backgroundDensity>0$.
	\begin{definition}[Logarithmic density and specific vorticity]
		\label{D:LOGDENSITYANDSPECIFICFORTICITY}
		We define the logarithmic density $\logDensity$ and the specific vorticity vectorfield $\vortrenormalized$ as follows:
		\begin{align}\label{DE:vortrenomal}
			\logDensity&:=\ln\left(\frac{\density}{\backgroundDensity}\right),&
			\vortrenormalized&:=e^{-\logDensity}\vort.
		\end{align}
	\end{definition}
	
	\begin{definition}[Entropy gradient]\label{D:entgradient}
		We define the entropy gradient vectorfield $\gradEnt$ as follows:
			\begin{align}\label{DE:ent}
				\gradEnt^i&=\updelta^{ia}\p_a\Ent.
		\end{align}
	\end{definition}
	
	\subsection{First-order Formulation of the Equations}\label{SS:firstorderequations} We denote the material derivative vectorfield $\materialderivative:=\p_t+v^a\p_a$. Relative to the standard Cartesian coordinates $(t,x^1,x^2,x^3)$ on $\mathbb{R}^{1+3}$, the compressible Euler equations can be expressed in the following standard first-order form:
	\begin{subequations}\label{E:firstorderformulation}
		\begin{align}
			\label{E:Brho}\materialderivative\logDensity&=-\dive v,
				\\
			\label{E:Bv}\materialderivative v^i&=-\speed^2\updelta^{ia}\p_a\logDensity-e^{-\logDensity}\frac{\pressure_{;\Ent}}{\backgroundDensity}\updelta^{ia}\p_a\Ent,
				\\
			\label{E:Bs}\materialderivative s&=0,
		\end{align}
	\end{subequations}
	where $\delta^{ab}$ is the Kronecker delta.
	
	We define the array of fluid variables as follows:
	\begin{align}\label{DE:vvariables}
		\vvariables:=(\logDensity,v^1,v^2,v^3,\Ent)
	\end{align}
	
	The above system can be expressed in the following form:
	\begin{align}\label{FOEuler}
		A^\alpha(\vvariables)\p_\alpha\vvariables=0,
	\end{align}
	where, for each $\alpha=0,1,2,3$, $A^\alpha(\vvariables)$ is a solution-dependent matrix whose precise form we do not record here.
	
	\section{The Geometric Formulation of the Compressible Euler Equations}\label{S:Geometricformulation}
	\begin{definition}\cite[Definition 1.3. Modified fluid variables]{ANewFormulationofthe3DCompressibleEulerEquationswithDynamicEntropy}.\label{D:modifiedfluidvariables} We define the modified fluid variables as follows:
		\begin{subequations}\label{DE:CD}
			\begin{align}
				\label{Def:C}\modivort^i:=&e^{-\logDensity}(\curl\vortrenormalized)^i+e^{-3\logDensity}\speed^{-2}\frac{\pressure_{;\Ent}}{\backgroundDensity}\gradEnt^a\p_av^i-e^{-3\logDensity}\speed^{-2}\frac{\pressure_{;\Ent}}{\backgroundDensity}(\p_av^a)\gradEnt^i,\\
				\label{Def:D}\DivGradEnt:=&e^{-2\logDensity}\dive\gradEnt-e^{-2\logDensity}\gradEnt^a\p_a\logDensity.
			\end{align}
		\end{subequations}
	\end{definition}
	
	\begin{definition}[The acoustical metric and its inverse]\label{D:acousticmetric}	
		We define the acoustical metric $\gfour$ and its inverse\footnote{For convenience, we write $\gfour^{\alpha\beta}$ instead of $(\gfour^{-1})^{\alpha\beta}$ in this paper.} $\gfour^{-1}$ relative to the Cartesian coordinates as follows:
		\begin{subequations}
			\begin{align}
				\gfour&:=-\diff t \otimes \diff t + \speed^{-2}\sum\limits^3_{a=1}(\diff x^a-v^a\diff t) \otimes (\diff x^a-v^a\diff t),
				\label{E:ACOUSTICALMETRIC} \\
				\gfour^{-1}&:=-\materialderivative \otimes \materialderivative+\speed^2\sum\limits^3_{a=1}\p_a \otimes\p_a.
				\label{E:INVERSEACOUSTICALMETRIC}
			\end{align}					
		\end{subequations}
	\end{definition}
	\begin{definition}[Differential operators defined by $\gfour$]\label{D:LeviCivita}
		$\Dfour$ denotes the Levi-Civita connection
		of $\gfour$ and $\boxg:=(\gfour^{-1})^{\alpha\beta}\Dfour_\alpha\Dfour_\beta$ denotes the corresponding covariant wave operator. In arbitrary coordinates,\footnote{The operator $\boxg$ is coordinate invariant.} 
		for a scalar function $\varphi$, $\boxg\varphi=\frac{1}{\sqrt{-\abs{\gfour}}}\p_\alpha\left(\sqrt{-\abs{\gfour}}(\gfour^{-1})^{\alpha\beta}\p_\beta\varphi\right)$, 
		where $\abs{\gfour}$ is the determinant of $\gfour$.
	\end{definition}
	\begin{definition}[Arrays of variables] For convenience, in presenting the formulations and analysis, we define the following arrays of solution variables:\label{D:cartesiancomponentsofvariables}
		\begin{align}
			\vec{v}:&=(v^1,v^2,v^3),& \vvort:&=(\vortrenormalized^1,\vortrenormalized^2,\vortrenormalized^3),& \vgradientEnt:&=(\gradEnt^1,\gradEnt^2,\gradEnt^3),& \vectvort:=(\modivort^1,\modivort^2,\modivort^3).
		\end{align}
	\end{definition}
	\begin{definition}[Definition of $\dfour{\alpha}{\beta}$]\label{DE:dfour}
		Recall (\ref{DE:vvariables}) for the definition of $\vvariables$. We define $\dfour{\alpha}{\beta}$ as follows:
		\begin{align}\label{D:G}
			\dfour{\alpha}{\beta}:=\frac{\pfour\gfour_{\alpha\beta}(\vvariables)}{\pfour\vvariables}.
		\end{align}
	\end{definition}
	\begin{definition}[Christoffel symbols]\label{D:Christoffelsymbols}
		We define the Christoffel symbols $\Chfour_{\alpha\kappa\lambda}$, $\Chfour^\beta_{\kappa\lambda}$ and $\Chfour_\alpha$ of the acoustical metric $\gfour$ as follows:
		\begin{subequations}
			\begin{align}
				\Chfour_{\alpha\kappa\lambda}:=&\frac{1}{2}\left(\p_\kappa\gfour_{\alpha\lambda}+\p_\lambda\gfour_{\alpha\kappa}-\p_\alpha\gfour_{\kappa\lambda}\right)=\frac{1}{2}(G_{\alpha\lambda}\p_\kappa\vvariables+G_{\alpha\kappa}\p_\lambda\vvariables-G_{\kappa\lambda}\p_\alpha\vvariables),\\
				\Chfour^\beta_{\kappa\lambda}:=&(\gfour^{-1})^{\alpha\beta}\Chfour_{\alpha\kappa\lambda},\\
				\Chfour^\beta:=&\left(\gfour^{-1}\right)^{\kappa\lambda}\Chfour^\beta_{\kappa\lambda}.
			\end{align}
		\end{subequations}
	\end{definition}
	
	\begin{theorem}\cite[Theorem 1. The geometric formulation of the compressible Euler equations]{ANewFormulationofthe3DCompressibleEulerEquationswithDynamicEntropy}.\label{T:geometricformulation}
		If $(v^1,v^2,v^3,\logDensity,\Ent)$ solves the first-order formulation of the compressible Euler equations (\ref{E:firstorderformulation}), the following formulation holds for $(v^i,\logDensity,\Ent)$ and $(\vortrenormalized^i,\gradEnt^i,\modivort^i,\DivGradEnt)$:\\
		
		\textbf{Wave equations}
		\begin{subequations}\label{E:WaveEQ}
			With $\boxg$ denoting the covariant wave operator\footnote{On LHS~\eqref{E:boxgv}, the Cartesian components $v^i$ are treated as scalar functions when being
			differentiated by $\boxg$.} 
			from Def.\,\ref{D:LeviCivita}, we have:
			\begin{align}
				\label{E:boxgv}\boxg v^i&=-\speed^2e^{2\logDensity}\modivort^i+\mathfrak{Q}_{(v)}^i+\mathfrak{L}_{(v)}^i,\\
				\label{E:boxgdensity}\boxg\logDensity&
				=
				-e^{\logDensity}\frac{\pressure_{;\Ent}}{\backgroundDensity}\DivGradEnt+\mathfrak{Q}_{(\logDensity)}+\mathfrak{L}_{(\logDensity)},\\
				\label{E:boxgent}\boxg\Ent&=\speed^2e^{2\logDensity}\DivGradEnt+\mathfrak{L}_{(\Ent)}.
			\end{align}
		\end{subequations}
		
		\textbf{Transport equations}
		\begin{subequations}
			\begin{align}
				\materialderivative\vortrenormalized^i&=\mathfrak{L}_{(\vortrenormalized)}^i,\\
				\materialderivative\Ent&=0,\\
				\label{E:BS}\materialderivative\gradEnt^i&=\mathfrak{L}_{(\gradEnt)}^i,
			\end{align}
		\end{subequations}
		
		\textbf{Transport-div-curl system}
		\begin{subequations}\label{E:divcurlsystem}
			\begin{align}
				\label{E:divvort}\dive\vortrenormalized=&\mathfrak{L}_{(\dive\vortrenormalized)},\\
				\label{E:BC}\materialderivative\modivort^i=&-2\delta_{jk}\antisymmetic_{iab}e^{-\logDensity}(\p_av^j)\p_b\vortrenormalized^k+\antisymmetic_{ajk}e^{-\logDensity}(\p_av^j)\p_j\vortrenormalized^k\\
				\notag&+e^{-3\logDensity}\speed^{-2}\frac{\pressure_{;\Ent}}{\backgroundDensity}\left\{(\materialderivative\gradEnt^a)\p_av^i-(\materialderivative v^i)\p_a\gradEnt^a\right\}\\
				\notag&+e^{-3\logDensity}\speed^{-2}\frac{\pressure_{;\Ent}}{\backgroundDensity}\left\{(\materialderivative v^a)\p_a\gradEnt^i-(\p_a v^a)\materialderivative\gradEnt^i\right\}\\
				\notag&+\mathfrak{Q}^i_{(\modivort)}+\mathfrak{L}^i_{(\modivort)},\\
				\label{E:BD}\materialderivative\DivGradEnt=&2e^{-2\logDensity}\left\{(\p_av^a)\p_b\gradEnt^b-(\p_av^b)\p_b\gradEnt^a\right\}+e^{-\logDensity}\delta_{ab}(\curl\vortrenormalized)^a\gradEnt^b+\mathfrak{Q}_{(\DivGradEnt)},\\
				\label{E:curlS}\curl\gradEnt^i=&0.
			\end{align}
		\end{subequations}
		Above, $\mathfrak{Q}_{(v)}^i, \mathfrak{Q}_{(\logDensity)}, \mathfrak{Q}^i_{(\modivort)}$ and $\mathfrak{Q}_{(\DivGradEnt)}$ are the $\gfour$-null forms defined as follows:
		\begin{subequations}\label{E:Qequations}
			\begin{align}
				\mathfrak{Q}_{(v)}^i:=&-(1+\speed^{-1}\speed_{;\logDensity})(\gfour^{-1})^{\alpha\beta}(\p_\alpha\logDensity)\p_\beta v^i,\\
				\mathfrak{Q}_{(\logDensity)}:=&-3\speed^{-1}\speed_{;\logDensity}(\gfour^{-1})^{\alpha\beta}(\p_\alpha\logDensity)\p_\beta \logDensity+\left\{(\p_av^a)\p_bv^b-(\p_av^b)\p_bv^a\right\},\\
				\mathfrak{Q}^i_{(\modivort)}:=& e^{-3\logDensity}\speed^{-2}\frac{\pressure_{;\Ent}}{\backgroundDensity}\gradEnt^i\left\{(\p_av^b)\p_bv^a-(\p_av^a)\p_bv^b\right\}\\
				\notag&+e^{-3\logDensity}\speed^{-2}\frac{\pressure_{;\Ent}}{\backgroundDensity}\gradEnt^b\left\{(\p_av^a)\p_bv^i-(\p_av^i)\p_bv^a\right\}\\
				\notag&+2e^{-3\logDensity}\speed^{-2}\frac{\pressure_{;\Ent}}{\backgroundDensity}\gradEnt^a\left\{(\p_a\logDensity)\materialderivative v^i-(\p_av^i)\materialderivative \logDensity\right\}\\
				\notag&+2e^{-3\logDensity}\speed^{-3}\speed_{;\logDensity}\frac{\pressure_{;\Ent}}{\backgroundDensity}\gradEnt^a\left\{(\p_a\logDensity)\materialderivative v^i-(\p_av^i)\materialderivative \logDensity\right\}\\
				\notag&+e^{-3\logDensity}\speed^{-2}\frac{\pressure_{;\Ent;\logDensity}}{\backgroundDensity}\gradEnt^a\left\{(\p_av^i)\materialderivative \logDensity-(\p_a\logDensity)\materialderivative v^i\right\}\\
				\notag&+e^{-3\logDensity}\speed^{-2}\frac{\pressure_{;\Ent;\logDensity}}{\backgroundDensity}\gradEnt^i\left\{(\materialderivative v^a)\p_a  \logDensity-(\materialderivative\logDensity)\p_av^a\right\}\\			\notag&+2e^{-3\logDensity}\speed^{-2}\frac{\pressure_{;\Ent}}{\backgroundDensity}\gradEnt^i\left\{(\materialderivative\logDensity)\p_av^a-(\materialderivative v^a)\p_a\logDensity\right\}\\
				\notag&+2e^{-3\logDensity}\speed^{-3}\speed_{;\logDensity}\frac{\pressure_{;\Ent}}{\backgroundDensity}\gradEnt^i\left\{(\materialderivative\logDensity)\p_av^a-(\materialderivative v^a)\p_a\logDensity\right\},\\
				\mathfrak{Q}_{(\DivGradEnt)}:=&2e^{-2\logDensity}\gradEnt^a\left\{(\p_av^b)\p_b\logDensity-(\p_a\logDensity)\p_bv^b\right\}.
			\end{align}
		\end{subequations}
		Also, $\mathfrak{L}_{(v)}^i,\mathfrak{L}_{(\logDensity)},\mathfrak{L}_{(\Ent)},\mathfrak{L}_{(\vortrenormalized)}^i,\mathfrak{L}_{(\dive\vortrenormalized)},\mathfrak{L}^i_{(\modivort)}$ are linear in $(\pfour\vvariables)$ defined as follows:
		\begin{subequations}\label{E:Lequations}
			\begin{align}
				\mathfrak{L}_{(v)}^i:=&2e^{\logDensity}\antisymmetic_{iab}(\materialderivative v^a)\vortrenormalized^b-\frac{\pressure_{;\Ent}}{\backgroundDensity}\antisymmetic_{iab}\vortrenormalized^a\gradEnt^b-\frac{1}{2}e^{-\logDensity}\frac{\pressure_{;\Ent;\logDensity}}{\backgroundDensity}\gradEnt^a\p_av^i\\
				\notag&-2e^{-\logDensity}\speed^{-1}\speed_{;\logDensity}\frac{\pressure_{;\Ent}}{\backgroundDensity}(\materialderivative\logDensity)\gradEnt^i+e^{-\logDensity}\frac{\pressure_{;\Ent;\logDensity}}{\backgroundDensity}(\materialderivative\logDensity)\gradEnt^i,\\
				\mathfrak{L}_{(\logDensity)}:=&-\frac{5}{2}e^{-\logDensity}\frac{\pressure_{;\Ent;\logDensity}}{\backgroundDensity}\gradEnt^a\p_a\logDensity-e^{-\logDensity}\frac{\pressure_{;\Ent;\Ent}}{\backgroundDensity}\delta_{ab}\gradEnt^a\gradEnt^b,\\
				\mathfrak{L}_{(\Ent)}:=&\speed^2\gradEnt^a\p_a\logDensity-\speed\speed_{;\logDensity}\gradEnt^a\p_a\logDensity-\speed\speed_{;\Ent}\delta_{ab}\gradEnt^a\gradEnt^b,\\
				\label{E:LOmega}\mathfrak{L}_{(\vortrenormalized)}^i:=&\vortrenormalized^a\p_av^i-e^{-2\logDensity}\speed^{-2}\frac{\pressure_{;\Ent}}{\backgroundDensity}\antisymmetic_{iab}(\materialderivative v^a)\gradEnt^b,\\
				\label{E:LS}\mathfrak{L}_{(\gradEnt)}^i:=&-\delta_{ab}\gradEnt^a\p_iv^b+\antisymmetic_{iab}e^{\logDensity}\vortrenormalized^a\gradEnt^b,\\
				\mathfrak{L}_{(\dive\vortrenormalized)}:=&-\vortrenormalized^a\p_a\logDensity,\\
				\mathfrak{L}^i_{(\modivort)}:=&2e^{-3\logDensity}\speed^{-3}\speed_{;\Ent}\frac{\pressure_{;\Ent}}{\backgroundDensity}(\materialderivative v^i)\delta_{ab}\gradEnt^a\gradEnt^b-2e^{-3\logDensity}\speed^{-3}\speed_{;\Ent}\frac{\pressure_{;\Ent}}{\backgroundDensity}\delta_{ab}(\materialderivative v^a)\gradEnt^b\gradEnt^i,\\
				\notag&e^{-3\logDensity}\speed^{-2}\frac{\pressure_{;\Ent;\Ent}}{\backgroundDensity}(\materialderivative v^i)\delta_{ab}\gradEnt^a\gradEnt^b-e^{-3\logDensity}\speed^{-2}\frac{\pressure_{;\Ent}}{\backgroundDensity}\delta_{ab}(\materialderivative v^a)\gradEnt^b\gradEnt^i.
			\end{align}
		\end{subequations}
	\end{theorem}

	\section{Acoustic Geometry} \label{S:ACOUSTICGEOMETRY}
	In this section, we introduce the geometric frame for the compressible Euler equations. We emphasize that, \textbf{we are not allowed to use the following quantities and equations until we construct initial data and derive a local existence and uniqueness result in a neighborhood of the initial sphere}, since these quantities are not well-defined at this moment. Nevertheless, the equations in this section are crucial in constructing the initial data.
	
	\subsection{The Eikonal Functions and the Geometric Time and Radial Functions} \label{SS:EIKONALFUNCTIONANDTIMEFUNCTION}
	We will construct acoustical eikonal functions $u$ and $\ubar$. By definition, they are
	solutions to the eikonal equation, which is the following fully nonlinear hyperbolic PDE,
	where $\gfour$ is the acoustical metric from Def.\,\ref{D:acousticmetric}:
	\begin{align}\label{EQ:eikonalequations}
		(\gfour^{-1})^{\alpha\beta}\p_\alpha u\p_\beta u&=0,
		&
		(\gfour^{-1})^{\alpha\beta}\p_\alpha \ubar\p_\beta \ubar&=0.
	\end{align}
	The level sets of $u$ and $\ubar$ are \emph{characteristic surfaces}, i.e., sound cones in the context of the present paper.
	Our main existence result concerns regions that are foliated by the level sets of $u$ and $\ubar$, which will intersect transversally
	in surfaces that are diffeomorphic to $\mathbb{S}^2$. This is called a \emph{double-null foliation}.

	Note that the eikonal functions depend on the fluid solution since $\gfour^{-1}$ depends on the fluid.
	In our applications, $u, \ubar$ will have transversal gradients and will
	be future directed, i.e.,
	\begin{align}
		\materialderivative u & > 0,
		& \materialderivative \ubar & > 0,
	\end{align}
	where $\materialderivative$ is the material derivative vectorfield from \eqref{E:MATDERIV}.
	
	Our ensuing analysis will rely on the functions in the following definition.
	
	\begin{definition}[Geometric time function and radial variable] \label{DE:TIMEFUNCTIONANDRADIALCOORDINATE}
		In terms of the eikonal functions $u$ and $\ubar$,
		we define the geometric time function $\timefunction$ and the geometric radial variable $\rgeo$ as follows:	
		\begin{subequations}
			\begin{align}
				\timefunction 
				& 
				:
				=
				u
				+
				\ubar,
				\label{E:GEOMETRICTIMEFUNCTION} 
				\\
				\rgeo
				& 
				:=
				\ubar-u. 
				\label{E:GEOMETRICRADIALVARIABLE}
			\end{align}
		\end{subequations}
	\end{definition}

	\subsection{Subsets of Spacetime Tied to the Eikonal Functions}	
	\label{SS:SUBSETSOFSPACETIME}	
	Given a pair of real numbers $u',\ubar'$, we define:
	\begin{align} \label{E:DEFSUUBAR}
		S_{u',\ubar'}
		& := \lbrace
		(t,x^1,x^2,x^3) \in \mathbb{R}^{1+3}
		\ | \
		u(t,x^1,x^2,x^3) = u'
		\rbrace
		\cap
		\lbrace
		(t,x^1,x^2,x^3) \in \mathbb{R}^{1+3}
		\ | \
		\ubar(t,x^1,x^2,x^3) = \ubar'
		\rbrace.
	\end{align}
	In our applications, the $S_{u',\ubar'}$ will be co-dimension $2$ embedded submanifolds that are diffeomorphic to $\mathbb{S}^2$
	and spacelike with respect to $\gfour$.
	
	We define/compute our geometric null frame and its relative vectorfields based on $u$ and $\ubar$ as follows:\\

	\textbf{Geometric Radius $\rgeo$}
	\begin{align}\label{DE:Defr}
		\rgeo:=\ubar-u.
	\end{align}
	
	\textbf{Geodesic vectorfields $\Lgeo^\alpha$ and $\uLgeo^\alpha$}
	\begin{align}\label{D:geo}
		\Lgeo^\alpha&:=-(\gfour^{-1})^{\alpha\beta}\p_\beta u,&
		\uLgeo^\alpha&:=-(\gfour^{-1})^{\alpha\beta}\p_\beta \ubar.
	\end{align}
	
	\textbf{Inverse foliation densities $\invfolia$ and $\uinvfolia$}
	\begin{align}\label{P:invuinv}
		\invfolia&:=\frac{-1}{(\gfour^{-1})^{\alpha\beta}\p_\alpha t\p_\beta u},&
		\uinvfolia&:=\frac{-1}{(\gfour^{-1})^{\alpha\beta}\p_\alpha t\p_\beta \ubar}.
	\end{align}
	
	\textbf{Rescaled null vectorfields $\Lunit$ and $\uLunit$}
	\begin{align}\label{D:LunituLunit}
		\Lunit&:=\invfolia\Lgeo,&
		\uLunit&:=\uinvfolia\uLgeo.
	\end{align}
	
	Note that:
	\begin{align} \label{E:LUNIT0ANDULUNIT0AREUNITY}
		\Lunit^0
		& = \uLunit^0
		= 1.
	\end{align}

	\textbf{Null Lapse $\lapse$}\label{P:lapse}
	\begin{align}\label{D:nullapse}
		\lapse^2 & := -2 \gfour(\Lunit,\uLunit)^{-1}.
	\end{align}
	
	\textbf{Principal null vectorfields $e_3$ and $e_4$}
	\begin{align}\label{D:e3e4}
		e_4 &:=\lapse \Lunit,
		&
		e_3 &:=\lapse \uLunit.
	\end{align}
	
	Note that \eqref{E:LUNIT0ANDULUNIT0AREUNITY} implies the following identity:
	\begin{align} \label{E:E3ANDE40}
		e_4^0 & = e_3^0
		= \lapse.
	\end{align}
	
	We also note the following identities:
	\begin{align} \label{E:GFOUROFMATERIALDERIVATIVEWITHE4ANDE3}
		\gfour(\materialderivative,e_4) 
		& =
		\gfour(\materialderivative,e_3)
		= - \lapse, 
	\end{align}
	which follow from
	\eqref{E:LUNIT0ANDULUNIT0AREUNITY},
	Def.\,\ref{D:e3e4},
	\eqref{E:E3ANDE40},
	and the fact that for any vectorfield $V$,
	$\gfour(\materialderivative,V) = - V^0$.

	
	\textbf{Eikonal function normalized null vectorfields $\Lrenormal$ and $\uLrenormal$}
	\begin{align}\label{D:breveL}
		\Lrenormal&:=\uinvfolia\lapse e_4,&
		\uLrenormal&:=\invfolia\lapse e_3.
	\end{align}
	
	\textbf{Renormalized geometric-time derivative vectorfield $\Trenomal$}
	\begin{align} \label{E:RENORMALIZEDTIMEDERIVATIVEVECTORFIELDDEFINITION}
		\Trenomal:=\frac{1}{2}(\Lrenormal+\uLrenormal).
	\end{align}

	Note that:
	\begin{align} \label{E:COORDINATENULLVECTORFIELDSAPPLIEDTOEIKONALFUNCTIONIS2}
		\Lrenormal
		\ubar
		& = \uLrenormal u
		= - \uinvfolia \invfolia \lapse^2 \gfour(\Lgeo,\uLgeo)
		= - \lapse^2 \gfour(\Lunit,\uLunit)
		= 2,
	\end{align}
	and thus
	\begin{align} \label{E:COORDINATEIMELIKEVECTORFIELDAPPLIEDTOEIKONALFUNCTIONIS1}
		\Trenomal u
		& = \Trenomal \ubar
		= 1.
	\end{align}
	By \eqref{E:COORDINATEIMELIKEVECTORFIELDAPPLIEDTOEIKONALFUNCTIONIS1}, 
	the integral curves of $\Trenomal$ yield bijections between
	the spheres $S_{0,\ubar}$ and $S_{z,\ubar + z}$
	as well as the spheres 
	$S_{\ubar,0}$ and $S_{\ubar + z,z}$.
	In our applications, these bijections will be diffeomorphisms.

	\textbf{Geometric-time unit-derivative vectorfield $\Timelike$}
	\begin{align}\label{D:Timelike}
		\Timelike:=\frac{1}{2}\left(\sqrt{\frac{\uinvfolia}{\invfolia}}e_4+\sqrt{\frac{\invfolia}{\uinvfolia}}e_3\right).
	\end{align}
	
	\textbf{Geometric coordinates}
	In our applications, $\inihypersur \cup \inihypersurbar$
	will be equipped with functions
	$(\inisphere^1,\inisphere^2)$
	such that $(\inisphere^1,\inisphere^2)$ are local angular coordinates on
	every sphere $S_{0,\ubar}, S_{u,0} \in
	\inihypersur \cup \inihypersurbar$
	and such that on $\inihypersur$,
	there exists a null generator $\iniLunit$
	(see \eqref{E:INITIALLUNITONH}) with
	$\iniLunit \inisphere^A = 0$ for $A=1,2$.
	By the uniqueness of null generators of a null hypersurface up to scalar multiple, 
	we will also have: 
	\begin{align} \label{E:DATAANGULARCOORDINATESARETRANSPORTEDINLDIRECTION}
		\Lrenormal \inisphere^A = 0, 
		&&
		\mbox{on } \inihypersur,
		& & A = 1,2.
	\end{align}
	Moreover, the vectorfield $\Lrenormal$ will be everywhere transversal to  
	$\inihypersurbar$. Therefore, after we have constructed $\Lrenormal$ in 
	a future-neighborhood of $\inihypersur \cup \inihypersurbar$, 
	we can extend the angular coordinate functions to the same future-neighborhood
	by solving the following transport equation initial value problems:
	\begin{align} \label{E:TRANSPORTINGTHEANGULARCOORDINATEFUNCTIONS} 
		\Lrenormal \cartiasphere & = 0,
		&
		\cartiasphere \restriction_{\inihypersurbar}
		= \inisphere^A,
		&&
		A =1,2.
	\end{align}
	Note that the transport equation \eqref{E:TRANSPORTINGTHEANGULARCOORDINATEFUNCTIONS}  
	is consistent with \eqref{E:DATAANGULARCOORDINATESARETRANSPORTEDINLDIRECTION}.
	
	\begin{definition}[Geometric coordinates]
		\label{D:GEOMETRICCOORDINATES}
		We refer to $\lbrace u, \ubar, \sphere^1, \sphere^2 \rbrace$ as
		the \emph{geometric coordinates}. We denote the corresponding partial derivative vectorfields
		by $\left\lbrace \frac{\partial}{\partial u}, \frac{\partial}{\partial \ubar}, \frac{\partial}{\partial \sphere^1}, 
		\frac{\partial}{\partial \sphere^2}
		\right\rbrace$.
	\end{definition}

	\textbf{$\suubar$-tangent vectorfield $\Theta$}
	\begin{align}\label{D:Theta}
		\Theta:=\materialderivative-\frac{\lapse}{2}(e_4+e_3).
	\end{align}
	
	\textbf{$\Stau$-tangent $\suubar$-normal vectorfield $\spherenormal$}
	\begin{align}\label{DE:N}
		\spherenormal:=\frac{1}{2}(\sqrt{\frac{\uinvfolia}{\invfolia}}e_4-\sqrt{\frac{\uinvfolia}{\invfolia}}e_3).
	\end{align}

	\textbf{$\Stau$-projected version of $\materialderivative$}
	\begin{align}\label{Def:K}
		K:=\frac{2}{\lapse}\cdot\left(\sqrt{\frac{\uinvfolia}{\invfolia}}+\sqrt{\frac{\invfolia}{\uinvfolia}}\right)^{-1}\materialderivative-\Timelike
	\end{align}
	
	\textbf{Angular coordinate functions $\cartiasphere$}

	\textbf{$\suubar$ tangent vectorfield $\Lambda$}
	\begin{align}\label{E:breveL}
		\Lrenormal&=\frac{\p}{\p\ubar},&
		\uLrenormal&=\frac{\p}{\p u}+\Lambda.
	\end{align}
	
	We denote the level sets of $u$ by $\hypersur_u$ and the level sets of $\ubar$ by $\hypersurbar_{\ubar}$. 
	Notice that $\hypersur_u$ and $\hypersurbar_{\ubar}$ are null hypersurfaces with principal null generators $e_4$ and $e_3$ respectively.
	In view of \eqref{E:DEFSUUBAR}, we see that the $2$-sphere $\suubar$ satisfies:
	\begin{align}
		\suubar & = \hypersur_u \bigcap \hypersurbar_{\ubar}.
	\end{align}
	
	We define $\region=\bigcup\limits_{u}\hypersur_{u}=\bigcup\limits_{\ubar}\hypersurbar_{\ubar}$. 
	We refer readers to Section \ref{SS:Defregion} for rough definitions of $\region$, and Section \ref{SS:SUBSETSOFSPACETIME} for precise definition of $\hypersur,\hypersurbar$ and $\region$.
	
	We define $\gsphere$ to be the inverse of $\gsphere^{-1}$, and $\gsphere^{-1}$ as follows:
	\begin{align}\label{D:decompositionofmetric}
		(\gsphere^{-1})^{\alpha\beta}:=(\gfour^{-1})^{\alpha\beta}+\frac{1}{2}e_4^\alpha e_3^\beta+\frac{1}{2}e_3^\alpha e_4^\beta.
	\end{align}
	
	It is easy to check that $\gsphere$ is an induced metric of $\gfour$ on $\suubar$. We denote $\angnabla$ the Levi-Civita connection on $\suubar$ with respect to $\gsphere$. We fix a pair of locally defined unit orthogonal spherical vectorfields\footnote{In the rest of the paper, we automatically sum over $A$ if there are two $A$'s in the expression.} on $\suubar$ by $\{e_A\}_{A=1,2}$ such that $(\gsphere^{-1})^{\alpha\beta}=\sum\limits_{A=1,2}e_A^\alpha e_A^\beta$. We call $\{e_I\}_{I=1}^4$ a null frame for the geometry. For convenience, we call $\{e_J\}_{J=3}^4$ a null pair for the geometry.
	
	By the construction above, we have the following identities:
	\begin{subequations}\label{EQ:innerproduct}
		\begin{align}
			\gfour(e_4,e_4)=&0,&\gfour(e_3,e_3)=&0,\\
			\gfour(e_3,e_4)=&-2,&\gfour(e_A,e_B)=&\gsphere_{AB},\\
			\gfour(e_A,e_4)=&0,&\gfour(e_A,e_3)=&0.
		\end{align}
	\end{subequations}
	\subsection{First Fundamental Forms and Projections}\label{SS:fundamentalforms} In this subsection, we define various first fundamental forms and projections. 
	
	\begin{definition}[First fundamental forms]\label{D:firstfundamentalforms}
		$\materialderivative$ is the material derivative defined in Section \ref{SS:firstorderequations}, $\Timelike$ is the unit $\gfour$-normal to $\Stau$ defined in (\ref{D:Timelike}). $e_3,e_4$ is the principal null pair defined in (\ref{D:e3e4}). We define $\gt,\gtau,$ and $\gsphere$ the first fundamental form of $\St$, $\Stau$ and $\suubar$ respectively:
		\begin{subequations}\label{D:defg}
			\begin{align}
				\gt_{\alpha\beta}:=&\gfour_{\alpha\beta}+\materialderivative_\alpha\materialderivative_\beta,\\
				\gtau_{\alpha\beta}:=&\gfour_{\alpha\beta}+\Timelike_\alpha\Timelike_\beta,\\
				\label{DE:gsphere}\gsphere_{\alpha\beta}:=&\gfour_{\alpha\beta}+\frac{1}{2}(e_4)_\alpha (e_3)_\beta+\frac{1}{2}(e_3)_\alpha (e_4)_\beta.
			\end{align}
		\end{subequations}
	\end{definition}
	\begin{definition}[Projection operators]\label{D:projectionoperators}
		We define the $\St$-projection $\Sigmatproject$, $\Stau$-projection $\Sigmatauproject$ and $\suubar$-projection $\sphereproject$ as follows:
		\begin{subequations}
			\begin{align}
				\Sigmatproject^\alpha_\beta:=&\delta^\alpha_\beta+\materialderivative^\alpha\materialderivative_\beta,\\
				\Sigmatauproject^\alpha_\beta:=&\delta^\alpha_\beta+\Timelike^\alpha\Timelike_\beta,\\
				\sphereproject^\alpha_\beta:=&\delta^\alpha_\beta+\frac{1}{2}(e_4)^\alpha (e_3)_\beta+\frac{1}{2}(e_3)^\alpha (e_4)_\beta.
			\end{align}
		\end{subequations}
	\end{definition}
	\begin{definition}[Projection of Cartesian coordinate partial derivatives]\label{DE:Projections}
		We define the $\Stau$-projected vectorfields $\taup_\alpha$ and $\suubar$-projected vectorfields $\angp_\alpha$ as follows, where $\alpha=0,1,2,3$:
		\begin{subequations}\label{D:taupangp}
			\begin{align}
				\taup_\alpha:=&\Sigmatauproject^\beta_\alpha\p_\beta,&
				\taup^\alpha:=&(\gtau^{-1})^{\alpha\beta}\p_\beta,\\
				\angp_\alpha:=&\sphereproject^\beta_\alpha\p_\beta,&
				\angp^\alpha:=&(\gsphere^{-1})^{\alpha\beta}\p_\beta.
			\end{align}
		\end{subequations}
		Moreover, recalling Definition \ref{DE:dfour} of $\dfour{\alpha}{\beta}$, we define its projection $\angdfour$ as follows:
		\begin{align}\label{DE:Projecteddfour}
			\angdfour^\alpha_\beta:=&(\gsphere^{-1})^{\alpha\gamma}\dfour{\gamma}{\beta},&
			\angdfour_{\alpha\beta}:=&\sphereproject^\gamma_\alpha\dfour{\gamma}{\beta}.
		\end{align}
		In addition, we define the projection of the material derivative $\angB$ as follows:
		\begin{align}
			\label{D:Bdecomp0}\angB^\alpha:=&\sphereproject^\alpha_\gamma\materialderivative^\gamma.
		\end{align}
		Lastly, we define the projection of the covariant derivative associated to a vectorfield $V$ as follows:
		\begin{align}
			\angD_V\xi=\sphereproject\Dfour_V\xi
		\end{align}
		Note that $\angD_V\upxi=\angnabla_V\upxi$ when both $V$ and $\upxi$ are $\stu$-tangent.
	\end{definition}
	\subsection{Connection Coefficients and Commutator Formula}\label{SS:EquationsforData}
	In this section, we define the connection coefficients and prove commutator formulas.
	\begin{definition}[Connection coefficients]\label{D:Defconnectioncoeff}
		Given a 2-sphere $\suubar$ and $\{e_A\}_{A=1,2}$ an arbitrary frame tangent to it, in particular, for $\{e_A\}_{A=1,2}$ defined in Section \ref{SS:SUBSETSOFSPACETIME}, we define the connection coefficients $\ACC:=\{\upchi,\uchi,\upeta,\ueta,\upzeta,\omega,\uomega\}$ relative to the null frame in Section \ref{SS:SUBSETSOFSPACETIME} as follows:
		\begin{subequations}\label{DE:ConnectionCoefficients}
			\begin{align}
				\upchi_{AB}&=\gfour(\Dfour_Ae_4,e_B),&
				\uchi_{AB}&=\gfour(\Dfour_Ae_3,e_B),\\
				\upeta_A&=\frac{1}{2}\gfour(\Dfour_3e_4,e_A),&
				\ueta_A&=\frac{1}{2}\gfour(\Dfour_4e_3,e_A),\\
				\omega&=-\frac{1}{4}\gfour(\Dfour_4e_3,e_4),&
				\uomega&=-\frac{1}{4}\gfour(\Dfour_3e_4,e_3),\\
				\upzeta_A&=\frac{1}{2}\gfour(\Dfour_Ae_4,e_3).
			\end{align}
		\end{subequations}
		We denote $\gtr\upchi:=(\gsphere^{-1})^{AB}\upchi_{AB}$, $\gtr\uchi:=(\gsphere^{-1})^{AB}\uchi_{AB}$ the trace part of $\upchi_{AB},\uchi_{AB}$ respectively and $\hchi:=\upchi-\frac{1}{2}\gtr\upchi\gsphere,\huchi:=\uchi-\frac{1}{2}\gtr\uchi\gsphere$ the traceless part.
	\end{definition}
	\begin{proposition}[Connection coefficients and relationship between tensors]\label{P:PropConnection}
		By Definition \ref{D:Defconnectioncoeff}, the covariant derivatives of the geometric frame can be expressed as follows:
		\begin{subequations}\label{E:propconnection}
			\begin{align}
				\label{DLL}\Dfour_{4}e_4&=-2\omega e_4,&
				\Dfour_{3}e_3&=-2\uomega e_3,\\
				\Dfour_Ae_4&=\upchi_{AB}e_B-\upzeta_Ae_4,&
				\Dfour_Ae_3&=\uchi_{AB}e_B+\upzeta_Ae_3,\\
				\Dfour_{3}e_4&=2\upeta_A e_A+2\uomega e_4,&
				\label{DLbarL}\Dfour_{4}e_3&=2\ueta_Ae_A+2\omega e_3,\\		
				\Dfour_{4} e_A&=\angD_{4} e_A+\ueta_Ae_4,&
				\Dfour_{3}e_A&=\angD_{3}e_A+\upeta_Ae_3,\\
				\Dfour_Ae_B&=\angnabla_Ae_B+\frac{1}{2}\uchi_{AB}e_4+\frac{1}{2}\upchi_{AB}e_3,\\
				\Dfour_\Timelike\Timelike&=\angnabla_A(\ln\frac{\invfolia}{\uinvfolia}) e_A+\frac{1}{2}(\uomega-\omega)e_4+\frac{1}{2}(\omega-\uomega)e_3,\\
				\Dfour_\spherenormal\spherenormal&=-\angnabla_A(\ln\frac{\invfolia}{\uinvfolia}) e_A-\frac{1}{2}(\omega+\uomega)e_4+\frac{1}{2}(\omega+\uomega)e_3.
			\end{align}
		\end{subequations}
		Moreover,
		\begin{subequations}
			\begin{align}
				\label{omegauomega}\omega&=-\frac{1}{2}e_4\ln (\lapse\invfolia),&
				\uomega&=-\frac{1}{2}e_3\ln (\lapse\uinvfolia),\\
				\label{etaupeta}\upeta_A&=\upzeta_A+\angnabla_A(\ln \frac{\invfolia}{\uinvfolia}),&
				\ueta_A&=-\upzeta_A+\angnabla_A(\ln\frac{\invfolia}{\uinvfolia}).&
			\end{align}
		\end{subequations}
		Note that $\upchi$ and $\uchi$ are symmetric, and the following relations hold:
		\begin{subequations}\label{DE:SECONDFORMLIEDIFFERENTIATIONDEF}
			\begin{align}
				\upchi & = \frac{1}{2} \angLie_{4} \gsphere
				= \frac{1}{2} \angLie_{4} \gfour,
				& 
				\uchi
				& = \frac{1}{2} \angLie_{3} \gsphere
				= \frac{1}{2} \angLie_{3} \gfour,
			\end{align}
		\end{subequations}
		where $\Lie$ is the standard Lie derivative, and $\angLie=\sphereproject\Lie$. 
	\end{proposition}

	We need the following lemma in preparation for understanding the Euler constraint equation in Proposition \ref{P:constraint0}.
	\begin{lemma}[Decomposition of derivatives]\label{L:derivativesdecomp}
		Let $\angB,\angp$ be defined as in Definition \ref{DE:Projections}. Given the geometric frame in Section \ref{SS:SUBSETSOFSPACETIME}, we can decompose the material derivative $\materialderivative$ and partial derivatives as follows:
		\begin{subequations}
			\begin{align}
				\label{D:Bdecomp}\angB=&\materialderivative+\frac{1}{2}\gfour(\materialderivative,e_4) e_3+\frac{1}{2}\gfour(\materialderivative,e_3) e_4=(\gsphere^{-1})^{AB}\gfour(\materialderivative,e_A)e_B,\\
				\label{D:pdecomp}\angp_a=&\p_a+\frac{1}{2}\gfour(\p_a,e_4) e_3+\frac{1}{2}\gfour(\p_a,e_3) e_4=(\gsphere^{-1})^{AB}\gfour(\p_a,e_A)e_B.
			\end{align}
		\end{subequations}
	\end{lemma}
	\begin{proof}[Discussion of the proof of Lemma \ref{L:derivativesdecomp}]
		Lemma \ref{L:derivativesdecomp} can be proved by straightforward computation and definitions (\ref{D:LunituLunit})-(\ref{EQ:innerproduct}).
	\end{proof}
	\begin{lemma}[An expression for $\angD_3\angp_\alpha$ and $\angD_4\angp_\alpha$]\label{LE:ExpressionofD3angp}
		\begin{subequations}
			We can express $\angD_3\angp_\alpha$ and $\angD_4\angp_\alpha$ as follows:
			\begin{align}
				\label{EQ:ExpressionofD4angp}\angD_4\angp_\alpha^\beta=&\sphereproject^\beta_\gamma\Chfour_{\alpha 4}^\gamma+(\gsphere^{-1})^{\mu\beta}\ueta_{\mu}(e_4)_\alpha,\\
				\label{EQ:ExpressionofD3angp}\angD_3\angp_\alpha^\beta=&\sphereproject^\beta_\gamma\Chfour_{\alpha 3}^\gamma+(\gsphere^{-1})^{\mu\beta}\upeta_{\mu}(e_3)_\alpha.
			\end{align}
		\end{subequations}
	\end{lemma}
	\begin{proof}[Proof of Lemma \ref{LE:ExpressionofD3angp}]
		We only prove \eqref{EQ:ExpressionofD3angp}, and the proof of \eqref{EQ:ExpressionofD4angp} is similar. By Definition \ref{D:projectionoperators} and Definition \ref{DE:Projections}, we can express $\angD_3\angp_\alpha$ as follows:
		\begin{align}
			\angD_3\angp_\alpha^\beta=&\sphereproject^\beta_\gamma\Dfour_3\left\{\delta^\gamma_\alpha+\frac{1}{2}(e_4)^\gamma(e_3)_\alpha+\frac{1}{2}(e_3)^\gamma(e_4)_\alpha\right\}.
		\end{align}
		The desired result \eqref{EQ:ExpressionofD3angp} follows from straightforward computation and Proposition \ref{P:PropConnection}.
	\end{proof}
	\begin{lemma}[Commutator formula]\label{LE:CommutatorFormula}
		Let $e_J$ with $J=3,4$ be as defined in Section \ref{SS:SUBSETSOFSPACETIME}, $\sphereproject$ be as defined in Definition \ref{D:projectionoperators}, $\angp$ be as defined in Definition \ref{DE:Projections}, $\Chfour$ be as defined in Definition \ref{D:Christoffelsymbols}. Then we have the following commutator formulas:
		\begin{subequations}\label{EQ:commuteformula1}
			\begin{align}
				\label{E:commuteformula00}[e_4,\angp_\alpha]=&~\angD_{4}\angp_\alpha-(\gsphere^{-1})^{\beta\gamma}\upchi_{\angp_\alpha\angp_\beta}\angp_\gamma+(\ueta_{\angp_\alpha}+\upzeta_{\angp_\alpha})e_4,\\
				\notag=&~(\gsphere^{-1})^{\beta\gamma}\gfour(e_4\angp_\alpha,\angp_\beta)\angp_\gamma-(\gsphere^{-1})^{\beta\gamma}\gfour(\angp_\alpha e_4,\angp_\beta)\angp_\gamma+\frac{1}{2}\left[\gfour(e_4e_3,\angp_\alpha)+\gfour(\angp_\alpha e_4,e_3)\right]e_4\\
				\notag&+\frac{1}{2}(\dfour{3}{\angp_\alpha}e_4\variables) e_4,\\
				\label{E:commuteformula01}[e_3,\angp_\alpha]=&~\angD_{3}\angp_\alpha-(\gsphere^{-1})^{\beta\gamma}\uchi_{\angp_\alpha\angp_\beta}\angp_\gamma+(\upeta_{\angp_\alpha}-\upzeta_{\angp_\alpha})e_3,\\
				\notag=&~(\gsphere^{-1})^{\beta\gamma}\gfour(e_3\angp_\alpha,\angp_\beta)\angp_\gamma-(\gsphere^{-1})^{\beta\gamma}\gfour(\angp_\alpha e_3,\angp_\beta)\angp_\gamma+\frac{1}{2}\left[\gfour(e_3e_4,\angp_\alpha)+\gfour(\angp_\alpha e_3,e_4)\right]e_3\\
				\notag&+\frac{1}{2}(\dfour{4}{\angp_\alpha}e_3\variables) e_3,\\
				\label{E:commuteformula02}[e_3,e_4]=&~4\upzeta\cdot\angp+2\uomega e_4-2\omega e_3\\
				\notag=&~2(\gsphere^{-1})^{\alpha\beta}\gfour(\angp_\alpha e_4,e_3)\angp_\beta-\frac{1}{2}\gfour(e_3e_4,e_3)e_4+\frac{1}{2}\gfour(e_4e_3,e_4)e_3\\
				\notag&+2(\gsphere^{-1})^{\alpha\beta}\Chfour_{3\angp_\alpha4}\angp_\beta-\frac{1}{2}\Chfour_{334}e_4+\frac{1}{2}\Chfour_{443}e_3,\\
				\label{E:commuteformula03}[\angp_\alpha,\angp_\beta]=&~\angnabla_{\angp_\alpha}\angp_\beta-\angnabla_{\angp_\beta}\angp_\alpha\\
				\notag=&~(\gsphere^{-1})^{\gamma\delta}\left[\gfour(\angp_\alpha\angp_\beta,\angp_\gamma)\angp_\delta-\gfour(\angp_\beta\angp_\alpha,\angp_\gamma)\angp_\delta\right],
			\end{align}	
		\end{subequations}
		Moreover, for $\suubar$-tangent $\binom{0}{1}$-tensorfield $\xi$, we have the following commutator formula\footnote{The last line on the RHS of \eqref{EQ:commuteformula2.0} and \eqref{EQ:commuteformula2.1} are schematic expressions.}:
		\begin{subequations}\label{EQ:commuteformula2}
			\begin{align}
				\label{EQ:commuteformula2.0}[\angD_4,\angp_\alpha]\xi_{\angp_\beta}=&-\Chfour_{(\angp\xi_{\angp_\beta})4\angp_\alpha}-\Chfour_{(\angp_\alpha\xi)4\angp_\alpha}+\angp_\alpha\Chfour_{\xi4\angp_\beta}-\Chfour_{\xi(\angp_\alpha e_4)\angp_\beta}-\Chfour_{\xi 4(\angp_\alpha\angp_\beta)}\\
				\notag&-(\angD_{4}\angp_\alpha)\xi_{\angp_\beta}+(\gsphere^{-1})^{\delta\gamma}\upchi_{\angp_\alpha\angp_\delta}\angp_\gamma\xi_{\angp_\beta}-(\ueta_{\angp_\alpha}+\upzeta_{\angp_\alpha})e_4\xi_{\angp_\beta}\\
				\notag=&~\gensmoothfunction(\vvariables)(e_4,\angp)^2\vvariables\cdot\xi+\gensmoothfunction(\vvariables)\left\{(e_4,\angp)(\vvariables,e_J)\right\}\cdot(e_4,\angp)\xi,\\
				\label{EQ:commuteformula2.1}[\angD_3,\angp_\alpha]\xi_{\angp_\beta}=&-\Chfour_{(\angp\xi_{\angp_\beta})3\angp_\alpha}-\Chfour_{(\angp_\alpha\xi)3\angp_\alpha}+\angp_\alpha\Chfour_{\xi3\angp_\beta}-\Chfour_{\xi(\angp_\alpha e_3)\angp_\beta}-\Chfour_{\xi 3(\angp_\alpha\angp_\beta)}\\
				\notag&-(\angD_{3}\angp_\alpha)\xi_{\angp_\beta}+(\gsphere^{-1})^{\delta\gamma}\uchi_{\angp_\alpha\angp_\delta}\angp_\gamma\xi_{\angp_\beta}-(\upeta_{\angp_\alpha}-\upzeta_{\angp_\alpha})e_3\xi_{\angp_\beta}\\
				\notag=&~\gensmoothfunction(\vvariables)(e_3,\angp)^2\vvariables\cdot\xi+\gensmoothfunction(\vvariables)\left\{(e_4,\angp)(\vvariables,e_J)\right\}\cdot(e_3,\angp)\xi.
			\end{align}
		\end{subequations}
		
		In addition, for $\suubar$-tangent $\binom{0}{2}$-tensorfield $\xi$, we have the following commutator formula\footnote{The last line on the RHS of \eqref{EQ:commuteformula3.0} and \eqref{EQ:commuteformula3.1} are schematic expressions.}:
		\begin{subequations}\label{EQ:commuteformula3}
			\begin{align}
				\label{EQ:commuteformula3.0}[\angD_4,\angp_\alpha]\xi_{\angp_\mu\angp_\nu}=&-\sphereproject^\gamma_\delta\Chfour_{4\angp_\alpha}^\delta\p_\gamma\xi_{\angp_\mu\angp_\nu}-\angp_\alpha\xi_{\gamma\angp_\nu}\cdot\sphereproject^\gamma_\delta\Chfour_{4\angp_\mu}^\delta-\angp_\alpha\xi_{\gamma\angp_\mu}\cdot\sphereproject^\gamma_\delta\Chfour_{4\angp_\nu}^\delta\\
				\notag&-\angp_\alpha\left[-(\gsphere^{-1})^{\kappa\gamma}\Chfour_{\gamma3\angp_\nu}\xi_{\angp_\mu\kappa}-(\gsphere^{-1})^{\kappa\gamma}\Chfour_{\gamma4\angp_\mu}\xi_{\angp_\nu\kappa}\right]\\
				\notag&+(\gsphere^{-1})^{\kappa\gamma}\Chfour_{\gamma(\angp_\alpha e_4)\angp_\mu}\xi_{\angp_\nu\kappa}+(\gsphere^{-1})^{\kappa\gamma}\Chfour_{\gamma(\angp_\alpha e_4)\angp_\nu}\xi_{\angp_\mu\kappa}\\
				\notag&+(\gsphere^{-1})^{\kappa\gamma}\Chfour_{\gamma 4(\angp_\alpha\angp_\nu)}\xi_{\angp_\mu\kappa}+(\gsphere^{-1})^{\kappa\gamma}\Chfour_{\gamma 4(\angp_\alpha\angp_\mu)}\xi_{\angp_\nu\kappa}\\
				\notag&+(\gsphere^{-1})^{\kappa\gamma}\Chfour_{\gamma 4\angp_\mu}\xi_{(\angp_\alpha\angp_\nu)\kappa}+(\gsphere^{-1})^{\kappa\gamma}\Chfour_{\gamma 4\angp_\nu}\xi_{(\angp_\alpha\angp_\mu)\kappa}\\
				\notag&-(\angD_{4}\angp_\alpha)\xi_{\angp_\mu\angp_\nu}+(\gsphere^{-1})^{\delta\gamma}\uchi_{\angp_\alpha\angp_\delta}\angp_\gamma\xi_{\angp_\mu\angp_\nu}-(\ueta_{\angp_\alpha}+\upzeta_{\angp_\alpha})e_4\xi_{\angp_\mu\angp_\nu}\\
				\notag=&~\gensmoothfunction(\vvariables)(e_4,\angp)^2\vvariables\cdot\xi+\gensmoothfunction(\vvariables)\left\{(e_4,\angp)(\vvariables,e_J)\right\}\cdot(e_4,\angp)\xi,\\
				\label{EQ:commuteformula3.1}[\angD_3,\angp_\alpha]\xi_{\angp_\mu\angp_\nu}=&-\sphereproject^\gamma_\delta\Chfour_{3\angp_\alpha}^\delta\p_\gamma\xi_{\angp_\mu\angp_\nu}-\angp_\alpha\xi_{\gamma\angp_\nu}\cdot\sphereproject^\gamma_\delta\Chfour_{3\angp_\mu}^\delta-\angp_\alpha\xi_{\gamma\angp_\mu}\cdot\sphereproject^\gamma_\delta\Chfour_{3\angp_\nu}^\delta\\
				\notag&-\angp_\alpha\left[-(\gsphere^{-1})^{\kappa\gamma}\Chfour_{\gamma3\angp_\nu}\xi_{\angp_\mu\kappa}-(\gsphere^{-1})^{\kappa\gamma}\Chfour_{\gamma3\angp_\mu}\xi_{\angp_\nu\kappa}\right]\\
				\notag&+(\gsphere^{-1})^{\kappa\gamma}\Chfour_{\gamma(\angp_\alpha e_3)\angp_\mu}\xi_{\angp_\nu\kappa}+(\gsphere^{-1})^{\kappa\gamma}\Chfour_{\gamma(\angp_\alpha e_3)\angp_\nu}\xi_{\angp_\mu\kappa}\\
				\notag&+(\gsphere^{-1})^{\kappa\gamma}\Chfour_{\gamma 3(\angp_\alpha\angp_\nu)}\xi_{\angp_\mu\kappa}+(\gsphere^{-1})^{\kappa\gamma}\Chfour_{\gamma 3(\angp_\alpha\angp_\mu)}\xi_{\angp_\nu\kappa}\\
				\notag&+(\gsphere^{-1})^{\kappa\gamma}\Chfour_{\gamma 3\angp_\mu}\xi_{(\angp_\alpha\angp_\nu)\kappa}+(\gsphere^{-1})^{\kappa\gamma}\Chfour_{\gamma 3\angp_\nu}\xi_{(\angp_\alpha\angp_\mu)\kappa}\\
				\notag&-(\angD_{3}\angp_\alpha)\xi_{\angp_\mu\angp_\nu}+(\gsphere^{-1})^{\delta\gamma}\uchi_{\angp_\alpha\angp_\delta}\angp_\gamma\xi_{\angp_\mu\angp_\nu}-(\upeta_{\angp_\alpha}-\upzeta_{\angp_\alpha})e_3\xi_{\angp_\mu\angp_\nu}\\
				\notag=&~\gensmoothfunction(\vvariables)(e_3,\angp)^2\vvariables\cdot\xi+\gensmoothfunction(\vvariables)\left\{(e_4,\angp)(\vvariables,e_J)\right\}\cdot(e_3,\angp)\xi.
			\end{align}
		\end{subequations}
	\end{lemma}
	\begin{remark}
		Notice that $[e_4,\angp]=\gensmoothfunction\left((e_4,\angp) (e_{J},\vvariables)\right)(e_4,\angp)$, where $J=3,4$. $e_3$-derivative is not present in the commutator due to the fact that $e_4,\angp$ are tangent to manifold $\hypersur$. Similarly, $[e_3,\angp]=\gensmoothfunction\left((e_3,\angp) (e_{J},\vvariables)\right)(e_3,\angp)$.
	\end{remark}
	
	\begin{proof}[Proof of Lemma \ref{LE:CommutatorFormula}]
		\eqref{EQ:commuteformula1} follows from Definition \ref{D:Defconnectioncoeff}, Proposition \ref{P:PropConnection}, and the expression of connection coefficients in standard derivative as follows:
		\begin{align}
			\gfour(\Dfour_Z X,Y)=\gfour(ZX,Y)+\Chfour_{YZX},
		\end{align}
		where $\Chfour$ is defined in Definition \ref{D:Christoffelsymbols}.
		
		We now prove \eqref{EQ:commuteformula3.1}. First of all, by definition of the covariant derivative, we have:
		\begin{align}
			\angD_3\angp_\alpha^\gamma=&~e_3\angp_\alpha^\gamma+\sphereproject^\gamma_\beta\Chfour_{3\angp_\alpha}^\beta,\\
			\angD_3\xi_{\angp_\mu\angp_\nu}=&~e_3\xi_{\angp_\mu\angp_\nu}-(\gsphere^{-1})^{\kappa\gamma}\Chfour_{\gamma3\angp_\nu}\xi_{\angp_\mu\kappa}-(\gsphere^{-1})^{\kappa\gamma}\Chfour_{\gamma3\angp_\mu}\xi_{\angp_\nu\kappa}.
		\end{align}
		Then,
		\begin{subequations}\label{EQ:commuteformula2.1.1}
			\begin{align}
				\angD_3\angp_\alpha\xi_{\angp_\mu\angp_\nu}=&~e_3(\angp_\alpha\xi_{\angp_\mu\angp_\nu})-(\angD_3\angp_\alpha)\xi_{\angp_\mu\angp_\nu}-\angp_\alpha\xi(\angD_3\angp_\mu,\angp_\nu)-\angp_\alpha\xi(\angD_3\angp_\nu,\angp_\mu),\\
				e_3\angp_\alpha\xi_{\angp_\mu\angp_\nu}=&~e_3(\angp_\alpha\xi_{\angp_\mu\angp_\nu})-(e_3\angp_\alpha)\xi_{\angp_\mu\angp_\nu}-\angp_\alpha\xi(e_3\angp_\mu,\angp_\nu)-\angp_\alpha\xi(e_3\angp_\nu,\angp_\mu),
			\end{align}
		\end{subequations}
		and
		\begin{subequations}\label{EQ:commuteformula2.1.2}
			\begin{align}
				\angp_\alpha\angD_3\xi_{\angp_\mu\angp_\nu}=&~\angp_\alpha (\angD_3\xi_{\angp_\mu\angp_\nu})-\angD_{\angp_\alpha e_3}\xi_{\angp_\mu\angp_\nu}-\angD_3\xi(\angp_\alpha\angp_\mu,\angp_\nu)-\angD_3\xi(\angp_\alpha\angp_\nu,\angp_\mu),\\
				\angp_\alpha e_3\xi_{\angp_\mu\angp_\nu}=&~\angp_\alpha (e_3\xi_{\angp_\mu\angp_\nu})-(\angp_\alpha e_3)\xi_{\angp_\mu\angp_\nu}-e_3\xi(\angp_\alpha\angp_\mu,\angp_\nu)-e_3\xi(\angp_\alpha\angp_\nu,\angp_\mu).
			\end{align}
		\end{subequations}
		By  \eqref{EQ:commuteformula2.1.1} and \eqref{EQ:commuteformula2.1.2}, we have:
		\begin{subequations}\label{EQ:commuteformula2.1.3}
			\begin{align}
				\angD_3\angp_\alpha\xi_{\angp_\mu\angp_\nu}-e_3\angp_\alpha\xi_{\angp_\mu\angp_\nu}=&-\sphereproject^\gamma_\delta\Chfour_{3\angp_\alpha}^\delta\p_\gamma\xi_{\angp_\mu\angp_\nu}-\angp_\alpha\xi_{\gamma\angp_\nu}\cdot\sphereproject^\gamma_\delta\Chfour_{3\angp_\mu}^\delta-\angp_\alpha\xi_{\gamma\angp_\mu}\cdot\sphereproject^\gamma_\delta\Chfour_{3\angp_\nu}^\delta,\\
				\angp_\alpha\angD_3\xi_{\angp_\mu\angp_\nu}-\angp_\alpha e_3\xi_{\angp_\mu\angp_\nu}=&-\angp_\alpha\left[-(\gsphere^{-1})^{\kappa\gamma}\Chfour_{\gamma3\angp_\nu}\xi_{\angp_\mu\kappa}-(\gsphere^{-1})^{\kappa\gamma}\Chfour_{\gamma3\angp_\mu}\xi_{\angp_\nu\kappa}\right]\\
				\notag&+(\gsphere^{-1})^{\kappa\gamma}\Chfour_{\gamma(\angp_\alpha e_3)\angp_\mu}\xi_{\angp_\nu\kappa}+(\gsphere^{-1})^{\kappa\gamma}\Chfour_{\gamma(\angp_\alpha e_3)\angp_\nu}\xi_{\angp_\mu\kappa}\\
				\notag&+(\gsphere^{-1})^{\kappa\gamma}\Chfour_{\gamma 3(\angp_\alpha\angp_\nu)}\xi_{\angp_\mu\kappa}+(\gsphere^{-1})^{\kappa\gamma}\Chfour_{\gamma 3(\angp_\alpha\angp_\mu)}\xi_{\angp_\nu\kappa}\\
				\notag&+(\gsphere^{-1})^{\kappa\gamma}\Chfour_{\gamma 3\angp_\mu}\xi_{(\angp_\alpha\angp_\nu)\kappa}+(\gsphere^{-1})^{\kappa\gamma}\Chfour_{\gamma 3\angp_\nu}\xi_{(\angp_\alpha\angp_\mu)\kappa}.
			\end{align}
		\end{subequations}
		Combining \eqref{EQ:commuteformula2.1.3} and \eqref{E:commuteformula01}, we have the desired result \eqref{EQ:commuteformula3.1}. The proofs of \eqref{EQ:commuteformula3.0} and \eqref{EQ:commuteformula2} are of a similar fashion.
	\end{proof}
	\begin{remark}
		We note that, by Definition \ref{D:projectionoperators} and Definition \ref{DE:Projections}, for the $\alpha$-th component of $\angp_\beta$, we have:
		\begin{align}
			\angp_\beta^\alpha=\sphereproject^\alpha_\beta=\delta^\alpha_\beta+\frac{1}{2}(e_4)^\alpha (e_3)_\beta+\frac{1}{2}(e_3)^\alpha (e_4)_\beta.
		\end{align}
		Therefore, schematically, we can write $(e_J^\alpha,\angp^\alpha)=\gensmoothfunction(\vvariables,e_J^\alpha)$ with $J=3,4$. In the rest of the article, we avoid using vectorfield $e_A$ with $A=1,2$ in the regularity analysis. We instead use $\angp_\alpha$ with $\alpha=0,1,2,3$ since they are globally defined by Definition \ref{DE:Projections}.
	\end{remark}
	\begin{lemma}[Decomposition of $\angnabla$ and $\anglap$]\label{LE:Decompositionofanglap}
		For $\suubar$-tangent tensorfield $\xi$, we can write $\angnabla\xi$ as follows:
		\begin{align}\label{EQ:angnabla}
			\angnabla\xi=\angp\xi+\sphereproject\Chfour\cdot\xi,
		\end{align}
		where $\sphereproject\Chfour\cdot\xi$ schematically denotes the summation over terms of the same type. For example, if $\xi_{AB}$ is a type $\binom{0}{2}$ $\suubar$-tangent tensorfield, then $(\sphereproject\Chfour\cdot\xi)_{BCD}=\Chfour_{BD}^A\xi_{AC}+\Chfour_{CD}^A\xi_{BA}$.
		
		For function $\varphi$, we can write $\anglap \varphi$ as follows\footnote{The last line on \eqref{EQ:anglaplace} is a schematic expression.}:
		\begin{align}\label{EQ:anglaplace}
			\anglap \varphi=&~(\gsphere^{-1})^{\alpha\beta}\left[\angp_\alpha\angp_\beta \varphi-\sphereproject^\delta_\alpha\Chfour_{\delta\beta}^\kappa\angp_\kappa \varphi\right]\\
			\notag=&~\gensmoothfunction(\vvariables,e_J)\angp^2\varphi+\gensmoothfunction(\vvariables,e_J)\angp\vvariables\cdot\angp \varphi.
		\end{align}
	\end{lemma}
	\begin{proof}[Proof of Lemma \ref{LE:Decompositionofanglap}]
		\eqref{EQ:angnabla} follows from the definition of the Levi-Civita connection on $\suubar$. We now prove \eqref{EQ:anglaplace}.
		\begin{align}\label{EQ:anglaplace1}
			\anglap \varphi=&~(\gsphere^{-1})^{\alpha\beta}\angnabla_\alpha\angnabla_\beta \varphi\\
			\notag=&~(\gsphere^{-1})^{\alpha\beta}\sphereproject^\delta_\alpha\Dfour_\delta(\sphereproject^\gamma_\beta\Dfour_\gamma \varphi)\\
			\notag=&~(\gsphere^{-1})^{\alpha\beta}\sphereproject^\delta_\alpha\left[\p_\delta(\sphereproject^\gamma_\beta\p_\gamma \varphi)-\Chfour_{\delta\beta}^\kappa(\sphereproject^\gamma_\kappa\p_\gamma \varphi)\right]\\
			\notag=&~(\gsphere^{-1})^{\alpha\beta}\left[\angp_\alpha\angp_\beta \varphi-\sphereproject^\delta_\alpha\Chfour_{\delta\beta}^\kappa\angp_\kappa \varphi\right],
		\end{align}
		where $\Chfour$ is defined in Definition \ref{D:Christoffelsymbols}. Then the desired result \eqref{EQ:anglaplace} follows from \eqref{EQ:anglaplace1} and the fact that $\gfour(e_3,\angp)=\gfour(e_4,\angp)=0$.
	\end{proof}
	The following proposition provides a precise description of the characteristic initial data assumptions, that is, initial data $(\inivel^i,\inilogDensity,\inient)$ on the initial hypersurfaces $\inihypersur$ and $\inihypersurbar$ (see Section \ref{SS:initialembedding2} for definitions of $\inihypersur$ and $\inihypersurbar$) must satisfy the Euler equations (\ref{E:Brho})-(\ref{E:Bs}).

	\subsection{Notation for the Left-Hand Side of the Compressible Euler Equations}
	\label{SS:NOTAITONFOREQUATIONS}
	In this section, we introduce some notation that we find convenient for our proof that the compressible Euler equations
	are satisfied to all orders on portions of the initial null hypersurfaces $\inihypersur \cup \inihypersurbar$ 
	by the initial data and the corresponding transversal derivatives, which we construct in 
	Prop.\,\ref{P:DETERMININGTRANSVERALDERIVATIVESANDSATISFYINGEQUATIONSONINITIALNULLHYPERSURFACES}.
	Our notation is just abbreviations for the left-hand side of the compressible Euler equations \eqref{E:Brho}--\eqref{E:Bv}.
	Specifically, we define:
	\begin{subequations}\label{EQ:EquationThemselves}
		\begin{align} \label{E:CONSTRAINTDENSITYEVOLUTIONEQUATION} 
			\FR
			& := \materialderivative \logDensity + \dive v,
			\\
			\FV^i
			& 
			: =
			\materialderivative v^i
			+
			\speed^2\delta^{ia} \p_a \logDensity
			+
			e^{-\logDensity}\frac{\pressure_{;\Ent}}{\backgroundDensity}\delta^{ia}\p_a\Ent,
			\label{E:CONSTRAINTVELOCITYEVOLUTIONEQUATION} 
			\\
			\FS
			& : = \materialderivative \Ent.
			\label{E:CONSTRAINTENTROPYEVOLUTIONEQUATION} 
		\end{align}
	\end{subequations}
	
	Under this notation, the compressible Euler equations \eqref{E:Brho}--\eqref{E:Bv} 
	are precisely the system
	$\FR = 0, \FV^i = 0, \FS = 0$.
	
	For future use, we also define:
	\begin{align} \label{E:SPECKLUKWAVEEQUATIONTERMS}
		\FW
		& := \mbox{LHS~\eqref{E:boxgdensity} - RHS~\eqref{E:boxgdensity}},
	\end{align}
	so that the geometric wave equation \eqref{E:boxgdensity} for $\logDensity$ is precisely $\FW = 0$.

	\subsection{Decompositions of the First-order Compressible Euler Equations with respect to the Null Frame}
	\label{SS:DECOMPOSITIONSOFTHECOMPRESSIBLEEULEREQUATIONS}
	Given any spacetime vectorfield $V$,
	we let $V_{\perp}$ denote its $\gfour$-orthogonal projection onto $\Sigma_t$, i.e., with
	$\Sigmatproject$ denoting the $\Sigma_t$-projection from Def.\,\ref{D:projectionoperators},
	we set:
	\begin{align} \label{E:PROJECTIONOFVONTOSIGMAT}
		V_{\perp}^{\alpha} 
		& := \Sigmatproject^{\alpha}_\beta V^{\beta}.
	\end{align}
	Since $\materialderivative$ is $\gfour$-orthogonal to $\Sigma_t$ and satisfies
	$\gfour(\materialderivative,\materialderivative) = -1$,
	we can decompose:
	\begin{align} \label{E:DECOMPOSITIONOFVECTORFIELDINTOMATERIALDERIVATIVEPARTANDSPATIALPART}
		V 
		&
		=
		-
		\gfour(V,\materialderivative) 
		\materialderivative 
		+ 
		V_{\perp}.
	\end{align}

	\begin{lemma}[An identity tied to some components of the compressible Euler equations]
		\label{L:COMPRESSIBLEEULERSOLUTIONIDENTITYFORVDERIVATIVEOFDENSITY}
		Let $V$ be any spacetime vectorfield. Under the notation of Sect.\,\ref{SS:NOTAITONFOREQUATIONS},
		the following algebraic identity holds:
		\begin{align} \label{E:NOCOMPRESSIBLEEULERSOLUTIONIDENTITYFORVDERIVATIVEOFDENSITY}
			\begin{split}
				&
				\frac{1}{2}
				\left\lbrace
				-
				\gfour(V,\materialderivative)
				\gfour(\p_a,e_3) 
				+
				\gfour(\materialderivative,e_3)
				\speed^{-2} V_{\perp}^a
				\right\rbrace
				e_4 v^a
				+
				\frac{1}{2}
				\left\lbrace
				-
				\gfour(V,\materialderivative)
				\gfour(\p_a,e_4)
				+
				\gfour(\materialderivative,e_4)
				\speed^{-2} V_{\perp}^a
				\right\rbrace
				e_3 v^a
				\\
				& = V \logDensity
				-
				\gfour(V,\materialderivative) 
				\angp_a v^a
				+
				\speed^{-2} V_{\perp}^a \angB v^a
				+
				\speed^{-2} e^{-\logDensity}\frac{\pressure_{;\Ent}}{\backgroundDensity} V_{\perp}^a \p_a \Ent
				\\
				& \ \
				+
				\gfour(V,\materialderivative) \FR
				-
				(V_{\perp})_a \FV^a.
			\end{split}
		\end{align}
		In particular, for solutions to the compressible Euler equations 
		\eqref{E:Brho}--\eqref{E:Bv} (and thus $\FR = 0$ and $\FV^a = 0$),
		the following equation holds:
		\begin{align} \label{E:COMPRESSIBLEEULERSOLUTIONIDENTITYFORVDERIVATIVEOFDENSITY}
			\begin{split}
				&
				\frac{1}{2}
				\left\lbrace
				-
				\gfour(V,\materialderivative)
				\gfour(\p_a,e_3) 
				+
				\gfour(\materialderivative,e_3)
				\speed^{-2} V_{\perp}^a
				\right\rbrace
				e_4 v^a
				+
				\frac{1}{2}
				\left\lbrace
				-
				\gfour(V,\materialderivative)
				\gfour(\p_a,e_4)
				+
				\gfour(\materialderivative,e_4)
				\speed^{-2} V_{\perp}^a
				\right\rbrace
				e_3 v^a
				\\
				& = V \logDensity
				-
				\gfour(V,\materialderivative) 
				\angp_a v^a
				+
				\speed^{-2} V_{\perp}^a \angB v^a
				+
				\speed^{-2} e^{-\logDensity}\frac{\pressure_{;\Ent}}{\backgroundDensity} V_{\perp}^a \p_a \Ent.
			\end{split}
		\end{align}
		
		\begin{proof}
			First, using \eqref{E:DECOMPOSITIONOFVECTORFIELDINTOMATERIALDERIVATIVEPARTANDSPATIALPART},
			we decompose:
			\begin{align} \label{E:FIRSTPROOFSTEPCOMPRESSIBLEEULERSOLUTIONIDENTITYFORVDERIVATIVEOFDENSITY}
				V \logDensity 
				& = -
				\gfour(V,\materialderivative) 
				\materialderivative \logDensity
				+ 
				V_{\perp}^a \partial_a \logDensity.
			\end{align}
			Next, we use \eqref{E:CONSTRAINTDENSITYEVOLUTIONEQUATION} 
			to algebraically substitute for the factor $\materialderivative \logDensity$ on 
			RHS~\eqref{E:FIRSTPROOFSTEPCOMPRESSIBLEEULERSOLUTIONIDENTITYFORVDERIVATIVEOFDENSITY},
			and we use \eqref{E:CONSTRAINTVELOCITYEVOLUTIONEQUATION} 
			to algebraically substitute for the factors $\partial_a \logDensity$ on
			RHS~\eqref{E:FIRSTPROOFSTEPCOMPRESSIBLEEULERSOLUTIONIDENTITYFORVDERIVATIVEOFDENSITY},
			thereby obtaining the following identity:
			\begin{align} 
				\begin{split} \label{E:SECONDPROOFSTEPCOMPRESSIBLEEULERSOLUTIONIDENTITYFORVDERIVATIVEOFDENSITY}
					V \logDensity
					& = 
					\gfour(V,\materialderivative) \partial_a v^a
					-
					\speed^{-2} V_{\perp}^a \materialderivative v^a
					- 
					\speed^{-2} e^{-\logDensity}\frac{\pressure_{;\Ent}}{\backgroundDensity} V_{\perp}^a \p_a \Ent
					\\
					& \ \
					-
					\gfour(V,\materialderivative) \FR
					+
					(V_{\perp})_a \FV^a.
				\end{split}
			\end{align}
			From \eqref{E:SECONDPROOFSTEPCOMPRESSIBLEEULERSOLUTIONIDENTITYFORVDERIVATIVEOFDENSITY} and the decompositions
			of Lemma~\ref{L:derivativesdecomp}, we deduce:
			\begin{align} 
				\begin{split}
					\label{E:THIRDPROOFSTEPCOMPRESSIBLEEULERSOLUTIONIDENTITYFORVDERIVATIVEOFDENSITY}
					V \logDensity
					& = 
					\gfour(V,\materialderivative) 
					\left\lbrace
					\angp_a v^a
					-
					\frac{1}{2}\gfour(\p_a,e_4) e_3 v^a
					-
					\frac{1}{2}\gfour(\p_a,e_3) e_4 v^a
					\right\rbrace	
					-
					\speed^{-2} V_{\perp}^a \materialderivative v^a
					- 
					\speed^{-2} e^{-\logDensity}\frac{\pressure_{;\Ent}}{\backgroundDensity} V_{\perp}^a \p_a \Ent
					\\
					& \ \
					-
					\gfour(V,\materialderivative) \FR
					+
					(V_{\perp})_a \FV^a.
				\end{split}
			\end{align}	
			Finally, from \eqref{E:THIRDPROOFSTEPCOMPRESSIBLEEULERSOLUTIONIDENTITYFORVDERIVATIVEOFDENSITY}
			and straightforward algebraic computations, we conclude that 
			\eqref{E:THIRDPROOFSTEPCOMPRESSIBLEEULERSOLUTIONIDENTITYFORVDERIVATIVEOFDENSITY} is equivalent to
			\eqref{E:NOCOMPRESSIBLEEULERSOLUTIONIDENTITYFORVDERIVATIVEOFDENSITY}.

		\end{proof}

	\end{lemma}

	\begin{proposition}[Equivalent formulations of the compressible Euler equations involving the Euler constraints] \label{P:constraint0}
		\ \\
		
		\noindent \underline{\textbf{Equivalencies for individual equations}}.
		The following equation is equivalent to the velocity compressible Euler equation \eqref{E:Bv}:
		\begin{align} \label{E:e3v}
			e_3 v^i
			& =
			-
			\lapse^{-1}
			\left\{-\speed^2\delta^{ia}\left[2\angp_a-\gfour(\p_a,e_4)e_3-\gfour(\p_a,e_3)e_4\right]\logDensity
			-
			2e^{-\logDensity}\frac{\pressure_{;\Ent}}{\average{\density}}\delta^{ia}\p_a\Ent
			+
			\left(2 \angB + \lapse e_4 \right)v^i
			\right\}.
		\end{align}
		
		In addition, the following equation is equivalent to the entropy compressible Euler equation \eqref{E:Bs}:
		\begin{align} \label{E:ENTROPYEULEREQUATIONRELATIVETONULLFRAME}
			e_4 \Ent
			+
			e_3 \Ent
			& 
			= 
			- 
			2 \lapse^{-1} \angB \Ent.
		\end{align}

		
		\medskip
		
		\noindent \underline{\textbf{Equivalent formulation I for the entire first-order compressible Euler system}}.
		The full set of compressible Euler equations \eqref{E:Brho}--\eqref{E:Bs} is equivalent to
		the system comprising equation \eqref{E:e3v}, 
		equation \eqref{E:ENTROPYEULEREQUATIONRELATIVETONULLFRAME}, and
		either one of the following two equations, which we refer to as the \textbf{Euler constraints}:
		\begin{subequations}\label{E:constraint}
			\begin{align}
				\begin{split} \label{E:constraint1} 
					& e_4 \logDensity + 
					\frac{1}{2}\left\{\left[\gfour(\materialderivative,e_3)\gfour(\p_i,e_4)-\gfour(\materialderivative,e_4)\gfour(\p_i,e_3)\right] \right\}
					e_4 v^i
					\\
					& \ \ 
					= 
					\gfour(\materialderivative,e_4)\left(\angB\logDensity+\angp_iv^i\right)+\gfour(\p_i,e_4)\left(\angB v^i+\speed^2\delta^{ia} \angp_a \logDensity
					+
					e^{-\logDensity}\frac{\pressure_{;\Ent}}{\average{\density}}\delta^{ia}\p_a\Ent\right),
				\end{split}	
				\\
				\begin{split} \label{E:constraint2}
					& e_3 \logDensity
					+
					\frac{1}{2}\left\{\left[\gfour(\materialderivative,e_4)\gfour(\p_i,e_3)-\gfour(\materialderivative,e_3)\gfour(\p_i,e_4)\right] \right\}
					e_3 v^i
					\\
					& \ \ = \gfour(\materialderivative,e_3)\left(\angB\logDensity+\angp_iv^i\right)+\gfour(\p_i,e_3)\left(\angB v^i+\speed^2\delta^{ia}\angp_a\logDensity+e^{-\logDensity}\frac{\pressure_{;\Ent}}{\average{\density}}\delta^{ia}\p_a\Ent\right),
				\end{split}
			\end{align}
		\end{subequations}
		where we emphasize that \eqref{E:constraint1} involves only $\hypersur_u$-tangent derivatives
		and \eqref{E:constraint2} involves only $\hypersurbar_{\ubar}$-tangent derivatives.

		\medskip
		
		\noindent \underline{\textbf{Equivalent formulation II for the entire first-order compressible Euler system}}.
		The full set of compressible Euler equations \eqref{E:Brho}--\eqref{E:Bs} is equivalent to
		the system comprising equation \eqref{E:ENTROPYEULEREQUATIONRELATIVETONULLFRAME},
		the Euler constraints \eqref{E:constraint1}--\eqref{E:constraint2},
		and the following two equations:
		\begin{align} \label{E:EULEREQUATIONANGULARCOMPONENTSOFVELOCITY}
			\begin{split}
				&
				\frac{1}{2}
				\left\lbrace
				-
				\gfour(e_A,\materialderivative)
				(e_3)_a
				+
				\gfour(\materialderivative,e_3)
				(e_A)_a
				\right\rbrace
				e_4 v^a
				+
				\frac{1}{2}
				\left\lbrace
				-
				\gfour(e_A,\materialderivative)
				(e_4)_a
				+
				\gfour(\materialderivative,e_4)
				(e_A)_a
				\right\rbrace
				e_3 v^a
				\\
				& = e_A \logDensity
				-
				\gfour(e_A,\materialderivative) 
				\angp_a v^a
				+
				(e_A)_a \angB v^a
				+
				e^{-\logDensity}\frac{\pressure_{;\Ent}}{\backgroundDensity} (e_A)_a \p_a \Ent,
				\\ A & =1,2,
			\end{split}
		\end{align}
		where the coefficients
		$	\left\lbrace
		-
		\gfour(e_A,\materialderivative)
		(e_3)_a
		+
		\gfour(\materialderivative,e_3)
		(e_A)_a
		\right\rbrace$
		and $		\left\lbrace
		-
		\gfour(e_A,\materialderivative)
		(e_4)_a
		+
		(e_A)_a
		\right\rbrace$,
		viewed as ``spatial'' one-forms with Cartesian spatial index $a$,
		are everywhere non-vanishing.
	\end{proposition}	
	
	\begin{remark}[Equivalence with the Riemann invariant transport equations in the $1D$ isentropic case]
		\label{R:EQUIVALENCEWITHRIEMANNINVARIANTEQUATIONSIN1D}
		We note that in the 1D isentropic case, the Euler constraint equations \eqref{E:constraint1} and \eqref{E:constraint2} are equivalent to the transport equations \eqref{EQ:Riemann1}--\eqref{EQ:Riemann2}
		satisfied by the Riemann invariants.
	\end{remark}
	
	\begin{proof}[Proof of Proposition \ref{P:constraint0}]
		The equivalence of equations \eqref{E:Bv} and \eqref{E:e3v} follows from Lemma~\ref{L:derivativesdecomp}
		and the identity \eqref{E:GFOUROFMATERIALDERIVATIVEWITHE4ANDE3}.
		The equivalence of equations \eqref{E:Bs} and \eqref{E:ENTROPYEULEREQUATIONRELATIVETONULLFRAME} follows from the same argument.
		
		We now prove that \eqref{E:Brho}--\eqref{E:Bv} imply equation \eqref{E:constraint1}. 
		This follows from \eqref{E:DECOMPOSITIONOFVECTORFIELDINTOMATERIALDERIVATIVEPARTANDSPATIALPART} and 
		\eqref{E:COMPRESSIBLEEULERSOLUTIONIDENTITYFORVDERIVATIVEOFDENSITY} with
		$e_4$ in the role of $V$, together with the fact that $\materialderivative_a = 0$ for $a=1,2,3$.
		Similarly, 
		\eqref{E:DECOMPOSITIONOFVECTORFIELDINTOMATERIALDERIVATIVEPARTANDSPATIALPART} and 
		\eqref{E:COMPRESSIBLEEULERSOLUTIONIDENTITYFORVDERIVATIVEOFDENSITY} with
		$e_3$ in the role of $V$ imply \eqref{E:constraint2}.
		
		Conversely, 
		if equations \eqref{E:constraint1} and \eqref{E:Bv} are satisfied,
		then by considering the identity \eqref{E:NOCOMPRESSIBLEEULERSOLUTIONIDENTITYFORVDERIVATIVEOFDENSITY}
		with $e_4$ in the role of $V$ and using that $\FV^1=\FV^2=\FV^3 = 0$,
		we conclude that $\FR = 0$, i.e., that equation \eqref{E:Brho}
		is satisfied. Similarly, equations \eqref{E:constraint2} and \eqref{E:Bv} together imply equation
		\eqref{E:Brho}. We have therefore proved the equivalence of the two formulations of compressible Euler flow
		stated just above equation \eqref{E:constraint}.
		
		To see that \eqref{E:EULEREQUATIONANGULARCOMPONENTSOFVELOCITY}
		follows from the first-order compressible Euler equations 
		\eqref{E:Brho}--\eqref{E:Bs}, we simply use \eqref{E:COMPRESSIBLEEULERSOLUTIONIDENTITYFORVDERIVATIVEOFDENSITY} with 
		$e_A$ in the role of $V$.
		We next show that the spatial one-forms
		$	-
		\gfour(e_A,\materialderivative)
		(e_3)_a
		+
		\gfour(\materialderivative,e_3)
		(e_A)_a$
		and
		$		
		-
		\gfour(e_A,\materialderivative)
		(e_4)_a
		+
		(e_A)_a
		$
		are everywhere non-vanishing,
		as is claimed in the proposition.
		We give the proof for the first one since the second can be handled similarly.
		It suffices to show that its norm with respect to $\gfour$ is everywhere positive. By
		\eqref{E:DECOMPOSITIONOFVECTORFIELDINTOMATERIALDERIVATIVEPARTANDSPATIALPART} and the fact that
		$\materialderivative_a = 0$ for $a=1,2,3$, we have that $(e_A)_a = ((e_A)_{\perp})_a$ and
		$(e_3)_a = ((e_3)_{\perp})_a$, where $(e_A)_{\perp}$ and $(e_3)_{\perp}$ are $\Sigma_t$-tangent vectorfields
		(they have vanishing $0$ index when the indices are upstairs).
		Hence, taking the $\gfour$-inner product of the one-form with itself, we see that it suffices to show that:
		\begin{align} \label{E:ONEFORMNORMTHATWEWANTTOBEPOSITIVE}
			[\gfour(e_A,\materialderivative)]^2
			\gfour((e_3)_{\perp},(e_3)_{\perp})
			-
			2 \gfour(e_A,\materialderivative) \gfour(\materialderivative,e_3)
			\gfour((e_3)_{\perp},(e_A)_{\perp})
			+
			[\gfour(\materialderivative,e_3)]^2
			\gfour((e_A)_{\perp},(e_A)_{\perp})
			& > 0.
		\end{align}
		We will prove the following identity (no summation over $A$):
		\begin{align} \label{E:IDENTITYFORONEFORMNORMTHATWEWANTTOBEPOSITIVE}
			\mbox{LHS~\eqref{E:ONEFORMNORMTHATWEWANTTOBEPOSITIVE}}
			& = 
			[\gfour(\materialderivative,e_3)]^2
			\gfour(e_A,e_A),
		\end{align}
		which implies the desired inequality \eqref{E:ONEFORMNORMTHATWEWANTTOBEPOSITIVE} since $\gfour(e_A,e_A) > 0$ and
		$\gfour(\materialderivative,e_3) < 0$,
		where the former inequality holds because $e_A$ is $\gfour$-spacelike and the latter holds
		because $\materialderivative$ is $\gfour$-timelike and future-directed while
		$e_3$ is $\gfour$-null.
		To prove \eqref{E:IDENTITYFORONEFORMNORMTHATWEWANTTOBEPOSITIVE},
		we first use
		\eqref{E:DECOMPOSITIONOFVECTORFIELDINTOMATERIALDERIVATIVEPARTANDSPATIALPART}
		to decompose:
		\begin{align} \label{E:E3DECOMPOSITIONINTOE3PERPANDMATERIALDERIVATIVEPART}
			(e_3)_{\perp}
			& = e_3
			+ 
			\gfour(e_3,\materialderivative)
			\materialderivative,
			\\
			(e_A)_{\perp}
			& = e_A
			+ 
			\gfour(e_A,\materialderivative)
			\materialderivative.
			\label{E:EADECOMPOSITIONINTOEAPERPANDMATERIALDERIVATIVEPART} 
		\end{align}
		
		Next, we recall that:
		\begin{align} \label{E:SIMPLEFRAMEINNERPRODUCTRELATIONS}
			\gfour(e_3,e_3)
			& = 
			\gfour(e_3,e_A) 
			= 0,
			&
			\gfour(\materialderivative,\materialderivative)
			& = -1.
		\end{align}
		Using \eqref{E:EADECOMPOSITIONINTOEAPERPANDMATERIALDERIVATIVEPART}--\eqref{E:SIMPLEFRAMEINNERPRODUCTRELATIONS},
		we compute that:
		\begin{align}
			\gfour((e_3)_{\perp},(e_3)_{\perp})
			& = 
			[\gfour(e_3,\materialderivative)]^2,
			\label{E:E3PERPSQUARENORMIDENTITY} \\
			\gfour((e_3)_{\perp},(e_A)_{\perp})
			& = \gfour(e_3,\materialderivative) \gfour(e_A,\materialderivative),
			\label{E:EAPERPE3PERPINNERPRODUCTIDENTITYIDENTITY}
			\\
			\gfour((e_A)_{\perp},(e_A)_{\perp})
			& = \gfour(e_A,e_A)	
			+
			[\gfour(e_A,\materialderivative)]^2.
			\label{E:EAPERPSQUARENORMIDENTITY} 
		\end{align}
		Finally, we use \eqref{E:EAPERPSQUARENORMIDENTITY}--\eqref{E:EAPERPE3PERPINNERPRODUCTIDENTITYIDENTITY}
		to algebraically substitute for the second factor in each of the three
		terms on LHS~\eqref{E:ONEFORMNORMTHATWEWANTTOBEPOSITIVE}, thereby arriving at \eqref{E:IDENTITYFORONEFORMNORMTHATWEWANTTOBEPOSITIVE}
		after some straightforward computations.
		
		To complete the proof of the proposition, it remains for us to show that
		equations \eqref{E:ENTROPYEULEREQUATIONRELATIVETONULLFRAME},
		\eqref{E:constraint1}--\eqref{E:constraint2},
		and \eqref{E:EULEREQUATIONANGULARCOMPONENTSOFVELOCITY}
		collectively imply the compressible Euler equations \eqref{E:Brho}--\eqref{E:Bs}.
		We have already shown that \eqref{E:Bs} is equivalent to \eqref{E:ENTROPYEULEREQUATIONRELATIVETONULLFRAME}.
		The above arguments and Lemma~\ref{L:COMPRESSIBLEEULERSOLUTIONIDENTITYFORVDERIVATIVEOFDENSITY} 
		show that equations \eqref{E:constraint1}--\eqref{E:constraint2},
		and \eqref{E:EULEREQUATIONANGULARCOMPONENTSOFVELOCITY}
		imply that the terms on the last line of 
		\eqref{E:NOCOMPRESSIBLEEULERSOLUTIONIDENTITYFORVDERIVATIVEOFDENSITY}
		vanish for $V \in \lbrace e_4, e_3, e_1, e_2 \rbrace$,
		i.e.,
		\begin{align}
			-
			\gfour(e_4,\materialderivative) \FR
			+
			((e_4)_{\perp})_a \FV^a
			& = 0,
			\label{E:E4COMPRESSIBLEEULERSYSTEMCOMPONENTSVANISH} \\
			-
			\gfour(e_3,\materialderivative) \FR
			+
			((e_3)_{\perp})_a \FV^a
			& = 0,
			\label{E:E3COMPRESSIBLEEULERSYSTEMCOMPONENTSVANISH} \\
			-
			\gfour(e_1,\materialderivative) \FR
			+
			((e_1)_{\perp})_a \FV^a
			& = 0,
			\label{E:E1COMPRESSIBLEEULERSYSTEMCOMPONENTSVANISH} \\
			-
			\gfour(e_2,\materialderivative) \FR
			+
			((e_2)_{\perp})_a \FV^a
			& = 0.
			\label{E:E2COMPRESSIBLEEULERSYSTEMCOMPONENTSVANISH}
		\end{align}
		We view equations \eqref{E:E4COMPRESSIBLEEULERSYSTEMCOMPONENTSVANISH}--\eqref{E:E2COMPRESSIBLEEULERSYSTEMCOMPONENTSVANISH}
		as a $4 \times 4$ linear system in the unknowns $(\FR,\FV^1,\FV^2,\FV^3)$. We denote the coefficient matrix by $M$. From \eqref{E:DECOMPOSITIONOFVECTORFIELDINTOMATERIALDERIVATIVEPARTANDSPATIALPART},
		we see that the transpose of $M$ 
		can be viewed as the transition matrix between the bases 
		$\lbrace e_4, e_3, e_1, e_2 \rbrace$
		and $\lbrace \materialderivative, \partial_1, \partial_2, \partial_3 \rbrace$.
		Hence, $M$ is invertible, and thus the only solution to the system is
		$(\FR,\FV^1,\FV^2,\FV^3) = (0,0,0,0)$.
		Hence, in view of \eqref{E:CONSTRAINTDENSITYEVOLUTIONEQUATION}--\eqref{E:CONSTRAINTVELOCITYEVOLUTIONEQUATION},
		we conclude that the compressible Euler equations
		\eqref{E:Brho} and \eqref{E:Bv} are satisfied, as is desired.
		
	\end{proof}

	\subsection{Decomposition of Wave Operators with Respect to the Null Frame}\label{SS:DECOMPOSITIONOFWAVEEQUATIONS}
	In this section, we decompose the wave equations with respect to the null frame in Proposition \ref{P:Propdecompboxgf} followed by presenting the evolution equation for the transversal derivative on the null hypersurfaces in Lemma \ref{LE:Solvee3rho} and Corollary \ref{CO:e4e3k}.
	\subsubsection{An algebraic lemma connecting wave equations and first-order Euler equations}
	\begin{proposition}\cite[Lemma 7. Canonical decomposition of the wave operator]{ANewFormulationofthe3DCompressibleEulerEquationswithDynamicEntropy}.\label{P:Propdecompboxgf0}
		Let $\gfour$ be the acoustical metric
		from Definition \ref{D:acousticmetric}. For scalar functions $\varphi$, the expression $\boxg \varphi$ can be decomposed as follows:
		\begin{align}\label{EQ:Boxgf}
			\boxg \varphi=&-\materialderivative\materialderivative \varphi+\speed^2\delta^{ab}\p_a\p_b \varphi+2\speed^{-1}\speed_{;\logDensity}(\materialderivative\logDensity)\materialderivative \varphi-(\p_av^a)\materialderivative \varphi\\
			\notag&-\speed^{-1}\speed_{;\logDensity}(\gfour^{-1})^{\alpha\beta}\p_\alpha\logDensity\p_\beta \varphi-\speed\speed_{;\Ent}\gradEnt^a\p_a\varphi+3\speed^{-1}\speed_{;\Ent}(\materialderivative\Ent)\materialderivative \varphi.
		\end{align}
	\end{proposition}
	\begin{remark}
		We note that (\ref{EQ:Boxgf}) is derived by Definition \ref{D:acousticmetric} and the definition of wave operator in Definition \ref{D:LeviCivita} only, and that 
		no null structure is needed.
	\end{remark}
	\label{SSS:PRELIMINARYALGEBRAICLEMMAFORTRANSVERSALDERIVATIVES}
	The following lemma provides algebraic relationships between the first-order compressible Euler equations, 
	the geometric wave equation \eqref{E:boxgdensity}
	for $\logDensity$, and the reformulation of this equation
	that we established in \eqref{E:WAVEEQUATIONFORDENSITYINADAPATEDCOORDIANTESALONGINITIALNULL}
	relative to the adapted coordinates.

	\begin{lemma}[Algebraic relationship between various equations]
		\label{L:ALGEBRAICRELATIONSHIPESBETWEENSOMEEQUATIONS}
		Let $\FR,\FV,\FS,$ and $\FW$
		be as defined in \eqref{E:CONSTRAINTDENSITYEVOLUTIONEQUATION}--\eqref{E:CONSTRAINTENTROPYEVOLUTIONEQUATION} 
		and
		\eqref{E:SPECKLUKWAVEEQUATIONTERMS}.
		Then the following identity holds:
		\begin{align}
			- \materialderivative 
			\FR
			& = \FW
			-
			\partial_a \FV^a
			+
			(\p_a v^a) \FR
			-
			3\speed^{-1}\speed_{;\Ent} 
			(\materialderivative \logDensity) \FS.
			\label{E:DIFFERENTIATEDDENSITYCONSTRAINEDCOMPAREDTOLUKSPECKWAVE}
		\end{align}
	\end{lemma}

	\begin{proof}
		We first prove \eqref{E:DIFFERENTIATEDDENSITYCONSTRAINEDCOMPAREDTOLUKSPECKWAVE}. 
		We start by using \eqref{E:CONSTRAINTDENSITYEVOLUTIONEQUATION}--\eqref{E:CONSTRAINTVELOCITYEVOLUTIONEQUATION} 
		and Def.\,\ref{DE:ent}
		to compute the following identity, where to obtain the final equality, we use 
		definition~\eqref{Def:D}:
		\begin{align} 
			\begin{split} \label{E:MATERIALDERIVATIVEOFCONSTRAINTVELOCITYEVOLUTIONEQUATION}
				- \materialderivative 
				\FR
				& =
				-
				\materialderivative(\materialderivative \logDensity + \dive v)
				= 
				-\materialderivative \materialderivative 
				\logDensity
				-
				\p_a(\materialderivative v^a)
				+
				\p_a v^b \p_b v^a
				\\
				& = 
				-
				\materialderivative \materialderivative \logDensity
				+
				\speed^2 \delta^{ab} \p_a \p_b \logDensity
				+
				\delta^{ab}(\p_a\speed^2)\p_b \logDensity
				+
				\p_a\left[e^{-\logDensity}\frac{p_{;\Ent}}{\backgroundDensity}\delta^{ab}
				\p_b\Ent\right]
				+
				\p_a v^b \p_b v^a
				-
				\partial_a \FV^a
				\\
				& 
				= 
				-
				\materialderivative \materialderivative \logDensity
				+
				\speed^2 \delta^{ab} \p_a \p_b \logDensity
				+
				2 \speed \speed_{;\logDensity} \delta^{ab} (\p_a \logDensity) \p_b \logDensity
				+
				2 \speed \speed_{;\Ent} \gradEnt^a \p_a \logDensity
				\\
				& \ \
				+
				e^{-\logDensity}\frac{p_{;\Ent}}
				{\backgroundDensity}
				\p_a
				\gradEnt^a
				+
				e^{-\logDensity}
				\frac{\pressure_{;\Ent;\logDensity}}
				{\backgroundDensity}
				\gradEnt^a
				\p_a
				\logDensity
				+
				e^{-\logDensity}\frac{p_{;\Ent;\Ent}}
				{\backgroundDensity}\delta_{ab}
				\gradEnt^a
				\gradEnt^b
				-
				e^{-\logDensity} \frac{p_{;\Ent}}{\backgroundDensity} \gradEnt^a \p_a \logDensity
				+
				\p_a v^b \p_b v^a
				-
				\partial_a \FV^a
				\\
				& 
				= 
				-
				\materialderivative \materialderivative \logDensity
				+
				\speed^2 \delta^{ab} \p_a \p_b \logDensity
				+
				2 \speed \speed_{;\logDensity} \delta^{ab} (\p_a \logDensity) \p_b \logDensity
				+
				2 \speed \speed_{;\Ent} \gradEnt^a \p_a \logDensity
				\\
				& \ \
				+
				e^{\logDensity}\frac{\pressure_{;\Ent}}{\backgroundDensity}\DivGradEnt
				+
				e^{-\logDensity}
				\frac{\pressure_{;\Ent;\logDensity}}
				{\backgroundDensity}
				\gradEnt^a
				\p_a
				\logDensity
				+
				e^{-\logDensity}\frac{p_{;\Ent;\Ent}}
				{\backgroundDensity}\delta_{ab}
				\gradEnt^a
				\gradEnt^b
				+
				\p_a v^b \p_b v^a
				-
				\partial_a \FV^a.
			\end{split}
		\end{align}
		Next, using \eqref{EQ:Boxgf} with $\varphi := \logDensity$,
		we deduce the following identity, where to obtain the second equality, we use
		\eqref{E:CONSTRAINTDENSITYEVOLUTIONEQUATION}--\eqref{E:CONSTRAINTENTROPYEVOLUTIONEQUATION}:
		\begin{align}
			\begin{split} \label{EQ:BoxLOGDENSITYPROOFSTEP}
				\boxg \logDensity 
				& = 
				-
				\materialderivative\materialderivative \logDensity
				+
				\speed^2\delta^{ab}\p_a\p_b \logDensity
				+
				2\speed^{-1}\speed_{;\logDensity}(\materialderivative\logDensity)
				\materialderivative \logDensity-(\p_av^a)\materialderivative \logDensity
				\\
				&	\ \ 
				-\speed^{-1}\speed_{;\logDensity}(\gfour^{-1})^{\alpha\beta}\p_\alpha\logDensity\p_\beta \logDensity-\speed\speed_{;\Ent}\gradEnt^a
				\p_a\logDensity+3\speed^{-1}\speed_{;\Ent}(\materialderivative\Ent) \materialderivative \logDensity
				\\
				& 
				= 
				-
				\materialderivative\materialderivative \logDensity
				+
				\speed^2\delta^{ab}\p_a\p_b \logDensity
				+
				2\speed^{-1}\speed_{;\logDensity}(\materialderivative\logDensity)
				\materialderivative \logDensity
				+
				(\p_a v^a)
				\p_b v^b
				-
				(\p_av^a) \FR
				\\
				&	
				\ \ 
				-\speed^{-1}\speed_{;\logDensity}(\gfour^{-1})^{\alpha\beta}\p_\alpha\logDensity\p_\beta \logDensity
				-
				\speed
				\speed_{;\Ent}\gradEnt^a
				\p_a\logDensity
				+
				3\speed^{-1}\speed_{;\Ent} 
				(\materialderivative \logDensity) \FS.
			\end{split}
		\end{align}
		
		Comparing the right-hand sides of \eqref{E:MATERIALDERIVATIVEOFCONSTRAINTVELOCITYEVOLUTIONEQUATION}
		and \eqref{EQ:BoxLOGDENSITYPROOFSTEP}, we deduce, after some straightforward algebra, that
		takes into account the identity
		$ -
		2\speed^{-1}\speed_{;\logDensity}(\materialderivative\logDensity)
		\materialderivative \logDensity
		+
		2 \speed \speed_{;\logDensity} \delta^{ab} (\p_a \logDensity) \p_b \logDensity
		= 
		2\speed^{-1}\speed_{;\logDensity}
		(\gfour^{-1})^{\alpha\beta}\p_\alpha\logDensity\p_\beta \logDensity
		$, which follows from \eqref{E:INVERSEACOUSTICALMETRIC},
		we deduce:
		\begin{align} 
			\begin{split} \label{E:MATERIALDERIVATIVEOFCONSTRAINTVELOCITYEVOLUTIONEQUATIONEQUALTOBOXDENSITYPLUSSOMEERRORSPROOFSTEP}
				- \materialderivative 
				\FR
				& = \boxg \logDensity 
				+
				e^{\logDensity}\frac{\pressure_{;\Ent}}{\backgroundDensity}\DivGradEnt
				\\
				& \ \
				+
				3 \speed^{-1}\speed_{;\logDensity}(\gfour^{-1})^{\alpha\beta}\p_\alpha\logDensity\p_\beta \logDensity
				+
				\p_a v^b \p_b v^a
				-
				(\p_a v^a)
				\p_b v^b	
				\\
				& \ \
				+
				3 \speed \speed_{;\Ent} \gradEnt^a \p_a \logDensity
				+
				e^{-\logDensity}
				\frac{\pressure_{;\Ent;\logDensity}}
				{\backgroundDensity}
				\gradEnt^a
				\p_a
				\logDensity
				+
				e^{-\logDensity}\frac{p_{;\Ent;\Ent}}
				{\backgroundDensity}\delta_{ab}
				\gradEnt^a
				\gradEnt^b
				\\
				&	
				\ \ 
				-
				\partial_a \FV^a
				+
				(\p_a v^a) \FR
				-
				3\speed^{-1}\speed_{;\Ent} 
				(\materialderivative \logDensity) \FS.
			\end{split}
		\end{align}
		Next, we note that the chain rule yields the identity
		$
		\speed \speed_{;\Ent} 
		= \frac{1}{2} (\speed^2)_{;\Ent}
		= \frac{1}{2} (e^{-\logDensity}
		\frac{\pressure_{;\logDensity}}{\backgroundDensity})_{;\Ent}
		=
		\frac{1}{2} e^{-\logDensity}
		\frac{\pressure_{;\Ent;\logDensity}}{\backgroundDensity}
		$,
		and thus the first two products on the next-to-last line of
		RHS~\eqref{E:MATERIALDERIVATIVEOFCONSTRAINTVELOCITYEVOLUTIONEQUATIONEQUALTOBOXDENSITYPLUSSOMEERRORSPROOFSTEP} sum to
		$
		\frac{5}{2}
		e^{-\logDensity}
		\frac{\pressure_{;\Ent;\logDensity}}
		{\backgroundDensity}
		\gradEnt^a
		\p_a
		\logDensity
		$.
		Using this identity, we compare 
		\eqref{E:MATERIALDERIVATIVEOFCONSTRAINTVELOCITYEVOLUTIONEQUATIONEQUALTOBOXDENSITYPLUSSOMEERRORSPROOFSTEP}
		and \eqref{E:boxgdensity}
		to arrive at the desired identity \eqref{E:DIFFERENTIATEDDENSITYCONSTRAINEDCOMPAREDTOLUKSPECKWAVE}.
	\end{proof}

	\begin{proposition}[Decomposition of the wave operator under the null frame]\label{P:Propdecompboxgf} By \eqref{D:decompositionofmetric} and Definition \ref{D:Defconnectioncoeff}, we have the following decompositions for $\boxg \varphi$, where $\varphi$ is a scalar function.
		\begin{subequations}
			\begin{align}
				\label{E:boxgf}\boxg \varphi&=-e_4 e_3 \varphi+\anglap \varphi-\frac{1}{2}\gtr\upchi e_3 \varphi-\frac{1}{2}\gtr\uchi e_4 \varphi+(\upeta+\ueta)\cdot\angnabla \varphi+2\omega e_3 \varphi,\\
				\label{E:boxgf2}\boxg \varphi&=-e_3 e_4 \varphi+\anglap \varphi-\frac{1}{2}\gtr\upchi e_3 \varphi-\frac{1}{2}\gtr\uchi e_4 \varphi+(\upeta+\ueta)\cdot\angnabla \varphi+2\uomega e_4 \varphi.
			\end{align}
		\end{subequations}
	\end{proposition}
	\begin{proof}[Proof of Proposition \ref{P:Propdecompboxgf}]
		We use the fact that $\anglap \varphi=\angnabla_A\angnabla_A\varphi=\angnabla_A\left(\angnabla_A \varphi\right)-\angnabla_{\angnabla_Ae_A}\varphi$ and the decomposition of the metric $\gfour^{-1}$ in \eqref{D:decompositionofmetric} to obtain
		\begin{align}
			\boxg \varphi=\left(e_A^\alpha e_A^\beta-\frac{1}{2}e_4^\alpha e_3^\beta-\frac{1}{2} e_3^\alpha e_4^\beta\right)\Dfour_\alpha\Dfour_\beta \varphi,
		\end{align}
		where
		\begin{subequations}
			\begin{align}
				e_A^\alpha e_A^\alpha\Dfour_\alpha\Dfour_\beta \varphi&=\Dfour_A\left(\Dfour_A\varphi\right)-\Dfour_{\Dfour_Ae_A}\varphi=\anglap \varphi-\frac{1}{2}\gtr\upchi e_3 \varphi-\frac{1}{2}\gtr\uchi e_4 \varphi,\\
				-\frac{1}{2} e_4^\alpha e_3^\beta\Dfour_\alpha\Dfour_\beta \varphi&=-\frac{1}{2}\Dfour_ {e_4}(\Dfour_{3}\varphi)+\frac{1}{2}(\Dfour_ {4} e_3)\varphi=-\frac{1}{2} e_4 e_3 \varphi+\ueta\cdot\angnabla \varphi+\omega e_3 \varphi,\\
				-\frac{1}{2} e_4^\alpha e_3^\beta\Dfour_\alpha\Dfour_\beta \varphi&=-\frac{1}{2}\Dfour_{3}(\Dfour_{ 4}\varphi)+\frac{1}{2}(\Dfour_{3} e_4)\varphi=-\frac{1}{2} e_3 e_4 \varphi+\upeta\cdot\angnabla \varphi+\uomega e_4 \varphi\\
				\notag&=-\frac{1}{2} e_4 e_3 \varphi+\upeta\cdot\angnabla \varphi+\frac{1}{2}[ e_4, e_3]\varphi+\uomega e_4 \varphi\\
				\notag&=-\frac{1}{2} e_4 e_3 \varphi+\upeta\cdot\angnabla \varphi+\omega e_3 \varphi
				.\end{align}
		\end{subequations}
		Combining the above equations, we obtain (\ref{E:boxgf}). 
		(\ref{E:boxgf2}) can be obtained in a similar fashion.
	\end{proof}

	\begin{lemma}[Evolution equation for $e_3\logDensity$ on $\hypersur$, and for $e_4 \logDensity$ on $\hypersurbar$]
		\label{LE:Solvee3rho}
		For a ``smooth'' solution to the compressible Euler equations with initial data $\vvariables$ on $\hypersur$, 
		$e_3\logDensity$ must satisfy the following (schematic) evolution equation along $\hypersur$:
		\begin{align}\label{EQ:e4e3rho}
			e_4 e_3 \logDensity=&\gensmoothfunction(\vvariables,e_J^\alpha)\cdot \angp^2\logDensity+\gensmoothfunction(\vvariables,e_J^\alpha)\cdot(\angp,e_4)(\vvariables,\gradEnt^i,e_J^\alpha)\cdot(\angp,e_4)(\vvariables,\gradEnt^i,e_J^\alpha)\\
			\notag&+\gensmoothfunction(\vvariables,e_J^\alpha)\cdot \left\{(\angp,e_4)(\vvariables,\gradEnt^i,e_J^\alpha),e_3\logDensity\right\}\cdot e_3\logDensity,
		\end{align}
		where $e_J^\alpha$ is viewed as a component of vectorfield with $J=3,4$, and $\gensmoothfunction$ is the schematic notation of a smooth function (each $\gensmoothfunction$ can vary).
		Similarly, $e_4\logDensity$ must satisfy the following evolution equation along $\hypersurbar$:
		\begin{align}
			\label{EQ:e3e4rho}e_3 e_4 \logDensity=&\gensmoothfunction(\vvariables,e_J^\alpha)\cdot \angp^2\logDensity+\gensmoothfunction(\vvariables,e_J^\alpha)\cdot(\angp,e_3)(\vvariables,\gradEnt^i,e_J^\alpha)\cdot(\angp,e_3)(\vvariables,\gradEnt^i,e_J^\alpha)\\
			\notag&+\gensmoothfunction(\vvariables,e_J^\alpha)\cdot \left\{(\angp,e_3)(\vvariables,\gradEnt^i,e_J^\alpha),e_4\logDensity\right\}\cdot e_4\logDensity.
		\end{align}
	\end{lemma}
	
	\begin{proof}[Proof of Lemma \ref{LE:Solvee3rho}]
		We only prove (\ref{EQ:e4e3rho}); (\ref{EQ:e3e4rho}) can be proved in a similar fashion. First of all, by (\ref{E:boxgdensity}) and (\ref{E:boxgf}) for $e_3\logDensity$, we have:
		\begin{align}
			\label{E:reforboxgdensity}-e_4 e_3 \logDensity+\anglap\logDensity -\frac{1}{2}\gtr\upchi e_3 \logDensity-\frac{1}{2}\gtr\uchi e_4 \logDensity+(\upeta+\ueta)\cdot\angp\logDensity+2\omega e_3 \logDensity=-e^{\logDensity}\frac{\pressure_{;\Ent}}{\average{\density}}\DivGradEnt+\mathfrak{Q}_{(\logDensity)}+\mathfrak{L}_{(\logDensity)},
		\end{align}
		where $\upeta+\ueta=2(\ueta+\upzeta)$ and
		\begin{subequations}\label{D:Defconnect1}
			\begin{align}
				\label{E:omegainG}\omega=&-\frac{1}{4}\gfour(e_4e_3,e_4)-\frac{1}{4}\dfour{4}{4}e_3\variables,\\
				\label{E:chiinG}\gtr\upchi=&\angp_\alpha e_4^\alpha+\frac{1}{2}(\gsphere^{-1})^{\alpha\beta}\angdfour_{\alpha\beta}e_4\variables,\\
				\label{E:uchiinG}\gtr\uchi=&\angp_\alpha e_3^\alpha+\frac{1}{2}(\gsphere^{-1})^{\alpha\beta}\angdfour_{\alpha\beta}e_3\variables,\\
				\label{E:uetainG}\ueta=&\frac{1}{2}\gfour(e_4e_3,\cdot)+\frac{1}{4}\left(\angdfour_{\cdot4}e_3\variables+\angdfour_{\cdot3}e_4\variables-\dfour{4}{3}\angp\variables\right),\\
				\label{E:upzetainG}\upzeta=&\frac{1}{2}\gfour(\angp e_4,e_3)+\frac{1}{4}\left(\angdfour_{\cdot3}e_4\variables+\dfour{3}{4}\angp\variables-\angdfour_{\cdot4}e_3\variables\right).
			\end{align}
		\end{subequations}
		By Definition \ref{D:modifiedfluidvariables} and equation \eqref{E:BS}, using Lemma \ref{L:derivativesdecomp}, we can rewrite $\DivGradEnt$ schematically as follows:
		\begin{align}\label{EQ:SchematicDivGradEnt}
			\DivGradEnt=&\gensmoothfunction(\vvariables,e_J^\alpha)\cdot(e_4,\angp)(\gradEnt,v)+\gensmoothfunction(\vvariables,e_J^\alpha)\cdot e_3v.
		\end{align}
		By (\ref{E:Bv}) and Lemma \ref{L:derivativesdecomp}, we write $e_3v^i$ as the linear combination of $e_3\logDensity,e_4(\logDensity,v^i),\angp(\logDensity,v^i),\p\Ent$, that is, \eqref{E:e3v}.

		Combining (\ref{E:reforboxgdensity}) and (\ref{E:e3v}), we obtain the desired result (\ref{EQ:e4e3rho}).
	\end{proof}

	\begin{corollary}\label{CO:e4e3k}
		By taking $e_3^{k-1}$-derivative of (\ref{EQ:e4e3rho}) and $e_4^{k-1}$-derivative of (\ref{EQ:e3e4rho}) with $k\geq1$, we have the following evolution equations for $e_3^k\logDensity$ and $e_4^k\logDensity$.
		\begin{subequations}\label{EQ:e4e3krho}
			\begin{align}
				\label{EQ:e4e3krho1.1}	e_4e_3^k\logDensity=&\gensmoothfunction\left\{(e_4,\angp)^2e_3^{k-1}(v,\logDensity,\Ent),(e_4,\angp)e_3^{k-1}e_J^\alpha,e_3^k\logDensity\right\}+\mathcal{E}_{0,k},\\
				\label{EQ:e4e3krho1.2}	e_3e_4^k\logDensity=&\gensmoothfunction\left\{(e_3,\angp)^2e_4^{k-1}(v,\logDensity,\Ent),(e_3,\angp)e_4^{k-1}e_J^\alpha,e_4^k\logDensity\right\}+\underline{\mathcal{E}}_{0,k},
			\end{align}
		\end{subequations}
		where if $k\geq 2$, then $\gensmoothfunction$ is linear in $e_3^k\logDensity$ in \eqref{EQ:e4e3krho1.1}, and linear in $e_4^k\logDensity$ in \eqref{EQ:e4e3krho1.2}. In \eqref{EQ:e4e3krho}, the lower order terms $\mathcal{E}_{0,k},\underline{\mathcal{E}}_{0,k}$ can be expressed schematically as follows:
		\begin{subequations}
			\begin{align}
				\mathcal{E}_{0,k}=&\sum\limits_{\substack{k_i\leq k\\l<2}}\gensmoothfunction\left\{e_3^{k_1}\logDensity,(e_4,\angp)^{l}e_3^{k_2-1}(v,\logDensity,\Ent),(e_4,\angp)^{l-1}e_3^{k_4-1}e_J^\alpha\right\},\\
				\underline{\mathcal{E}}_{0,k}=&\sum\limits_{\substack{k_i\leq k\\l<2}}\gensmoothfunction\left\{e_4^{k_1}\logDensity,(e_3,\angp)^{l}e_4^{k_2-1}(v,\logDensity,\Ent),(e_3,\angp)^{l-1}e_4^{k_4-1}e_J^\alpha\right\}.
			\end{align}
		\end{subequations}
	\end{corollary}
	
	\begin{corollary}\label{CO:e4eAl}
		By taking $(e_4,\angp)^l$-derivative of (\ref{EQ:e4e3krho1.1}) and $(e_3,\angp)^l$-derivative of (\ref{EQ:e4e3krho1.2}) with $k\geq1$, we have the following evolution equations for $e_3^k\logDensity$ and $e_4^k\logDensity$.
		\begin{subequations}\label{EQ:e4eAlrho}
			\begin{align}
				\label{EQ:e4eAlrho1.1}	e_4(e_4,\angp)^le_3^k\logDensity=&\gensmoothfunction\left\{(e_4,\angp)^{l+2}e_3^{k-1}(v,\logDensity,\Ent),(e_4,\angp)^{l+1}e_3^{k-1}e_J^\alpha,(e_4,\angp)^le_3^k\logDensity\right\}+\mathcal{E}_{l,k},\\
				\label{EQ:e4eAlrho1.2}	e_3(e_3,\angp)^le_4^k\logDensity=&\gensmoothfunction\left\{(e_3,\angp)^{l+2}e_4^{k-1}(v,\logDensity,\Ent),(e_3,\angp)^{l+1}e_4^{k-1}e_J^\alpha,(e_3,\angp)^le_4^k\logDensity\right\}+\underline{\mathcal{E}}_{l,k},
			\end{align}
		\end{subequations}
		where if $k\geq 2$, then $\gensmoothfunction$ is linear in $(e_4,\angp)^le_3^k\logDensity$ in \eqref{EQ:e4eAlrho1.1}, and linear in $(e_3,\angp)^le_4^k\logDensity$ in \eqref{EQ:e4eAlrho1.2}. In \eqref{EQ:e4eAlrho}, the lower order terms $\mathcal{E}_{l,k},\underline{\mathcal{E}}_{l,k}$ can be expressed schematically as follows:
		\begin{subequations}
			\begin{align}
				\mathcal{E}_{l,k}=&\sum\limits_{\substack{k_i\leq k\\l_i<2+l}}\gensmoothfunction\left\{(e_4,\angp)^{l_1}e_3^{k_1}\logDensity,(e_4,\angp)^{l_2}e_3^{k_2-1}(v,\logDensity,\Ent),(e_4,\angp)^{l_3-1}e_3^{k_3-1}e_J^\alpha\right\},\\
				\underline{\mathcal{E}}_{l,k}=&\sum\limits_{\substack{k_i\leq k\\l_i<2+l}}\gensmoothfunction\left\{(e_3,\angp)^{l_1}e_4^{k_1}\logDensity,(e_3,\angp)^{l_2}e_4^{k_2-1}(v,\logDensity,\Ent),(e_3,\angp)^{l_3-1}e_4^{k_3-1}e_J^\alpha\right\}.
			\end{align}
		\end{subequations}
	\end{corollary}	
	\begin{proof}[Discussion of the proof of Corollary \ref{CO:e4e3k} and Corollary \ref{CO:e4eAl}]
		Corollary \ref{CO:e4e3k} and Corollary \ref{CO:e4eAl} are direct results of commuting $(e_4,\angp)^le_3^k$ ($(e_3,\angp)^le_4^k$ for proving \eqref{EQ:e4eAlrho1.2}) with (\ref{EQ:e4e3rho}) (\eqref{EQ:e3e4rho}  for proving \eqref{EQ:e4eAlrho1.2}) and Lemma \ref{LE:CommutatorFormula}. We note that by Lemma \ref{LE:CommutatorFormula}, all commutator terms are of lower-order. We omit the details of the proof.
	\end{proof}

	\section{Characteristic initial data, Characteristic Embeddings, and the Initial Geometry}\label{S:CID}
	In this section, we rigorously define what it means to have characteristic initial data for the compressible Euler equations,
	and we construct various geometric quantities out of the data.
	Roughly, the initial data comprise an ``initial'' pair of transversally intersecting null 
	(with respect to the acoustical metric $\gfour$)
	hypersurfaces $\inihypersur,\inihypersurbar$ embedded into spacetime $\mathbb{R}^{1+3}$, 
	as well as fluid data along $\inihypersur,\inihypersurbar$, which we denote by
	$(\mathring{v},\inilogDensity,\inient)$.
	A crucial point is that the null hypersurfaces and fluid data cannot be
	prescribed arbitrarily, but rather they must be specified in such a way that
	$\inihypersur,\inihypersurbar$
	are null as measured by $\gfour(\mathring{v},\inilogDensity,\inient)$.
	Another less obvious -- though equally important -- consideration is that the compressible
	Euler equations imply PDEs along each of $\inihypersur$, $\inihypersurbar$
	where all differentiations are tangential. This means that any admissible data 
	$(\mathring{v},\inilogDensity,\inient)$
	must satisfy these PDEs, which we refer to as the \emph{characteristic constraints}.

	\subsection{The characteristic Initial Data and Characteristic Embeddings}
	\label{SS:initialembedding2}
	In this section, we state some of our assumptions on the embeddings of the characteristic data surfaces into spacetime
	and on the fluid data.
	
	\begin{itemize}
		\item (\textbf{Initial embeddings}) 
		Throughout, $\inisphere$ denotes a point in $\mathbb{S}^2$ and $(\inisphere^1,\inisphere^2)$
		denote its coordinates in a local chart.
		
		We assume that there exist numbers $\ubar_* > 0$ and $u_* > 0$
		and maps 
		\begin{subequations}
			\begin{align}
				\embedding :[0,\ubar_*] \times \mathbb{S}^2 & \rightarrow \mathbb{R}^{1+3},
				&
				(\iniubar,\inisphere^1,\inisphere^2) & \rightarrow \embedding(\iniubar,\inisphere^1,\inisphere^2),
				\label{E:INITIALHEMBEDDING} \\
				\uembedding : [0,u_*] \times \mathbb{S}^2 & \rightarrow \mathbb{R}^{1+3},
				& 
				(\iniu,\inisphere^1,\inisphere^2) &\rightarrow \uembedding(\iniu,\inisphere^1,\inisphere^2),
				\label{E:INITIALHBAREMBEDDING}
			\end{align}
		\end{subequations}
		such that $\embedding$ and $\uembedding$ are $C^{\infty}$ embeddings, i.e.,  
		injective and with full rank differential (the rank is $3$).
		We denote the image sets in $\mathbb{R}^{1+3}$ as 
		\begin{align} \label{E:INITIALNULLHYPERSURFACEIMAGESETSINR1PLUS3}
			\inihypersur
			& := \embedding([0,u_*] \times \mathbb{S}^2),
			&
			\inihypersurbar 
			&
			:= \uembedding([0,\ubar_*] \times \mathbb{S}^2).
		\end{align}
		Thus, $\inihypersur,\inihypersurbar$ are $C^{\infty}$ manifolds-with-boundary embedded in 
		spacetime $\mathbb{R}^{1+3}$, each diffeomorphic to a compact interval crossed with $\mathbb{S}^2$.
		We assume that $\inihypersur$ intersects $\inihypersurbar$ transversally in a set that is diffeomorphic to 
		$\mathbb{S}^2$. We denote the intersection by $S_{0,0} := \inihypersur \cap \inihypersurbar$.
		More precisely, we assume that 
		$\inihypersur \cap \inihypersurbar = \uembedding(\lbrace 0 \rbrace \times \mathbb{S}^2)
		= \embedding(\lbrace 0 \rbrace \times \mathbb{S}^2)$,
		that for every point $\inisphere \in \mathbb{S}^2$ we have
		\begin{align} \label{E:EMBEDDINSPHERECOMPATIBILITYCONDITION}
			\embedding(0,\inisphere^1,\inisphere^2) & = \uembedding(0,\inisphere^1,\inisphere^2),
		\end{align}
		and that the (Euclidean) normals to $\inihypersur$ and $\inihypersurbar$ are \emph{not parallel}
		at any point in $S_{0,0}$.
		Later on, we will make precise assumptions on various norms of the embeddings.
		\item  (\textbf{Initial fluid data})
		We assume that: 
		\begin{align} \label{E:DATAONH}
			\inivvariables := (\inivel^1,\inivel^2,\inivel^3,\inilogDensity,\inient),
		\end{align}
		are $C^{\infty}$ fluid data
		on $\inihypersur$, that: 
		\begin{align} \label{E:DATAONHBAR}
			\inivvariablesbar & := (\inivelbar^1,\inivelbar^2,\inivelbar^3,\inilogDensitybar,\inientbar)
		\end{align}
		are $C^{\infty}$ fluid data on $\inihypersurbar$,
		and that: 
		\begin{align} \label{E:ORDER0compatibilityCONDITION}
			\inivvariables \restriction_{S_{0,0}}
			& = \inivvariablesbar \restriction_{S_{0,0}}.
		\end{align}
		\eqref{E:ORDER0compatibilityCONDITION} is a compatibility condition enforcing
		that the data functions agree on their ``common domain'' $S_{0,0}$.
		In Def.\,\ref{D:CORNERCOMPATIBILITYCONDITIONSONDATA}, we make higher-order compatibility assumptions
		on the data. Later on, we will make precise assumptions on various norms of the data.
		We often view $\inivvariables$ as a function of $(\iniubar,\inisphere^1,\inisphere^2)$ 
		defined on
		the domain $[0,\ubar_*] \times \mathbb{S}^2$, and 
		$\inivvariablesbar$
		as a function of $(\iniu,\inisphere^1,\inisphere^2)$ defined on
		the domain $[0,u_*] \times \mathbb{S}^2$.
		
		\begin{notation}[Tensors along $\inihypersur \cup \inihypersurbar$]	
			\label{N:TENSORSONDATAHYPERSURFACES}
			For tensors $Q$ that are functions of $\vvariables$ 
			(i.e., $Q = Q(\vvariables)$),
			we use the notation $\mathring{Q}$ to denote 
			either $Q(\inivvariables)$ or $Q(\inivvariablesbar)$,
			i.e., the restriction of $Q$ to  $\inihypersur \cup \inihypersurbar$.
			The precise meaning will always be clear from context.
		\end{notation}
		\item (\textbf{$\inihypersur$ and $\inihypersurbar$ are $\gfour$-null}).
		Recall that $\gfour$ is the acoustical metric, and that its rectangular components $\gfour_{\alpha \beta}$
		are functions of $\vvariables$, i.e.,  $\gfour_{\alpha \beta} = \gfour_{\alpha \beta}(\vvariables)$.
		On $\inihypersur$, we assume that:
		\begin{subequations} \label{E:FIRSTINITIALEMBEDDINGISNULL}
			\begin{align}
				\inigfour_{\alpha \beta}
				\frac{\p\embedding^{\alpha}}{\p\iniubar}
				\frac{\p\embedding^\beta}{\p\iniubar}
				& =0,	
				& \frac{\p\embedding^0}{\p\iniubar} 
				& =1,\\
				\inigfour_{\alpha \beta} 
				\frac{\p\embedding^{\alpha}}{\p \iniubar}
				\frac{\p\embedding^\beta}{\p\inicartiasphere}
				& = 0,
				& A &=1,2.
				\label{E:FIRSTINITIALEMBEDDINGNULLDIRECTIONORTHOGONALTOANGULARDIRECTIONS}
			\end{align}
		\end{subequations}
		\eqref{E:FIRSTINITIALEMBEDDINGISNULL}-\eqref{E:FIRSTINITIALEMBEDDINGNULLDIRECTIONORTHOGONALTOANGULARDIRECTIONS} imply that
		$\inihypersur$ is a $\inigfour$-null hypersurface ruled by integral curves of the future-directed null generator 
		$\frac{\p\embedding^{\alpha}}{\p\iniubar}$, which is normalized by $\frac{\p\embedding^0}{\p\iniubar} = 1$.
		In particular, for $A=1,2$, because $\inigfour_{\alpha \beta}$ is Lorentzian, one can show that
		$\inigfour_{\alpha \beta} 
		\frac{\p\embedding^{\alpha}}{\p \inicartiasphere}
		\frac{\p\embedding^\beta}{\p \inicartiasphere}
		> 0
		$,
		i.e., $\frac{\p\embedding^{\alpha}}{\p \inicartiasphere}$ is $\inigfour$-spacelike.
		Note that if $\frac{\p\embedding^0}{\p\iniubar}$ were assumed only to be non-vanishing instead of equal to unity, then
		by solving a transport equation along the integral curves of $\frac{\p\embedding^{\alpha}}{\p\iniubar}$,
		one could construct a new coordinate, also denoted $\iniubar$ (by abuse of notation),
		such that with respect to it, the condition $\frac{\p\embedding^0}{\p\iniubar} = 1$ holds.
		
		Similarly, on  $\inihypersurbar$, we assume that:
		\begin{subequations} \label{E:SECONDINITIALEMBEDDINGISNULL}
			\begin{align}
				\inigfour_{\alpha \beta}
				\frac{\p\uembedding^{\alpha}}{\p \iniu}
				\frac{\p\uembedding^\beta}{\p \iniu}
				& =0,
				& \frac{\p\uembedding^0}{\p \iniu} 
				& = 1,
				\\
				\inigfour_{\alpha \beta} 
				\frac{\p\uembedding^{\alpha}}{\p \iniu}
				\frac{\p\uembedding^\beta}{\p\inicartiasphere}
				& = 0,
				& A&=1,2.
				\label{E:SECONDINITIALEMBEDDINGNULLDIRECTIONORTHOGONALTOANGULARDIRECTIONS}
			\end{align}
		\end{subequations}
		\eqref{E:SECONDINITIALEMBEDDINGISNULL}-\eqref{E:SECONDINITIALEMBEDDINGNULLDIRECTIONORTHOGONALTOANGULARDIRECTIONS} imply that
		$\inihypersurbar$ is a $\inigfour$-null hypersurface ruled by integral curves of the future-directed null generator 
		$\frac{\p\embedding^{\alpha}}{\p\iniu}$. As above,
		for $A=1,2$, $\frac{\p\uembedding^{\alpha}}{\p\inicartiasphere}$ is $\inigfour$-spacelike.
		
		\item (\textbf{The Euler constraint equations are satisfied on $\inihypersur\bigcup\inihypersurbar$}).	On $\inihypersur$, we assume that:
		\begin{align}\label{EQ:EulerConstraintsinCoordinates1}
			\frac{\p}{\p\iniubar}\inilogDensity+\sum\limits_{i=1}^3\gensmoothfunction_i(\frac{\p\embedding}{\p\iniubar},\frac{\p\embedding}{\p\inisphere^A},\inivel)\frac{\p}{\p\iniubar}\inivel^i=\gensmoothfunction\left(\frac{\p\embedding}{\p\iniubar},\frac{\p\embedding}{\p\inisphere^A},\frac{\p}{\p\inisphere^A}(\inivel,\inilogDensity,\inient)\right),
		\end{align}
		where the precise equation is provided in \eqref{E:constraint1}. \eqref{EQ:EulerConstraintsinCoordinates1} is the projection of the first-order Euler equations \eqref{E:firstorderformulation} onto $\inihypersur$. Notice that \eqref{EQ:EulerConstraintsinCoordinates1} only contains $\inihypersur$-tangential derivatives of the fluid variables.
		Similarly, on $\inihypersurbar$, we assume that:
		\begin{align}\label{EQ:EulerConstraintsinCoordinates2}
			\frac{\p}{\p\iniu}\inilogDensity+\sum\limits_{i=1}^3\gensmoothfunction_i(\frac{\p\uembedding}{\p\iniu},\frac{\p\uembedding}{\p\inisphere^A},\inivel)\frac{\p}{\p\iniu}\inivel^i=\gensmoothfunction\left(\p\uembedding,\frac{\p}{\p\inisphere^A}(\inivel,\inilogDensity,\inient)\right),
		\end{align}
		where the precise equation is provided in \eqref{E:constraint2}. \eqref{EQ:EulerConstraintsinCoordinates2} is the projection of the first-order Euler equations \eqref{E:firstorderformulation} onto $\inihypersurbar$. Notice that \eqref{EQ:EulerConstraintsinCoordinates2} only contains $\inihypersurbar$-tangential derivatives of the fluid variables.
	\end{itemize}

	\begin{remark}[Surfaces of constant $\iniu > 0$ or constant $\iniubar > 0$ are not generally null]
		\label{R:NOTGENERALLYNULL}
		Although $\inihypersur = 
		\embed(\lbrace 0 \rbrace \times [0,u_*] \times \mathbb{S}^2)$ is a $\inigfour$-null hypersurface,
		for the solutions we will construct, for fixed $\iniu > 0$,
		$\embed(\lbrace \iniu \rbrace \times [0,u_*] \times \mathbb{S}^2)$ will \emph{not} generally be 
		a $\gfour$-null hypersurface. Analogous remarks apply for
		$\inihypersurbar = \embed([0,\ubar_*] \times \lbrace 0 \rbrace \times \mathbb{S}^2)$
		and the surfaces $\embed([0,\ubar_*] \times \lbrace \iniubar \rbrace \times \mathbb{S}^2)$
		of constant $\iniubar > 0$.
		Instead, our $\gfour$-null surfaces will be the level sets of the eikonal functions, i.e., the solutions to
		\eqref{EQ:eikonalequations}.
	\end{remark}
	
	\begin{remark}
		In the flat case with a single characteristic hypersurface, where $\embed(\lbrace 0 \rbrace \times\lbrace u \rbrace \times \mathbb{S}^2)\subset\Sigma_t$, the null frame condition \eqref{E:FIRSTINITIALEMBEDDINGISNULL} and the Euler constraint equation \eqref{EQ:EulerConstraintsinCoordinates1} can be solved; see \cite{WangYuYu25}. 
	\end{remark}

	\subsection{Geometry along the Initial Characteristic Hypersurfaces}
	\label{SS:GEOMETRYONDATAHYPERSURFACES}
	In this section, we define and construct various geometric quantities along the initial null hypersurfaces.
	
	\subsubsection{Null vectorfields along the initial characteristic hypersurfaces}
	\label{SSS:NULLVECTORFIELDSONDATAHYPERSURFACES}
	We define the vectorfields $\iniLunit$ and $\uiniuLunit$ respectively as follows:
	\begin{subequations}
		\begin{align}
			\iniLunit^{\alpha}&:=\frac{\p\embedding^{\alpha}}{\p\iniubar},&
			\text{on}\ \inihypersur,
			\label{E:INITIALLUNITONH} 
			\\
			\uiniuLunit^{\alpha}&:=\frac{\p\uembedding^{\alpha}}{\p \iniu},&
			\text{on}\ \inihypersurbar.
			\label{E:INITIALULUNITONHBAR}
		\end{align}
	\end{subequations}
	\eqref{E:FIRSTINITIALEMBEDDINGISNULL}
	implies that	
	$\iniLunit$ is a null generator of $\inihypersur$, i.e., 
	$\inigfour_{\alpha \beta} \iniLunit^{\alpha} \iniLunit^{\beta} = 0$ along 
	$\inihypersur$, $\iniLunit$ is tangent to $\inihypersur$, 
	$\inigfour$-normal to $\inihypersur$,
	and $\inihypersur$ is ruled by integral curves of $\iniLunit$.
	Similarly, \eqref{E:SECONDINITIALEMBEDDINGISNULL}-\eqref{E:SECONDINITIALEMBEDDINGNULLDIRECTIONORTHOGONALTOANGULARDIRECTIONS}
	imply that $\uiniuLunit$ is a null generator of $\inihypersurbar$.
	
	In our subsequent constructions, we will use the vectorfields $\mathring{\mathcal{T}}^{\alpha}$ and 
	$\mathring{\mathcal{S}}^{\alpha}$, which are defined to be
	the $\inigfour$-duals of the following one-forms,
	where we recall our notation \ref{N:TENSORSONDATAHYPERSURFACES}:
	\begin{subequations}
		\begin{align}
			\mathring{\mathcal{S}}_\alpha
			&:=\antisymmetic_{\alpha\beta\gamma\delta}\inimaterialderivative^\beta\frac{\p\embedding^\gamma}{\p\inisphere^1}\frac{\p\embedding^\delta}{\p\inisphere^2},&
			\mathring{\mathcal{T}}_\alpha
			&:=\frac{-\antisymmetic_{\alpha\beta\gamma\delta}\mathring{\mathcal{S}}^\beta\frac{\p\embedding^\gamma}{\p\inisphere^1}\frac{\p\embedding^\delta}{\p\inisphere^2}}{\antisymmetic_{\alpha\beta\gamma\delta}\inimaterialderivative^\alpha\mathring{\mathcal{S}}^\beta\frac{\p\embedding^\gamma}{\p\inisphere^1}\frac{\p\embedding^\delta}{\p\inisphere^2}},&
			\text{on}\ \inihypersur,
			\label{E:INITIALSANDTVECTORFIELDSNEEDEDFORNULLFRAME1}\\
			\mathring{\mathcal{S}}_\alpha
			&:=\antisymmetic_{\alpha\beta\gamma\delta}\inimaterialderivative^\beta\frac{\p\uembedding^\gamma}{\p\inisphere^1}\frac{\p\uembedding^\delta}{\p\inisphere^2},&
			\mathring{\mathcal{T}}_\alpha
			&:=\frac{-\antisymmetic_{\alpha\beta\gamma\delta}\mathring{\mathcal{S}}^\beta\frac{\p\uembedding^\gamma}{\p\inisphere^1}\frac{\p\uembedding^\delta}{\p\inisphere^2}}{\antisymmetic_{\alpha\beta\gamma\delta}\inimaterialderivative^\alpha\mathring{\mathcal{S}}^\beta\frac{\p\uembedding^\gamma}{\p\inisphere^1}\frac{\p\uembedding^\delta}{\p\inisphere^2}},&
			\text{on}\ \inihypersurbar.
			\label{E:INITIALSANDTVECTORFIELDSNEEDEDFORNULLFRAME2}
		\end{align}
	\end{subequations}
	Above, $\antisymmetic_{\alpha\beta\gamma\delta}$ is the fully antisymmetric symbol normalized by 
	$\antisymmetic_{0123} = 1$.

	\begin{lemma}[Properties of $\iniuLunit^{\alpha}$,
		$\uiniLunit^{\alpha}$, 
		$\mathring{\mathcal{S}}^{\alpha}$, 
		and $\mathring{\mathcal{T}}^{\alpha}$]
		\label{L:PROPERTIESOFSOMEINITIALVECTORFIELDS}
		The following results hold, 
		where for vectorfields $V$, we denote $\inigfabs{V}:=\sqrt{\abs{\inigfour_{\alpha\beta}V^\alpha V^\beta}}$.
		\ \\
		\begin{itemize}
			\item The vectorfield $\mathring{\mathcal{S}}^{\alpha} := (\inigfourinverse)^{\alpha \beta} \mathring{\mathcal{S}}_{\beta}$
			is $\Sigma_t$-tangent (i.e., $\mathring{\mathcal{S}}^0 = 0$),
			$\inigfour$-spacelike (i.e., $\inigfour(\mathring{\mathcal{S}},\mathring{\mathcal{S}}) > 0$),
			$\inigfour$-orthogonal to $\frac{\p\embedding^{\alpha}}{\p\inisphere^1}$
			and $\frac{\p\embedding^{\alpha}}{\p\inisphere^2}$ on  $\inihypersur$,
			and
			$\inigfour$-orthogonal to $\frac{\p\uembedding^{\alpha}}{\p\inisphere^1}$
			and $\frac{\p\uembedding^{\alpha}}{\p\inisphere^2}$ on  $\inihypersurbar$.
			\item The vectorfield $\mathring{\mathcal{T}}^{\alpha} := (\inigfourinverse)^{\alpha \beta} \mathring{\mathcal{S}}_{\beta}$
			is non-vanishing, $\inigfour$-orthogonal to $\mathring{\mathcal{S}}$, 
			$\inigfour$-timelike (i.e., $\inigfour(\mathring{\mathcal{T}},\mathring{\mathcal{T}}) < 0$),
			and satisfies $\inimaterialderivative_{\alpha} \mathring{\mathcal{T}}^{\alpha} = -1$,
			i.e., $\mathring{\mathcal{T}}^0 = 1$.
			\item The vectorfield $\iniLunit^{\alpha}$ defined in \eqref{E:INITIALLUNITONH}
			satisfies the following identity:
			\begin{align} \label{E:IDENTITYFORINITIALLONH}
				\iniLunit^{\alpha}
				& 
				= 
				\mathring{\mathcal{T}}^{\alpha}
				+ 
				\mypmH
				\frac{\inigfabs{\mathring{\mathcal{T}}}}{\inigfabs{\mathring{\mathcal{S}}}}
				\mathring{\mathcal{S}}^{\alpha},
				&&
				\text{on}\ \inihypersur,
			\end{align}
			where $\mypmH$ is a constant always equal to $1$ or always equal to $-1$,
			depending on whether $\mathring{\mathcal{S}}^{\alpha}$ is inwards-pointing or outwards pointing.
			
			Similarly, with $\mypmH$ as above,
			the vectorfield $\uiniuLunit^{\alpha}$ defined in \eqref{E:INITIALULUNITONHBAR}
			satisfies the following identity:
			\begin{align} \label{E:IDENTITYFORINITIALLONHBAR}
				\uiniuLunit^{\alpha}
				& 
				= 
				\mathring{\mathcal{T}}^{\alpha}
				- 
				\mypmH
				\frac{\inigfabs{\mathring{\mathcal{T}}}}{\inigfabs{\mathring{\mathcal{S}}}}
				\mathring{\mathcal{S}}^{\alpha},
				&&
				\text{on} \ \inihypersurbar,
			\end{align}
			\item With $\mypmH$ as above, the vectorfield $\iniuLunit^{\alpha}$ defined by:
			\begin{align} \label{E:INITIALULUNITONH}
				\iniuLunit^{\alpha}
				&
				:= 
				\mathring{\mathcal{T}}^{\alpha}
				- 
				\mypmH
				\frac{\inigfabs{\mathring{\mathcal{T}}}}{\inigfabs{\mathring{\mathcal{S}}}}
				\mathring{\mathcal{S}}^{\alpha},
				&
				&\text{on}\ \inihypersur,
			\end{align}
			is $\inigfour$-null (i.e., $\inigfour(\iniuLunit,\iniuLunit) = 0$),
			$\inigfour$-orthogonal to the ``angular'' vectorfields $\frac{\p\embedding^{\alpha}}{\p\inisphere^1}$
			and $\frac{\p\embedding^{\alpha}}{\p\inisphere^2}$,
			normalized by $\iniuLunit^0 = 1$, and linearly independent from $\iniLunit^{\alpha}$.
			\item 
			Similarly, with $\mypmH$ as above, 
			the vectorfield $\uiniLunit^{\alpha}$ defined by:
			\begin{align} \label{E:INITIALLUNITONHBAR}
				\uiniLunit^{\alpha}
				&
				:= 
				\mathring{\mathcal{T}}^{\alpha}
				+ 
				\mypmH
				\frac{\inigfabs{\mathring{\mathcal{T}}}}{\inigfabs{\mathring{\mathcal{S}}}}
				\mathring{\mathcal{S}}^{\alpha},
				&
				&\text{on}\ \inihypersurbar,
			\end{align}
			is $\inigfour$-null,
			$\inigfour$-orthogonal to the ``angular'' vectorfields $\frac{\p\uembedding^{\alpha}}{\p\inisphere^1}$
			and $\frac{\p\uembedding^{\alpha}}{\p\inisphere^2}$,
			normalized by $\uiniLunit^0 = 1$, and linearly independent from $\uiniuLunit^{\alpha}$.
			\item The following identities hold:
			\begin{subequations}
				\begin{align} \label{E:LSAGREEONINITIALSPHERE}
					\iniLunit^{\alpha} & = \uiniLunit^{\alpha}
					&
					& \text{on}\ \inihypersur \cap \inihypersurbar,
					\\
					\uiniuLunit^{\alpha} & = \iniuLunit^{\alpha}
					&
					& \text{on}\ \inihypersur \cap \inihypersurbar.
					\label{E:LBARSAGREEONINITIALSPHERE}
				\end{align}
			\end{subequations}
		\end{itemize}
	\end{lemma}
	
	\begin{proof}
		\eqref{E:INITIALSANDTVECTORFIELDSNEEDEDFORNULLFRAME1}--\eqref{E:INITIALSANDTVECTORFIELDSNEEDEDFORNULLFRAME2} 
		and the antisymmetry of $\antisymmetic_{\alpha\beta\gamma\delta}$
		imply that
		$\mathring{\mathcal{S}}_{\alpha} \inimaterialderivative^{\alpha} = 0$,
		and since $\materialderivative$ is $\gfour$-normal to $\Sigma_t$,
		we conclude that $\mathring{\mathcal{S}}$ is $\Sigma_t$-tangent.
		$\mathring{\mathcal{S}}$ cannot vanish since
		$
		\inimaterialderivative^{\alpha}
		$,
		$
		\frac{\p\embedding^{\alpha}}{\p\inisphere^1}
		$,
		and
		$
		\frac{\p\embedding^{\alpha}}{\p\inisphere^2}
		$
		are a linearly independent triplet of vectors
		($\inimaterialderivative^{\alpha}$ is $\inigfour$-timelike while the latter two vectors are linearly independent and spacelike).
		The fact that $\mathring{\mathcal{S}}^{\alpha}$
		is $\inigfour$-orthogonal to 
		$\frac{\p\embedding^{\alpha}}{\p\inisphere^1}$
		and $\frac{\p\embedding^{\alpha}}{\p\inisphere^2}$ on 
		$\inihypersur$
		and
		$\inigfour$-orthogonal to $\frac{\p\uembedding^{\alpha}}{\p\inisphere^1}$
		and $\frac{\p\uembedding^{\alpha}}{\p\inisphere^2}$ on  $\inihypersurbar$			
		follows from \eqref{E:INITIALSANDTVECTORFIELDSNEEDEDFORNULLFRAME1}--\eqref{E:INITIALSANDTVECTORFIELDSNEEDEDFORNULLFRAME2}  
		and the antisymmetry of $\antisymmetic_{\alpha\beta\gamma\delta}$.
		
		The properties of $\mathring{\mathcal{T}}^{\alpha}$ follow from similar considerations,
		where we note that
		\eqref{E:INITIALSANDTVECTORFIELDSNEEDEDFORNULLFRAME1}--\eqref{E:INITIALSANDTVECTORFIELDSNEEDEDFORNULLFRAME2} 
		imply that $\mathring{\mathcal{T}}^{\alpha}$ is $\inigfour$-orthogonal to three linearly independent $\inigfour$-spacelike
		vectors and $\inimaterialderivative^{\alpha} \mathring{\mathcal{T}}_{\alpha} = -1$,
		which implies in particular that $\mathring{\mathcal{T}}^{\alpha}$ is $\inigfour$-timelike.
		
		To prove \eqref{E:IDENTITYFORINITIALLONH}, we note that both LHS~\eqref{E:IDENTITYFORINITIALLONH}
		and RHS~\eqref{E:IDENTITYFORINITIALLONH} are 
		$\inigfour$-orthogonal to $\frac{\p\embedding^{\alpha}}{\p\inisphere^1}$
		and $\frac{\p\embedding^{\alpha}}{\p\inisphere^2}$ on  $\inihypersur$
		and have $0$ component equal to $1$. 
		Since $\mathring{\mathcal{T}}^0=1$, and  $\mathring{\mathcal{S}}^0=0$, 
		this implies that $\iniLunit^{\alpha}$ is a linear combination of 
		$\mathring{\mathcal{T}}^{\alpha}$ and $\mathring{\mathcal{S}}^{\alpha}$,
		where the coefficient of $\mathring{\mathcal{T}}^{\alpha}$ 
		in the combination is unity. This yields \eqref{E:IDENTITYFORINITIALLONH},
		up to understanding the coefficient $\mypmH$.
		Next, we note that LHS~\eqref{E:IDENTITYFORINITIALLONH} is
		$\inigfour$-null by \eqref{E:FIRSTINITIALEMBEDDINGISNULL},
		and a simple computation based on the
		$\inigfour$-orthogonality of $\mathring{\mathcal{T}}^{\alpha}$ and $\mathring{\mathcal{S}}^{\alpha}$,
		the $\inigfour$-timelike nature of $\mathring{\mathcal{T}}^{\alpha}$,
		and the $\inigfour$-spacelike nature of $\mathring{\mathcal{S}}^{\alpha}$,
		implies that RHS~\eqref{E:IDENTITYFORINITIALLONH} is
		$\inigfour$-null precisely if $\mypmH = - 1$ or $\mypmH = 1$.
		One of these choices of $\mypmH$ must yield the identity \eqref{E:IDENTITYFORINITIALLONH},
		and the choice with the opposite sign yields the vectorfield $\iniuLunit^{\alpha}$ 
		defined in \eqref{E:INITIALULUNITONH} and its claimed properties.
		
		The identity \eqref{E:IDENTITYFORINITIALLONHBAR} as well as 
		the claimed properties surrounding
		\eqref{E:IDENTITYFORINITIALLONHBAR} and \eqref{E:INITIALLUNITONHBAR}
		follow from similar arguments, and we omit the details.	
		
		To prove \eqref{E:LSAGREEONINITIALSPHERE}, we note that in view of the assumptions stated in Sect.\,\ref{SS:initialembedding2}, 
		both vectorfields in \eqref{E:LSAGREEONINITIALSPHERE}
		satisfy the properties stated just below equation \eqref{E:INITIALLUNITONHBAR},
		and these properties suffice to uniquely determine a vectorfield.
		\eqref{E:LBARSAGREEONINITIALSPHERE} follows from a similar argument.
	\end{proof}

	\subsection{Adapted Coordinates in the Future of the Data Null Hypersurfaces}
	\label{SS:ADAPTEDCOORDINATES}
	Recall that in Sect.\,\ref{SS:initialembedding2}, we assumed the existence of coordinates 
	$(\iniubar,\inisphere^1,\inisphere^2)$ on $\inihypersur$ and coordinates
	$(\iniu,\inisphere^1,\inisphere^2)$ on $\inihypersurbar$. 
	To carry out the ``soft'' (i.e., non-quantitative) part of our proof of local existence and uniqueness, 
	we will extend $(\iniu,\iniubar,\inisphere^1,\inisphere^2)$ to coordinates defined
	on a small future-neighborhood of $\inihypersur\bigcup\inihypersurbar$, and by an abuse of notation, we will
	denote the extended coordinates by the same notation 
	$(\iniu,\iniubar,\inisphere^1,\inisphere^2)$.
	
	To construct the adapted coordinates, 
	we will foliate a future-neighborhood of $\inihypersur \bigcup \inihypersurbar$ with two families of transversal hypersurfaces 
	such that each of $\inihypersur$ and $\inihypersurbar$ is a member of one of the families.
	In total, this will provide us with coordinates $(\iniu,\iniubar,\inisphere^1,\inisphere^2)$
	on a small future-neighborhood of $\inihypersur\bigcup\inihypersurbar$ 
	such that the two hypersurface families are level sets of $\iniu$ and $\iniubar$ respectively, such that
	$\inihypersur$ is equal to the level set $\lbrace \iniu = 0 \rbrace$ 
	and $\inihypersurbar$ is equal to the level set $\lbrace \iniubar = 0 \rbrace$.
	
	\begin{definition}[Adapted coordinates]
		\label{D:ADAPTEDCOORDINATES}
		We refer to the coordinates $(\iniu,\iniubar,\inisphere^1,\inisphere^2)$,
		which we again emphasize will be defined 
		on a future-neighborhood of $\inihypersur\bigcup\inihypersurbar$,
		as \emph{adapted coordinates}.
		We denote the corresponding partial derivative vectorfields
		by $\lbrace \frac{\partial}{\partial \iniu}, \frac{\partial}{\partial \iniubar}, \frac{\partial}{\partial \inisphere^1}, 
		\frac{\partial}{\partial \inisphere^2}
		\rbrace$. 
		Moreover, for any function $f$,
		\begin{align} \label{E:GENERICADAPTEDCOORDINATEPARTIALDERIVATIVES}
			\inip f
			& := \left\lbrace \frac{\partial}{\partial \iniu} f, \frac{\partial}{\partial \iniubar} f, \frac{\partial}{\partial \inisphere^1 }f, 
			\frac{\partial}{\partial \inisphere^2 }f
			\right\rbrace
		\end{align}
		denotes the array of partial derivatives of $f$ with respect to the adapted coordinates.
		
	\end{definition}
	
	\begin{remark}[Schematic notation for angular partial derivatives]
		\label{R:SCHEMATICFORANGULARPARTIALDERIVATIVES}
		When the precise details do not matter, we use the notation $\frac{\partial}{\partial \inisphere}$ 
		to denote either $\frac{\partial}{\partial \inisphere^1}$ or $\frac{\partial}{\partial \inisphere^2}$.
	\end{remark}
	
	\begin{remark}[$\iniu$ and $\iniubar$ are not eikonal functions]
		The above remarks imply that 
		the zero level sets of $\iniu$ and $\iniubar$ are $\gfour$-null hypersurfaces,
		but the same is not generally true for their other level sets. That is, 
		$\iniu$ and $\iniubar$ are not solutions of the acoustic eikonal equation.
	\end{remark}
	
	To initiate the construction of the adapted coordinates, we fix a $C^{\infty}$ cutoff function $\psi : \mathbb{R} \rightarrow [0,\infty)$ such that:
	\begin{subequations}
		\begin{align}
			\psi(y)=& 1,
			& y \leq 1,
			\\
			\psi(y)=& 0,
			& y \geq 2.
		\end{align}
	\end{subequations}
	
	In what follows, $\tilde{\embedkappa}\geq\embedkappa>0$ denote 
	small constants that depend on the fluid data and null hypersurface embeddings
	$\embedding$ and $\uembedding$;  $\tilde{\embedkappa}$ and $\embedkappa$ will be determined later. 
	
	Recall that $\uiniuLunit^{\alpha}(\iniu,\inisphere^1,\inisphere^2)
	:= \frac{\p\uembedding^{\alpha}}{\p\iniu}(\iniu,\inisphere^1,\inisphere^2)$ 
	is the null vectorfield on $\inihypersurbar$ defined in \eqref{E:INITIALULUNITONHBAR}
	and that $\iniuLunit^{\alpha}(\iniubar,\inisphere^1,\inisphere^2)$ 
	is the null vectorfield on $\inihypersur$ defined in \eqref{E:INITIALULUNITONH}.
	Recall that $\uiniuLunit^{\alpha}(\iniu,\inisphere^1,\inisphere^2)$ is everywhere tangent to $\inihypersurbar$ and transversal to 
	$\inihypersur$ at $S_{0,0}$, i.e., when $\iniu=0$. Recall also that
	$\iniuLunit^{\alpha}(\iniubar,\inisphere^1,\inisphere^2)$ is everywhere transversal to $\inihypersur$.

	For each point $p=\embedding(\iniubar,\inisphere^2,\inisphere^2)$ on $\inihypersur$, 
	we define the $\inihypersur$-transversal curve  
	$\iniLbargeo_{\iniubar,\inisphere^2,\inisphere^2} : [0,\tilde{\embedkappa}] \rightarrow \mathbb{R}^{1+3}$,
	$\iniu \rightarrow \iniLbargeo_{\iniubar,\inisphere^2,\inisphere^2}(\iniu)$
	emanating from $p$
	to be the solution to the following ODE initial value problem:
	\begin{subequations}\label{DE:iniUpsilonbar}
		\begin{align}
			\frac{d}{d \iniu} \iniLbargeo^{\alpha}_{\iniubar,\inisphere^2,\inisphere^2}(\iniu)=&\begin{cases}
				\uiniuLunit^{\alpha}(\iniu,\inisphere^1,\inisphere^2),
				& 
				\iniu\in[0,\tilde{\embedkappa}),
				\\
				\psi(\iniubar/ \tilde{\embedkappa})
				\left\{
				\uiniuLunit^{\alpha}(\iniu,\inisphere^1,\inisphere^2)
				-
				\iniuLunit^{\alpha}(\iniubar,\inisphere^1,\inisphere^2)
				\right\}
				+
				\iniuLunit^\alpha(\iniubar,\inisphere^1,\inisphere^2),
				& 
				\iniubar\in [\tilde{\embedkappa},2\tilde{\embedkappa}],
				\\
				\iniuLunit^{\alpha}(\iniubar,\inisphere^1,\inisphere^2),
				& 
				\iniubar \in (2 \tilde{\embedkappa},\ubar_*].
			\end{cases} \\
			\iniLbargeo^\alpha_{\iniubar,\inisphere^2,\inisphere^2}(0)=&\embedding^{\alpha}(\iniubar,\inisphere^2,\inisphere^2).
		\end{align}
	\end{subequations} 
	Note that by ODE uniqueness,
	\eqref{E:INITIALULUNITONHBAR},
	\eqref{DE:iniUpsilonbar},
	and \eqref{E:EMBEDDINSPHERECOMPATIBILITYCONDITION}, 
	we have:
	\begin{align} \label{E:DATAHBAREMBEDDINGAGREESWITHFLOW}
		\iniLbargeo^{\alpha}_{0,\inisphere^2,\inisphere^2}(\iniu) = \uembedding^{\alpha}(\iniu,\inisphere^1,\inisphere^2).
	\end{align}

	Similarly, with $\iniLunit^{\alpha}(\iniubar,\inisphere^1,\inisphere^2):=\frac{\p\embedding^{\alpha}}{\p\iniubar}(\iniubar,\inisphere^1,\inisphere^2)$
	denoting the null vectorfield on $\inihypersur$ defined in \eqref{E:INITIALLUNITONH}
	and 
	$\uiniuLunit^{\alpha}(\iniu,\inisphere^1,\inisphere^2)$ denoting
	the null vectorfield on $\inihypersurbar$ defined in \eqref{E:IDENTITYFORINITIALLONHBAR},
	for each point $q=\uembedding(\iniu,\inisphere^2,\inisphere^2)$ on $\inihypersurbar$, we define the $\inihypersurbar$-transversal curve 
	$\iniLgeo_{\iniu,\inisphere^2,\inisphere^2} : [0,\tilde{\embedkappa}] \rightarrow \mathbb{R}^{1+3}$,
	$\iniubar \rightarrow \iniLgeo_{\iniu,\inisphere^2,\inisphere^2}(\iniubar)$ emanating from $q$ to be the solution to the following
	ODE initial value problem:
	\begin{subequations}\label{DE:iniUpsilon}
		\begin{align}
			\frac{d}{d \iniubar} \iniLgeo^{\alpha}_{\iniu,\inisphere^2,\inisphere^2}(\iniubar)
			& = 
			\begin{cases}
				\iniLunit^{\alpha}(\iniubar,\inisphere^1,\inisphere^2), 
				& 
				\iniu\in[0,\tilde{\embedkappa}),
				\\
				\psi(\iniu/\embedkappa)
				\left\{
				\iniLunit^{\alpha}(\iniubar,\inisphere^1,\inisphere^2)
				-
				\uiniLunit^\alpha(\iniu,\inisphere^1,\inisphere^2)
				\right\}
				+
				\uiniLunit^\alpha(\iniu,\inisphere^1,\inisphere^2),
				&
				\iniu\in[\tilde{\embedkappa},2\tilde{\embedkappa}],
				\\
				\uiniLunit^\alpha(\iniu,\inisphere^1,\inisphere^2),
				& 
				\iniu\in(2\tilde{\embedkappa},u_*].
			\end{cases} \\
			\iniLgeo^\alpha_{\iniu,\inisphere^2,\inisphere^2}(0)=&\uembedding^{\alpha}(\iniu,\inisphere^2,\inisphere^2).
		\end{align}
	\end{subequations} 
	Note that by ODE uniqueness,
	\eqref{E:INITIALLUNITONH},
	\eqref{DE:iniUpsilon},
	and \eqref{E:EMBEDDINSPHERECOMPATIBILITYCONDITION}, 
	we have: 
	\begin{align} \label{E:DATAHEMBEDDINGAGREESWITHFLOW}
		\iniLgeo^{\alpha}_{0,\inisphere^2,\inisphere^2}(\iniubar) = \embedding^{\alpha}(\iniubar,\inisphere^1,\inisphere^2).
	\end{align}
	
	Also observe that \eqref{E:EMBEDDINSPHERECOMPATIBILITYCONDITION} and the fundamental theorem of calculus imply that
	for any $(\iniu,\iniubar,\inisphere^1,\inisphere^2)\in[0,\tilde{\embedkappa}]\times[0,\tilde{\embedkappa}]\times\Stwo$, 
	we have:
	\begin{align}\label{EQ:DangerofMultiDefining}
		\iniLbargeo_{\iniubar,\inisphere^2,\inisphere^2}(\iniu)
		& = 
		\iniLgeo_{\iniu,\inisphere^2,\inisphere^2}(\iniubar)
		=
		\embedding^{\alpha}(0,\inisphere^1,\inisphere^2)
		+
		\int_{0}^{\iniu}\frac{\p\uembedding^{\alpha}}{\p\iniu}(\iniu^\prime,\inisphere^1,\inisphere^2)\diff\iniu^\prime
		+
		\int_{0}^{\iniubar}\frac{\p\embedding^{\alpha}}{\p\iniubar}(\iniubar^\prime,\inisphere^1,\inisphere^2)\diff\iniubar^\prime.
	\end{align}

	\begin{definition}[A diffeomorphism constructed out of the initial embeddings]
		\label{D:LDIFFEOCONSTRUCTEDOUTOFEMBEDDINGS}
		Recall the definitions of $\iniLgeo,\iniLbargeo$ in \eqref{DE:iniUpsilonbar}-\eqref{DE:iniUpsilon} and the identity 
		\eqref{EQ:DangerofMultiDefining}. 
		We define the domain:
		\begin{align} \label{E:DIFFEOMORPHISMDOMAIN}
			\embeddomain_{\tilde{\embedkappa},\ubar_*,u_*}
			:=
			([0,\ubar_*] \times [0,\tilde{\embedkappa}] \times \mathbb{S}^2)
			\bigcup
			([0,\tilde{\embedkappa}] \times [0,u_*] \times \mathbb{S}^2),
		\end{align}
		and we define the map
		$
		\embed
		:
		\embeddomain_{\tilde{\embedkappa},\ubar_*,u_*}
		\rightarrow \mathbb{R}^{1+3}$ as follows:
		\begin{align} \label{E:LOCALDIFFEOCONSTRUCTEDOUTOFEMBEDDINGS}
			\embed(\iniu,\iniubar,\inisphere^1,\inisphere^2)
			& 
			:=\begin{cases} 
				\iniLgeo_{\iniu,\inisphere^2,\inisphere^2}(\iniubar), & (\iniu,\iniubar,\inisphere^1,\inisphere^2)\in [0,u_*] \times[0,\tilde{\embedkappa}] \times \mathbb{S}^2, \\
				\iniLbargeo_{\iniubar,\inisphere^2,\inisphere^2}(\iniu), & (\iniu,\iniubar,\inisphere^1,\inisphere^2)
				\in [0,\tilde{\embedkappa}] \times[0,u_*] \times \mathbb{S}^2.
			\end{cases}
		\end{align}
	\end{definition}
	In this paper, $\embed$ will play the role of the change of variables map from adapted coordinates to geometric coordinates.
	In particular, we have:
	\begin{align} \label{E:EMBEDDINGASCHOV}
		\embed(\iniu,\iniubar,\inisphere^1,\inisphere^2)	
		& = \left(
		x^0(\iniu,\iniubar,\inisphere^1,\inisphere^2),
		x^1(\iniu,\iniubar,\inisphere^1,\inisphere^2),
		x^2(\iniu,\iniubar,\inisphere^1,\inisphere^2),
		x^3(\iniu,\iniubar,\inisphere^1,\inisphere^2) 
		\right),
	\end{align}
	where on RHS~\eqref{E:EMBEDDINGASCHOV}, we are viewing the Cartesian coordinates as functions of the adapted coordinates.
	
	At first glance at RHS~\eqref{E:LOCALDIFFEOCONSTRUCTEDOUTOFEMBEDDINGS}, one might think that $\embed$
	is not well-defined on 
	$(\iniu,\iniubar,\inisphere^1,\inisphere^2) \in [0,\tilde{\embedkappa}] \times[0,\tilde{\embedkappa}] \times \mathbb{S}^2$
	due to multi-valuedness. However, \eqref{EQ:DangerofMultiDefining} guarantees that $\embed$ is well-defined.
	
\begin{figure}
		\centering
		\includegraphics[width=10cm]{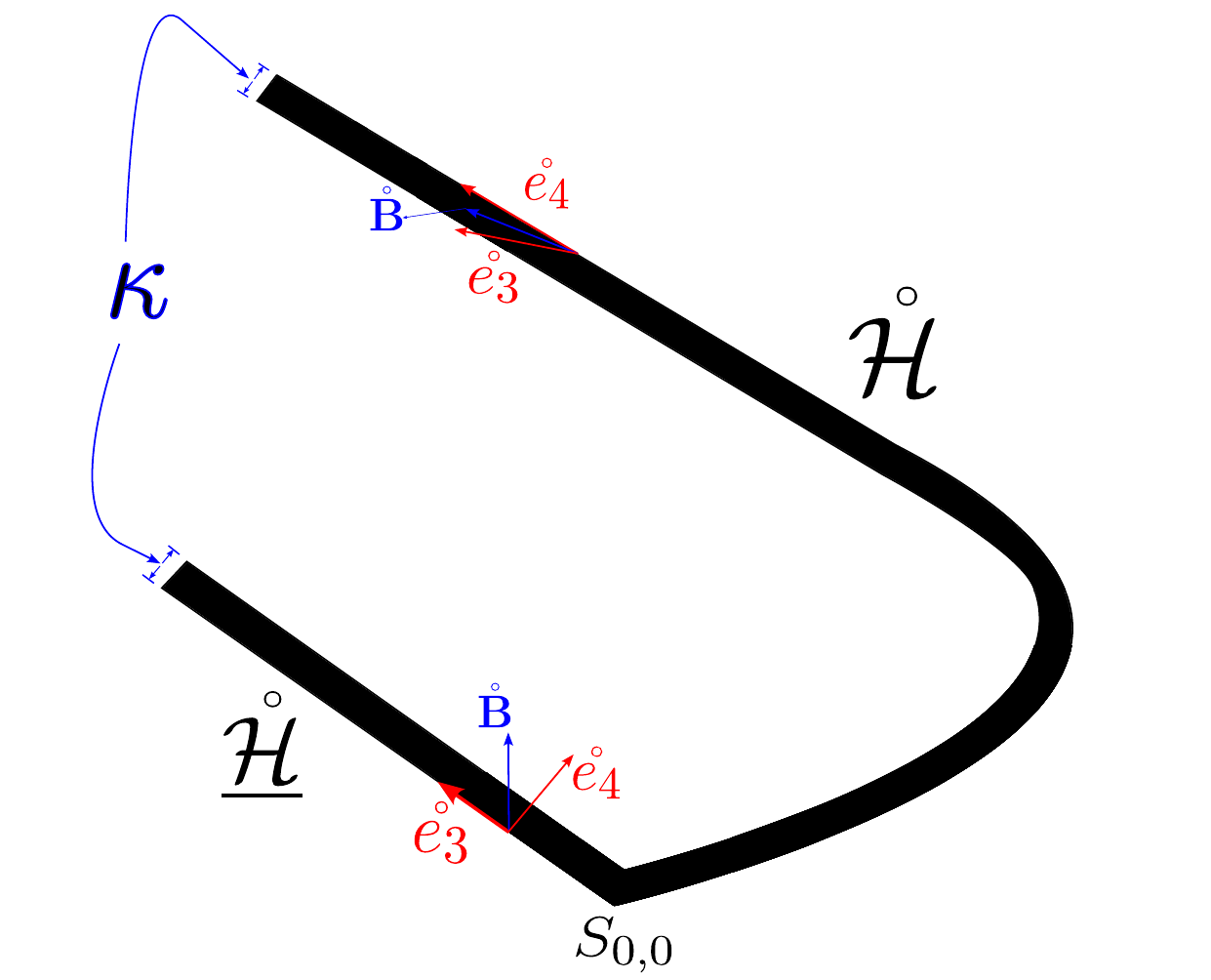}
		\caption{A representative image of $\embed$ from $\embeddomain_{\embedkappa,\ubar_*,u_*}$ for small $\embedkappa > 0$}
		\label{fig2}
	\end{figure}

	Since $\inihypersur$ intersects $\inihypersurbar$ transversally by assumption, $\embed$ is a diffeomorphism whenever
	$\tilde{\embedkappa}$ is sufficiently small; we give a proof of this fact in the next lemma.
	See also Fig.\,\ref{fig2}.
	
	\begin{lemma}[$\embed$ is a diffeomorphism on $\embeddomain_{\embedkappa,\ubar_*,u_*}$ for small $\embedkappa > 0$]
		\label{L:EXISTENCEOFALOCALDIFFEOMORPHISM}
		There is a number $\embedkappa > 0$ such that the map $\embed$ from Def.\,\ref{D:LDIFFEOCONSTRUCTEDOUTOFEMBEDDINGS}
		is a diffeomorphism from the domain $\embeddomain_{\embedkappa,\ubar_*,u_*}$
		(see definition~\eqref{E:DIFFEOMORPHISMDOMAIN})
		onto its image.
		
		Moreover, the following identities hold: 
		\begin{subequations}
			\begin{align} \label{E:EMBEDDINGRESPECTSINITIALHBAR}
				\embed(0,\iniubar,\inisphere^1,\inisphere^2) 
				& = \embedding(\iniubar,\inisphere^1,\inisphere^2),
				&& (\iniubar,\inisphere^1,\inisphere^2) \in [0,\ubar_*] \times \mathbb{S}^2 
				\\
				\embed(\iniu,0,\inisphere^1,\inisphere^2) & = \uembedding(\iniu,\inisphere^1,\inisphere^2),
				&& (\iniu,\inisphere^1,\inisphere^2) \in [0,u_*] \times \mathbb{S}^2.
				\label{E:EMBEDDINGRESPECTSINITIALH}
			\end{align}	
		\end{subequations}

	\end{lemma}
	
	\begin{proof}	
		\eqref{E:EMBEDDINGRESPECTSINITIALHBAR}--\eqref{E:EMBEDDINGRESPECTSINITIALH} follow from
		\eqref{E:DATAHBAREMBEDDINGAGREESWITHFLOW} and \eqref{E:DATAHEMBEDDINGAGREESWITHFLOW}.
		
		Using the definition of $\embed$ given in \eqref{E:LOCALDIFFEOCONSTRUCTEDOUTOFEMBEDDINGS}, 
		\eqref{E:INITIALLUNITONH}--\eqref{E:INITIALULUNITONHBAR}, 
		and \eqref{DE:iniUpsilonbar}--\eqref{E:DATAHEMBEDDINGAGREESWITHFLOW}, 
		we compute the following identities:
		\\
		
		\noindent\textbf{On $\inihypersur = \embed(\lbrace 0 \rbrace \times [0,\ubar_*]  \times \mathbb{S}^2)$:}
		\begin{subequations}\label{EQ:JacobianofEonH}
			\begin{align}
				\frac{\p\embed^\alpha}{\p \iniu}
				(0,\iniubar,\inisphere^1,\inisphere^2) 
				& =
				\begin{cases}
					\uiniuLunit^{\alpha}(0,\inisphere^1,\inisphere^2), & \iniubar \in [0,\tilde{\embedkappa}),
					\\
					\psi(\iniubar/\tilde{\embedkappa}) 
					\left\{
					\uiniuLunit^{\alpha}(0,\inisphere^1,\inisphere^2)
					-
					\iniuLunit^\alpha(\iniubar,\inisphere^1,\inisphere^2)
					\right\}
					+
					\iniuLunit^\alpha(\iniubar,\inisphere^1,\inisphere^2),&\iniubar\in[\tilde{\embedkappa},2\tilde{\embedkappa}],
					\\
					\iniuLunit^\alpha(\iniubar,\inisphere^1,\inisphere^2),&\iniubar\in(2\tilde{\embedkappa},\ubar_*],
				\end{cases}
				\\
				\frac{\p\embed^\alpha}{\p \iniubar}(0,\iniubar,\inisphere^1,\inisphere^2)
				& = \iniLunit^{\alpha}(\iniubar,\inisphere^1,\inisphere^2),
				\\
				\frac{\p\embed^\alpha}{\p \inisphere^A}
				(0,\iniubar,\inisphere^1,\inisphere^2)
				& = \frac{\p\embedding^{\alpha}}{\p \inisphere^A}(\iniubar,\inisphere^1,\inisphere^2),
			\end{align}
		\end{subequations}
		\textbf{On $\inihypersurbar = \embed([0,u_*] \times \lbrace 0 \rbrace \times \mathbb{S}^2)$:}
		\begin{subequations}\label{EQ:JacobianofEonHbar}
			\begin{align}
				\frac{\p\embed^\alpha}{\p \iniubar}
				(\iniu,0,\inisphere^1,\inisphere^2)
				& = 
				\begin{cases}
					\iniLunit^{\alpha}(0,\inisphere^1,\inisphere^2),&\iniu\in[0,\tilde{\embedkappa}]
					\\
					\psi(\iniu/\tilde{\embedkappa})
					\left\{
					\iniLunit^{\alpha}(0,\inisphere^1,\inisphere^2)
					-
					\uiniLunit^\alpha(\iniu,\inisphere^1,\inisphere^2)
					\right\}
					+
					\uiniLunit^\alpha(\iniu,\inisphere^1,\inisphere^2),&\iniu\in[\tilde{\embedkappa},2\tilde{\embedkappa}],
					\\
					\uiniLunit^\alpha(\iniu,\inisphere^1,\inisphere^2),&\iniu\in(2\tilde{\embedkappa},\ubar_*],
				\end{cases}\\
				\frac{\p\embed^\alpha}{\p \iniu}
				(\iniu,0,\inisphere^1,\inisphere^2)
				& = \uiniuLunit^{\alpha}(\iniu,\inisphere^1,\inisphere^2),\\
				\frac{\p\embed^\alpha}{\p \inisphere^A}
				(\iniu,0,\inisphere^1,\inisphere^2)
				& = \frac{\p\uembedding^{\alpha}}{\p \inisphere^A}(\iniu,\inisphere^1,\inisphere^2).
			\end{align}
		\end{subequations}
		From \eqref{EQ:JacobianofEonH}--\eqref{EQ:JacobianofEonHbar},
		the assumptions on the embeddings $\embedding^{\alpha}$ and $\uembedding^{\alpha}$ stated in
		Sect.\,\ref{SS:initialembedding2} 
		(which imply in particular that 
		$\lbrace \iniLunit^{\alpha}(0,0,\theta), \uiniuLunit^{\alpha}(0,0,\theta) \rbrace$ is linearly independent for all
		$\theta \in \mathbb{S}^2$), 
		and the results of Lemma~\ref{L:PROPERTIESOFSOMEINITIALVECTORFIELDS},
		including the linear independence of 
		$\lbrace \iniLunit^{\alpha}, \iniuLunit^{\alpha} \rbrace$ 
		and linear independence of $\lbrace \uiniLunit^{\alpha}, \uiniuLunit^{\alpha} \rbrace$,
		the identities \eqref{E:LSAGREEONINITIALSPHERE}--\eqref{E:LBARSAGREEONINITIALSPHERE},
		and the inverse function theorem (see, e.g. \cite{GuilleminPollack74}),
		it follows that if $\tilde{\embedkappa} > 0$ is sufficiently small, then on the compact set
		$\embeddomain_{\tilde{\embedkappa},\ubar_*,u_*}
		:=
		[0,\ubar_*] \times [0,\tilde{\embedkappa}] \times \mathbb{S}^2
		\bigcup
		[0,\embedkappa] \times [0,u_*] \times \mathbb{S}^2$,
		$\embed(\iniu,\iniubar,\inisphere^1,\inisphere^2)$ 
		is a smooth injective map with full rank differential (the rank is $4$), i.e., a 
		diffeomorphism from $\embeddomain_{\tilde{\embedkappa},\ubar_*,u_*}$ onto its image.
		We fix such a choice of $\tilde{\embedkappa}$ and denote it by $\embedkappa$, which
		finishes the proof of the lemma.
	\end{proof}
	
	\subsection{Some First-order Compressible Euler Equations in Adapted Coordinates}
	\label{SS:SOMEFIRSTORDERCOMPRESSIBLEEULERINADAPATEDCOORDIANTES}

	\begin{lemma}[Some compressible Euler equations in adapted coordinates]
		\label{E:SOMECOMPRESSIBLEEULEREQUATIONSINADPATEDCOORDINATES}
		Let $E_3$ denote the vectorfield on the domain $\embeddomain_{\tilde{\embedkappa},\ubar_*,u_*}$
		from Lemma~\ref{L:TRANSVERALANDINITIALLYNULLCOMMUTATORSFORRENDALLMETHOD}, and let $\vec{E}_3$ denote
		its Cartesian components.
		In adapted coordinates, in a neighborhood of $\inihypersur$, the following hold:
		
		The velocity compressible Euler equation \eqref{E:Bv} 
		(which by \eqref{E:CONSTRAINTVELOCITYEVOLUTIONEQUATION} is equivalent to $\FV^i = 0$)
		can be expressed as follows:
		\begin{align} \label{E:e3vADAPTEDCOORDINATES}
			E_3 v^i
			& =
			\gensmoothfunction_{\star}(\inip \embed,\vec{E}_4,\vvariables)
			E_3 \logDensity
			+
			\gensmoothfunction
			\left(\inip \embed,\vec{E}_4,\vvariables,\frac{\p}{\p \iniubar} \vvariables, \frac{\p}{\p \inisphere} \vvariables,  
			\frac{\partial}{\partial \iniu} \Ent \right).
		\end{align}
		
		In particular, with $\FV^i$ as defined in \eqref{E:CONSTRAINTVELOCITYEVOLUTIONEQUATION},
		we have:
		\begin{align} \label{E:E3VINADPATEDCOORDINATESISAMULTIPLEOFVELOCITYCONSTRAINT}
			\mbox{LHS~\eqref{E:e3vADAPTEDCOORDINATES}} 
			-
			\mbox{RHS~\eqref{E:e3vADAPTEDCOORDINATES}}
			& = 
			\gensmoothfunction_{\star}(\inip \embed,\vec{E}_4,\vvariables)
			\FV^i.
		\end{align}
		
		The entropy compressible Euler equation \eqref{E:Bs} 
		can be expressed in any of the following forms:
		\begin{subequations}
			\begin{align} 
				\label{E:ENTROPYEULEREQUATIONADAPTEDCOORDINATES}
				\frac{\partial}{\partial \iniu} \Ent
				& 
				= 
				\gensmoothfunction_{\star}(\inip \embed,v^i)
				\cdot
				\frac{\partial}{\partial \iniubar} \Ent
				+
				\gensmoothfunction^A(\inip \embed,v^i)
				\cdot
				\frac{\partial}{\partial \inisphere^A} \Ent,
				\\
				\frac{\partial}{\partial \iniubar} \Ent
				& 
				= 
				\gensmoothfunction_{\star}(\inip \embed,v^i)
				\cdot
				\frac{\partial}{\partial \iniu} \Ent
				+
				\gensmoothfunction^A(\inip \embed,v^i)
				\cdot
				\frac{\partial}{\partial \inisphere^A} \Ent,
				\label{E:ALTERNATEENTROPYEULEREQUATIONADAPTEDCOORDINATES}
				\\
				E_3 \Ent
				& 
				= 
				\gensmoothfunction_{\star}(\inip \embed,\vec{E}_4,v^i)
				\cdot
				\frac{\partial}{\partial \iniubar} \Ent
				+
				\gensmoothfunction^A(\inip \embed,\vec{E}_4,v^i)
				\cdot
				\frac{\partial}{\partial \inisphere^A} \Ent.
				\label{E:ENTROPYEULEREQUATIONWITHE3DERIVATIVE}
			\end{align}
		\end{subequations}
		
		In particular, with $\FS$ as defined in \eqref{E:CONSTRAINTENTROPYEVOLUTIONEQUATION},
		we have:
		\begin{align} \label{E:ENTROPYEULEREQUATIONWITHE3DERIVATIVEISAMULTIPLEOFENTROPYCONSTRAINT}
			\mbox{LHS~\eqref{E:ENTROPYEULEREQUATIONWITHE3DERIVATIVE}} 
			-
			\mbox{RHS~\eqref{E:ENTROPYEULEREQUATIONWITHE3DERIVATIVE}}
			& = 
			\gensmoothfunction_{\star}(\inip \embed,\vec{E}_4,\vvariables)
			\FS.
		\end{align}

	\end{lemma}

	\subsubsection{Null lapses and null frames along the initial characteristic hypersurfaces}
	\label{SSS:LAPSESANDNULLFRAMESONDATAHYPERSURFACES}
	In this section, we define various geometric quantities on the 
	initial characteristic hypersurfaces and use them to construct a $\gfour$-null frame along them.
	
	\begin{definition}[Null lapses along the initial characteristic hypersurfaces]
		\label{D:NULLLAPSESDATA}
		Let $\iniLunit$, $\iniuLunit$, $\uiniLunit$, and $\uiniuLunit$ be the $\gfour$-null vectorfields
		from \eqref{E:INITIALLUNITONH}--\eqref{E:INITIALULUNITONHBAR}
		and Lemma~\ref{L:PROPERTIESOFSOMEINITIALVECTORFIELDS}.
		We define the initial null lapse $\mathring{\nullapse}$ as follows:
		\begin{subequations}
			\begin{align}
				\mathring{\nullapse}:=&-\frac{1}{2}\inigfour(\iniLunit,\iniuLunit),&
				\text{on } &\inihypersur,\\
				\mathring{\nullapse}:=&-\frac{1}{2}\inigfour(\uiniLunit,\uiniuLunit),&
				\text{on } &\inihypersurbar.
			\end{align}	
		\end{subequations}
	\end{definition}	
	
	Since $\iniLunit,\iniuLunit,\uiniLunit,\uiniuLunit$ are future-directed null vectors, it follows that there exists 
	a constant $\boundinilapse \geq 1$ (depending on the initial data and initial embeddings) such that:
	\begin{align} \label{E:INITIALLAPSEBOUNDS}
		\frac{1}{\boundinilapse^2}
		\leq
		\mathring{\nullapse}
		\leq \boundinilapse^2.
	\end{align}
	
	\begin{definition}[Null frames along the initial characteristic hypersurfaces]	
		\label{D:NULLFRAMEDATA}
		We define the initial null frame $\{\mathring{e_1},\mathring{e_2},\mathring{e}_3,\mathring{e}_4\}$ as follows:
		
		\textbf{On $\inihypersur$:}
		\begin{subequations}
			\begin{align}
				\mathring{e}_A&:=\frac{\p\embedding}{\p\cartiasphere},&
				A&=1,2.\\
				\label{D:e3e42}\mathring{e}_3&:=\mathring{\nullapse}^{-1/2}\iniuLunit,&
				\mathring{e}_4&:=\mathring{\nullapse}^{-1/2}\iniLunit. 
			\end{align}
		\end{subequations}
		
		\textbf{On $\inihypersurbar$:}
		\begin{subequations}
			\begin{align}
				\mathring{e}_A&:=\frac{\p\uembedding}{\p \cartiasphere},&
				A&=1,2.\\
				\mathring{e}_3&:=\mathring{\nullapse}^{-1/2}\uiniuLunit,&
				\mathring{e}_4&:=\mathring{\nullapse}^{-1/2}\uiniLunit.
			\end{align}
		\end{subequations}
	\end{definition}
	
	By construction, the following relations hold:
	\begin{subequations}\label{D:key}
		\begin{align}
			\inigfour(\mathring{e}_4,\mathring{e}_4)=&0,&\inigfour(\mathring{e}_3,\mathring{e}_3)=&0,&\inigfour(\mathring{e}_3,\mathring{e}_4)=&-2,\\
			\inigfour(\mathring{e}_A,\mathring{e}_B):=&\inigsphere_{AB},& \inigfour(\mathring{e}_4,\mathring{e}_A)=&0,&\inigfour(\mathring{e}_3,\mathring{e}_A)=&0,
		\end{align}\
	\end{subequations}
	where we clarify that $\inigsphere_{AB}$ is defined by the above equation.
	
	By construction, we have $\iniLunit^0 = \iniuLunit^0 = \uiniLunit^0 = \uiniuLunit^0 = 1$.
	It follows that:
	\begin{align} \label{E:E40ANDE30ARENITIALLYUNITY}
		\mathring{e}_4^0 = \mathring{e}_3^0 = \mathring{\nullapse}^{-1/2}, &&  \mbox{on } \inihypersur \cup \inihypersurbar,
	\end{align}
	which justifies the equalities in the following definition.
	
	\begin{definition}[The function $\inilapse$]
		\label{D:lapsequantities}
		We define $\inilapse$ to be the following scalar function:
		\begin{align} \label{D:lapse}
			\inilapse
			&:= 
			\mathring{e}_4^0
			=
			\mathring{e}_3^0 =
			\mathring{\nullapse}^{-1/2},
			&&
			\text{on } \inihypersur \cup \inihypersurbar.
		\end{align}
	\end{definition}
	
	Note that by construction, $\mathring{e}_3$ and $\mathring{e}_4$ are future-directed,
	and by \eqref{E:INITIALLAPSEBOUNDS}, we have:
	\begin{align} \label{E:E40ANDE30INITIALBOUNDS}		
		\frac{1}{\boundinilapse}
		&
		\leq
		\inilapse
		\leq 
		\boundinilapse.
	\end{align}

	\subsection{Transversal and Initially Null Commutators for Rendall's Method}
	\label{SS:TRANSVERSALCOMMUTATORSFORTHERENDALLMETHOD}
	
	\begin{lemma}[Transversal and initially null commutators for Rendall's Method]
		\label{L:TRANSVERALANDINITIALLYNULLCOMMUTATORSFORRENDALLMETHOD}
		On the domain $\embeddomain_{\tilde{\embedkappa},\ubar_*,u_*}$
		from \eqref{E:DIFFEOMORPHISMDOMAIN}, there exists a non-vanishing, $C^{\infty}$ vectorfield
		$E_3$
		such that:
		\begin{align} \label{E:HTRANSVERALCOMMUTATORFORRENDALLMETHOD}
			E_3 \restriction_{\inihypersur}
			& = \mathring{e}_3,
		\end{align}
		where $\mathring{e}_3$ is the null vectorfield from Def.\,\ref{D:NULLFRAMEDATA}.
		We denote the array of Cartesian components of $E_3$ by $\vec{E}_3$, i.e.,
		\begin{align} \label{E:CARTESIANCOMPONENTSOFHTRANSVERALCOMMUTATORFORRENDALLMETHOD}
			\vec{E}_3
			& : = (E_3^0,E_3^1,E_3^2,E_3^2).
		\end{align}
		
		Similarly, on
		$\embeddomain_{\tilde{\embedkappa},\ubar_*,u_*}$,
		there exists a non-vanishing $C^{\infty}$ vectorfield
		$E_4$
		such that:
		\begin{align} \label{E:HBARTRANSVERALCOMMUTATORFORRENDALLMETHOD}
			E_4 \restriction_{\inihypersurbar}
			& = \mathring{e}_4,
		\end{align}
		where $\mathring{e}_4$ is the null vectorfield from Def.\,\ref{D:NULLFRAMEDATA},
		and we denote the array of Cartesian components of $E_4$ by $\vec{E}_4$, i.e.,
		\begin{align} \label{E:CARTESIANCOMPONENTSOFHBARTRANSVERALCOMMUTATORFORRENDALLMETHOD}
			\vec{E}_4
			& : = (E_4^0,E_4^1,E_4^2,E_4^3).
		\end{align}
		
		Moreover, we can decompose the inverse acoustical metric from \eqref{E:INVERSEACOUSTICALMETRIC} as follows,
		where $\gensmoothfunction$ schematically denotes a $C^{\infty}$ function that is free to vary from line to line,
		$\gensmoothfunction_{\star}$ schematically denotes a $C^{\infty}$ function that is free to vary from line to line
		and is non-vanishing on $\inihypersur$,
		and $\gensmoothfunction_0$ schematically denotes a $C^{\infty}$ function that is free to vary from line to line, 
		such that it vanishes on $\inihypersur$ (though its spacetime derivatives do not necessarily vanish):
		\begin{align} 
			\begin{split} \label{E:INVERSEACOUSTICALMETRICWITHNULLFRAMEONH}
				\gfour^{-1}
				& = 
				-
				\gensmoothfunction_{\star}(\inip \embed,\vec{E}_3,\vvariables)
				\left\lbrace
				\frac{\p}{\p \iniubar} 
				\otimes
				E_3
				+
				E_3
				\otimes
				\frac{\p }{\p \iniubar} 
				\right\rbrace
				+
				\gensmoothfunction^{AB}(\inip \embed,\vec{E}_3,\vvariables)
				\cdot 
				\frac{\p}{\p \inisphere^A} \otimes \frac{\p}{\p \inisphere^B} 
				\\
				& \ \
				+
				\gensmoothfunction_0(\inip \embed,\vec{E}_3,\vvariables) 
				\frac{\p}{\p \iniubar} \otimes \frac{\p}{\p \iniubar}
				+
				\gensmoothfunction_0(\inip \embed, \vvariables) 
				E_3 \otimes E_3
				\\
				& \ \ +
				\gensmoothfunction_0^A(\inip \embed,\vec{E}_3,\vvariables)
				\left\lbrace
				\frac{\p}{\p \inisphere^A} \otimes \frac{\p}{\p \iniu}
				+
				\frac{\p}{\p \iniu} \otimes \frac{\p}{\p \inisphere^A}
				\right\rbrace
				+
				\gensmoothfunction_0^A(\inip \embed,\vec{E}_3,\vvariables) 
				\left\lbrace
				\frac{\p}{\p \inisphere^A} \otimes E_3
				+
				E_3 \otimes \frac{\p}{\p \inisphere^A}
				\right\rbrace.
			\end{split}	
		\end{align}
		
		Similarly, we can decompose the inverse acoustical metric as follows,
		where $\gensmoothfunction$ schematically denotes a $C^{\infty}$ function that is free to vary from line to line,
		$\gensmoothfunction_{\star}$ schematically denotes a $C^{\infty}$ function that is free to vary from line to line
		and is non-vanishing on $\inihypersurbar$,
		and $\gensmoothfunction_0$ schematically denotes a $C^{\infty}$ function that is free to vary from line to line, 
		such that it vanishes on $\inihypersurbar$ (though its spacetime derivatives do not necessarily vanish):
		\begin{align} 
			\begin{split} \label{E:INVERSEACOUSTICALMETRICWITHNULLFRAMEONHBAR}
				\gfour^{-1}
				& = 
				-
				\gensmoothfunction_{\star}(\inip \embed,\vec{E}_4,\vvariables)
				\left\lbrace
				\frac{\p}{\p \iniu} 
				\otimes
				E_4
				+
				E_4
				\otimes
				\frac{\p }{\p \iniu} 
				\right\rbrace
				+
				\gensmoothfunction^{AB}(\inip \embed,\vec{E}_4,\vvariables)
				\cdot 
				\frac{\p}{\p \inisphere^A} \otimes \frac{\p}{\p \inisphere^B} 
				\\
				&  \ \
				+
				\gensmoothfunction_0(\inip \embed,\vec{E}_4,\vvariables)
				\frac{\p}{\p \iniubar} \otimes \frac{\p}{\p \iniubar}
				+
				\gensmoothfunction_0(\inip \embed,\vec{E}_4,\vvariables) 
				E_4 \otimes E_4
				\\
				& \ \
				+
				\gensmoothfunction_0^A(\inip \embed,\vec{E}_4,\vvariables)
				\left\lbrace
				\frac{\p}{\p \inisphere^A} \otimes \frac{\p}{\p \iniubar}
				+
				\frac{\p}{\p \iniubar} \otimes \frac{\p}{\p \inisphere^A}
				\right\rbrace
				+
				\gensmoothfunction_0^A(\inip \embed,\vec{E}_4,\vvariables)
				\left\lbrace
				\frac{\p}{\p \inisphere^A} \otimes E_4
				+
				E_4 \otimes \frac{\p}{\p \inisphere^A}
				\right\rbrace.
			\end{split}	
		\end{align}
		
	\end{lemma}
	
	\begin{proof}
		On $\inihypersur$, $\mathring{e}_3$ is a $C^{\infty}$ function of $(\iniubar,\inisphere)$, i.e.,
		$\mathring{e}_3 = \mathring{e}_3(\iniubar,\inisphere)$.
		We can therefore define $E_3(\iniu,\iniubar,\inisphere) := \mathring{e}_3(\iniubar,\inisphere)$,
		and $E_3$ has the desired properties.
		Similarly, on $\inihypersurbar$, $\mathring{e}_4$ is a $C^{\infty}$ function of $(u,\inisphere)$,
		so we can define $E_4(\iniu,\iniubar,\inisphere) := \mathring{e}_4(u,\inisphere)$, and then $E_4$
		has the desired properties.
		
		To prove \eqref{E:INVERSEACOUSTICALMETRICWITHNULLFRAMEONH}, we first note that \eqref{E:INVERSEACOUSTICALMETRIC}
		implies that the Cartesian components of $\gfour^{-1}$ are smooth functions of $\vvariables$.
		We then express $\gfour^{-1}$ in adapted coordinates and note that the components are smooth functions of
		$\vvariables$ and the Jacobian matrix entries $\inip \embed$. 
		Next, we express $\frac{\p}{\p \iniu}$ as a linear combination of the elements of
		$\lbrace E_3, \frac{\p}{\p \iniubar}, \frac{\p}{\p \inisphere^1}, \frac{\p}{\p \inisphere^2} \rbrace$ as follows: 
		\begin{align} \label{E:PARTIALPARTIALADAPATEDUONHINTERMSOFFRAME}
			\frac{\p}{\p \iniu}
			& =
			\gensmoothfunction_{\star}(\inip \embed,\vec{E}_3)
			E_3
			+
			\gensmoothfunction(\inip \embed,\vec{E}_3) \frac{\p}{\p \iniubar}
			+
			\gensmoothfunction^A(\inip \embed,\vec{E}_3) \frac{\p}{\p \inisphere^A},
		\end{align}
		where $\gensmoothfunction_{\star}$ is non-vanishing on $\inihypersur$ because
		$\frac{\p}{\p \iniu}$ and $E_3$ are transversal to $\inihypersur$ 
		while $\frac{\p}{\p \iniubar}, \frac{\p}{\p \inisphere^1}, \frac{\p}{\p \inisphere^1}$
		are tangent to $\inihypersur$.
		From these facts, the null frame
		relations \eqref{D:key}, the  $\inihypersur$-restriction identity \eqref{E:HBARTRANSVERALCOMMUTATORFORRENDALLMETHOD},
		and the fact that $\frac{\p}{\p \iniubar}$ is a positive multiple of
		$\mathring{e}_4$ on $\inihypersur$, we conclude \eqref{E:INVERSEACOUSTICALMETRICWITHNULLFRAMEONH}.
		The identity \eqref{E:INVERSEACOUSTICALMETRICWITHNULLFRAMEONHBAR} can be proved using similar arguments, 
		and we omit the details.

	\end{proof}

	\begin{lemma}[Some identities on the initial null hypersurfaces]
		\label{L:SOMEIDENTITIESONTHEINITIALNULLHYPERSURFACES}
		The following identities hold:
		\begin{subequations}
			\begin{align}
				\label{D:ddudecomp}\frac{\p\embed}{\p \iniu}+\frac{1}{2}\inigfour(\frac{\p\embed}{\p \iniu},\mathring{e}_4) \mathring{e}_3+\frac{1}{2}\inigfour(\frac{\p\embed}{\p \iniu},\mathring{e}_3) \mathring{	e_4}=&(\inigsphere^{-1})^{AB}\inigfour(\frac{\p\embed}{\p \iniu},\mathring{e}_A)\mathring{e}_B,& \text{on }& \inihypersur,\\
				\label{D:ddubardecomp}\frac{\p\embed}{\p\iniubar}+\frac{1}{2}\inigfour(\frac{\p\embed}{\p\iniubar},\mathring{e}_4) \mathring{	e_3}+\frac{1}{2}\inigfour(\frac{\p\embed}{\p\iniubar},\mathring{e}_3) \mathring{e}_4=&(\inigsphere^{-1})^{AB}\inigfour(\frac{\p\embed}{\p\iniubar},\mathring{e}_A)\mathring{e}_B,& \text{on }& \inihypersurbar.
			\end{align}
		\end{subequations}
	\end{lemma}
	
	\begin{proof}
		\eqref{D:ddudecomp}--\eqref{D:ddubardecomp}
		follow from
		the identity
		\eqref{D:decompositionofmetric}, 
		\eqref{EQ:JacobianofEonH}--\eqref{EQ:JacobianofEonHbar},
		Def.\,\ref{D:NULLFRAMEDATA},
	\end{proof}

	\section{Transversal Derivatives along the Initial Null Hypersurfaces, Corner Compatibility Conditions, 
	and the Euler Constraints} \label{S:FLUIDCHARDATAFURTHERASPECTS}
	Recall that $\inivvariables := (\inivel^1,\inivel^2,\inivel^3,\inilogDensity,\inient)$ are $C^{\infty}$ fluid data
	on $\inihypersur$, that $\inivvariablesbar := (\inivelbar^1,\inivelbar^2,\inivelbar^3,\inilogDensitybar,\inientbar)$ 
	are $C^{\infty}$ fluid data on $\inihypersurbar$.	
	In this section, we formally assume that the compressible Euler equations hold on a neighborhood of 
	$\inihypersur \cup \inihypersurbar$, and we 
	construct the solution jet (i.e., the array comprising all of the solution's partial derivatives to all orders) 
	along a \emph{portion} of $\inihypersur$ and a \emph{portion} of $\inihypersurbar$. Our main result in this direction is 
	Prop.\,\ref{P:DETERMININGTRANSVERALDERIVATIVESANDSATISFYINGEQUATIONSONINITIALNULLHYPERSURFACES}. 
	See Remark~\ref{R:BLOWUPONDATAHYPERSURFACES}
	for elaboration on the portion of $\inihypersur \cup \inihypersurbar$ on which the solution jet exists.
	More precisely, in Prop.\,\ref{P:DETERMININGTRANSVERALDERIVATIVESANDSATISFYINGEQUATIONSONINITIALNULLHYPERSURFACES}, 
	we \emph{separately} construct the solution jets on portions of
	$\inihypersur$ and $\inihypersurbar$. There is therefore potential for the solution jet to be
	``multiply defined'' (and thus ill-defined) on $S_{0,0} = \inihypersur \cap \inihypersurbar$.
	However, in Prop.\,\ref{P:CORNERCOMPATIBILITYCONDITIONS}, we show that the solution jet is well-defined on $S_{0,0}$
	as long as a set of ``corner compatibility conditions''
	that we state in Def.\,\ref{D:CORNERCOMPATIBILITYCONDITIONSONDATA}
	are satisfied by the data at $S_{0,0}$.

	\begin{remark}[Blowup of transversal derivatives can occur on the data hypersurfaces]
		\label{R:BLOWUPONDATAHYPERSURFACES}
		In general, the portion of $\inihypersur \cup \inihypersurbar$ on which the solution jet is finite is a strict subset of
		$\inihypersur \cup \inihypersurbar$, corresponding to our choice of the numbers $\ubarfinite$ and $\ufinite$ in 
		Prop.\,\ref{P:DETERMININGTRANSVERALDERIVATIVESANDSATISFYINGEQUATIONSONINITIALNULLHYPERSURFACES}.
		The reason is that equation 
		\eqref{E:NEEDEDFORMWAVEEQUATIONFORDENSITYINADAPATEDCOORDIANTESALONGINITIALNULL} 
		can be viewed as a transport equation for $E_3 \logDensity$ on $\inihypersur$ (with data on $S_{0,0}$),
		and generally, there are terms that are quadratic in $E_3 \logDensity$ on the right-hand side.
		Such quadratic ``Riccati-terms'' can cause the blowup of $E_3 \logDensity$ on $\inihypersur$, which would correspond to the
		formation of the shock on the data hypersurface. In fact, for every fluid equation of state except for that of the
		Chaplygin gas, such Riccati-type terms will be present. This is tied to the fact that the compressible Euler equations fail to satisfy
		the null condition under every equation of state except that of the Chaplygin gas; these issues are discussed in more detail in 
		the introduction of \cite{AbbresciaSpeckYu2026}. 
	\end{remark}

	\subsection{The Wave Equations for the Logarithmic Density in Adapted Coordinates along the Initial Null Hypersurfaces}
	\label{E:WAVEEQUATIONFORDENSITYINADPATEDCOORDINATESALONGINITIALNULL}
	To derive the transversal derivatives of the solution along the initial null hypersurfaces,
	we will use the wave equations in the following lemma. 
	
	\begin{lemma} \label{L:WAVEEQUATIONFORDENSITYINADPATEDCOORDINATESALONGINITIALNULL}
		Let $E_3$ denote the vectorfield on the domain $\embeddomain_{\tilde{\embedkappa},\ubar_*,u_*}$
		from Lemma~\ref{L:TRANSVERALANDINITIALLYNULLCOMMUTATORSFORRENDALLMETHOD}, and let $\vec{E}_3$ denote
		its Cartesian components.
		Then in adapted coordinates, equation \eqref{E:boxgdensity} can be expressed as follows,
		where $\gensmoothfunction$ schematically denotes a $C^{\infty}$ function that is free to vary from line to line
		and $\gensmoothfunction_0$ schematically denotes a $C^{\infty}$ function that is free to vary from line to line
		such that when evaluated at the initial data on 
		$\inihypersur$, it vanishes (though its spacetime derivatives do not necessarily vanish):
		\begin{subequations}
		\begin{align} \label{E:WAVEEQUATIONFORDENSITYINADAPATEDCOORDIANTESALONGINITIALNULL}
			\begin{split}
				\frac{\p}{\p \iniubar}
				E_3 \logDensity
				& = 
				\gensmoothfunction^{AB}(\inip \embed, \vec{E}_3, \vvariables) 
				\cdot 
				\frac{\p}{\p \inisphere^A} \frac{\p}{\p \inisphere^B} \logDensity
				\\
				& 
				+
				\gensmoothfunction_0(\inip \embed, \vec{E}_3, \vvariables) 
				\cdot 
				E_3 E_3 \logDensity
				+
				\gensmoothfunction_0(\inip \embed, \vec{E}_3, \vvariables) 
				\cdot 
				\frac{\p^2}{\p \iniubar^2} \logDensity
				+
				\gensmoothfunction_0^A(\inip \embed, \vec{E}_3, \vvariables)  
				\cdot 
				\frac{\p}{\p \inisphere^A} E_3 \logDensity
				+
				\gensmoothfunction_0^A(\inip \embed, \vec{E}_3, \vvariables) 
				\cdot 
				\frac{\p}{\p \inisphere^A}  \frac{\p}{\p \iniubar}  \logDensity		
				\\
				& 
				+
				\mathscr{L}\left(\inip^{\leq 2},\inip^{\leq 1} \vec{E}_3,\vvariables,
				\frac{\p}{\p \iniubar} \vvariables, \frac{\p}{\p \inisphere} \vvariables, E_3 \vvariables \right)
				[\frac{\p}{\p \iniubar} \vvariables, \frac{\p}{\p \inisphere} \vvariables,E_3 \vvariables]
				+
				\mathscr{L}(\inip \embed,\inip^{\leq 1} \vec{E}_3,\vvariables) 
				[\frac{\p}{\p \iniubar} \vvariables, \frac{\p}{\p \inisphere} \vvariables,E_3 \gradEnt],
			\end{split}	
		\end{align}
		where throughout this lemma, $\mathscr{L}(\cdot)[\cdot]$ denotes a linear function of the terms in square brackets
		with coefficients that are $C^{\infty}$ functions of the terms in parentheses.
		
		Moreover, by using the first-order compressible Euler equations for substitution, 
		we can eliminate all $E_3$ derivatives of $\vvariables$ on
		RHS~\eqref{E:WAVEEQUATIONFORDENSITYINADAPATEDCOORDIANTESALONGINITIALNULL} 
		except for $E_3 \logDensity$, thereby
		expressing the equation as follows:
		\begin{align} \label{E:NEEDEDFORMWAVEEQUATIONFORDENSITYINADAPATEDCOORDIANTESALONGINITIALNULL}
			\begin{split}
				\frac{\p}{\p \iniubar}
				E_3 \logDensity
				& = 
				\gensmoothfunction^{AB}(\inip \embed, \vec{E}_3, \vvariables) 
				\cdot 
				\frac{\p}{\p \inisphere^A} \frac{\p}{\p \inisphere^B} \logDensity
				\\
				& 
				+
				\gensmoothfunction_0(\inip \embed, \vec{E}_3, \vvariables) 
				\cdot 
				E_3 E_3 \logDensity
				+
				\gensmoothfunction_0(\inip \embed, \vec{E}_3, \vvariables) 
				\cdot 
				\frac{\p^2}{\p \iniubar^2} \logDensity
				+
				\gensmoothfunction_0^A(\inip \embed, \vec{E}_3, \vvariables)  
				\cdot 
				\frac{\p}{\p \inisphere^A} E_3 \logDensity
				+
				\gensmoothfunction_0^A(\inip \embed, \vec{E}_3, \vvariables) 
				\cdot 
				\frac{\p}{\p \inisphere^A}  \frac{\p}{\p \iniubar} \logDensity		
				\\
				& 
				+
				\mathscr{L}\left(\inip^{\leq 2} \embed,\inip^{\leq 1} \vec{E}_3,\vvariables,
				\frac{\p}{\p \iniubar} \vvariables, \frac{\p}{\p \inisphere} \vvariables, E_3 \logDensity \right)
				\left[\frac{\p}{\p \iniubar} \vvariables, \frac{\p}{\p \inisphere} \vvariables,E_3 \logDensity \right]
				+
				\mathscr{L}(\inip^{\leq 2} \embed,\vec{E}_3,\vvariables) 
				\left[(\frac{\p}{\p \iniubar} \vvariables, \frac{\p}{\p \inisphere})^2 \Ent \right].
			\end{split}	
		\end{align}
		
		In particular, the following algebraic identity holds independent of whether the compressible Euler equations are satisfied:
		\begin{align} \label{E:ERRORTERMSDEPENDINGONFIRSTORDEREULER}
			\mbox{RHS\eqref{E:WAVEEQUATIONFORDENSITYINADAPATEDCOORDIANTESALONGINITIALNULL} - RHS~\eqref{E:NEEDEDFORMWAVEEQUATIONFORDENSITYINADAPATEDCOORDIANTESALONGINITIALNULL}}
			& = 
			\mathscr{L}\left(\inip^{\leq 2} \embed,\inip^{\leq 1} \vec{E}_3,\vvariables,
			\frac{\p}{\p \iniubar} \vvariables, \frac{\p}{\p \inisphere} \vvariables, E_3 \logDensity, 
			\FR, \FV,\FS,E_3 \FS
			\right)
			[\FR,\FV,\FS,E_3 \FS],
		\end{align}
		\end{subequations}
		and thus $\mbox{RHS~\eqref{E:ERRORTERMSDEPENDINGONFIRSTORDEREULER}} = 0$
		when $\FR= \FV = \FS = E_3 \FS = 0$.
		
		Similarly, let $E_4$ denote the vectorfield on the domain $\embeddomain_{\tilde{\embedkappa},\ubar_*,u_*}$
		from Lemma~\ref{L:TRANSVERALANDINITIALLYNULLCOMMUTATORSFORRENDALLMETHOD}, and let $\vec{E}_4$ denote
		its Cartesian components.
		Then in adapted coordinates, equation \eqref{E:boxgdensity} can be expressed as follows,
		where $\gensmoothfunction$ schematically denotes a $C^{\infty}$ function that is free to vary from line to line
		and $\gensmoothfunction_0$ schematically denotes a $C^{\infty}$ function that is free to vary from line to line
		such that when evaluated at the initial data on 
		$\inihypersurbar$, it vanishes (though its spacetime derivatives do not necessarily vanish):
		\begin{subequations}
		\begin{align} \label{E:WAVEEQUATIONFORDENSITYINADAPATEDCOORDIANTESALONGINITIALNULLHBAR}
			\begin{split}
				\frac{\p}{\p \iniu}
				E_4 \logDensity
				& = 
				\gensmoothfunction^{AB}(\inip \embed, \vec{E}_4, \vvariables) 
				\cdot 
				\frac{\p}{\p \inisphere^A} \frac{\p}{\p \inisphere^B} \logDensity
				\\
				& 
				+
				\gensmoothfunction_0(\inip \embed, \vec{E}_4, \vvariables) 
				\cdot 
				E_4 E_4 \logDensity
				+
				\gensmoothfunction_0(\inip \embed, \vec{E}_4, \vvariables) 
				\cdot 
				\frac{\p^2}{\p \iniu} \logDensity
				+
				\gensmoothfunction_0^A(\inip \embed, \vec{E}_4, \vvariables)  
				\cdot 
				\frac{\p}{\p \inisphere^A} E_4 \logDensity
				+
				\gensmoothfunction_0^A(\inip \embed, \vec{E}_4, \vvariables) 
				\cdot 
				\frac{\p}{\p \inisphere^A}  \frac{\p}{\p \iniubar}  \logDensity		
				\\
				& 
				+
				\mathscr{L}\left(\inip^{\leq 2} \embed,\inip^{\leq 1} \vec{E}_3,\vvariables,
				\frac{\p}{\p \iniubar} \vvariables, \frac{\p}{\p \inisphere} \vvariables, E_4 \vvariables \right)
				[\frac{\p}{\p \iniubar} \vvariables, \frac{\p}{\p \inisphere} \vvariables,E_4 \vvariables]
				+
				\mathscr{L}(\inip \embed,\inip^{\leq 1} \vec{E}_3,\vvariables) 
				[\frac{\p}{\p \iniubar} \vvariables, \frac{\p}{\p \inisphere} \vvariables,E_4 \gradEnt].
			\end{split}	
		\end{align}
		 
		Moreover, by using the first-order compressible Euler equations for substitution, 
		we can eliminate all $E_4$ derivatives of $\vvariables$ on
		RHS~\eqref{E:WAVEEQUATIONFORDENSITYINADAPATEDCOORDIANTESALONGINITIALNULLHBAR} 
		except for $E_4 \logDensity$, thereby
		expressing the equation as follows:
		\begin{align} \label{E:NEEDEDFORMWAVEEQUATIONFORDENSITYINADAPATEDCOORDIANTESALONGINITIALNULLHBAR}
			\begin{split}
				\frac{\p}{\p \iniu}
				E_4 \logDensity
				& = 
				\gensmoothfunction^{AB}(\inip \embed, \vec{E}_4, \vvariables) 
				\cdot 
				\frac{\p}{\p \inisphere^A} \frac{\p}{\p \inisphere^B} \logDensity
				\\
				& 
				+
				\gensmoothfunction_0(\inip \embed, \vec{E}_4, \vvariables) 
				\cdot 
				E_4 E_4 \logDensity
				+
				\gensmoothfunction_0(\inip \embed, \vec{E}_4, \vvariables) 
				\cdot 
				\frac{\p^2}{\p \iniu^2} \logDensity
				+
				\gensmoothfunction_0^A(\inip \embed, \vec{E}_4, \vvariables)  
				\cdot 
				\frac{\p}{\p \inisphere^A} E_4 \logDensity
				+
				\gensmoothfunction_0^A(\inip \embed, \vec{E}_4, \vvariables) 
				\cdot 
				\frac{\p}{\p \inisphere^A}  \frac{\p}{\p \iniubar} \logDensity		
				\\
				& 
				+
				\mathscr{L}\left(\inip^{\leq 2} \embed,\inip^{\leq 1} \vec{E}_4,\vvariables,
				\frac{\p}{\p \iniubar} \vvariables, \frac{\p}{\p \inisphere} \vvariables, E_4 \logDensity \right)
				\left[\frac{\p}{\p \iniubar} \vvariables, \frac{\p}{\p \inisphere} \vvariables,E_4 \logDensity \right]
				+
				\mathscr{L}(\inip^{\leq 2} \embed,\vec{E}_4,\vvariables) 
				\left[(\frac{\p}{\p \iniubar} \vvariables, \frac{\p}{\p \inisphere})^2 \Ent \right].
			\end{split}	
		\end{align}
		
		In particular, the following algebraic identity holds independent of whether the compressible Euler equations are satisfied:
		\begin{align} \label{E:ERRORTERMSDEPENDINGONFIRSTORDEREULERHBAR}
			\mbox{RHS\eqref{E:WAVEEQUATIONFORDENSITYINADAPATEDCOORDIANTESALONGINITIALNULLHBAR} - RHS~\eqref{E:NEEDEDFORMWAVEEQUATIONFORDENSITYINADAPATEDCOORDIANTESALONGINITIALNULLHBAR}}
			& = 
			\mathscr{L}\left(\inip^{\leq 2} \embed,\inip^{\leq 1} \vec{E}_4,\vvariables,
			\frac{\p}{\p \iniubar} \vvariables, \frac{\p}{\p \inisphere} \vvariables, E_4 \logDensity, 
			\FR, \FV,\FS,E_4 \FS
			\right)
			[\FR,\FV,\FS,E_4 \FS],
		\end{align}
		\end{subequations}
		and thus $\mbox{RHS~\eqref{E:ERRORTERMSDEPENDINGONFIRSTORDEREULERHBAR}} = 0$
		when $\FR= \FV = \FS = E_4 \FS  = 0$.

	\end{lemma}
	
	\begin{proof}
		\eqref{E:WAVEEQUATIONFORDENSITYINADAPATEDCOORDIANTESALONGINITIALNULL} follows from
		expressing equation \eqref{E:boxgdensity} in adapted coordinates
		and using the expression
		\eqref{E:INVERSEACOUSTICALMETRICWITHNULLFRAMEONH}
		for the acoustical metric.
	
		To pass from \eqref{E:WAVEEQUATIONFORDENSITYINADAPATEDCOORDIANTESALONGINITIALNULL}
		to \eqref{E:NEEDEDFORMWAVEEQUATIONFORDENSITYINADAPATEDCOORDIANTESALONGINITIALNULL},
		we simply use Lemma~\ref{E:SOMECOMPRESSIBLEEULEREQUATIONSINADPATEDCOORDINATES}
		to substitute $E_3 \Ent$ with error terms involving 
		$\inihypersur$-tangential derivatives of $\Ent$
		and to substitute $E_3 v^i$ with $E_3 \logDensity$ plus error terms involving 
		$\inihypersur$-tangential derivatives.
		
		The results stated in \eqref{E:WAVEEQUATIONFORDENSITYINADAPATEDCOORDIANTESALONGINITIALNULLHBAR}--\eqref{E:ERRORTERMSDEPENDINGONFIRSTORDEREULERHBAR}
		follow from similar arguments based on \eqref{E:INVERSEACOUSTICALMETRICWITHNULLFRAMEONHBAR}, 
		and we omit the details.
	\end{proof}

	For future use, we define:
	\begin{subequations}
	\begin{align} \label{E:ADAPTEDCOORDINATESVERSIONOFSPECKLUKWAVEEQUATIONTERMS}
		\widetilde{\FW}
		& := \mbox{LHS~\eqref{E:NEEDEDFORMWAVEEQUATIONFORDENSITYINADAPATEDCOORDIANTESALONGINITIALNULL} - RHS~\eqref{E:NEEDEDFORMWAVEEQUATIONFORDENSITYINADAPATEDCOORDIANTESALONGINITIALNULL}}, 
			\\
	\widehat{\FW}
		& := \mbox{LHS~\eqref{E:NEEDEDFORMWAVEEQUATIONFORDENSITYINADAPATEDCOORDIANTESALONGINITIALNULLHBAR} - RHS~\eqref{E:NEEDEDFORMWAVEEQUATIONFORDENSITYINADAPATEDCOORDIANTESALONGINITIALNULLHBAR}}, 
		\label{E:ADAPTEDCOORDINATESVERSIONOFSPECKLUKWAVEEQUATIONTERMSHBAR}
	\end{align}
	\end{subequations}
	i.e., equation \eqref{E:NEEDEDFORMWAVEEQUATIONFORDENSITYINADAPATEDCOORDIANTESALONGINITIALNULL} is precisely 
	$\widetilde{\FW} = 0$
	and equation \eqref{E:NEEDEDFORMWAVEEQUATIONFORDENSITYINADAPATEDCOORDIANTESALONGINITIALNULLHBAR} is precisely 
	$\widehat{\FW} = 0$.

	\subsection{Determining the Transversal Derivatives and Satisfying the Equations to All Orders on the Initial Null Hypersurfaces}
	\label{SS:TRANSVERALDERIVATIVESONINITIALNULLHYPERSURFACES}
	We are almost ready to prove the main result of this section,
	Prop.\,\ref{P:DETERMININGTRANSVERALDERIVATIVESANDSATISFYINGEQUATIONSONINITIALNULLHYPERSURFACES}, 
	in which we construct the solution jet on a portion of $\inihypersur$ and separately on a portion of $\inihypersurbar$
	in such a way that the compressible Euler equations are satisfied to all orders on the portions of
	$\inihypersur$ and on $\inihypersurbar$.

	\subsubsection{A preliminary lemma}
	\begin{lemma}[Further algebraic relationship between various equations]
		\label{L:ALGEBRAICRELATIONSHIPESBETWEENSOMEEQUATIONS2}
		Let $\FR,\FV,\FS,$ and $\FW$
		be as defined in \eqref{E:CONSTRAINTDENSITYEVOLUTIONEQUATION}--\eqref{E:CONSTRAINTENTROPYEVOLUTIONEQUATION} 
		and
		\eqref{E:SPECKLUKWAVEEQUATIONTERMS}.
		With $\widetilde{\FW}$ 
		as defined in \eqref{E:ADAPTEDCOORDINATESVERSIONOFSPECKLUKWAVEEQUATIONTERMS},
		the following identity holds relative to the adapted coordinates,\footnote{In particular, it should be understood 
		that the term $\FW$
		on LHS~\eqref{E:EQUIVALENCEOFTWOWAVEEQUATIONSUPTOFIRSTORDEREULER2} and LHS~\eqref{E:EQUIVALENCEOFTWOWAVEEQUATIONSUPTOFIRSTORDEREULER2HBAR} is expressed in adapted coordinates.}
		where
		$E_3$ denotes the vectorfield on the domain $\embeddomain_{\tilde{\embedkappa},\ubar_*,u_*}$
		from Lemma~\ref{L:TRANSVERALANDINITIALLYNULLCOMMUTATORSFORRENDALLMETHOD} and $\vec{E}_3$ denotes
		its Cartesian components:
		\begin{subequations}
		\begin{align} \label{E:EQUIVALENCEOFTWOWAVEEQUATIONSUPTOFIRSTORDEREULER2}
			\FW
			& 
			=
			\gensmoothfunction_{\star}(\inip \embed, \vvariables)
			\widetilde{\FW}
			+
			\mathscr{L}(\inip^{\leq 2} \embed,\inip^{\leq 1} \vec{E}_3, \vvariables, \FR, \FV,\FS,E_3 \FS)
			[\FR,\FV,\FS,E_3 \FS],
		\end{align}
		where $\gensmoothfunction_{\star}$ is a $C^{\infty}$ non-vanishing function, 
		and $\mathscr{L}(\cdot)[\cdot]$ is a linear function of the terms in square brackets
		with coefficients that are $C^{\infty}$ functions of the terms in parentheses.
	
		Similarly, with $\widehat{\FW}$ 
		as defined in \eqref{E:ADAPTEDCOORDINATESVERSIONOFSPECKLUKWAVEEQUATIONTERMSHBAR},
		the following identity holds relative to the adapted coordinates, where
		$E_4$ denotes the vectorfield on the domain $\embeddomain_{\tilde{\embedkappa},\ubar_*,u_*}$
		from Lemma~\ref{L:TRANSVERALANDINITIALLYNULLCOMMUTATORSFORRENDALLMETHOD} and $\vec{E}_4$ denotes
		its Cartesian components:
		\begin{align} \label{E:EQUIVALENCEOFTWOWAVEEQUATIONSUPTOFIRSTORDEREULER2HBAR}
			\FW
			& 
			=
			\gensmoothfunction_{\star}(\inip \embed, \vvariables)
			\widehat{\FW}
			+
			\mathscr{L}(\inip^{\leq 2} \embed,\inip^{\leq 1} \vec{E}_4, \vvariables, \FR, \FV,\FS,E_4 \FS)
			[\FR,\FV,\FS,E_4 \FS],
		\end{align}
		where $\gensmoothfunction_{\star}$ is a $C^{\infty}$ non-vanishing function, 
		and $\mathscr{L}(\cdot)[\cdot]$ is a linear function of the terms in square brackets
		with coefficients that are $C^{\infty}$ functions of the terms in parentheses.
		
		Finally, we have:
		\begin{align} \label{E:TWOWAVEEQUATIONEXPRESSIONSARESAMEUPTOLOWERORDERTERMS}
			\widetilde{\FW}
			& = 
			\gensmoothfunction_{\star}(\inip \embed, \vvariables)
			\widehat{\FW}
			+
			\mathscr{L}(\inip^{\leq 2} \embed,\inip^{\leq 1} \vec{E}_3,\inip^{\leq 1} \vec{E}_4,\vvariables, \FR, \FV,\FS,E_3 \FS,E_4 \FS)
			[\FR,\FV,\FS,E_3 \FS,E_4 \FS],
		\end{align}
		\end{subequations}
		where on RHS~\eqref{E:TWOWAVEEQUATIONEXPRESSIONSARESAMEUPTOLOWERORDERTERMS},
		\begin{align} \label{E:TWOWAVEQUANTITIESHAVESOMEPRINCIPALPARTATS00}
			\gensmoothfunction_{\star}(\inip \embed, \vvariables) \restriction_{S_{0,0}}  = 1. 
		\end{align}
	\end{lemma}
	
	\begin{proof}[Proof of Lemma \ref{L:ALGEBRAICRELATIONSHIPESBETWEENSOMEEQUATIONS2}]
		The identity \eqref{E:EQUIVALENCEOFTWOWAVEEQUATIONSUPTOFIRSTORDEREULER2}
		follows from Lemma \ref{L:ALGEBRAICRELATIONSHIPESBETWEENSOMEEQUATIONS} and revisiting the proof of \eqref{E:NEEDEDFORMWAVEEQUATIONFORDENSITYINADAPATEDCOORDIANTESALONGINITIALNULL}
		and keeping track of the steps in which the first-order 
		compressible Euler equations \eqref{E:Brho}--\eqref{E:Bs}
		were used for substitution; since we are not assuming that these equations are satisfied, 
		every time we ``use the equations for substitution,''
		we generate terms that depend on $\FR,\FV,\FS,E_3 \FS$,
		which is why these terms appear in \eqref{E:EQUIVALENCEOFTWOWAVEEQUATIONSUPTOFIRSTORDEREULER2}.
		
		The identity \eqref{E:EQUIVALENCEOFTWOWAVEEQUATIONSUPTOFIRSTORDEREULER2HBAR} follows from similar arguments
		based on revisiting the proof of
		\eqref{E:NEEDEDFORMWAVEEQUATIONFORDENSITYINADAPATEDCOORDIANTESALONGINITIALNULLHBAR}.
		
		\eqref{E:TWOWAVEEQUATIONEXPRESSIONSARESAMEUPTOLOWERORDERTERMS}
		follows easily from
		\eqref{E:EQUIVALENCEOFTWOWAVEEQUATIONSUPTOFIRSTORDEREULER2}
		and
		\eqref{E:EQUIVALENCEOFTWOWAVEEQUATIONSUPTOFIRSTORDEREULER2HBAR}.
	\end{proof}

	\subsubsection{Determining the transversal derivatives and satisfying the equations 
		to all orders on portions of the initial null hypersurfaces}
	\label{SSS:DETERMININGTRANSVERALDERIVATIVESANDSATISFYINGEQUATIONSONINITIALNULLHYPERSURFACES}
	We now state and prove the main result of this section.
	
	\begin{proposition}[Determining the transversal derivatives and satisfying the equations to all orders on portions of the initial null hypersurfaces]
		\label{P:DETERMININGTRANSVERALDERIVATIVESANDSATISFYINGEQUATIONSONINITIALNULLHYPERSURFACES}
		Recall that $\inivvariables$ are given $C^{\infty}$ fluid data
		on $\inihypersur$ (see \eqref{E:DATAONH}) and that
		$\inivvariablesbar$ are given $C^{\infty}$ fluid data on $\inihypersurbar$ (see \eqref{E:DATAONHBAR}). 
		Assume that the Euler constraint equation \eqref{E:constraint1} is satisfied on
		$\inihypersur$, that the Euler constraint equation \eqref{E:constraint2} is satisfied on $\inihypersurbar$.
		Recall also that the first-order compressible Euler equations are the system
		$\FR = 0, \FS = 0, \FV = 0$ (see \eqref{EQ:EquationThemselves}), that the equations
		$\widetilde{\FW} = 0$ and $\widehat{\FW} = 0$ are ``auxiliary'' wave equations 
		(see \eqref{E:ADAPTEDCOORDINATESVERSIONOFSPECKLUKWAVEEQUATIONTERMS}--\eqref{E:ADAPTEDCOORDINATESVERSIONOFSPECKLUKWAVEEQUATIONTERMSHBAR})
		that must be satisfied by any sufficiently
		smooth solution to the first-order equations,
		and that $E_3$ and $E_4$ are the vectorfields on the domain $\embeddomain_{\tilde{\embedkappa},\ubar_*,u_*}$
		from Lemma~\ref{L:TRANSVERALANDINITIALLYNULLCOMMUTATORSFORRENDALLMETHOD}, where
		$E_3$ is transversal to $\inihypersur$ and $E_4$ is transversal to $\inihypersurbar$.
		
		\medskip
		
		\noindent \underline{\textbf{Conclusions on $\inihypersur$}}.
		Then there exists a number $\ubarfinite \in (0,\ubar_*]$ such that for every integer $K \geq 0$, these initial data, together with 
		the following assumed equations: 	
 		\begin{align} \label{E:KTHORDERTRANSVERALDERIVATIVESOFMOSTEQUATIONSSATISFIEDONH}
			E_3^{\leq K} \widetilde{\FW} \restriction_{\inihypersur_{[0,\ubarfinite]}} & = 0,
			&
			E_3^{\leq K} \FS \restriction_{\inihypersur_{[0,\ubarfinite]}} & = 0,
			&
			E_3^{\leq K} \FV^i \restriction_{\inihypersur_{[0,\ubarfinite]}} & = 0,
		\end{align}
		uniquely determine the solution's $K+1$ transversal derivatives on: 
		\begin{align} \label{E:PORTIONOFHDATASURFACEWITHFINITEJET}
			\inihypersur_{[0,\ubarfinite]} 
			:= \embed\left(\lbrace (0,\iniubar,\inisphere)  \ | \ 
			0 \leq \iniubar \leq \ubarfinite, \inisphere \in \mathbb{S}^2 \rbrace \right),
		\end{align}
		i.e., they uniquely determine:
		\begin{align} \label{E:TRANSVERALDERIVATIVESDETERMINEDONH}
			\inilogDensity_{(K+1)}
			& := E_3^{K+1} \logDensity \restriction_{\inihypersur_{[0,\ubarfinite]}},
			& 
			\inient_{(K+1)} & := E_3^{K+1} \Ent \restriction_{\inihypersur_{[0,\ubarfinite]}},
			&
			\inivel_{(K+1)}^i & := E_3^{K+1} v^i \restriction_{\inihypersur_{[0,\ubarfinite]}}.
		\end{align}
		In fact, $\ubarfinite$ can be chosen to be any number in $[0,u_*]$ such that the solution
		$E_3 \logDensity$ to equation \eqref{E:NEEDEDFORMWAVEEQUATIONFORDENSITYINADAPATEDCOORDIANTESALONGINITIALNULL} 
		is finite on $\inihypersur_{[0,\ubarfinite]}$.
		Moreover, all of the quantities in \eqref{E:TRANSVERALDERIVATIVESDETERMINEDONH} are $C^{\infty}$ functions on 
		$\inihypersur_{[0,\ubarfinite]}$.
		
		In addition, for any integer $K \geq 0$, the following equation is satisfied:
		\begin{align}  \label{E:KTHORDERTRANSVERALDERIVATIVESOFDENSITYEQUATION}
			E_3^K \FR \restriction_{\inihypersur_{[0,\ubarfinite]}} & = 0.
		\end{align}
		In particular, \eqref{E:KTHORDERTRANSVERALDERIVATIVESOFMOSTEQUATIONSSATISFIEDONH} and \eqref{E:KTHORDERTRANSVERALDERIVATIVESOFDENSITYEQUATION}
		imply that the first-order compressible Euler equations and their transversal derivatives to all orders are satisfied on 
		$\inihypersur_{[0,\ubarfinite]}$.
		
		\medskip
		\noindent \underline{\textbf{Conclusions on $\inihypersurbar$}}.
		Similarly, there exists a number $\ufinite \in [0,u_*]$ 
		such that for every integer $K \geq 0$, these initial data, together with 
		the following assumed equations: 	
		\begin{align} \label{E:KTHORDERTRANSVERALDERIVATIVESOFMOSTEQUATIONSSATISFIEDONHBAR}
			E_4^{\leq K} \widehat{\FW} \restriction_{\inihypersurbar_{[0,\ufinite]}} & = 0,
			&
			E_4^{\leq K} \FS \restriction_{\inihypersurbar_{[0,\ufinite]}} & = 0,
			&
			E_4^{\leq K} \FV^i \restriction_{\inihypersurbar_{[0,\ufinite]}} & = 0,
		\end{align}
		uniquely determine the solution's $K+1$ transversal derivatives on 
		\begin{align} \label{E:PORTIONOFHBARDATASURFACEWITHFINITEJET}
			\inihypersurbar_{[0,\ufinite]}
			:= \embed\left(\lbrace (\iniu,0,\inisphere) \ | \ 0 \leq \iniu \leq \ufinite, \inisphere \in \mathbb{S}^2 \rbrace \right),
		\end{align}
		i.e., they determine:
		\begin{align} \label{E:TRANSVERALDERIVATIVESDETERMINEDONHBAR}
			\inilogDensitybar_{(K+1)}
			& := E_4^{K+1} \logDensity \restriction_{\inihypersurbar_{[0,\ufinite]}},
			& 
			\inientbar_{(K+1)} & := E_4^{K+1} \Ent \restriction_{\inihypersurbar_{[0,\ufinite]}},
			&
			\inivelbar_{(K+1)}^i & := E_4^{K+1} v^i \restriction_{\inihypersurbar_{[0,\ufinite]}}.
		\end{align}
		In fact, $\ufinite$ can be chosen to be any number in $[0,u_*]$ such that 
		$E_4 \logDensity$ is finite on $\inihypersurbar_{[0,\ufinite]}$.
		Moreover, all of the quantities in \eqref{E:TRANSVERALDERIVATIVESDETERMINEDONHBAR} are $C^{\infty}$ functions on 
		$\inihypersurbar_{[0,\ufinite]}$.
		
		In addition, for any integer $K \geq 0$, the following equation is satisfied:
		\begin{align}  \label{E:KTHORDERTRANSVERALDERIVATIVESOFDENSITYEQUATIONONHBAR}
			E_4^K \FR \restriction_{\inihypersurbar_{[0,\ufinite]}} & = 0.
		\end{align}
		In particular, \eqref{E:KTHORDERTRANSVERALDERIVATIVESOFMOSTEQUATIONSSATISFIEDONHBAR} and \eqref{E:KTHORDERTRANSVERALDERIVATIVESOFDENSITYEQUATIONONHBAR}
		imply that the first-order compressible Euler equations and their transversal derivatives to all orders are satisfied on 
		$\inihypersurbar_{[0,\ufinite]}$.
		
	\end{proposition}
	
	\begin{proof}
		We carry out the analysis only on $\inihypersur$ since the analysis on $\inihypersurbar$ can be carried out using similar arguments.
		We will inductively solve the equations \eqref{E:KTHORDERTRANSVERALDERIVATIVESOFMOSTEQUATIONSSATISFIEDONH} and solve for the 
		quantities in \eqref{E:TRANSVERALDERIVATIVESDETERMINEDONH}, starting from the base case $K=0$.
		
		Recall that $\widetilde{\FW} = 0$ is a nonlinear transport equation for $E_3 \logDensity$.
		Solving $E_3^K \widetilde{\FW} = 0$ is equivalent 
		to commuting equation \eqref{E:NEEDEDFORMWAVEEQUATIONFORDENSITYINADAPATEDCOORDIANTESALONGINITIALNULL}
		with $E_3^K$ and solving the resulting transport equation for 
		$E_3^{K+1} \logDensity$, with initial data on $S_{0,0}$.
		For this procedure, it is important that the terms $\gensmoothfunction_0$ on
		RHS~\eqref{E:NEEDEDFORMWAVEEQUATIONFORDENSITYINADAPATEDCOORDIANTESALONGINITIALNULL} vanish to order $0$.
		The initial data for $E_3^{K+1} \logDensity$ on $S_{0,0}$ is computable because  
		the fluid data on $\inihypersurbar$ are prescribed, and $E_3 \restriction_{S_{0,0}}$ is tangent to 
		$\inihypersurbar$. 
		We therefore commute \eqref{E:NEEDEDFORMWAVEEQUATIONFORDENSITYINADAPATEDCOORDIANTESALONGINITIALNULL} with
		$E_3^K$.
		All of the terms appearing on the 
		RHS of the commuted equation are $C^{\infty}(\inihypersur)$-terms that are 
		determined directly by the initial data or the induction hypotheses, except for $E_3^{K+1} \logDensity$, 
		which appears nonlinearly on the RHS of the transport equation when $K=0$ and, when $K \geq 1$, appears  
		linearly on the RHS of the transport equation.  Due to the nonlinearity in the case $K=0$,
		$E_3 \logDensity$ exists only on $\inihypersur_{[0,\ubarfinite]}$
		for some number $\ubarfinite \in (0,\ubar_*]$. On the other hand, 
		for $K \geq 1$, the $E_3^K$-commuted equation is linear in 
		$E_3^{K+1} \logDensity$, and the solution to the linear equation therefore exists on the same portion
		$\inihypersur_{[0,\ubarfinite]}$ as in the base case.
		We have therefore determined $\inilogDensity_{(K+1)}$
		on $\inihypersur_{[0,\ubarfinite]}$. Standard ODE theory implies that solutions 
		to transport equations with $C^{\infty}$ coefficients and source terms/data depending in a $C^{\infty}$ fashion on parameters are
		$C^{\infty}$ solutions that depend in a $C^{\infty}$ fashion on the parameters. This implies that
		$\inilogDensity_{(K+1)}$ is a $C^{\infty}$ function on $\inihypersur_{[0,\ubarfinite]}$.
		
		We then differentiate equation \eqref{E:ENTROPYEULEREQUATIONWITHE3DERIVATIVE} with 
		$E_3^K$, which allows us to algebraically express
		$E_3^{K+1} \Ent$ in terms of already determined quantities. 
		This yields $\inient_{(K+1)} = E_3^{K+1} \Ent \restriction_{\inihypersur_{[0,\ubarfinite]}}$ as 
		a $C^{\infty}$ function on $\inihypersur_{[0,\ubarfinite]}$,
		and by \eqref{E:ENTROPYEULEREQUATIONWITHE3DERIVATIVEISAMULTIPLEOFENTROPYCONSTRAINT} and the induction hypotheses, 
		it ensures that $E_3^K \FS = 0$ on $\inihypersur_{[0,\ubarfinite]}$.
		
		We then differentiate equation \eqref{E:e3vADAPTEDCOORDINATES} with 
		$E_3^K$, which allows us to algebraically express
		$E_3^{K+1} v^i$ in terms of already determined quantities. 
		This yields $\inivel_{(K+1)} = E_3^{K+1} v \restriction_{\inihypersur_{[0,\ubarfinite]}}$
		as a $C^{\infty}$ function on $\inihypersur_{[0,\ubarfinite]}$, 
		and by \eqref{E:E3VINADPATEDCOORDINATESISAMULTIPLEOFVELOCITYCONSTRAINT} and the induction hypotheses,  
		it ensures that
		$E_3^K \FV \restriction_{\inihypersur_{[0,\ubarfinite]}} = 0$.

		To summarize, we have thus far proved \eqref{E:TRANSVERALDERIVATIVESDETERMINEDONH} and 
		\eqref{E:KTHORDERTRANSVERALDERIVATIVESOFMOSTEQUATIONSSATISFIEDONH} for all integers $K \geq 0$.
		
		It remains for us to prove \eqref{E:KTHORDERTRANSVERALDERIVATIVESOFDENSITYEQUATION}. We will again use induction in $K$,
		starting with $K=0$. First, we note that Prop.\,\ref{P:constraint0} implies that since the Euler constraint \eqref{E:constraint1},
		the equation $\FS = 0$, and the equation $\FV = 0$ are all satisfied on
		$\inihypersur_{[0,\ubarfinite]}$, the compressible Euler density equation $\FR = 0$ is also satisfied on 
		$\inihypersur_{[0,\ubarfinite]}$, i.e., \eqref{E:KTHORDERTRANSVERALDERIVATIVESOFDENSITYEQUATION} holds in the case $K=0$.
		For the induction step, we assume that \eqref{E:KTHORDERTRANSVERALDERIVATIVESOFDENSITYEQUATION} has been shown for $0,1,\cdots,K$.
		We aim to show that \eqref{E:KTHORDERTRANSVERALDERIVATIVESOFDENSITYEQUATION} holds with $K+1$ in the role of $K$.
		From equation \eqref{E:EQUIVALENCEOFTWOWAVEEQUATIONSUPTOFIRSTORDEREULER2}, the induction hypotheses, and
		the results stated in \eqref{E:KTHORDERTRANSVERALDERIVATIVESOFMOSTEQUATIONSSATISFIEDONH} 
		(which we have already shown to hold for all $K$),
		we find that $E_3^K \FW \restriction_{\inihypersur_{[0,\ubarfinite]}} = 0$.
		Also using equation \eqref{E:DIFFERENTIATEDDENSITYCONSTRAINEDCOMPAREDTOLUKSPECKWAVE},
		we further deduce that
		\begin{align} \label{E:E3KDERIVATIVEOFMATERIALDERIVATIVEOFEULERDENSITYCONSTRAINTVANISHES}
			E_3^K (\materialderivative \FR) \restriction_{\inihypersur_{[0,\ubarfinite]}} = 0.
		\end{align}
		The key point is that $\materialderivative$, being timelike, is transversal to the null hypersurface $\inihypersur_{[0,\ubarfinite]}$
		and therefore, when expanded relative to the frame 
		$\lbrace E_3, \frac{\p}{\p \iniubar},  \frac{\p}{\p \inisphere^1},\frac{\p}{\p \inisphere^2} \rbrace$, 
		must have a non-zero component in the direction of $E_3$ 
		(since the other three vectorfields in the frame are $\inihypersur_{[0,\ubarfinite]}$-tangent). 
		From this fact, \eqref{E:E3KDERIVATIVEOFMATERIALDERIVATIVEOFEULERDENSITYCONSTRAINTVANISHES}, 
		and the induction hypotheses, we conclude that
		$E_3^{K+1} \FR \restriction_{\inihypersur_{[0,\ubarfinite]}} = 0$,
		which completes the induction hypothesis step and finishes the proof of the proposition.
		
	\end{proof}

	\subsection{Corner Compatibility Conditions on $S_{0,0}$}
	\label{SS:CORNERcompatibilityCONDITIONS}
	Recall that in Prop.\,\ref{P:DETERMININGTRANSVERALDERIVATIVESANDSATISFYINGEQUATIONSONINITIALNULLHYPERSURFACES}, 
	we separately constructed the solution jets on portions of $\inihypersur$ and $\inihypersurbar$.
	On $S_{0,0} = \inihypersur \cap \inihypersurbar$, there is potential for the solution's derivatives to be multiply defined:
	in one way from the data on $\inihypersur$ + the compressible Euler equations and in a second way from the data on 
	$\inihypersurbar$ + the compressible Euler equations. In this section, we show that such inconsistencies can be 
	avoided as long as the initial data satisfy 
	corner compatibility conditions that we state in Def.\,\ref{D:CORNERCOMPATIBILITYCONDITIONSONDATA}.
	In Prop.\,\ref{P:CORNERCOMPATIBILITYCONDITIONS}, we show that these corner compatibility conditions
	are sufficient to yield the well-definedness of the solution jet at $S_{0,0}$.

	\subsubsection{Various solution-jet equations}
	\label{SSS:SOMESOLUTIONJETEQUATIONS}
	Before stating the corner compatibility conditions, we first derive various PDEs that must be satisfied 
	by any $C^{\infty}$ compressible Euler solution at $S_{0,0}$.
	For convenience, we will use the adapted coordinates $(\iniu,\iniubar,\inisphere^1,\inisphere^2)$
	and the diffeomorphism $\embed$ from \eqref{E:LOCALDIFFEOCONSTRUCTEDOUTOFEMBEDDINGS}
	when stating the equations.

	\begin{lemma}[Various solution-jet equations at $S_{0,0}$]
	\label{L:SOMESOLUTIONJETEQUATIONS}
	Recall that $(\iniu,\iniubar,\inisphere^1,\inisphere^2)$ are the adapted coordinates from Def.\,\ref{D:ADAPTEDCOORDINATES},
	that $\inip$ denotes partial derivatives with respect to the adapted coordinates,
	and that $\embed$ is the diffeomorphism from \eqref{E:LOCALDIFFEOCONSTRUCTEDOUTOFEMBEDDINGS}.
	Recall also that the first-order compressible Euler equations are the system
	$\FR = 0, \FV^i = 0, \FS = 0$ (see \eqref{EQ:EquationThemselves}).
	Then the following hold at the initial sphere $S_{0,0}$,
	where $\embedkappa > 0$ is the number from Lemma~\ref{L:EXISTENCEOFALOCALDIFFEOMORPHISM},
	where in what follows, $\gensmoothfunction_*$, $\gensmoothfunction_0$, and $\gensmoothfunction$ denote smooth functions
	that can vary from line to line, $\gensmoothfunction_*$ is non-vanishing on $S_{0,0}$, and
	$\gensmoothfunction_0$ is vanishing on $S_{0,0}$ (though its derivatives are not necessarily vanishing there).

	\medskip
	
	\noindent \underline{\textbf{Base-level equations for $\Ent$ and $v$}}.
	$\frac{\partial^{J_3}}{\partial \inisphere^{J_3}} \FS = 0$ holds
	for $0 \leq J_3 \leq K_3$
	if and only if the following equations, expressed in adapted coordinates,
	hold for $0 \leq J_3 \leq K_3$:
	\begin{align} 
	\begin{split} \label{E:ENTROPYEULEREQUATIONRELATIVETOALMOSTNULLCOORDINATES}
		\frac{\partial}{\partial \iniu} 
		\frac{\partial^{J_3}}{\partial \inisphere^{J_3}}
		\Ent 
		& 
		= 
		\gensmoothfunction_*(\inip \embed,v^i)
		\cdot
		\frac{\partial}{\partial \iniubar} 
		\frac{\partial^{J_3}}{\partial \inisphere^{J_3}} \Ent
		+
		\gensmoothfunction^A(\inip \embed,v^i)
		\cdot
		\frac{\partial}{\partial \inisphere^A} 
		\frac{\partial^{J_3}}{\partial \inisphere^{J_3}}
		\Ent
			\\
		& \ \
			+
		\gensmoothfunction\left(\inip^{\leq J_3 + 1} \embed, 
		\frac{\partial}{\partial \inisphere}^{\leq J_3} v,
		(\frac{\partial}{\partial \iniubar},\frac{\partial}{\partial \inisphere})^{\leq J_3}
		\Ent
		\right).
	\end{split}
	\end{align}
	
	Assume that $\frac{\partial^{J_3}}{\partial \inisphere^{J_3}} \FS = 0$ holds
	for $0 \leq J_3 \leq K_3$. 
	Then $\frac{\partial^{J_3}}{\partial \inisphere^{J_3}} \FV = 0$ holds
	for $0 \leq J_3 \leq K_3$ if and only if the following equations, 
	expressed in adapted coordinates, hold for $0 \leq J_3 \leq K_3$:
	\begin{align} 
	\begin{split} \label{E:VELOCITYEULEREQUATIONRELATIVETOALMOSTNULLCOORDINATES}
		\frac{\partial}{\partial \iniu} 
		\frac{\partial^{J_3}}{\partial \inisphere^{J_3}}
		v^i
		& 
		= 
		\gensmoothfunction_*(\inip \embed,v^i)
		\cdot
		\frac{\partial}{\partial \iniubar} 
		\frac{\partial^{J_3}}{\partial \inisphere^{J_3}} v^i
		+ 
		\gensmoothfunction^A(\inip \embed,v^i)
		\cdot
		\frac{\partial}{\partial \inisphere^A} 
		\frac{\partial^{J_3}}{\partial \inisphere^{J_3}} v^i
		\\
		& \ \
		+
		\gensmoothfunction\left(\inip^{\leq J_3 + 1} \embed, \inip^{\leq 1} \frac{\partial^{J_3}}{\partial \inisphere^{J_3}} \logDensity, 
		(\frac{\partial}{\partial \iniubar},\frac{\partial}{\partial \inisphere})^{\leq J_3} v, 
		(\frac{\partial}{\partial \iniubar},\frac{\partial}{\partial \inisphere})^{\leq J_3+1} \Ent 
		\right).
	\end{split}
	\end{align}
	
	\medskip

	\noindent \underline{\textbf{Base-level wave equations for $\logDensity$}}.
	Let $K_3 \geq 0$ be an integer.
	Assume that 
	$\frac{\partial^{\leq K_3 + 1}}{\partial \inisphere^{\leq K_3 + 1}} \FS = 0$,
	$\frac{\p}{\p \iniu} \frac{\partial^{\leq K_3}}{\partial \inisphere^{\leq K_3}} \FS = 0$,
	and
	$\frac{\p}{\p \iniubar} \frac{\partial^{\leq K_3}}{\partial \inisphere^{\leq K_3}} \FS = 0$,
	and $\frac{\partial^{\leq K_3}}{\partial \inisphere^{\leq K_3}} \FV^i = 0$.
	Then $\frac{\partial^{J_3}}{\partial \inisphere^{J_3}} \FW = 0$ holds
	for $0 \leq J_3 \leq K_3$
	if and only if the following equations, expressed in adapted coordinates,
	hold for $0 \leq J_3 \leq K_3$:
	\begin{align} \label{E:ADAPTEDCOORDINATESNEEDEDFORMWAVEEQUATIONFORDENSITYSALONGINITIALNULL}
			\begin{split}
				\frac{\p}{\p \iniubar}
				\frac{\p}{\p \iniu} 
				\logDensity
				& = 
				\gensmoothfunction^{AB}(\inip \embed, \vvariables) 
				\cdot 
				\frac{\p}{\p \inisphere^A} \frac{\p}{\p \inisphere^B} \logDensity
				\\
				& 
				+
				\gensmoothfunction_0(\inip \embed, \vvariables) 
				\cdot 
				\frac{\p^2}{\p \iniu^2}
				\logDensity
				+
				\gensmoothfunction_0(\inip \embed,\vvariables) 
				\cdot 
				\frac{\p^2}{\p \iniubar^2} \logDensity
				+
				\gensmoothfunction_0^A(\inip \embed, \vvariables)  
				\cdot 
				\frac{\p}{\p \inisphere^A} \frac{\p}{\p \iniu}\logDensity
				+
				\gensmoothfunction_0^A(\inip \embed, \vvariables) 
				\cdot 
				\frac{\p}{\p \inisphere^A}  \frac{\p}{\p \iniubar} \logDensity		
				\\
				& 
				+
				\mathscr{L}\left(\inip^{\leq 2} \embed,\vvariables,
				\frac{\p}{\p \iniubar} \vvariables, \frac{\p}{\p \inisphere} \vvariables, \frac{\p}{\p \iniu} \logDensity \right)
				\left[\frac{\p}{\p \iniubar} \vvariables, \frac{\p}{\p \inisphere} \vvariables, \frac{\p}{\p \iniu} \logDensity \right]
				+
				\mathscr{L}(\inip^{\leq 2} \embed,\vvariables) 
				\left[(\frac{\p}{\p \iniubar} \vvariables, \frac{\p}{\p \inisphere})^2 \Ent \right].
			\end{split}	
		\end{align}
		
	Let $K_3 \geq 0$ be an integer. 
	Assume that 
	$\frac{\partial^{\leq K_3 + 1}}{\partial \inisphere^{\leq K_3 + 1}} \FS = 0$,
	$\frac{\p}{\p \iniu} \frac{\partial^{\leq K_3}}{\partial \inisphere^{\leq K_3}} \FS = 0$,
	and
	$\frac{\p}{\p \iniubar} \frac{\partial^{\leq K_3}}{\partial \inisphere^{\leq K_3}} \FS = 0$,
	and $\frac{\partial^{\leq K_3}}{\partial \inisphere^{\leq K_3}} \FV^i = 0$.
	Then $\frac{\partial^{J_3}}{\partial \inisphere^{J_3}} \FW = 0$ holds
	for $0 \leq J_3 \leq K_3$
	if and only if the following equations, expressed in adapted coordinates,
	hold for $0 \leq J_3 \leq K_3$:
		\begin{align} \label{E:ADAPTEDCOORDINATESMWAVEEQUATIONFORDENSITYALONGINITIALNULLHBAR}
			\begin{split}
				\frac{\p}{\p \iniu}
				 \frac{\p}{\p \iniubar} 
				\logDensity
				& = 
				\gensmoothfunction^{AB}(\inip \embed, \vvariables) 
				\cdot 
				\frac{\p}{\p \inisphere^A} \frac{\p}{\p \inisphere^B} \logDensity
				\\
				& 
				+
				\gensmoothfunction_0(\inip \embed, \vvariables) 
				\cdot 
				 \frac{\p^2}{\p \iniubar^2}
				\logDensity
				+
				\gensmoothfunction_0(\inip \embed,\vvariables) 
				\cdot 
				\frac{\p^2}{\p \iniu^2} \logDensity
				+
				\gensmoothfunction_0^A(\inip \embed, \vvariables)  
				\cdot 
				\frac{\p}{\p \inisphere^A}  \frac{\p}{\p \iniubar} \logDensity
				+
				\gensmoothfunction_0^A(\inip \embed, \vvariables) 
				\cdot 
				\frac{\p}{\p \inisphere^A}  \frac{\p}{\p \iniu} \logDensity		
				\\
				& 
				+
				\mathscr{L}\left(\inip^{\leq 2} \embed,\vvariables,
				\frac{\p}{\p \iniubar} \vvariables, \frac{\p}{\p \inisphere} \vvariables,  \frac{\p}{\p \iniubar} \logDensity \right)
				\left[\frac{\p}{\p \iniu} \vvariables, \frac{\p}{\p \inisphere} \vvariables, \frac{\p}{\p \iniubar} \logDensity \right]
				+
				\mathscr{L}(\inip^{\leq 2} \embed,\vvariables) 
				\left[(\frac{\p}{\p \iniu} \vvariables, \frac{\p}{\p \inisphere})^2 \Ent \right].
			\end{split}	
		\end{align}

	\medskip
	\noindent \underline{\textbf{Higher-order equations for $\Ent$ and $v$ involving a term with a $\frac{\partial}{\partial \iniubar}$-differentiation}}.
	Assume that the following equations hold:
	\begin{itemize}
	 \item $
	\frac{\partial^{\leq K_1}}{\partial \iniu^{\leq K_1}}
	\frac{\partial^{\leq K_2}}{\partial \iniubar^{\leq K_2}}
	\frac{\partial^{\leq K_3 + 1}}{\partial \inisphere^{\leq K_3 + 1}} 
	\FS = 0$
		\item
	$
	\frac{\partial^{\leq K_1}}{\partial \iniu^{\leq K_1}}
	\frac{\partial^{\leq K_2}}{\partial \iniubar^{\leq K_2}}
	\frac{\partial^{\leq K_3}}{\partial \inisphere^{\leq K_3}} 
	\FV = 0$
	\end{itemize}
	Then 
	$\frac{\partial^{J_1}}{\partial \iniu^{J_1}} \frac{\partial^{J_2+1}}{\partial \iniubar^{J_2+1}} 
	\frac{\partial^{J_3}}{\partial \inisphere^{J_3}} \FS = 0$ 
	holds for $0 \leq J_1 \leq K_1$, $0 \leq J_2 \leq K_2$, and $0 \leq J_3 \leq K_3$ if and only if
	the following equations, 
	expressed in adapted coordinates, hold for $0 \leq J_1 \leq K_1$, $0 \leq J_2 \leq K_2$, and $0 \leq J_3 \leq K_3$:
	\begin{align} 
	\begin{split}
	\label{E:HIGHERORDERMIXEDENTROPYEULEREQUATIONRELATIVETOALMOSTNULLCOORDINATES}
		\frac{\partial^{J_1 + 1}}{\partial \iniu^{J_1 + 1}} 
		\frac{\partial^{J_2 + 1}}{\partial \iniubar^{J_2 + 1}} 
		\frac{\partial^{J_3}}{\partial \inisphere^{J_3}}
		\Ent 
		& 
		= 
		\gensmoothfunction_*(\inip \embed,v^i)
		\cdot
		\frac{\partial^{J_1 + J_2 + 2}}{\partial \iniubar^{J_1 + J_2 + 2}}
		\frac{\partial^{J_3}}{\partial \inisphere^{J_3}}
		\Ent
		+
		\gensmoothfunction(\inip \embed,v^i)
		\cdot
		\frac{\partial^{J_1 + J_2 + 1}}{\partial \iniubar^{J_1 + J_2 + 1}}
		\frac{\partial^{J_3 + 1}}{\partial \inisphere^{J_3 + 1}}
		\Ent
			\\
	& \ \
		+
		\gensmoothfunction\left(\inip^{\leq J_1 + J_2 + J_3 + 2} \embed,
		\underbrace{\frac{\partial^{\leq J_1}}{\partial \iniu^{\leq J_1}} 
		(\frac{\partial}{\partial \iniubar},\frac{\partial}{\partial \inisphere})^{\leq J_2 + J_3 + 1}
		\logDensity}_{\mbox{absent if $J_1 = 0$}},
		(\frac{\partial}{\partial \iniubar},\frac{\partial}{\partial \inisphere})^{\leq J_1 + J_2 + J_3 + 1}
		(v,\Ent)
		\right).
	\end{split}
	\end{align}

	Assume that the following equations hold:
	\begin{itemize}
	\item 
	$
	\frac{\partial^{\leq K_1}}{\partial \iniu^{\leq K_1}}
	\frac{\partial^{\leq K_2}}{\partial \iniubar^{\leq K_2}}
	\frac{\partial^{\leq K_3 + 1}}{\partial \inisphere^{\leq K_3 + 1}} 
	\FS = 0$
	\item
	 $
	\frac{\partial^{\leq K_1}}{\partial \iniu^{\leq K_1}}
	\frac{\partial^{\leq K_2}}{\partial \iniubar^{\leq K_2}}
	\frac{\partial^{\leq K_3 + 1}}{\partial \inisphere^{\leq K_3 + 1}} 
	\FV = 0$
	\item
		$
	\frac{\partial^{\leq K_1}}{\partial \iniu^{\leq K_1}}
	\frac{\partial^{\leq K_2 + 1}}{\partial \iniubar^{\leq K_2 + 1}}
	\frac{\partial^{\leq K_3}}{\partial \inisphere^{\leq K_3}} 
	\FS = 0$	\end{itemize}
	Then $
	\frac{\partial^{J_1}}{\partial \iniu^{J_1}} \frac{\partial^{J_2 + 1}}{\partial \iniubar^{J_2 + 1}} 
	\frac{\partial^{J_3}}{\partial \inisphere^{J_3}} \FV = 0$ 
	holds for $0 \leq J_1 \leq K_1$, $0 \leq J_2 \leq K_2$, and $0 \leq J_3 \leq K_3$ if and only if
	the following equations, 
	expressed in adapted coordinates, hold for $0 \leq J_1 \leq K_1$, $0 \leq J_2 \leq K_2$, and $0 \leq J_3 \leq K_3$:
	\begin{align} 
	\begin{split} \label{E:HIGHERORDERMIXEDVELOCITYEULEREQUATIONRELATIVETOALMOSTNULLCOORDINATES}
		\frac{\partial^{J_1+1}}{\partial \iniu^{J_1+1}} \frac{\partial^{J_2 + 1}}{\partial \iniubar^{J_2 + 1}} \frac{\partial^{J_3}}{\partial \inisphere^{J_3}}
		v^i 
		& 
		= 
		\gensmoothfunction_*(\inip \embed,v^i)
		\cdot
		\frac{\partial^{J_1 + J_2 + 2}}{\partial \iniubar^{J_1 + J_2 + 2}}
		\frac{\partial^{J_3}}{\partial \inisphere^{J_3}}
		v^i
		+
		\gensmoothfunction(\inip \embed,v^i)
		\cdot
		\frac{\partial^{J_1 + J_2 + 1}}{\partial \iniubar^{J_1 + J_2 +1}}
		\frac{\partial^{J_3 + 1}}{\partial \inisphere^{J_3 + 1}}
		v^i
			\\
	& \ \
		+
	\gensmoothfunction\left(
	\begin{aligned}
		&\inip^{\leq J_1 + J_2 + J_3 + 1}\embed,\;
		\frac{\partial^{\leq J_1 + 1}}{\partial \iniu^{\leq J_1 + 1}}
		\left(\frac{\partial}{\partial \iniubar},
		\frac{\partial}{\partial \inisphere}\right)^{\leq J_2 + J_3 + 1}
		\logDensity,
		\\
		&\left(\frac{\partial}{\partial \iniubar},
		\frac{\partial}{\partial \inisphere}\right)^{\leq J_1 + J_2 + J_3 + 1}
		v,\;
		\left(\frac{\partial}{\partial \iniubar},
		\frac{\partial}{\partial \inisphere}\right)^{\leq J_1 + J_2 + J_3 + 2}
		\Ent
	\end{aligned}
	\right).
	\end{split}
	\end{align}
	
	\medskip
	\noindent \underline{\textbf{Higher-order equations for $\Ent$ and $v$ involving a term with a $\frac{\partial}{\partial \iniu}$-differentiation}}.
	Assume that the following equations hold:
	\begin{itemize}
	 \item $
	\frac{\partial^{\leq K_1}}{\partial \iniu^{\leq K_1}}
	\frac{\partial^{\leq K_2}}{\partial \iniubar^{\leq K_2}}
	\frac{\partial^{\leq K_3 + 1}}{\partial \inisphere^{\leq K_3 + 1}} 
	\FS = 0$
		\item
	$
	\frac{\partial^{\leq K_1}}{\partial \iniu^{\leq K_1}}
	\frac{\partial^{\leq K_2}}{\partial \iniubar^{\leq K_2}}
	\frac{\partial^{\leq K_3}}{\partial \inisphere^{\leq K_3}} 
	\FV = 0$
	\end{itemize}
	Then 
	$\frac{\partial^{J_1 + 1}}{\partial \iniu^{J_1 + 1}} \frac{\partial^{J_2}}{\partial \iniubar^{J_2}} 
	\frac{\partial^{J_3}}{\partial \inisphere^{J_3}} \FS = 0$ 
	holds for $0 \leq J_1 \leq K_1$, $0 \leq J_2 \leq K_2$, and $0 \leq J_3 \leq K_3$ if and only if
	the following equations, 
	expressed in adapted coordinates, hold for $0 \leq J_1 \leq K_1$, $0 \leq J_2 \leq K_2$, and $0 \leq J_3 \leq K_3$:
	\begin{align} 
	\begin{split}
	\label{E:HIGHERORDERMIXEDENTROPYEULEREQUATIONRELATIVETOALMOSTNULLCOORDINATESHBAR}
		\gensmoothfunction_*(\inip \embed,v^i)
		\cdot
		\frac{\partial^{J_1 + 1}}{\partial \iniu^{J_1 + 1}} 
		\frac{\partial^{J_2 + 1}}{\partial \iniubar^{J_2 + 1}} 
		\frac{\partial^{J_3}}{\partial \inisphere^{J_3}}
		\Ent 
		& 
		= 
		\frac{\partial^{J_1 + J_2 + 2}}{\partial \iniu^{J_1 + J_2 + 2}}
		\frac{\partial^{J_3}}{\partial \inisphere^{J_3}}
		\Ent
		+
		\gensmoothfunction(\inip \embed,v^i)
		\cdot
		\frac{\partial^{J_1 + J_2 + 1}}{\partial \iniu^{J_1 + J_2 + 1}}
		\frac{\partial^{J_3 + 1}}{\partial \inisphere^{J_3 + 1}}
		\Ent
			\\
	& \ \
		+
		\gensmoothfunction\left(\begin{aligned}
			&\inip^{\leq J_1 + J_2 + J_3 + 2} \embed,
			(\frac{\partial}{\partial \iniu},\frac{\partial}{\partial \inisphere})^{\leq J_1 + J_2 + J_3 + 1}
			(v,\Ent),
			\\
		&\underbrace{\frac{\partial^{\leq J_1}}{\partial \iniu^{\leq J_1}} 
		(\frac{\partial}{\partial \iniu},\frac{\partial}{\partial \inisphere})^{\leq J_2 + J_3 + 1}
		\logDensity}_{\mbox{absent if $J_2 = 0$}}
		\end{aligned}
		\right).
	\end{split}
	\end{align}

	Assume that the following equations hold:
	\begin{itemize}
	\item 
	$
	\frac{\partial^{\leq K_1}}{\partial \iniu^{\leq K_1}}
	\frac{\partial^{\leq K_2}}{\partial \iniubar^{\leq K_2}}
	\frac{\partial^{\leq K_3 + 1}}{\partial \inisphere^{\leq K_3 + 1}} 
	\FS = 0$
	\item
	 $
	\frac{\partial^{\leq K_1}}{\partial \iniu^{\leq K_1}}
	\frac{\partial^{\leq K_2}}{\partial \iniubar^{\leq K_2}}
	\frac{\partial^{\leq K_3 + 1}}{\partial \inisphere^{\leq K_3 + 1}} 
	\FV = 0$
	\item
		$
	\frac{\partial^{\leq K_1+1}}{\partial \iniu^{\leq K_1+1}}
	\frac{\partial^{\leq K_2}}{\partial \iniubar^{\leq K_2}}
	\frac{\partial^{\leq K_3}}{\partial \inisphere^{\leq K_3}} 
	\FS = 0$	\end{itemize}
	Then $
	\frac{\partial^{J_1+1}}{\partial \iniu^{J_1+1}} \frac{\partial^{J_2}}{\partial \iniubar^{J_2}}
	\frac{\partial^{J_3}}{\partial \inisphere^{J_3}} \FV = 0$ 
	holds for $0 \leq J_1 \leq K_1$, $0 \leq J_2 \leq K_2$, and $0 \leq J_3 \leq K_3$ if and only if
	the following equations, 
	expressed in adapted coordinates, hold for $0 \leq J_1 \leq K_1$, $0 \leq J_2 \leq K_2$, and $0 \leq J_3 \leq K_3$:
	\begin{align} 
	\begin{split} \label{E:HIGHERORDERMIXEDVELOCITYEULEREQUATIONRELATIVETOALMOSTNULLCOORDINATESHBAR}
		\frac{\partial^{J_1+1}}{\partial \iniu^{J_1+1}} \frac{\partial^{J_2 + 1}}{\partial \iniubar^{J_2 + 1}} 
		\frac{\partial^{J_3}}{\partial \inisphere^{J_3}}
		v^i 
		& 
		= 
		\gensmoothfunction_*(\inip \embed,v^i)
		\cdot
		\frac{\partial^{J_1 + J_2 + 2}}{\partial \iniubar^{J_1 + J_2 + 2}}
		\frac{\partial^{J_3}}{\partial \inisphere^{J_3}}
		v^i
		+
		\gensmoothfunction(\inip \embed,v^i)
		\cdot
		\frac{\partial^{J_1 + J_2 + 1}}{\partial \iniu^{J_1 + J_2 +1}}
		\frac{\partial^{J_3 + 1}}{\partial \inisphere^{J_3 + 1}}
		v^i
			\\
	& \ \
		+
		\gensmoothfunction\left(
		\begin{aligned}
			\inip^{\leq J_1 + J_2 + J_3 + 1} \embed,
		\frac{\partial^{\leq J_1 + 1}}{\partial \iniubar^{\leq J_1 + 1}} 
		(\frac{\partial}{\partial \iniu},\frac{\partial}{\partial \inisphere})^{\leq J_2 + J_3 + 1}
		\logDensity,
		\\
		(\frac{\partial}{\partial \iniu},\frac{\partial}{\partial \inisphere})^{\leq J_1 + J_2 + J_3 + 1}
		v,
		(\frac{\partial}{\partial \iniu},\frac{\partial}{\partial \inisphere})^{\leq J_1 + J_2 + J_3 + 2}
		\Ent
		\end{aligned}
		\right).
	\end{split}
	\end{align}
	
	\medskip
	\noindent \underline{\textbf{Higher-order equations for $\Ent$ and $v$ involving a term with no $\frac{\partial}{\partial \iniubar}$-differentiation}}.
	Let $K_1 \geq 1$ and $K_3 \geq 0$ be integers.
	Assume that the following equations hold:
	\begin{itemize}
	\item 
	$
	\frac{\partial^{J_1}}{\partial \iniu^{J_1}}
	\frac{\partial^{J_2}}{\partial \iniubar^{J_2}}
	\frac{\partial^{J_3}}{\partial \inisphere^{J_3}}
	\FS = 0$
for all non-negative integers $J_1$, $J_2$, $J_3$ satisfying 
	$J_1 \leq K_1 - 1$,
	$J_2 \leq K_1$,
	and
	$J_1 + J_2 + J_3 \leq K_1 + K_3$
	\item 
	$
	\frac{\partial^{J_1}}{\partial \iniu^{J_1}}
	\frac{\partial^{J_2}}{\partial \iniubar^{J_2}}
	\frac{\partial^{J_3}}{\partial \inisphere^{J_3}}
	\FV = 0$
	for all $J_1$, $J_2$, $J_3$ satisfying 
	$J_1 \leq K_1 - 1$,
	$J_2 \leq K_1$,
	and
	$J_1 + J_2 + J_3 \leq K_1 + K_3 - 1$.
	\end{itemize}
	Then 
	$\frac{\partial^{J_1}}{\partial \iniu^{J_1}}
	\frac{\partial^{J_3}}{\partial \inisphere^{J_3}} 
	\FS = 0$ 
	holds for $1 \leq J_1 \leq K_1$ and $0 \leq J_3 \leq K_3$ if and only if
	the following equations, 
	expressed in adapted coordinates, hold for $1 \leq J_1 \leq K_1$ and $0 \leq J_3 \leq K_3$:
	\begin{align} 
		\begin{split} \label{E:HIGHERORDERPUREENTROPYEULEREQUATIONRELATIVETOALMOSTNULLCOORDINATES}
		\frac{\partial^{J_1 + 1}}{\partial \iniu^{J_1 + 1}} \frac{\partial^{J_3}}{\partial \inisphere^{J_3}}
		\Ent 
		& 
		= 
		\gensmoothfunction_*(\inip \embed,v^i)
		\cdot
		\frac{\partial^{J_1 + 1}}{\partial \iniubar^{J_1+1}}
		\frac{\partial^{J_3}}{\partial \inisphere^{J_3}}
		\Ent
		+
		\gensmoothfunction(\inip \embed,v^i)
		\cdot
		(\frac{\partial}{\partial \iniubar},\frac{\partial}{\partial \inisphere})^{J_1}
		\frac{\partial^{J_3 + 1}}{\partial \inisphere^{J_3 +1}}
		\Ent
			\\
	& \ \
		+
		\gensmoothfunction\left(\inip^{\leq J_1 + J_3 + 1} \embed,
		\inip^{\leq J_1} 
		(\frac{\partial}{\partial \iniubar},\frac{\partial}{\partial \inisphere})^{\leq J_3 }
		\logDensity,
		(\frac{\partial}{\partial \iniubar},\frac{\partial}{\partial \inisphere})^{\leq J_1 + J_3}
		(v,\Ent)
		\right).
	\end{split}
	\end{align}
		
	Let $K_1 \geq 1$ and $K_3 \geq 0$ be integers.
	Assume that the following equations hold:
	\begin{itemize}
	\item
	$
	\frac{\partial^{J_1}}{\partial \iniu^{J_1}}
	\frac{\partial^{J_2}}{\partial \iniubar^{J_2}}
	\frac{\partial^{J_3}}{\partial \inisphere^{J_3}}
	\FS = 0$
	for all non-negative integers $J_1$, $J_2$, $J_3$ satisfying 
	$J_1 \leq K_1$,
	$J_2 \leq K_1$,
	and
	$J_1 + J_2 + J_3 \leq K_1 + K_3$
	\item
	$
	\frac{\partial^{J_1}}{\partial \iniu^{J_1}}
	\frac{\partial^{J_2}}{\partial \iniubar^{J_2}}
	\frac{\partial^{J_3}}{\partial \inisphere^{J_3}}
	\FV = 0$
	holds for all non-negative integers $J_1$, $J_2$, $J_3$ satisfying 
	$J_1 \leq K_1 - 1$,
	$J_2 \leq K_1$,
	and
	$J_1 + J_2 + J_3 \leq K_1 + K_3$.
	\end{itemize}
	Then 
	$\frac{\partial^{J_1}}{\partial \iniu^{J_1}}
	\frac{\partial^{J_3}}{\partial \inisphere^{J_3}} 
	\FV = 0$ 
	holds for $1 \leq J_1 \leq K_1$ and $0 \leq J_3 \leq K_3$ if and only if
	the following equations, 
	expressed in adapted coordinates, hold for $1 \leq J_1 \leq K_1$ and $0 \leq J_3 \leq K_3$:
	\begin{align} 
	\begin{split} \label{E:HIGHERORDERPUREVELOCITYEULEREQUATIONRELATIVETOALMOSTNULLCOORDINATES}
		\frac{\partial^{J_1 + 1}}{\partial \iniu^{J_1 + 1}} 
		\frac{\partial^{J_3}}{\partial \inisphere^{J_3}}
		v^i 
		& 
		= 
		\gensmoothfunction_*(\inip \embed,v^i)
		\cdot
		\frac{\partial^{J_1 + 1}}{\partial \iniubar^{J_1 + 1}}
		\frac{\partial^{J_3}}{\partial \inisphere^{J_3}}
		v^i
		+
		\gensmoothfunction(\inip \embed,v^i)
		\cdot
		(\frac{\partial}{\partial \iniubar},\frac{\partial}{\partial \inisphere})^{J_1}
		\frac{\partial^{J_3 + 1}}{\partial \inisphere^{J_3 +1}}
		v^i
			\\
	& \ \
		+
		\gensmoothfunction\left(\inip^{\leq J_1 + J_3 + 1} \embed,
		\inip^{\leq J_1+1} 
		(\frac{\partial}{\partial \iniubar},\frac{\partial}{\partial \inisphere})^{\leq J_3}
		\logDensity,
		(\frac{\partial}{\partial \iniubar},\frac{\partial}{\partial \inisphere})^{\leq J_1 + J_3}
		v,
		(\frac{\partial}{\partial \iniubar},\frac{\partial}{\partial \inisphere})^{\leq J_1 + J_3 + 1}
		\Ent
		\right).
	\end{split}
	\end{align}

\medskip
	\noindent \underline{\textbf{Higher-order wave equations for $\logDensity$}}.
	Let $K_1, K_2, K_3$ be non-negative integers.
	Assume that the following equations hold:
	\begin{itemize}
	\item 
	$
	\frac{\partial^{\leq K_1}}{\partial \iniu^{\leq K_1}}
	\frac{\partial^{\leq K_2}}{\partial \iniubar^{\leq K_2}}
	\frac{\partial^{\leq K_3 + 1}}{\partial \inisphere^{\leq K_3 + 1}} 
	\FS = 0$
	\item
		$
	\frac{\partial^{\leq K_1+1}}{\partial \iniu^{\leq K_1+1}}
	\frac{\partial^{\leq K_2}}{\partial \iniubar^{\leq K_2}}
	\frac{\partial^{\leq K_3}}{\partial \inisphere^{\leq K_3}} 
	\FS = 0$
	\item
	$
	\frac{\partial^{\leq K_1}}{\partial \iniu^{\leq K_1}}
	\frac{\partial^{\leq K_2+2}}{\partial \iniubar^{\leq K_2+2}}
	\frac{\partial^{\leq K_3}}{\partial \inisphere^{\leq K_3}}
	$
	\item
	 $
	\frac{\partial^{\leq K_1}}{\partial \iniu^{\leq K_1}}
	\frac{\partial^{\leq K_2}}{\partial \iniubar^{\leq K_2}}
	\frac{\partial^{\leq K_3}}{\partial \inisphere^{\leq K_3}} 
	\FV = 0$
	\end{itemize}
	Then 
	$
	\frac{\partial^{J_1}}{\partial \iniu^{J_1}}
	\frac{\partial^{J_2}}{\partial \iniubar^{J_2}}
	\frac{\partial^{J_3}}{\partial \inisphere^{J_3}} \FW = 0$ 
	holds for $0 \leq J_1 \leq K_1$, $0 \leq J_2 \leq K_2$, and $0 \leq J_3 \leq K_3$ if and only if
	the following equations, 
	expressed in adapted coordinates, hold for $0 \leq J_1 \leq K_1$, $0 \leq J_2 \leq K_2$, and $0 \leq J_3 \leq K_3$:
\begin{align} \label{E:HIGHERORDERADAPTEDCOORDINATESNEEDEDFORMWAVEEQUATIONFORDENSITYSALONGINITIALNULL}
			\begin{split}
				\frac{\p^{J_2+1}}{\p \iniubar^{J_2+1}}
				\frac{\p^{J_1+1}}{\p \iniu^{J_1+1}} 
				\frac{\partial^{J_3}}{\partial \inisphere^{J_3}}
				\logDensity
				& = 
				\gensmoothfunction^{AB}(\inip \embed, \vvariables) 
				\cdot 
				\frac{\p}{\p \inisphere^A} \frac{\p}{\p \inisphere^B} 
				\frac{\partial^{J_1}}{\partial \iniu^{J_1}}
				\frac{\partial^{J_2}}{\partial \iniubar^{J_2}}
				\frac{\partial^{J_3}}{\partial \inisphere^{J_3}}
				\logDensity
				\\
				& 
				+
				\gensmoothfunction_0(\inip \embed, \vvariables) 
				\cdot 
				\frac{\partial^{J_1 + 2}}{\partial \iniu^{J_1 + 2}}
				\frac{\partial^{J_2}}{\partial \iniubar^{J_2}}
				\frac{\partial^{J_3}}{\partial \inisphere^{J_3}}
				\logDensity
				+
				\gensmoothfunction_0(\inip \embed,\vvariables) 
				\cdot 
				\frac{\partial^{J_1}}{\partial \iniu^{J_1}}
				\frac{\partial^{J_2+2}}{\partial \iniubar^{J_2+2}}
				\frac{\partial^{J_3}}{\partial \inisphere^{J_3}}
				\logDensity
					\\
			& \ \
				+
				\gensmoothfunction_0^A(\inip \embed, \vvariables)  
				\cdot 
				\frac{\p}{\p \inisphere^A}
				\frac{\partial^{J_1 + 1}}{\partial \iniu^{J_1 + 1}}
				\frac{\partial^{J_2}}{\partial \iniubar^{J_2}}
				\frac{\partial^{J_3}}{\partial \inisphere^{J_3}}
				\logDensity
				+
				\gensmoothfunction_0^A(\inip \embed, \vvariables) 
				\cdot 
				\frac{\p}{\p \inisphere^A} 
				\frac{\partial^{J_1}}{\partial \iniu^{J_1}}
				\frac{\partial^{J_2+1}}{\partial \iniubar^{J_2+1}}
				\frac{\partial^{J_3}}{\partial \inisphere^{J_3}}
				\logDensity		
				\\
				& 
				+
				\gensmoothfunction\left(
				\begin{aligned}
					&\inip^{\leq J_1 + J_2 + J_3 + 2} \embed,
		\inip^{\leq J_1+1} 
		(\frac{\partial}{\partial \iniubar},\frac{\partial}{\partial \inisphere})^{\leq J_2 + J_3}
		\logDensity,
		\\
		&(\frac{\partial}{\partial \iniubar},\frac{\partial}{\partial \inisphere})^{\leq J_1 + J_2 + J_3 + 1}
		v,
		(\frac{\partial}{\partial \iniubar},\frac{\partial}{\partial \inisphere})^{\leq J_1 + J_2 + J_3 + 2}
		\Ent
				\end{aligned}
		\right).
			\end{split}	 
		\end{align}
		
	Moreover, under the same assumptions,
	the following holds for $0 \leq J_1 \leq K_1$, $0 \leq J_2 \leq K_2$, and $0 \leq J_3 \leq K_3$:
		\begin{align} \label{E:HIGHERORDERADAPTEDCOORDINATESMWAVEEQUATIONFORDENSITYALONGINITIALNULLHBAR}
			\begin{split}
				\mbox{RHS~\eqref{E:HIGHERORDERADAPTEDCOORDINATESNEEDEDFORMWAVEEQUATIONFORDENSITYSALONGINITIALNULL}}
				& = 
				\gensmoothfunction^{AB}(\inip \embed, \vvariables) 
				\cdot 
				\frac{\p}{\p \inisphere^A} \frac{\p}{\p \inisphere^B} 
				\frac{\partial^{J_1}}{\partial \iniu^{J_1}}
				\frac{\partial^{J_2}}{\partial \iniubar^{J_2}}
				\frac{\partial^{J_3}}{\partial \inisphere^{J_3}}
				\logDensity
				\\
				& 
				+
				\gensmoothfunction_0(\inip \embed, \vvariables) 
				\cdot 
				\frac{\partial^{J_1}}{\partial \iniu^{J_1}}
				\frac{\partial^{J_2 + 2}}{\partial \iniubar^{J_2 + 2}}
				\frac{\partial^{J_3}}{\partial \inisphere^{J_3}}
				\logDensity
				+
				\gensmoothfunction_0(\inip \embed,\vvariables) 
				\cdot 
				\frac{\partial^{J_1 + 2}}{\partial \iniu^{J_1 + 2}}
				\frac{\partial^{J_2}}{\partial \iniubar^{J_2}}
				\frac{\partial^{J_3}}{\partial \inisphere^{J_3}} 
				\logDensity
					\\
			& \ \
				+
				\gensmoothfunction_0^A(\inip \embed, \vvariables)  
				\cdot 
				\frac{\p}{\p \inisphere^A}  
				\frac{\partial^{J_1}}{\partial \iniu^{J_1}}
				\frac{\partial^{J_2 + 1}}{\partial \iniubar^{J_2 + 1}}
				\frac{\partial^{J_3}}{\partial \inisphere^{J_3}}
				\logDensity
				+
				\gensmoothfunction_0^A(\inip \embed, \vvariables) 
				\cdot 
				\frac{\p}{\p \inisphere^A}  
				\frac{\partial^{J_1}+1}{\partial \iniu^{J_1+1}}
				\frac{\partial^{J_2}}{\partial \iniubar^{J_2}}
				\frac{\partial^{J_3}}{\partial \inisphere^{J_3}}
				\logDensity		
				\\
				& 
				+
				\gensmoothfunction\left(\inip^{\leq J_1 + J_2 + J_3 + 2} \embed,
		\inip^{\leq J_1+1} 
		(\frac{\partial}{\partial \iniu},\frac{\partial}{\partial \inisphere})^{\leq J_2 + J_3}
		\logDensity,
		(\frac{\partial}{\partial \iniu},\frac{\partial}{\partial \inisphere})^{\leq J_1 + J_2 + J_3 + 1}
		v,
		(\frac{\partial}{\partial \iniu},\frac{\partial}{\partial \inisphere})^{\leq J_1 + J_2 + J_3 + 2}
		\Ent
		\right).
			\end{split}	
		\end{align}

	\end{lemma}
	
	\begin{proof}
	
	\noindent \underline{\textbf{Proof of \eqref{E:ENTROPYEULEREQUATIONRELATIVETOALMOSTNULLCOORDINATES} and \eqref{E:VELOCITYEULEREQUATIONRELATIVETOALMOSTNULLCOORDINATES}}}:
	Since $\materialderivative = \partial_t + v^a \partial_a$, 
	we can express the partial derivatives in equation \eqref{E:Bs}
	with respect to $(\iniu,\iniubar,\inisphere^1,\inisphere^2)$ 
	and then schematically express the resulting equation as follows
	(cf.\ equation \eqref{E:ENTROPYEULEREQUATIONRELATIVETONULLFRAME}):
	\begin{align} \label{E:PROOFSTEPENTROPYEULEREQUATIONRELATIVETOALMOSTNULLCOORDINATES}
		\frac{\partial}{\partial \iniu} \Ent
		& 
		= 
		\gensmoothfunction_*(\inip \embed,v^i)
		\cdot
		\frac{\partial}{\partial \iniubar} \Ent
		+
		\gensmoothfunction^A(\inip \embed,v^i)
		\cdot
		\frac{\partial}{\partial \inisphere^A} \Ent.
	\end{align}
	The function $\gensmoothfunction_*$ is non-vanishing on the initial null hypersurfaces because
	$\materialderivative$ is $\gfour$-timelike and hence transversal to all null hypersurfaces.
	Since equation \eqref{E:Bs} is equivalent to $\FS = 0$,
	\eqref{E:ENTROPYEULEREQUATIONRELATIVETOALMOSTNULLCOORDINATES} now follows from differentiating
	\eqref{E:PROOFSTEPENTROPYEULEREQUATIONRELATIVETOALMOSTNULLCOORDINATES} with $\frac{\partial^{J_3}}{\partial \inisphere^{J_3}}$.
	
	Similarly, in adapted coordinates, we can express \eqref{E:Bv} as follows:
	\begin{align} 
	\begin{split} \label{E:PROOFSTEPVELOCITYEULEREQUATIONRELATIVETOALMOSTNULLCOORDINATES}
		\frac{\partial}{\partial \iniu} v^i
		& 
		= 
		\gensmoothfunction_*(\inip \embed,v)
		\cdot
		\frac{\partial}{\partial \iniubar} v^i
		+
		\gensmoothfunction^A(\inip \embed,v)
		\cdot
		\frac{\partial}{\partial \inisphere^A} v^i
			\\
	& \ \
		+ \gensmoothfunction(\inip \embed,v,\Ent) \cdot \inip \logDensity 
		+ \gensmoothfunction(\inip \embed,v,\Ent) \cdot \inip \Ent.
	\end{split}
	\end{align}
	From the definition of the quantity $\FS$, we can express (with the help of
	equation \eqref{E:PROOFSTEPENTROPYEULEREQUATIONRELATIVETOALMOSTNULLCOORDINATES}) the last 
	factor $\inip \Ent$ on RHS~\eqref{E:PROOFSTEPVELOCITYEULEREQUATIONRELATIVETOALMOSTNULLCOORDINATES}
	as $\inip \Ent 
	= 
	\gensmoothfunction(\inip \embed,v,\Ent) \cdot (\frac{\partial}{\partial \iniubar},\frac{\partial}{\partial \inisphere}) \Ent
	+ 
	\gensmoothfunction(\inip \embed,v,\Ent) \FS$, i.e., we can eliminate the derivative $\frac{\partial}{\partial \iniu} \Ent$
	up to a factor of $\FS$.
	We then use this identity to substitute for the factor $\inip \Ent$ on RHS~\eqref{E:PROOFSTEPVELOCITYEULEREQUATIONRELATIVETOALMOSTNULLCOORDINATES}.
	Differentiating the resulting equation with $\frac{\partial^{J_3}}{\partial \inisphere^{J_3}}$
	and using our assumption that
	 $\frac{\partial^{J_3}}{\partial \inisphere^{J_3}} \FS = 0$ holds
	for $0 \leq J_3 \leq K_3$, we conclude the desired equation \eqref{E:VELOCITYEULEREQUATIONRELATIVETOALMOSTNULLCOORDINATES}.
	
	\medskip
	
	\noindent \underline{\textbf{Proof of \eqref{E:HIGHERORDERMIXEDENTROPYEULEREQUATIONRELATIVETOALMOSTNULLCOORDINATES}, 
	\eqref{E:HIGHERORDERMIXEDVELOCITYEULEREQUATIONRELATIVETOALMOSTNULLCOORDINATES},
	\eqref{E:HIGHERORDERPUREENTROPYEULEREQUATIONRELATIVETOALMOSTNULLCOORDINATES},
	and \eqref{E:HIGHERORDERPUREVELOCITYEULEREQUATIONRELATIVETOALMOSTNULLCOORDINATES}}}:
	These identities follow from differentiating 
	\eqref{E:ENTROPYEULEREQUATIONRELATIVETOALMOSTNULLCOORDINATES} and \eqref{E:VELOCITYEULEREQUATIONRELATIVETOALMOSTNULLCOORDINATES}
	and eliminating all $\frac{\partial}{\partial \iniu}$-derivatives of $\Ent$ and $v$ on the RHSs of the resulting equations
	by replacing them with $(\frac{\partial}{\partial \iniubar},\frac{\partial}{\partial \inisphere})$-derivatives and factors
	of $\FS$ and $\FV$. Our assumptions on $\FS$ and $\FV$ (and their derivatives) ensure that all $\FS$- and $\FV$-involving
	terms vanish.
	
	\medskip
	\noindent \underline{\textbf{Proof of 
	\eqref{E:HIGHERORDERMIXEDENTROPYEULEREQUATIONRELATIVETOALMOSTNULLCOORDINATESHBAR}
	and \eqref{E:HIGHERORDERMIXEDVELOCITYEULEREQUATIONRELATIVETOALMOSTNULLCOORDINATESHBAR}}}:
	These follow as in the proofs 
	of \eqref{E:HIGHERORDERMIXEDENTROPYEULEREQUATIONRELATIVETOALMOSTNULLCOORDINATES}--\eqref{E:HIGHERORDERMIXEDVELOCITYEULEREQUATIONRELATIVETOALMOSTNULLCOORDINATES},,
	except we now eliminate $\frac{\partial}{\partial \iniu}$-derivatives 
	in favor of $(\frac{\partial}{\partial \iniu},\frac{\partial}{\partial \inisphere})$-derivatives.
	
	\medskip
	\noindent \underline{\textbf{Proof of 
	\eqref{E:ADAPTEDCOORDINATESNEEDEDFORMWAVEEQUATIONFORDENSITYSALONGINITIALNULL}
	and \eqref{E:ADAPTEDCOORDINATESMWAVEEQUATIONFORDENSITYALONGINITIALNULLHBAR}}}:
	\eqref{E:ADAPTEDCOORDINATESNEEDEDFORMWAVEEQUATIONFORDENSITYSALONGINITIALNULL}
	follows as in the proof of
	\eqref{E:WAVEEQUATIONFORDENSITYINADAPATEDCOORDIANTESALONGINITIALNULL},
	except we express the equation relative to the adapted coordinates rather than using the frame vectorfield $E_3$
	and use that by construction
	$\frac{\partial}{\partial \iniu}$ and $E_3$ are parallel on $S_{0,0}$.
	Similarly, \eqref{E:ADAPTEDCOORDINATESMWAVEEQUATIONFORDENSITYALONGINITIALNULLHBAR}
	follows as in the proof of \eqref{E:WAVEEQUATIONFORDENSITYINADAPATEDCOORDIANTESALONGINITIALNULLHBAR}.

\medskip
	\noindent \underline{\textbf{Proof of 
	\eqref{E:HIGHERORDERADAPTEDCOORDINATESNEEDEDFORMWAVEEQUATIONFORDENSITYSALONGINITIALNULL}
	and \eqref{E:HIGHERORDERADAPTEDCOORDINATESMWAVEEQUATIONFORDENSITYALONGINITIALNULLHBAR}}}:
	\eqref{E:HIGHERORDERADAPTEDCOORDINATESNEEDEDFORMWAVEEQUATIONFORDENSITYSALONGINITIALNULL} follows from differentiating 
	\eqref{E:ADAPTEDCOORDINATESNEEDEDFORMWAVEEQUATIONFORDENSITYSALONGINITIALNULL} and \eqref{E:ADAPTEDCOORDINATESMWAVEEQUATIONFORDENSITYALONGINITIALNULLHBAR} 
	and eliminating all $\frac{\partial}{\partial \iniu}$-derivatives of $\Ent$ and $v$ on the RHSs of the resulting equations
	by replacing them with $(\frac{\partial}{\partial \iniubar},\frac{\partial}{\partial \inisphere})$-derivatives and factors
	of $\FS$ and $\FV$. Our assumptions on $\FS$ and $\FV$ (and their derivatives) ensure that all $\FS$- and $\FV$-involving
	terms vanish. Similarly, \eqref{E:HIGHERORDERADAPTEDCOORDINATESMWAVEEQUATIONFORDENSITYALONGINITIALNULLHBAR}
	follows from differentiating \eqref{E:ADAPTEDCOORDINATESMWAVEEQUATIONFORDENSITYALONGINITIALNULLHBAR}
	and eliminating all $\frac{\partial}{\partial \iniubar}$-derivatives of $\Ent$ and $v$ on the RHSs of the resulting equations
	by replacing them with $(\frac{\partial}{\partial \iniu},\frac{\partial}{\partial \inisphere})$-derivatives.

\end{proof}

	\subsubsection{Statement of the corner compatibility conditions on the data}
	\label{SSS:CORNERCOMPATIBLITYCONDITIONSONDATA}
	In the next definition, we state the corner compatibility conditions on the data.

	\begin{remark}[No explicit corner compatibility conditions for $\logDensity$]
	\label{R:SOLUTIONJETFORDENSITYDATS00WITHOUTCOMPATIBILITY}
	In Prop.\,\ref{P:CORNERCOMPATIBILITYCONDITIONS},
	we will show that for any integers $J_1, J_2, J_3 \geq 0$,
	the quantities
	$\frac{\partial^{J_1}}{\partial \iniu^{J_1}}
	\frac{\partial^{J_2}}{\partial \iniubar^{J_2}}
	\frac{\partial^{J_3}}{\partial \inisphere^{J_3}}
	\logDensity \restriction_{S_{0,0}}$
	are already consistently determined by the proof of
	Prop.\,\ref{P:DETERMININGTRANSVERALDERIVATIVESANDSATISFYINGEQUATIONSONINITIALNULLHYPERSURFACES},
	without the need for additional $\logDensity$-corner compatibility conditions,
	as long as the $\Ent$ and $v$-corner compatibility conditions from
	Def.\,\ref{D:CORNERCOMPATIBILITYCONDITIONSONDATA} are satisfied.
	That is, we do not assume any explicit corner compatibility conditions for $\logDensity$. On the other hand, our 
	corner compatibility conditions for $v$ and $\Ent$ (see Def.\,\ref{D:CORNERCOMPATIBILITYCONDITIONSONDATA}), in conjunction
	with our assumed Euler constraint equations
	\eqref{E:constraint1}--\eqref{E:constraint2},
	enforce ``implicit'' corner compatibility conditions for $\logDensity$.
	\end{remark}

	\begin{definition}[Corner compatibility conditions on the data at $S_{0,0}$]
	\label{D:CORNERCOMPATIBILITYCONDITIONSONDATA}
	\ \\
	In this definition, it is understood that the solution jet of $\logDensity$ at $S_{0,0}$ is determined by
	Prop.\,\ref{P:DETERMININGTRANSVERALDERIVATIVESANDSATISFYINGEQUATIONSONINITIALNULLHYPERSURFACES}
	and Remark~\ref{R:SOLUTIONJETFORDENSITYDATS00WITHOUTCOMPATIBILITY}.
	
	\noindent \underline{\textbf{Base-level equations}}.
	Consider \eqref{E:ENTROPYEULEREQUATIONRELATIVETOALMOSTNULLCOORDINATES} as an equation for the initial data, and note
	that the LHS depends only on terms on $\inihypersur$ while the RHS depends only on terms on $\inihypersurbar$.
	We say that the base-level entropy corner compatibility condition is satisfied if, along $S_{0,0}$,
	the LHS evaluated at the $\inihypersur$-data agrees with the RHS evaluated at the $\inihypersurbar$-data,
	i.e., if for all integers $J_3 \geq 0$, the following equation holds on $S_{0,0}$:
	\begin{align} 
	\begin{split} \label{E:CORNERCOMPATIBILITYENTROPYEULEREQUATIONRELATIVETOALMOSTNULLCOORDINATES}
		\frac{\partial}{\partial \iniu} 
		\frac{\partial^{J_3}}{\partial \inisphere^{J_3}}
		\inient 
		& 
		= 
		\gensmoothfunction_*(\inip \embed,\inivelbar)
		\cdot
		\frac{\partial}{\partial \iniubar} 
		\frac{\partial^{J_3}}{\partial \inisphere^{J_3}} \inientbar
		+
		\gensmoothfunction^A(\inip \embed,\inivelbar)
		\cdot
		\frac{\partial}{\partial \inisphere^A} 
		\frac{\partial^{J_3}}{\partial \inisphere^{J_3}}
		\inientbar
			\\
		& \ \
			+
		\gensmoothfunction\left(\inip^{\leq J_3 + 1} \embed, 
		\frac{\partial}{\partial \inisphere}^{\leq J_3} \inivelbar,
		(\frac{\partial}{\partial \iniubar},\frac{\partial}{\partial \inisphere})^{\leq J_3}
		\inientbar
		\right).
	\end{split}
	\end{align}
	
	Similarly, 
	we say that the base-level velocity corner compatibility condition is satisfied if, along $S_{0,0}$,
	LHS~\eqref{E:VELOCITYEULEREQUATIONRELATIVETOALMOSTNULLCOORDINATES} evaluated at the $\inihypersur$-data agrees with 
	RHS~\eqref{E:VELOCITYEULEREQUATIONRELATIVETOALMOSTNULLCOORDINATES} evaluated at the $\inihypersurbar$-data,
	i.e., if for all integers $J_3 \geq 0$, the following equation holds on $S_{0,0}$:
	\begin{align} 
	\begin{split} \label{E:CORNERCOMPATIBILIITYVELOCITYEULEREQUATIONRELATIVETOALMOSTNULLCOORDINATES}
		\frac{\partial}{\partial \iniu} 
		\frac{\partial^{J_3}}{\partial \inisphere^{J_3}}
		\inivel^i
		& 
		= 
		\gensmoothfunction_*(\inip \embed,v^i)
		\cdot
		\frac{\partial}{\partial \iniubar} 
		\frac{\partial^{J_3}}{\partial \inisphere^{J_3}} \inivelbar^i
		+ 
		\gensmoothfunction^A(\inip \embed,\inivelbar^i)
		\cdot
		\frac{\partial}{\partial \inisphere^A} 
		\frac{\partial^{J_3}}{\partial \inisphere^{J_3}} \inivelbar^i
		\\
		& \ \
		+
		\gensmoothfunction\left(\inip^{\leq J_3 + 1} \embed, \inip^{\leq 1} \frac{\partial^{J_3}}{\partial \inisphere^{J_3}} \logDensity, 
		(\frac{\partial}{\partial \iniubar},\frac{\partial}{\partial \inisphere})^{\leq J_3} \inivelbar, 
		(\frac{\partial}{\partial \iniubar},\frac{\partial}{\partial \inisphere})^{\leq J_3+1} \inientbar 
		\right).
	\end{split}
	\end{align}
	
	\medskip
	
	\noindent \underline{\textbf{Higher-order equations}}.
	Similarly, we say that the higher-order entropy corner compatibility conditions are satisfied if, along $S_{0,0}$,
	LHS~\eqref{E:HIGHERORDERPUREENTROPYEULEREQUATIONRELATIVETOALMOSTNULLCOORDINATES} 
	evaluated at the $\inihypersur$-data agrees with RHS~\eqref{E:HIGHERORDERPUREENTROPYEULEREQUATIONRELATIVETOALMOSTNULLCOORDINATES} 
	evaluated at the $\inihypersurbar$-data,
	i.e., if for all integers $M, J_3 \geq 0$, the following equation holds on $S_{0,0}$:
	\begin{align} 
		\begin{split} \label{E:CORNERCOMPATIBILITYHIGHERORDERPUREENTROPYEULEREQUATIONRELATIVETOALMOSTNULLCOORDINATES}
		\frac{\partial^{M + 1}}{\partial \iniu^{M + 1}} \frac{\partial^{J_3}}{\partial \inisphere^{J_3}}
		\inient 
		& 
		= 
		\gensmoothfunction_*(\inip \embed,\inivelbar)
		\cdot
		\frac{\partial^{M + 1}}{\partial \iniubar^{M+1}}
		\frac{\partial^{J_3}}{\partial \inisphere^{J_3}}
		\inientbar
		+
		\gensmoothfunction(\inip \embed,\inivelbar)
		\cdot
		(\frac{\partial}{\partial \iniubar},\frac{\partial}{\partial \inisphere})^{M}
		\frac{\partial^{J_3 + 1}}{\partial \inisphere^{J_3 +1}}
		\inientbar
			\\
	& \ \
		+
		\gensmoothfunction\left(\inip^{\leq M + J_3 + 1} \embed,
		\inip^{\leq M} 
		(\frac{\partial}{\partial \iniubar},\frac{\partial}{\partial \inisphere})^{\leq J_3 }
		\logDensity,
		(\frac{\partial}{\partial \iniubar},\frac{\partial}{\partial \inisphere})^{\leq M + J_3}
		(\inivelbar,\inientbar)
		\right).
	\end{split}
	\end{align}
	
	Similarly, we say that the higher-order velocity corner compatibility conditions are satisfied if, along $S_{0,0}$,
	LHS~\eqref{E:HIGHERORDERPUREVELOCITYEULEREQUATIONRELATIVETOALMOSTNULLCOORDINATES} 
	evaluated at the $\inihypersur$-data agrees with RHS~\eqref{E:HIGHERORDERPUREVELOCITYEULEREQUATIONRELATIVETOALMOSTNULLCOORDINATES}  
	evaluated at the $\inihypersurbar$-data,
	i.e., if for all integers $M, J_3 \geq 0$, the following equation holds on $S_{0,0}$:
	\begin{align} 
	\begin{split} \label{E:CORNERCOMPATIBILITYHIGHERORDERPUREVELOCITYEULEREQUATIONRELATIVETOALMOSTNULLCOORDINATES}
		\frac{\partial^{M + 1}}{\partial \iniu^{M + 1}} 
		\frac{\partial^{J_3}}{\partial \inisphere^{J_3}}
		\inivel^i 
		& 
		= 
		\gensmoothfunction_*(\inip \embed,\inivelbar)
		\cdot
		\frac{\partial^{M + 1}}{\partial \iniubar^{M + 1}}
		\frac{\partial^{J_3}}{\partial \inisphere^{J_3}}
		\inivelbar^i
		+
		\gensmoothfunction(\inip \embed,v^i)
		\cdot
		(\frac{\partial}{\partial \iniubar},\frac{\partial}{\partial \inisphere})^{M}
		\frac{\partial^{J_3 + 1}}{\partial \inisphere^{J_3 +1}}
		\inivelbar^i
			\\
	& \ \
		+
		\gensmoothfunction\left(\inip^{\leq M + J_3 + 1} \embed,
		\inip^{\leq M+1} 
		(\frac{\partial}{\partial \iniubar},\frac{\partial}{\partial \inisphere})^{\leq J_3}
		\logDensity,
		(\frac{\partial}{\partial \iniubar},\frac{\partial}{\partial \inisphere})^{\leq M + J_3}
		\inivelbar,
		(\frac{\partial}{\partial \iniubar},\frac{\partial}{\partial \inisphere})^{\leq M + J_3 + 1}
		\inientbar
		\right).
	\end{split}
	\end{align}
			
	\end{definition}

	\subsubsection{The solution-jet is well-defined at $S_{0,0}$}
	\label{SSS:SOLUTIONJETWELLDEFINED}
	In the next proposition, we show that the
	corner compatibility conditions of Def.\,\ref{D:CORNERCOMPATIBILITYCONDITIONSONDATA}
	are sufficient to yield a consistently defined solution-jet at $S_{0,0}$.
	
	\begin{proposition}[The corner compatibility conditions of Def.\,\ref{D:CORNERCOMPATIBILITYCONDITIONSONDATA} yield a consistently defined solution-jet at $S_{0,0}$]
	\label{P:CORNERCOMPATIBILITYCONDITIONS}
	Assume the hypotheses and conclusions of Prop.\,\ref{P:DETERMININGTRANSVERALDERIVATIVESANDSATISFYINGEQUATIONSONINITIALNULLHYPERSURFACES}.
	In particular, by the proposition, 
	the solution jet (i.e., the solution and all of its derivatives to all orders in all directions) is uniquely determined on
	$\inihypersur_{[0,\ubarfinite]}$ and uniquely determined on $\inihypersurbar_{[0,\ufinite]}$.
	Assume that the corner compatibility conditions of Def.\,\ref{D:CORNERCOMPATIBILITYCONDITIONSONDATA} are satisfied by the initial data at 
	$S_{0,0} = \inihypersur_{[0,\ubarfinite]} \cap \inihypersurbar_{[0,\ufinite]}$.
	
	For any integer $K \geq 0$, let
	$\inip^K \inivvariables := \left(\inip^K \inient, \inip^K \inivel^i, \inip^K \inilogDensity \right)
	$ denote the array of all order $K$ partial derivatives of the solution on $\inihypersur_{[0,\ubarfinite]}$
	with respect to the adapted coordinates
	(including derivatives in the $\inihypersur_{[0,\ubarfinite]}$-transversal direction  - see \eqref{E:GENERICADAPTEDCOORDINATEPARTIALDERIVATIVES}),
	and let $\inip^K \inivvariablesbar := := \left(\inip^K \inientbar, \inip^K \inivelbar^i, \inip^K \inilogDensitybar \right)$ 
	denote the array of all order $K$ partial derivatives of the solution on $\inihypersurbar_{[0,\ufinite]}$
	with respect to the adapted coordinates
	(including derivatives in the $\inihypersurbar_{[0,\ufinite]}$-transversal direction).
	In particular, the solution jet on $\inihypersur_{[0,\ubarfinite]}$ 
	yielded by Prop.\,\ref{P:DETERMININGTRANSVERALDERIVATIVESANDSATISFYINGEQUATIONSONINITIALNULLHYPERSURFACES}
	is $\lbrace \inip^K \inivvariables \rbrace_{K \in \mathbb{Z}_{\geq 0}}$,
	and the solution jet on $\inihypersurbar_{[0,\ufinite]}$ is
	$\lbrace \inip^K \inivvariablesbar \rbrace_{K \in \mathbb{Z}_{\geq 0}}$.
	
	Let $M_1, M_2, M_3 \geq 0$ be integers.
	Then:
	\begin{subequations}
	\begin{align}
	\frac{\partial^{M_1}}{\partial \iniu^{M_1}} \frac{\partial^{M_2}}{\partial \iniubar^{M_2}} \frac{\partial^{M_3}}{\partial \inisphere^{M_3}}
	\inient \restriction_{S_{0,0}} 
		& = 
		\frac{\partial^{M_1}}{\partial \iniu^{M_1}} \frac{\partial^{M_2}}{\partial \iniubar^{M_2}} 
		\frac{\partial^{M_3}}{\partial \inisphere^{M_3}} \inientbar \restriction_{S_{0,0}},
			\label{E:ENTROPIESSAGREETOALLORDERSATS00} 
			\\
		\frac{\partial^{M_1}}{\partial \iniu^{M_1}} \frac{\partial^{M_2}}{\partial \iniubar^{M_2}} \frac{\partial^{M_3}}{\partial \inisphere^{M_3}}
		\inilogDensity \restriction_{S_{0,0}} 
		& = 
		\frac{\partial^{M_1}}{\partial \iniu^{M_1}} \frac{\partial^{M_2}}{\partial \iniubar^{M_2}} 
		\frac{\partial^{M_3}}{\partial \inisphere^{M_3}} \inilogDensitybar \restriction_{S_{0,0}},
			\label{E:DENSITIESAGREETOALLORDERSATS00} 
				\\
		\frac{\partial^{M_1}}{\partial \iniu^{M_1}} \frac{\partial^{M_2}}{\partial \iniubar^{M_2}} 
		\frac{\partial^{M_3}}{\partial \inisphere^{M_3}} \inivel^i \restriction_{S_{0,0}}
		& =
		\frac{\partial^{M_1}}{\partial \iniu^{M_1}} \frac{\partial^{M_2}}{\partial \iniubar^{M_2}} 
		\frac{\partial^{M_3}}{\partial \inisphere^{M_3}} \inivelbar^i \restriction_{S_{0,0}}.
			\label{E:VELOCITIESSAGREETOALLORDERSATS00}
	\end{align}
		\end{subequations}

	\end{proposition}
	
	\begin{proof}
		\medskip
		\noindent \underline{\textbf{Step 1: No derivatives}}.
		The cases in which $M_1 = M_2 = 0$ and $M_3 \geq 0$ follow from our data-assumption \eqref{E:ORDER0compatibilityCONDITION}.	
		
		\medskip
		\noindent \underline{\textbf{Step 2: A single null derivative of $\logDensity$}}.
		Recall that in the proof of Prop.\,\ref{P:DETERMININGTRANSVERALDERIVATIVESANDSATISFYINGEQUATIONSONINITIALNULLHYPERSURFACES},
		we derived the transversal derivatives of $\logDensity$ on $\inihypersur_{[0,\ubarfinite]}$ 
		by solving the transport equation 
		\eqref{E:ADAPTEDCOORDINATESNEEDEDFORMWAVEEQUATIONFORDENSITYSALONGINITIALNULL}
		starting from $S_{0,0}$, where the data on $S_{0,0}$ were determined
		by the data on the other hypersurface $\inihypersurbar_{[0,\ufinite]}$. 
		Similarly, we derived the transversal derivatives of
		$\logDensity$ on $\inihypersurbar_{[0,\ufinite]}$ by solving the transport equation 
		\eqref{E:ADAPTEDCOORDINATESMWAVEEQUATIONFORDENSITYALONGINITIALNULLHBAR}
		starting from $S_{0,0}$, 
		where the data were determined by the data on $\inihypersur_{[0,\ubarfinite]}$. 
		In particular, by construction, \eqref{E:DENSITIESAGREETOALLORDERSATS00} holds in the cases in which
		$M_1 = 1, M_2 = 0, M_3 \geq 0$ or $M_1 = 0, M_2 = 1, M_3 \geq 0$.
		
		\medskip
		\noindent \underline{\textbf{Step 3: A single null derivative of $\Ent$}}.
		We now consider \eqref{E:ENTROPYEULEREQUATIONRELATIVETOALMOSTNULLCOORDINATES}. 
		In the proof of Prop.\,\ref{P:DETERMININGTRANSVERALDERIVATIVESANDSATISFYINGEQUATIONSONINITIALNULLHYPERSURFACES},
		the $\inihypersur_{[0,\ubarfinite]}$-transversal term $\frac{\partial}{\partial \iniu} \Ent$ on the LHS was \emph{defined} on $\inihypersur_{[0,\ubarfinite]}$ 
		in terms of the terms on the RHS, which involve only $\inihypersur_{[0,\ubarfinite]}$-tangent derivatives.
		The compatibility equation \eqref{E:CORNERCOMPATIBILITYENTROPYEULEREQUATIONRELATIVETOALMOSTNULLCOORDINATES} ensures that at $S_{0,0}$,
		this ``definition'' of $\frac{\partial}{\partial \iniu} \Ent$ agrees with the quantity that can be computed by 
		differentiating the $\Ent$-data on $\inihypersurbar_{[0,\ubarfinite]}$ with respect to $\frac{\partial}{\partial \iniu}$.
		Similarly, the $\inihypersur_{[0,\ubarfinite]}$-transversal 
		term $\frac{\partial}{\partial \iniubar} \Ent$ on RHS~\eqref{E:ENTROPYEULEREQUATIONRELATIVETOALMOSTNULLCOORDINATES} 
		was \emph{defined} on $\inihypersur_{[0,\ubarfinite]}$
		in terms of the remaining terms in the equation, which involve only $\inihypersurbar_{[0,\ubarfinite]}$-tangent derivatives.
		The same compatibility equation \eqref{E:CORNERCOMPATIBILITYENTROPYEULEREQUATIONRELATIVETOALMOSTNULLCOORDINATES} again ensures
		that at $S_{0,0}$, this definition of $\frac{\partial}{\partial \iniubar} \Ent$ agrees with the quantity that can be computed by 
		differentiating the $\Ent$-data on $\inihypersur_{[0,\ufinite]}$ with respect to $\frac{\partial}{\partial \iniubar}$.
		We have shown that
		\eqref{E:ENTROPIESSAGREETOALLORDERSATS00} holds in the cases in which
		$M_1 = 1, M_2 = 0, M_3 \geq 0$ or $M_1 = 0, M_2 = 1, M_3 \geq 0$.
		
		\medskip
		\noindent \underline{\textbf{Step 4: A single null derivative of $v$}}.
		In our proof of Prop.\,\ref{P:DETERMININGTRANSVERALDERIVATIVESANDSATISFYINGEQUATIONSONINITIALNULLHYPERSURFACES},
		the $\inihypersur_{[0,\ubarfinite]}$-transversal term $\frac{\partial}{\partial \iniu} v^i$ on LHS~ \eqref{E:VELOCITYEULEREQUATIONRELATIVETOALMOSTNULLCOORDINATES} 
		was \emph{defined} on $\inihypersur_{[0,\ubarfinite]}$ 
		in terms of the terms on the RHS, which involve only $\inihypersur_{[0,\ubarfinite]}$-tangent derivatives
		as well as the terms $\frac{\partial}{\partial \iniubar} \Ent$
		and $\inip \logDensity$,
		which have already been shown to satisfy the corner compatibility conditions at 
		$S_{0,0}$. Analogous remarks hold for the term $\frac{\partial}{\partial \iniu} v^i$ on
		$\inihypersurbar_{[0,\ubarfinite]}$.
		Hence, as in Step 3, the velocity corner compatibility equation
		\eqref{E:CORNERCOMPATIBILIITYVELOCITYEULEREQUATIONRELATIVETOALMOSTNULLCOORDINATES}
		ensures that \eqref{E:VELOCITIESSAGREETOALLORDERSATS00}
		holds in the cases in which
		$M_1 = 1, M_2 = 0, M_3 \geq 0$ or $M_1 = 0, M_2 = 1, M_3 \geq 0$.

		\medskip
		\noindent \underline{\textbf{Step 5: Two null derivatives of $\Ent$}}.
		In the proof of Prop.\,\ref{P:DETERMININGTRANSVERALDERIVATIVESANDSATISFYINGEQUATIONSONINITIALNULLHYPERSURFACES},
		we used equation \eqref{E:HIGHERORDERMIXEDENTROPYEULEREQUATIONRELATIVETOALMOSTNULLCOORDINATES} in the case $J_1 = 0$, $J_2 = 0$ 
		to define
		$\frac{\partial}{\partial \iniu} 
		\frac{\partial}{\partial \iniubar} 
		\frac{\partial^{J_3}}{\partial \inisphere^{J_3}}
		\Ent$ 
		on $\inihypersur_{[0,\ubarfinite]}$ in terms
		of $\inihypersur_{[0,\ubarfinite]}$-tangent derivatives.
		Similarly, we used
		equation \eqref{E:HIGHERORDERMIXEDENTROPYEULEREQUATIONRELATIVETOALMOSTNULLCOORDINATESHBAR}
		in the case $J_1 = 0$, $J_2 = 0$ 
		to define
		$
		\frac{\partial}{\partial \iniu} 
		\frac{\partial}{\partial \iniubar} 
		\frac{\partial^{J_3}}{\partial \inisphere^{J_3}}
		\Ent 
		$
		on $\inihypersurbar_{[0,\ufinite]}$
		in terms of $\inihypersurbar_{[0,\ufinite]}$-tangent derivatives.
		These two expressions for $\frac{\partial}{\partial \iniu} 
		\frac{\partial}{\partial \iniubar} 
		\frac{\partial^{J_3}}{\partial \inisphere^{J_3}}$ 
		agree on $S_{0,0}$
		precisely when 
		the compatibility equation \eqref{E:CORNERCOMPATIBILITYHIGHERORDERPUREENTROPYEULEREQUATIONRELATIVETOALMOSTNULLCOORDINATES}
		in the case $M := 1$ holds.
		This yields \eqref{E:ENTROPIESSAGREETOALLORDERSATS00} in the cases in which
		$M_1 + M_2 = 2, M_3 \geq 0$.

		\medskip
		\noindent \underline{\textbf{Step 6: Two mixed null derivatives of $\logDensity$}}.
		In the proof of Prop.\,\ref{P:DETERMININGTRANSVERALDERIVATIVESANDSATISFYINGEQUATIONSONINITIALNULLHYPERSURFACES},
		when solving for transversal derivatives along $\inihypersur_{[0,\ubarfinite]}$,
		we used the equation \eqref{E:ADAPTEDCOORDINATESNEEDEDFORMWAVEEQUATIONFORDENSITYSALONGINITIALNULL}
		along $\inihypersur_{[0,\ubarfinite]}$ .
		At $S_{0,0}$, this yields an 
		expression for $\frac{\p}{\p \iniubar} \frac{\p}{\p \iniu} \logDensity$ in terms of $\inihypersur_{[0,\ubarfinite]}$-tangent
		derivatives plus terms that have already been shown to satisfy the compatibility conditions on $S_{0,0}$. 
		Similarly, when solving for transversal derivatives along $\inihypersurbar_{[0,\ufinite]}$,
		we used the equation \eqref{E:ADAPTEDCOORDINATESNEEDEDFORMWAVEEQUATIONFORDENSITYSALONGINITIALNULL}
		along $\inihypersur_{[0,\ubarfinite]}$.
		This yields an ``alternative''
		expression for $\frac{\p}{\p \iniu} \frac{\p}{\p \iniubar} \logDensity$ in terms of $\inihypersur_{[0,\ubarfinite]}$-tangent derivatives
		plus terms that have already been shown to satisfy the compatibility conditions on $S_{0,0}$. 
		The term
		$\mathscr{L}(\inip^{\leq 2} \embed,\inip^{\leq 1} \vec{E}_3,\inip^{\leq 1} \vec{E}_4,\vvariables, \FR, \FV,\FS,E_3 \FS,E_4 \FS)
			[\FR,\FV,\FS,E_3 \FS,E_4 \FS]$
		on RHS~\eqref{E:TWOWAVEEQUATIONEXPRESSIONSARESAMEUPTOLOWERORDERTERMS} vanishes on $S_{0,0}$
		because it is linear in the quantities $[\FR,\FV,\FS,E_3 \FS,E_4 \FS]$,
		which vanish at $S_{0,0}$ as a consequence of the compatibility conditions that we have shown thus far.
		Hence, 
		by \eqref{E:TWOWAVEEQUATIONEXPRESSIONSARESAMEUPTOLOWERORDERTERMS}--\eqref{E:TWOWAVEQUANTITIESHAVESOMEPRINCIPALPARTATS00},
		the two expressions for $\frac{\p}{\p \iniu} \frac{\p}{\p \iniubar} \logDensity$ agree on $S_{0,0}$.
		This yields \eqref{E:DENSITIESAGREETOALLORDERSATS00} in the case $M_1 = 1, M_2 = 1, M_3 \geq 0$.
		
		\medskip
		\noindent \underline{\textbf{Step 7: Two pure null derivatives of $\logDensity$}}.
		For the same reasons as in Step 2, our proof of Prop.\,\ref{P:DETERMININGTRANSVERALDERIVATIVESANDSATISFYINGEQUATIONSONINITIALNULLHYPERSURFACES},
		yields \eqref{E:DENSITIESAGREETOALLORDERSATS00} in the cases in which
		$M_1 = 2, M_2 = 0, M_3 \geq 0$ or $M_1 = 0, M_2 = 2, M_3 \geq 0$.
		
		\medskip
		\noindent \underline{\textbf{Step 8: Two null derivatives of $v$}}.
		Using similar arguments as in Step 5, based on equations
		\eqref{E:HIGHERORDERMIXEDVELOCITYEULEREQUATIONRELATIVETOALMOSTNULLCOORDINATES}
		and \eqref{E:HIGHERORDERMIXEDVELOCITYEULEREQUATIONRELATIVETOALMOSTNULLCOORDINATESHBAR}
		and the corner compatibility equations \eqref{E:CORNERCOMPATIBILITYHIGHERORDERPUREVELOCITYEULEREQUATIONRELATIVETOALMOSTNULLCOORDINATES},
		we can show that
		\eqref{E:VELOCITIESSAGREETOALLORDERSATS00}
		holds in the cases in which
		$M_1 + M_2 = 2, M_3 \geq 0$.
		
		\medskip
		\noindent \underline{\textbf{Step 9: Induction hypotheses for $M$ null derivatives}}.
		Let $M \geq 2$ be an integer. We make the following induction hypotheses:
		\begin{itemize}
			\item \eqref{E:ENTROPIESSAGREETOALLORDERSATS00} holds in the cases $M_1 + M_2 \leq M, M_3 \geq 0$.
			\item \eqref{E:DENSITIESAGREETOALLORDERSATS00} holds in the cases $M_1 + M_2 \leq M, M_3 \geq 0$.
			\item \eqref{E:VELOCITIESSAGREETOALLORDERSATS00}
				holds in the cases $M_1 + M_2 \leq M, M_3 \geq 0$.
		\end{itemize}
		
		\medskip
		\noindent \underline{\textbf{Step 10: $M+1$ null derivatives of $\Ent$ via induction}}.
		In the proof of Prop.\,\ref{P:DETERMININGTRANSVERALDERIVATIVESANDSATISFYINGEQUATIONSONINITIALNULLHYPERSURFACES},
		we used equation \eqref{E:HIGHERORDERMIXEDENTROPYEULEREQUATIONRELATIVETOALMOSTNULLCOORDINATES} in the cases 
		$J_1 + J_2 + 2 = M + 1$
		to define
		$
		\frac{\partial^{J_1 + 1}}{\partial \iniu^{J_1 + 1}} 
		\frac{\partial^{J_2 + 1}}{\partial \iniubar^{J_2 + 1}} 
		\frac{\partial^{J_3}}{\partial \inisphere^{J_3}}
		\Ent$ 
		on $\inihypersur_{[0,\ubarfinite]}$ in terms
		of $\inihypersur_{[0,\ubarfinite]}$-tangent derivatives.
		Similarly, we used
		equation \eqref{E:HIGHERORDERMIXEDENTROPYEULEREQUATIONRELATIVETOALMOSTNULLCOORDINATESHBAR}
		in the cases 
		$J_1 + J_2 + 2 = M + 1$
		to define
		$
		\frac{\partial^{J_1 + 1}}{\partial \iniu^{J_1 + 1}} 
		\frac{\partial^{J_2 + 1}}{\partial \iniubar^{J_2 + 1}} 
		\frac{\partial^{J_3}}{\partial \inisphere^{J_3}}
		$
		on $\inihypersurbar_{[0,\ufinite]}$
		in terms of $\inihypersurbar_{[0,\ufinite]}$-tangent derivatives.
		Considering especially the first term on RHS~\eqref{E:HIGHERORDERMIXEDENTROPYEULEREQUATIONRELATIVETOALMOSTNULLCOORDINATES} 
		and LHS~\eqref{E:HIGHERORDERMIXEDENTROPYEULEREQUATIONRELATIVETOALMOSTNULLCOORDINATESHBAR},
		we see that these two expressions for
		$\frac{\partial^{J_1 + 1}}{\partial \iniu^{J_1 + 1}} 
		\frac{\partial^{J_2 + 1}}{\partial \iniubar^{J_2 + 1}} 
		\frac{\partial^{J_3}}{\partial \inisphere^{J_3}}$ 
		agree on $S_{0,0}$
		precisely when 
		the compatibility equation \eqref{E:CORNERCOMPATIBILITYHIGHERORDERPUREENTROPYEULEREQUATIONRELATIVETOALMOSTNULLCOORDINATES}
		in the case $M$ holds.
		This yields \eqref{E:ENTROPIESSAGREETOALLORDERSATS00} in the cases in which
		$M_1 + M_2 = M+1, M_3 \geq 0$.
		
		\medskip
		\noindent \underline{\textbf{Step 11: $M+1$ mixed null derivatives of $\logDensity$ via induction}}.
		In the proof of Prop.\,\ref{P:DETERMININGTRANSVERALDERIVATIVESANDSATISFYINGEQUATIONSONINITIALNULLHYPERSURFACES},
		when solving for transversal derivatives along $\inihypersur_{[0,\ubarfinite]}$,
		we used the equation \eqref{E:HIGHERORDERADAPTEDCOORDINATESNEEDEDFORMWAVEEQUATIONFORDENSITYSALONGINITIALNULL}
		in the cases 
		$J_1 + J_2 + 2 = M + 1$
		along $\inihypersur_{[0,\ubarfinite]}$ .
		At $S_{0,0}$, this yields an 
		expression for $\frac{\p^{J_2+1}}{\p \iniubar^{J_2+1}}
				\frac{\p^{J_1+1}}{\p \iniu^{J_1+1}} 
				\frac{\partial^{J_3}}{\partial \inisphere^{J_3}}
				\logDensity$ in terms of $\inihypersur_{[0,\ubarfinite]}$-tangent
		derivatives plus terms that have already been shown to satisfy the compatibility conditions on $S_{0,0}$. 
		Similarly, when solving for transversal derivatives along $\inihypersurbar_{[0,\ufinite]}$,
		we replaced RHS~\eqref{E:HIGHERORDERADAPTEDCOORDINATESNEEDEDFORMWAVEEQUATIONFORDENSITYSALONGINITIALNULL} with
		RHS~\eqref{E:HIGHERORDERADAPTEDCOORDINATESMWAVEEQUATIONFORDENSITYALONGINITIALNULLHBAR}
		(these are equal by the induction hypotheses and \eqref{E:HIGHERORDERADAPTEDCOORDINATESMWAVEEQUATIONFORDENSITYALONGINITIALNULLHBAR}),
		which yielded an alternate equation for
		$\frac{\p^{J_2+1}}{\p \iniubar^{J_2+1}}
				\frac{\p^{J_1+1}}{\p \iniu^{J_1+1}} 
				\frac{\partial^{J_3}}{\partial \inisphere^{J_3}}
				\logDensity$
		along $\inihypersur_{[0,\ubarfinite]}$ involving $\inihypersur_{[0,\ubarfinite]}$-tangent derivatives
		plus terms that have already been shown to satisfy the compatibility conditions on $S_{0,0}$.
		Because  
		$\mbox{RHS~\eqref{E:HIGHERORDERADAPTEDCOORDINATESNEEDEDFORMWAVEEQUATIONFORDENSITYSALONGINITIALNULL}}
		=
		\mbox{RHS~\eqref{E:HIGHERORDERADAPTEDCOORDINATESMWAVEEQUATIONFORDENSITYALONGINITIALNULLHBAR}}$,
		these two expressions for
		$\frac{\p^{J_2+1}}{\p \iniubar^{J_2+1}}
				\frac{\p^{J_1+1}}{\p \iniu^{J_1+1}} 
				\frac{\partial^{J_3}}{\partial \inisphere^{J_3}}
				\logDensity$ in terms of $\inihypersur_{[0,\ubarfinite]}
		$
		agree, we conclude \eqref{E:DENSITIESAGREETOALLORDERSATS00} in the case $M_1 + M_2 = M+1, M_3 \geq 0$.

		\medskip
		\noindent \underline{\textbf{Step 12: $M+1$ pure null derivatives of $\logDensity$ via induction}}.
		For the same reasons as in Step 2, our proof of Prop.\,\ref{P:DETERMININGTRANSVERALDERIVATIVESANDSATISFYINGEQUATIONSONINITIALNULLHYPERSURFACES},
		yields \eqref{E:DENSITIESAGREETOALLORDERSATS00} in the cases in which
		$M_1 = M+1, M_2 = 0, M_3 \geq 0$ or $M_1 = 0, M_2 = M+1, M_3 \geq 0$.

		\medskip
		\noindent \underline{\textbf{Step 13: $M+1$ null derivatives of $v$ via induction}}.
		We argue as in Step 10, but now using 
		equations
		\eqref{E:HIGHERORDERMIXEDVELOCITYEULEREQUATIONRELATIVETOALMOSTNULLCOORDINATES}
		and
		\eqref{E:HIGHERORDERMIXEDVELOCITYEULEREQUATIONRELATIVETOALMOSTNULLCOORDINATESHBAR}
		and the corner compatibility equation \eqref{E:CORNERCOMPATIBILITYHIGHERORDERPUREVELOCITYEULEREQUATIONRELATIVETOALMOSTNULLCOORDINATES}.
		This yields \eqref{E:VELOCITIESSAGREETOALLORDERSATS00} in the cases in which
		$M_1 + M_2 = M+1, M_3 \geq 0$
		and completes the induction argument.

	\end{proof}

	\subsubsection{Initial data of the fluid variable $\inimodivort$}\label{SSS:inimodivort}
	Recall that we defined the fluid variable $\modivort$ in Def.\ref{Def:C}.
	To solve for the data $\inimodivort$ on $\inihypersur$, we consider the following equation:
	\begin{align}
		\label{E:reforboxgv}-\mathring{e}_4 \mathring{v}_{(1)}^i+\mathring{\anglap}\mathring{	v}^i-\frac{1}{2}\inigtr\mathring{\upchi} \mathring{v}^i_{(1)}-\frac{1}{2}\inigtr\mathring{\uchi} \mathring{e}_4 \mathring{v}^i+(\mathring{\upeta}+\mathring{\ueta})\cdot\mathring{\angnabla}\mathring{v}^i+2\mathring{\omega} \mathring{v}^i_{(1)}=-\inispeed^2e^{2\inilogDensity}\inimodivort^i+\mathfrak{Q}_{(\mathring{v})}^i+\mathfrak{L}_{(\mathring{v})}^i.
	\end{align}
	
	In the above work in this section, we have prescribed/solved for (the candidates for) $\p \mathring{v}^i$, in particular, $\mathring{v}^i_{(1)}, \mathring{e}_4\mathring{v}^i$ and $\iniangp \mathring{v}^i$, along $\inihypersur$. Therefore, $\inihypersur$-tangent derivatives of $\mathring{v}^i_{(1)},\mathring{e}_4\mathring{v}^i$ and $\iniangp \mathring{v}^i$ are also known, in particular, $\mathring{e}_4\mathring{v}_{(1)}^i$ and $\inigsphere_{AB}\frac{\p^2\mathring{v}^i}{\p \inisphere^A\p \inisphere^B}$. Therefore, all the terms in (\ref{E:reforboxgv}) except $-\inispeed^2e^{2\inilogDensity}\inimodivort^i$ are known. We now simply algebraically solve for $\inimodivort^i$ 
	via the equation \eqref{E:reforboxgv}.

	\begin{remark}\label{rmke3rho}
		At this point, $\mathring{\omega},\mathring{\upchi},\mathring{\uchi},\mathring{\ueta},\mathring{\upzeta}$ are just some quantities that can be expressed in terms of fluid data and embedding using (\ref{D:Defconnect1}). They will coincide with the connection coefficients after we finish constructing our geometric frames. In addition, we aim to solve for $\inilogDensity_{(1)}$ using (\ref{E:reforboxgdensity}). The solution $\inilogDensity_{(1)}$ serves as a suitable candidate for $\mathring{e}_3\inilogDensity$ satisfying the Euler equations (\ref{E:boxgdensity}). We emphasize that (\ref{E:reforboxgdensity}) coincides with (\ref{E:boxgdensity}) only after we prescribe the geometric coordinates and define connection coefficients. However, the very definition of the connection coefficients requires the first derivatives of the fluid variables, including $\mathring{e}_3\inilogDensity$. Therefore, at present, (\ref{E:reforboxgdensity}) should be viewed merely as a ``good equation'' that we aim to solve for $\mathring{e}_3\inilogDensity$. We will discuss its consistency with the Euler equations and geometry in Section \ref{S:geometry}.
	\end{remark}

	\subsection{Regularity Assumptions on the Initial Data in Sobolev Norms} \label{SS:SOBOLEVREGULARITYOFDATA}
	Although we have assumed that the initial data are\footnote{Our analysis yields that the corresponding solutions are 
		also $C^{\infty}$. \label{FN:SOLUTIONSARECINFINITYWHENDATAARE}} $C^{\infty}$
	for technical convenience, all of our estimates are in terms of Sobolev norms. 
	In this section, we state the assumptions on the Sobolev regularity of the initial data that we use in our basic proof of existence and uniqueness. Let $\inihypersur,\inihypersurbar$ be as in \eqref{E:INITIALNULLHYPERSURFACEIMAGESETSINR1PLUS3}, $(\mathring{v}_{(k)}^i,\inilogDensity_{(k)},\inient_{(k)})$ be as in Proposition \ref{P:DETERMININGTRANSVERALDERIVATIVESANDSATISFYINGEQUATIONSONINITIALNULLHYPERSURFACES}, and $\embeddomain_{\tilde{\embedkappa},\ubar_*,u_*}$ be as in Definition \ref{D:LDIFFEOCONSTRUCTEDOUTOFEMBEDDINGS}. Then, for some integer $N\geq 17$ and real number $\iniLtwoS > 0$,
	we assume that the following bounds hold, where we refer to 
	\eqref{DE:mixednorms} for precise definitions of the mixed-norms:
	\begin{subequations}\label{A:regularity}
		\begin{align}
			&\embed\in C^\infty(\embeddomain_{\tilde{\embedkappa},\ubar_*,u_*}),\\
			&\twonorms{\left(\frac{\p }{\p \iniubar},\iniangp\right)^k(\mathring{v}^i,\inilogDensity,\inient),\left(\frac{\p }{\p \iniubar},\iniangp\right)^{k-1}\inient_{(1)}}{\ubar}{\infty}{\gsphere}{2}{\inihypersur}\leq\iniLtwoS&
			\text{ for }& 0\leq k\leq N,\\
			&\twonorms{\left(\frac{\p }{\p \iniu},\iniangp\right)^k(\mathring{v}^i,\inilogDensity,\inient),\left(\frac{\p }{\p \iniu},\iniangp\right)^{k-1}\inient_{(1)}}{u}{\infty}{\gsphere}{2}{\inihypersurbar}\leq\iniLtwoS&
			\text{ for }& 0\leq k\leq  N,\\	
			&\twonorms{\left(\frac{\p }{\p \iniubar},\iniangp\right)^l(\mathring{v}^i,\inilogDensity,\inient),\left(\frac{\p }{\p \iniubar},\iniangp\right)^{k-1}\inient_{(1)}}{\ubar}{\infty}{\gsphere}{\infty}{\inihypersur}\leq\iniLtwoS&
			\text{ for }& 0\leq l\leq \frac{N}{2}+1,\\
			&\twonorms{\left(\frac{\p }{\p \iniu},\iniangp\right)^l(\mathring{v}^i,\inilogDensity,\inient),\left(\frac{\p }{\p \iniu},\iniangp\right)^{k-1}\inient_{(1)}}{u}{\infty}{\gsphere}{\infty}{\inihypersurbar}\leq\iniLtwoS&
			\text{ for }& 0\leq l\leq  \frac{N}{2}+1.
		\end{align}
	\end{subequations}
	Moreover, since \eqref{A:regularity} involves statements about derivatives in all directions along
	$S_{0,0} = \inihypersur \cap \inihypersurbar$, it is consistent that we make the following additional assumptions:
	\begin{subequations}
		\begin{align}
			&\norm{\pfour^k\left(\mathring{v}^i,\inilogDensity,\inient\right)}_{L_{\gsphere}^2(S_{0,0})}\leq\iniLtwoS,&
			\text{ for }& 0\leq k\leq  N,
			\\
			&\norm{\pfour^l\left(\mathring{v}^i,\inilogDensity,\inient\right)}_{L_{\gsphere}^\infty(S_{0,0})}\leq\iniLtwoS,&
			\text{ for }& 0\leq l\leq  \frac{N}{2}+1.
		\end{align}
	\end{subequations}
	
	In addition, letting $\inihypersur_{[0,\ubarfinite]},\inihypersurbar_{[0,\ufinite]}$ be a portion of $\inihypersur,\inihypersurbar$ respectively as in Proposition \ref{P:DETERMININGTRANSVERALDERIVATIVESANDSATISFYINGEQUATIONSONINITIALNULLHYPERSURFACES}, we assume the following bound for the transversal derivative of $(\mathring{v}^i,\inilogDensity,\inient)$ on $\inihypersur_{[0,\ubarfinite]}\bigcup\inihypersurbar_{[0,\ufinite]}$:
	\begin{subequations}\label{A:regularity2}
		\begin{align}
			&\twonorms{\mathring{v}_{(1)}^i,\inilogDensity_{(1)},\inient_{(2)}}{\ubar}{\infty}{\gsphere}{\infty}{\inihypersur_{[0,\ubarfinite]}}\leq\iniLtwoS,\\
			&\twonorms{\mathring{v}_{(1)}^i,\inilogDensity_{(1)},\inient_{(2)}}{u}{\infty}{\gsphere}{\infty}{\inihypersurbar_{[0,\ufinite]}}\leq\iniLtwoS.
		\end{align}
	\end{subequations}.
	
	
	\begin{remark}
		By \eqref{E:ENTROPYEULEREQUATIONWITHE3DERIVATIVE}, schematically, we can express $\pfour\inient$ in terms of tangential derivatives of $\inient$ and $(\mathring{v}^i,\inilogDensity,\inient,\pfour\embed)$ algebraically. Therefore,  \eqref{A:regularity} for $(\inivel,\inilogDensity,\inient)$ implies \eqref{A:regularity} for $\inient_{(1)}$.
	\end{remark}
	
	\begin{remark}\label{RE:singularity}
		We note that \eqref{E:NEEDEDFORMWAVEEQUATIONFORDENSITYINADAPATEDCOORDIANTESALONGINITIALNULL} could form a Riccati-type blow-up solution. In this article, we only consider the region of initial hypersurfaces $\inihypersur,\inihypersurbar\subseteq\mathbb{R}\times\mathbb{S}^2$ (with the exception of a zero-measure set) such that $(\mathring{v}^i,\inilogDensity,\inient,\mathring{v}_{(k)}^i,\inilogDensity_{(k)},\inient_{(k)})$ remain regular for all $k\in \mathbb{N}$; see Proposition \ref{P:DETERMININGTRANSVERALDERIVATIVESANDSATISFYINGEQUATIONSONINITIALNULLHYPERSURFACES} for existence of $(\mathring{v}_{(k)}^i,\inilogDensity_{(k)},\inient_{(k)})$ on a portion of $\inihypersur\bigcup\inihypersurbar$ and Section \ref{SS:regularityoffluid} for discussion of the regularities of the initial data. 
	\end{remark}
	
	\subsection{The region $\region$}\label{SS:Defregion}
	We define $\region$ as the causal future (with respect to the characteristic initial hypersurface) region, wherever the solutions to (\ref{E:Brho})-(\ref{E:Bs}) exist. In particular, in Section \ref{S:localexistence}, $\region$ contains a neighborhood of the initial sphere $S_{0,0}$. We are unable to give a precise definition of $\region$ at the moment since the acoustical eikonal functions are solely defined on $\inihypersur,\inihypersurbar$. In Section \ref{SS:SUBSETSOFSPACETIME}, we present a detailed definition of $\region$.
	\subsection{The Acoustical Quantities on the Initial Characteristic Hypersurfaces}\label{SS:InitialCharData}
	In this section, we construct acoustical eikonal functions $u,\ubar$ on $\inihypersur$ and $\inihypersurbar$, and compute tangential derivatives of these functions $u,\ubar$ with the assistance of the characteristic initial embedding and data assumed in Section \ref{SS:initialembedding2}. $u,\ubar$ exhibit the characteristics of the Euler flow. Our goal, in this section and Section \ref{S:localexistence}, is to extend acoustical functions $u,\ubar$ to the region $\region$ bounded from below by $\inihypersur$ and $\inihypersurbar$, such that they satisfy the following acoustical eikonal equations:
	\begin{align}\label{E:eikonaleq}
		(\gfour^{-1})^{\alpha\beta}\p_\alpha u\p_\beta u&=0,&
		(\gfour^{-1})^{\alpha\beta}\p_\alpha \ubar\p_\beta \ubar&=0,
	\end{align}
	and such that $u,\ubar$ coincide with the adapted coordinates 
	$(\iniu,\iniubar,\inisphere^1,\inisphere^2)$ on $\inihypersur,\inihypersurbar$ (defined in Sect.\,\ref{SS:ADAPTEDCOORDINATES}). 
	
	In the following paragraph, recalling the definition of the null lapse $\inilapse$ (\ref{D:lapse}), we define the foliation densities $\iniinvfolia$ on $\inihypersurbar$, $\iniuinvfolia$ on $\inihypersur$ under the assumptions in Section \ref{SS:initialembedding2}. $\iniinvfolia,\iniuinvfolia$ are crucial in measuring the ``closeness'' between the characteristic hypersurfaces. We need $\iniinvfolia,\iniuinvfolia,\inilapse$ to extend $u,\ubar$ from $\inihypersur$, $\inihypersurbar$ to the region $\region$.
	
	First of all, consistent with Sect.\,\ref{SS:initialembedding2}, we assume that:
	\begin{subequations}
		\begin{align}
			u:=&0,&  &\text{on }\inihypersur,\\
			\ubar:=&0,&  &\text{on }\inihypersurbar.
		\end{align}
	\end{subequations}
	Moreover, by construction in Sect.\,\ref{SS:initialembedding2},  we have the following:\\
	
	\textbf{on} $\inihypersur$:
	\begin{subequations}\label{D:nullframe01}
		\begin{align}
			\mathring{e}_4 u=&0,\\
			\mathring{e}_A u=&0,& A=&1,2.
		\end{align}
	\end{subequations}
	
	\textbf{on} $\inihypersurbar$:
	\begin{subequations}\label{D:nullframe02}
		\begin{align}
			\mathring{e}_3 \ubar=&0,\\
			\mathring{e}_A \ubar=&0,& A=&1,2.
		\end{align}
	\end{subequations}
	In addition, we define the tangential derivatives of geometric coordinates $\ubar,\theta^A$ on $\inihypersur$ and $u,\theta^A$ on $\inihypersurbar$ as follows: \\
	
	\textbf{on} $\inihypersur$:
	\begin{subequations}\label{D:nullframe1}
		\begin{align}
			\mathring{e}_4 \ubar=&2\inilapse,\\
			\mathring{e}_A\ubar=&0, &\mathring{e}_4\theta^A=&0,& \mathring{e}_A\theta^B=&\inigsphere_{AB},& A,B=&1,2.
		\end{align}
	\end{subequations}
	
	\textbf{on} $\inihypersurbar$:
	\begin{subequations}\label{D:nullframe2}
		\begin{align}
			\mathring{e}_3 u=&2\inilapse,\\
			\mathring{e}_Au=&0, &\mathring{e}_3\theta^A=&0,& \mathring{e}_A\theta^B=&\inigsphere_{AB},& A=&1,2.
		\end{align}
	\end{subequations}
	In Section \ref{SS:CONSTRUCTIONOFACOUSTICEIKONALFUNCTIONSANDGEOMETRICANGULARCOORDINATESINRENDALLREGION}, we extend $u,\ubar$ to the region $\region$, and establish the level sets of $u,\ubar$ by solving the acoustical eikonal equations (\ref{E:eikonaleq}). We need the following two definitions for extending $u,\ubar$ to $\region$.
	
	\begin{definition}[Inverse foliation densities on $\inihypersur$ and $\inihypersurbar$]\label{D:inversefoliation}	
		Rooted in (\ref{P:invuinv}), (\ref{D:decompositionofmetric}) and (\ref{D:lapse}), we define the inverse foliation densities $\iniinvfolia$ on $\inihypersurbar$, $\iniuinvfolia$ on $\inihypersur$ as follows:
		\begin{subequations}\label{D:invfoliauinvfolia}
			\begin{align}
				\label{D:iniuinvfolia}\iniuinvfolia:=&\frac{2}{\mathring{e}_3^0\mathring{e}_4\ubar}=\inilapse^{-2},&
				\text{on } &\inihypersur,\\
				\iniinvfolia:=&\frac{2}{\mathring{e}_4^0\mathring{e}_3u}=\inilapse^{-2},&
				\text{on } &\inihypersurbar.
			\end{align}
		\end{subequations}
	\end{definition}
	
	\begin{definition}
		On $\inihypersur$ and $\inihypersurbar$, we define the null geodesic vectorfields $\Lgeo$ and $\uLgeo$ as follows:
		\begin{subequations}\label{D:uLgeoLgeo}
			\begin{align}
				\label{D:Lgeo}\Lgeo^{\alpha}:=&\iniinvfolia^{-1}\inilapse^{-1}\mathring{e}_4^\alpha,&
				\text{on }&\inihypersurbar\\
				\label{D:uLgeo}\uLgeo^{\alpha}:=&\iniuinvfolia^{-1}\inilapse^{-1}\mathring{e}_3^\alpha,&
				\text{on }&\inihypersur.
			\end{align}
		\end{subequations}
		where $\inilapse$ is defined in Definition \ref{D:lapsequantities}, $\iniuinvfolia$ and $\iniinvfolia$ are defined in Definition \ref{D:inversefoliation}.
	\end{definition}
	\begin{remark}
		$\Lgeo$ on $\inihypersurbar$, as defined in (\ref{D:Lgeo}), will serve as the initial condition for the geodesic equation (\ref{E:geodesiceq2}) to extend $u$ from $\inihypersurbar$ to the region $\region$. Similarly, $\uLgeo$ on $\inihypersur$, as defined in (\ref{D:uLgeo}), will serve as the initial condition for the geodesic equation (\ref{E:LBARGEODESICINRENDALLREGION}) to extend $\ubar$ from $\inihypersur$ to the region $\region$.
	\end{remark}

	\section{Local Existence and Uniqueness of the Fluid Variables and Eikonal Functions in a Neighborhood of the Initial Sphere $S_{0,0}$}\label{S:localexistence}
	Our first main goal in this section is to prove Theorem~\ref{T:Rendall'sTheorem},
	which yields a unique solution in a future-neighborhood of $S_{0,0}$. This theorem is an analog of the main result of 
	Rendall's paper \cite{ReductionoftheCharacteristicInitialValueProblem}, which
	provided a similar theorem for quasilinear wave equations.
	As a precursor, we prove several results tied to Whitney's extension theorem; these results are needed for the proof of Theorem~\ref{T:Rendall'sTheorem}.
	Our second main goal is to prove Prop.\,\ref{P:CONSTRUCTIONOFEIKONALANDGEOMETRICANGULARCOORDINATESINRENDALLREGION}, in which we construct the acoustic geometry
	(i.e., solutions $u$ and $\ubar$ to the eikonal equations \eqref{E:eikonaleq}) in a perhaps smaller future-neighborhood of $S_{0,0}$.
	We will use all of these results in our proof of Theorem~\ref{TH:localwellposedness},
	in which we extend the region of existence to a future-neighborhood of the initial null hypersurfaces 
	(more precisely, the portion of them on which the data solution-jet is finite).
	
	\subsection{Whitney extensions}
	\label{SS:WHITNEY}
	To implement Rendall's approach, we will use the Whitney extension theorem to extend our $C^{\infty}$ 
	initial data to $C^{\infty}$ functions defined on the 
	entire spacetime $\mathbb{R}^{1+3}$. Specifically, the version of the theorem that we will use is provided by the next lemma.
	
	\begin{lemma}\cite[Appendix A. Whitney extension theorem]{transversalmappingsandflows}.\label{LE:Whitneyextension}
		Let $A \subset \mathbb{R}^n$ be a closed subset of $\mathbb{R}^{1+3}$.
		Suppose that there exists a collection $\lbrace f_{\vec{\alpha}} \rbrace_{\vec{\alpha} \in \mathbb{N}^3}$ 
		of ``Cartesian coordinate Whitney seed functions''
		$f_{\vec{\alpha}}: A \rightarrow \mathbb{R}$,
		indexed by 
		coordinate partial derivative multi-indices $\vec{\alpha}$ with length denoted by $|\vec{\alpha}|$,
		such that
		for every $\vec{\alpha}$ and every integer $m \geq 0$,
		the following ``Taylor-type'' condition is satisfied:
		\begin{align}\label{EQ:Whitney1}
			f_{\vec{\alpha}}(y)
			& =
			\sum \limits_{\abs{\vec{\beta}} \leq m - |\vec{\alpha}|}
			\frac{f_{\vec{\alpha}+\vec{\beta}}(x)}{\abs{\vec{\beta}}!}(y-x)^{\vec{\beta}}
			+
			R_{\vec{\alpha}}(x,y),
		\end{align}
		where the remainder function $R_{\vec{\alpha}} : A \times A \rightarrow \mathbb{R}$ 
		satisfies:
		\begin{align} \label{whitney0} 
			\frac{\abs{R_{\vec{\alpha}}(x,y)}}{\norm{x-y}^{m - |\vec{\alpha}|}} 
			& \rightarrow 0,
			&&
			\mbox{uniformly in $x,y$ for every fixed $z \in A$ as $x,y \rightarrow z$}.
		\end{align}
		Then $f := f_{\vec{0}}$ extends to a $C^{\infty}$ function $F: \mathbb{R}^{1+3} \rightarrow F$ 
		such that for every $\vec{\alpha}$, we have $\left.\p^{\vec{\alpha}}F\right|_A = f_{\vec{\alpha}}$.
	\end{lemma}

	\begin{lemma}[$C^{\infty}$ extension of the fluid data and Whitney seeds to $\mathbb{R}^{1+3}$]
		\label{L:Whitneyextension}
		Under the hypotheses and conclusions of Props.\,\ref{P:DETERMININGTRANSVERALDERIVATIVESANDSATISFYINGEQUATIONSONINITIALNULLHYPERSURFACES} and
		\ref{P:CORNERCOMPATIBILITYCONDITIONS},
		the $C^{\infty}$ fluid data $\inivvariables, \inivvariablesbar$ (defined on $\inihypersur, \inihypersurbar$ respectively, 
		see \eqref{E:DATAONH} and \eqref{E:DATAONHBAR})
		extend to $C^{\infty}$ functions 
		$(\mathring{v}_{(\text{ext})}^i,\inilogDensity_{(\text{ext})},\inient_{(\text{ext})})$
		defined on all of $\mathbb{R}^{1+3}$
		such that the extended functions and 
		their Cartesian coordinate partial derivatives to all orders in all directions agree with the derivatives of the data on
		$\inihypersur_{[0,\ubarfinite]} \bigcup \inihypersurbar_{[0,\ufinite]}$, where the transversal derivatives on 
		the hypersurface portions
		$\inihypersur_{[0,\ubarfinite]} \bigcup \inihypersurbar_{[0,\ufinite]}$ are
		the formal ones provided by Prop.\,\ref{P:DETERMININGTRANSVERALDERIVATIVESANDSATISFYINGEQUATIONSONINITIALNULLHYPERSURFACES}.
	\end{lemma}

	\begin{proof}[Proof of Lemma \ref{L:Whitneyextension}]
		Prop.\,\ref{P:DETERMININGTRANSVERALDERIVATIVESANDSATISFYINGEQUATIONSONINITIALNULLHYPERSURFACES} 
		implies that for any integer $K \geq 0$, the Whitney seed 
		$(\inilogDensity_{(K)},\inivel_{(K)},\inient_{(K)})$ 
		(i.e., the formal order $K$ $E_3$-derivative of the solution on $\inihypersur$) 
		is $C^{\infty}$ on the $C^{\infty}$ manifold-with-boundary $\inihypersur_{[0,\ubarfinite]}$.
		Similarly, the Whitney seed 
		$(\inilogDensitybar_{(K)},\inivelbar_{(K)},\inientbar_{(K)})$
		(i.e., the formal order $K$ $E_4$-derivative of the solution on $\inihypersurbar$) 
		is $C^{\infty}$ on the $C^{\infty}$ manifold-with-boundary $\inihypersurbar_{[0,\ufinite]}$.
		Moreover, using \eqref{D:pdecomp} and \eqref{E:PARTIALPARTIALADAPATEDUONHINTERMSOFFRAME},
		on $\inihypersur_{[0,\ubarfinite]}$, 
		we can express any Cartesian coordinate partial derivative $\partial_{\alpha}$
		as a linear combination of elements of the frame 
		$\lbrace E_3, \frac{\p}{\p \iniubar}, \frac{\p}{\p \inisphere^1}, \frac{\p}{\p \inisphere^2} \rbrace$,
		where the coefficients in the expansion are elements of $C^{\infty}(\inihypersur_{[0,\ubarfinite]})$.
		Similarly, on $\inihypersurbar_{[0,\ufinite]}$, 
		we can express any Cartesian coordinate partial derivative $\partial_{\alpha}$
		as a linear combination of elements of the frame 
		$\lbrace E_4, \frac{\p}{\p \iniu}, \frac{\p}{\p \inisphere^1}, \frac{\p}{\p \inisphere^2} \rbrace$,
		where the coefficients in the expansion are elements of $C^{\infty}($.
		By Prop.\,\ref{P:CORNERCOMPATIBILITYCONDITIONS}, these seeds are consistently defined on
		the intersection $S_{0,0} = \inihypersur_{[0,\ubarfinite]} \cap \inihypersurbar_{[0,\ufinite]})$.
		Therefore, for every Cartesian coordinate multi-index $\vec{\alpha}$, we can construct
		the Cartesian coordinate Whitney seed functions
		$(\inilogDensity_{\vec{\alpha}},\inivel_{\vec{\alpha}}^i,\inient_{\vec{\alpha}})$
		on $\inihypersur_{[0,\ubarfinite]} \bigcup \inihypersurbar_{[0,\ufinite]}$,
		i.e., the formal Cartesian coordinate $\partial^{\vec{\alpha}}$-derivative of the solution on 
		$\inihypersur_{[0,\ubarfinite]} \bigcup \inihypersurbar_{[0,\ufinite]}$.
		Since Lemma~\ref{L:EXISTENCEOFALOCALDIFFEOMORPHISM} implies that 
		$\embed$ is a diffeomorphism from adapted coordinates to Cartesian coordinates
		on the domain $\embeddomain_{\embedkappa,\ubar_*,u_*}$ such that 
		$\inihypersur_{[0,\ubarfinite]} \bigcup \inihypersurbar_{[0,\ufinite]} \subset 
		\embed(\embeddomain_{\embedkappa,\ubar_*,u_*})$,
		applications of Taylor's theorem yield that
		the conditions \eqref{EQ:Whitney1}--\eqref{whitney0} are satisfied with 
		$f_{\vec{\alpha}} := (\inilogDensity_{\vec{\alpha}},\inivel_{\vec{\alpha}}^i,\inient_{\vec{\alpha}})$.
		Hence, by Lemma~\ref{LE:Whitneyextension} with the closed set $A$ defined to be 
		$A := \inihypersur_{[0,\ubarfinite]} \bigcup \inihypersurbar_{[0,\ufinite]}$,
		we conclude the existence of the extended $C^{\infty}$ functions.
		
	\end{proof}

	\subsection{The Extended Functions Satisfy the Compressible Euler Equations at All Orders on the Initial Hypersurfaces}\label{SS:DataSatisfiesEuler}
	In the next proposition, we prove that the extended functions provided by Lemma~\ref{L:Whitneyextension} solve the first-order compressible Euler equations
	to all orders \eqref{E:firstorderformulation} on $\inihypersur \bigcup \inihypersurbar$. We need this result for the proof of
	Theorem~\ref{T:Rendall'sTheorem}.
	
	\begin{proposition}[Extended initial data is a solution to all orders on $\inihypersur\bigcup\inihypersurbar$]\label{PR:datasatisifesEuler}
		Let $(\mathring{v}_{(\text{ext})}^i,\inilogDensity_{(\text{ext})},\inient_{(\text{ext})})$ be the extended initial data functions
		(which are defined on all of $\mathbb{R}^{1+3}$)
		provided by Lemma~\ref{L:Whitneyextension}.
		Then for all integers $K \geq 0$, on $\inihypersur_{[0,\ubarfinite]} \bigcup \inihypersurbar_{[0,\ufinite]}$
		we have:
		\begin{subequations}\label{E:firstorderformulation3}
			\begin{align}
				\label{E:Brho3}\pfour^K\left(\inimaterialderivative_{(\text{ext})} \inilogDensity_{(\text{ext})} 
				+ 
				\dive \mathring{v}_{(\text{ext})} \right) & = 0,
				\\
				\label{E:Bv3}\pfour^K\left(\inimaterialderivative_{(\text{ext})} \mathring{v}_{(\text{ext})}^i
				+
				\inispeed_{(\text{ext})}^2\delta^{ia}
				\p_a \inilogDensity_{(\text{ext})}
				+
				e^{-\inilogDensity_{(\text{ext})}}\frac{\inipressure_{;\Ent}}{\backgroundDensity}\delta^{ia}\p_a \inient_{(\text{ext})} \right)&=0,\\
				\label{E:Bs3}\pfour^K\left(\inimaterialderivative_{(\text{ext})} \inient_{(\text{ext})} \right)&=0,
			\end{align}
		\end{subequations}
		where $\inimaterialderivative_{(\text{ext})}:=\p_t+\mathring{v}_{(\text{ext})}^a \p_a$,
		$\inispeed_{(\text{ext})}$ is the speed of sound evaluated at the extended initial data,
		and $\inipressure_{;\Ent}$ is the partial derivative of the equation of state with respect to entropy evaluated
		at the extended initial data.
	\end{proposition}
	
	\begin{proof}[Proof of Proposition \ref{PR:datasatisifesEuler}]
		The $C^{\infty}$ extended data $(\mathring{v}_{(\text{ext})}^i,\inilogDensity_{(\text{ext})},\inient_{(\text{ext})})$ 
		provided by Lemma~\ref{L:Whitneyextension}
		have the property that along $\inihypersur \bigcup \inihypersurbar$,
		their derivatives to all orders in all directions agree with the ones formally determined
		by the Whitney seeds we have constructed. The
		proposition therefore follows from Prop.\,\ref{P:DETERMININGTRANSVERALDERIVATIVESANDSATISFYINGEQUATIONSONINITIALNULLHYPERSURFACES},
		in particular equations  
		\eqref{E:KTHORDERTRANSVERALDERIVATIVESOFMOSTEQUATIONSSATISFIEDONH}--\eqref{E:KTHORDERTRANSVERALDERIVATIVESOFDENSITYEQUATIONONHBAR}.
		
	\end{proof}

	\subsection{Local Existence and Uniqueness of $(v^i,\logDensity,\Ent)$ in a Neighborhood of $S_{0,0}$}\label{SS:Rendallmethod}
	\begin{figure}
		\centering
		\includegraphics[width=7cm]{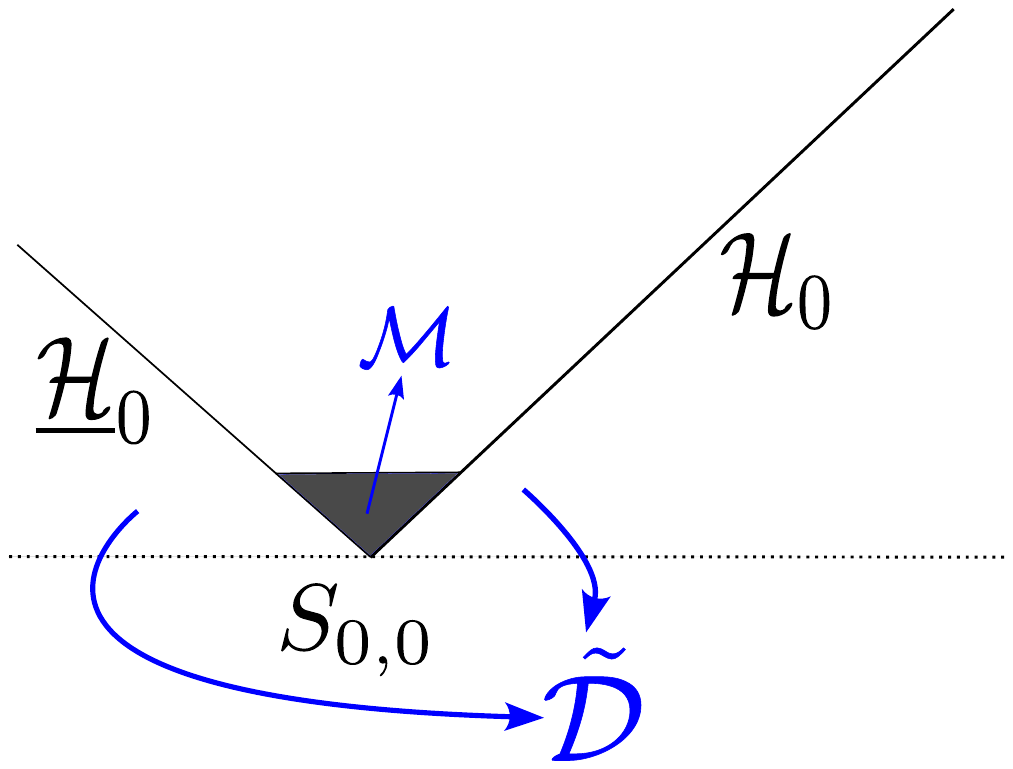}
		\caption{Region of Existence in Theorem~\ref{T:Rendall'sTheorem}}
		\label{Figure1}
	\end{figure}
	In this section, we combine the framework of \cite{ReductionoftheCharacteristicInitialValueProblem} 
	with our derivation of the solution jet on $\inihypersur_{[0,\ubarfinite]} \bigcup \inihypersurbar_{[0,\ufinite]}$
	in order 
	to prove Theorem~\ref{T:Rendall'sTheorem}, which yields 
	local existence and uniqueness of the solution $(v^i,\logDensity,\Ent)$ 
	to the characteristic initial value problem for the compressible Euler equations with initial data prescribed on $\inihypersur\bigcup\inihypersurbar$ in a future-neighborhood of the initial sphere $S_{0,0}$.
	
	\begin{theorem}[Local existence and uniqueness of the characteristic initial value problem in a future-neighborhood of $S_{0,0}$]
		\label{T:Rendall'sTheorem}
		Under the assumptions of Prop.\,\ref{P:CORNERCOMPATIBILITYCONDITIONS},
		there exists a data-dependent number 
		$\Rendallexistencetime$ satisfying $0 < \Rendallexistencetime \leq \frac{1}{2} (\ufinite + \ubarfinite)$
		and a corresponding Cartesian spacetime region $\region_{(\text{in})}$ defined by:
		\begin{align} \label{E:RENDALLREGION}
			\region_{(\text{in})}
			& := \embed 
			\left( 
			\left\lbrace 
			(\iniu,\iniubar,\inisphere) 
			\ | \ \iniu \geq 0, \iniubar \geq 0, \inisphere \in \mathbb{S}^2, 
			\frac{1}{2}(\iniu + \iniubar) \leq \Rendallexistencetime 
			\right\rbrace 
			\right),
		\end{align}
		which is a future-neighborhood 
		of the initial sphere $S_{0,0}$,
		such that on $\region_{(\text{in})}$, there exists a unique
		solution of the compressible Euler equations \eqref{E:firstorderformulation} 
		that agrees with the characteristic initial data on the initial null hypersurfaces.
	\end{theorem}
	
	\begin{proof}[Proof of Theorem~\ref{T:Rendall'sTheorem}]
		Recall that $\embed(\iniu,\iniubar,\inisphere^1,\inisphere^2)$ is the $C^{\infty}$-diffeomorphism from Lemma~\ref{L:EXISTENCEOFALOCALDIFFEOMORPHISM}
		defined on the domain $\embeddomain_{\embedkappa,\ubar_*,u_*} = 	
		[0,\ubar_*] \times [0,\tilde{\embedkappa}] \times \mathbb{S}^2
		\bigcup
		[0,\tilde{\embedkappa}] \times [0,u_*] \times \mathbb{S}^2$ from definition~\eqref{E:DIFFEOMORPHISMDOMAIN}.
		We now define the functions:
		\begin{align} \label{E:NEWCOORDINATESFORUSEINRENDALLMETHOD}
			\mathring{\uptau} 
			& := \frac{1}{2}(\iniu + \iniubar), 
			&
			\mathring{r} 
			& := \frac{1}{2}(\iniu - \iniubar),
		\end{align}
		and change coordinates from $(\iniu,\iniubar,\inisphere^1,\inisphere^2)$ to 
		$(\mathring{\uptau},\mathring{r},\inisphere^1,\inisphere^2)$.
		We view
		$(\mathring{\uptau},\mathring{r},\inisphere^1,\inisphere^2)$ to be spherical coordinates covering a compact subset $\mathfrak{K}$ 
		of $\mathbb{R}^{1+3}$ (in spherical coordinates, 
		$\mathfrak{K}$ is the image of $\embeddomain_{\embedkappa,\ubar_*,u_*}$ under the change of coordinates
		$
		(\iniu,\iniubar,\inisphere^1,\inisphere^2)
		\rightarrow
		(\mathring{\uptau},\mathring{r},\inisphere^1,\inisphere^2)
		$). 
		Suppressing changes of variables, we think of 
		$\embed$ as a $C^{\infty}$ diffeomorphism from $\mathfrak{K} \subset \mathbb{R}^{1+3}$ to $\mathbb{R}^{1+3}$,
		where the target is equipped with our standard Cartesian coordinates $(t,x^1,x^2,x^3)$ and the 
		domain is equipped with Cartesian coordinates $(\mathring{\uptau},\mathring{x}^1,\mathring{x}^2,\mathring{x}^3)$,
		where $(\mathring{x}^1,\mathring{x}^2,\mathring{x}^3)$ are standard Cartesian spatial 
		coordinates corresponding to the spherical coordinates $(\mathring{r},\inisphere^1,\inisphere^2)$.
		By Lemma~\ref{LE:Whitneyextension}, $\embed(\mathring{\uptau},\mathring{x}^1,\mathring{x}^2,\mathring{x}^3)$ 
		(we are again suppressing composition with the change of variables maps on the domain of $\embed$)
		extends to a $C^{\infty}$ function 
		$\embed_{(ext)}$
		defined on all of $\mathbb{R}^{1+3}$. By compactness, there is an open subset $\mathfrak{O} \subset \mathbb{R}^{1+3}$ containing $\mathfrak{K}$
		such that $\embed_{(ext)}$ is a diffeomorphism on $\mathfrak{O}$.
		Let $\widetilde{\Sigma}_0 : = \embed(\lbrace \mathring{\uptau} = 0\rbrace)$. 
		Note that $S_{0,0} \subset \widetilde{\Sigma}_0$,
		and that since $\embed$ is a diffeomorphism on $\mathfrak{O}$,
		there is a connected neighborhood $\mathcal{N}$ of $S_{0,0}$ such that 
		\begin{align} \label{E:SPACELIKEHYPERSURFACEPORTIONFORRENDALLMETHOD}
			\widetilde{\Sigma}_0^{\mathcal{N}} 
			& : = \mathcal{N} \cap \widetilde{\Sigma}_0
		\end{align}
		is an embedded hypersurface. 
		
		We will show that, after possibly shrinking $\mathcal{N}$ while maintaining 
		$S_{0,0} \subset \widetilde{\Sigma}_0^{\mathcal{N}}$, we can ensure that
		$\widetilde{\Sigma}_0^{\mathcal{N}}$ is $\gfour$-spacelike.
		To this end, we note that since $\mathring{\uptau} : = \frac{1}{2}(\iniu + \iniubar)$, 
		the Cartesian coordinate spacetime one-form
		$\upxi_{\alpha} 
		:= 
		\frac{1}{2}(\partial_{\alpha} \iniu + \partial_{\alpha} \iniubar) \restriction_{\widetilde{\Sigma}_0^{\mathcal{N}}}$
		is co-normal to $\widetilde{\Sigma}_0^{\mathcal{N}}$.
		Since $\inihypersur$ is a null hypersurface portion along which $\iniu = 0$
		and $\inihypersurbar$ is a null hypersurface portion along which $\iniubar = 0$,
		it follows that $(\gfour^{-1})^{\alpha \beta} \partial_{\beta} \iniu \restriction_{\inihypersur}$ 
		is proportional to the null generator $\iniLunit$ defined in \eqref{E:INITIALLUNITONH} 
		and that $(\gfour^{-1})^{\alpha \beta} \partial_{\beta} \iniubar \restriction_{\inihypersurbar}$ 
		is proportional to the null generator $\uiniuLunit$ defined in \eqref{E:INITIALULUNITONHBAR}.
		Both of these proportionality factors must be positive on a neighborhood of 
		$S_{0,0} = \inihypersur \cap \inihypersurbar$ because $\iniLunit,\uiniuLunit$ are future-directed
		and $\iniu$ and $\iniubar$ increase to the future (i.e., since 
		\eqref{EQ:JacobianofEonH}--\eqref{EQ:JacobianofEonHbar} imply that on
		$S_{0,0}$, we have $1 = \frac{\p}{\p \iniu} u = \uiniuLunit^{\alpha} \partial_{\alpha} u$
		and  $1 = \frac{\p }{\p \iniubar} \iniubar = \iniLunit^{\alpha} \partial_{\alpha} \iniubar$).
		It follows that on $S_{0,0}$, the vectorfield 
		$N^{\alpha} := (\gfour^{-1})^{\alpha \beta} \upxi_{\alpha}$, which is $\gfour$-normal to
		$\widetilde{\Sigma}_0^{\mathcal{N}}$, satisfies
		$N^{\alpha} = \frac{1}{2} (p_1 \iniLunit^{\alpha} + p_2 \uiniuLunit^{\alpha})$,
		where $p_1$ and $p_2$ are strictly positive functions. Also using that $\iniLunit$ and $\uiniuLunit$ are future-directed and null on $S_{0,0}$,
		we deduce that on $S_{0,0}$, we have $\gfour(N,N) = 2 p_1 p_2 \gfour(\iniLunit,\uiniuLunit) < 0$, i.e., that 
		$N \restriction_{S_{0,0}}$ 
		is $\gfour$-timelike.
		By continuity, $N$ is $\gfour$-timelike in a neighborhood of $S_{0,0}$. We have therefore shown that,
		after possibly shrinking $\mathcal{N}$, we can ensure that
		$\widetilde{\Sigma}_0^{\mathcal{N}}$ is $\gfour$-spacelike.
		For future use, we define the following three subsets of Cartesian coordinate space:
		\begin{subequations}
			\begin{align}
				\extendedregion_{(\text{in})} 
				& := \embed\left(\lbrace \mathring{\uptau} \geq 0 \rbrace \cap \lbrace  \iniu \geq 0 \rbrace \cap 
				\lbrace \iniubar \geq 0 \rbrace \right) \cap \mathcal{N},
				\label{E:INTERIORRENDALLREGION} \\
				\extendedregion_{(\text{out})} 
				& := \embed\left( \lbrace \mathring{\uptau} \geq 0 \rbrace 
				\cap 
				(\lbrace  \iniu < 0 \rbrace 
				\cup 
				\lbrace \iniubar < 0 \rbrace) 
				\right) \cap \mathcal{N},
				\label{E:EXTERIORRENDALLREGION} \\
				\extendedregion 
				& := 
				\extendedregion_{(\text{in})}  \cup \extendedregion_{(\text{out})}  
				=
				\embed\left( \lbrace \mathring{\uptau} \geq 0 \rbrace \right) \cap \mathcal{N}.
				\label{E:TOTALRENDALLREGION}
			\end{align}
		\end{subequations}
		If necessary, we can further shrink $\mathcal{N}$ so that 
		$\mathcal{N} \cap (\inihypersur \bigcup \inihypersurbar)
		\subset \inihypersur_{[0,\ubarfinite]} \bigcup \inihypersurbar_{[0,\ufinite]}$,
		where the data hypersurface portions $\inihypersur_{[0,\ubarfinite]}$ and $\inihypersurbar_{[0,\ufinite]}$
		are as in Lemma~\ref{L:Whitneyextension}.
		
		We now consider 
		the matrix-vector form 
		$A^{\alpha}[\vvariables]\p_\alpha\vvariables=0$
		of the first-order compressible Euler equations, as in \eqref{FOEuler}.
		Let 
		$\vvariables_{(\text{ext})} 
		:= 
		(\inilogDensity_{(\text{ext})},\mathring{v}_{(\text{ext})}^1,\mathring{v}_{(\text{ext})}^2,\mathring{v}_{(\text{ext})}^3,
		\inient_{(\text{ext})})$ denote the $C^{\infty}$ spacetime extension of the data furnished 
		by Lemma~\ref{L:Whitneyextension}.
		We aim to decompose the desired solution $\vvariables$ to \eqref{FOEuler} as follows:
		\begin{align} \label{E:SOLUTIONSPLITSINTOEXTENDEDDATAPLUSREMAINDER}
			\vvariables
			& = 
			\vvariables_{(\text{ext})}
			+
			\tilde{\vvariables}.
		\end{align}
		To this end, we consider the following initial value problem in the unknowns $\tilde{\vvariables}$:
		\begin{subequations}
			\begin{align} \label{E:SPLITFIRSTORDEREULERINTOEXTENDEDDATAPLUSREMAINDER}
				A^{\alpha}[\vvariables_{(\text{ext})}+\tilde{\vvariables}]
				\p_{\alpha}
				\left
				(\vvariables_{(\text{ext})}
				+
				\tilde{\vvariables}
				\right)
				& = F,
				\\
				\tilde{\vvariables} \restriction_{\widetilde{\Sigma}_0^{\mathcal{N}}}
				& = 0,
				\label{E:0SPACELIKEDATAFORSPLITFIRSTORDEREULERINTOEXTENDEDDATAPLUSREMAINDER}
			\end{align}
		\end{subequations}
		where
		\begin{align} \label{E:INHOMOGENEOUSTERMFORSPLITFIRSTORDEREULERINTOEXTENDEDDATAPLUSREMAINDER}
			F
			& := 
			\begin{cases}
				0, & \mbox{in } \extendedregion_{(\text{in})},
				\\
				A^{\alpha}[\vvariables_{(\text{ext})}]
				\p_{\alpha}
				\vvariables_{(\text{ext})},
				&  \mbox{in } \extendedregion_{(\text{out})}.
			\end{cases}
		\end{align}
		Crucially, from Prop.\,\ref{PR:datasatisifesEuler} and 
		\eqref{E:INHOMOGENEOUSTERMFORSPLITFIRSTORDEREULERINTOEXTENDEDDATAPLUSREMAINDER}, 
		it follows that $F$ is a $C^{\infty}$ function on $\extendedregion$ that vanishes in 
		$\extendedregion_{(\text{in})}$.
		Since $\widetilde{\Sigma}_0^{\mathcal{N}}$ is spacelike with respect to $\gfour(\vvariables_{(\text{ext})})$,
		standard local well-posedness for quasilinear hyperbolic systems yields that 
		the initial value problem 
		\eqref{E:SPLITFIRSTORDEREULERINTOEXTENDEDDATAPLUSREMAINDER}--\eqref{E:0SPACELIKEDATAFORSPLITFIRSTORDEREULERINTOEXTENDEDDATAPLUSREMAINDER}
		has a unique $C^{\infty}$ solution $\tilde{\vvariables}$ on a neighborhood $\mathcal{P}$ of $S_{0,0}$ in
		$\extendedregion$, where $\widetilde{\Sigma}_0^{\mathcal{N}} \subset \mathcal{P}$.
		Since the principal matrix on LHS~\eqref{E:SPLITFIRSTORDEREULERINTOEXTENDEDDATAPLUSREMAINDER} 
		is the compressible Euler matrix $A^{\alpha}$ evaluated at $\vvariables_{(\text{ext})}+\tilde{\vvariables}$,
		$\widetilde{\Sigma}_0^{\mathcal{N}}$
		is a $\gfour(\vvariables_{(\text{ext})} + \tilde{\vvariables})$-spacelike Cauchy hypersurface for the 
		solution domain $\mathcal{P}$.
		We also define:
		\begin{subequations}
			\begin{align}
				\region_{(\text{in})} & := \mathcal{P} \cap \extendedregion_{(\text{in})},
				\label{E:RENDALLREGIONOFLOCALEXISTENCEINSIDECONES} \\
				\region_{(\text{out})} & : = \mathcal{P} \cap \extendedregion_{(\text{out})},
				\label{E:RENDALLREGIONOFLOCALEXISTENCEOUTSIDECONES}
				\\
				\region
				& := \region_{(\text{in})} \cup \region_{(\text{out})}.
				\label{E:TOTALRENDALLREGIONINCLUDINGEXTERIOR}
			\end{align}
		\end{subequations}
		$\region_{(\text{in})}$ is the main region of existence in the theorem; see Fig.\,\ref{Figure1} for a depiction of it.
		
		By observation,
		we see that the trivial array $\tilde{\vvariables} := 0$ is a solution to the initial value problem
		\eqref{E:SPLITFIRSTORDEREULERINTOEXTENDEDDATAPLUSREMAINDER}--\eqref{E:0SPACELIKEDATAFORSPLITFIRSTORDEREULERINTOEXTENDEDDATAPLUSREMAINDER} 
		in $\region_{(\text{out})}$; by standard uniqueness in a domain of dependence for a hyperbolic PDE,
		$\tilde{\vvariables} := 0$ is in fact the unique such solution.
		Also using that the restriction of $\vvariables_{(\text{ext})}$
		to $\inihypersur_{[0,\ubarfinite]} \bigcup \inihypersurbar_{[0,\ufinite]}$ is equal to the initial data on
		$\inihypersur_{[0,\ubarfinite]} \bigcup \inihypersurbar_{[0,\ufinite]}$
		(see Lemma~\ref{L:Whitneyextension}),
		we further deduce that:
		\begin{align} \label{E:RENDALLSOLUTIONHASCORRECTDATA}
			\vvariables_{(\text{ext})}+\tilde{\vvariables}
			\restriction_{\inihypersur \cap \region}
			& = \inivvariables,
			&
			\vvariables_{(\text{ext})}+\tilde{\vvariables}
			\restriction_{\inihypersurbar \cap \region}
			& = \inivvariablesbar,
		\end{align}
		i.e., that $\vvariables_{(\text{ext})}+\tilde{\vvariables}$ is equal to the desired initial data on
		$(\inihypersur \bigcup \inihypersurbar) \cap \region_{(\text{in})}$.
		Moreover, since \eqref{E:INHOMOGENEOUSTERMFORSPLITFIRSTORDEREULERINTOEXTENDEDDATAPLUSREMAINDER} 
		and \eqref{E:RENDALLREGIONOFLOCALEXISTENCEINSIDECONES}
		imply that
		$F \restriction_{\region_{(\text{in})}} = 0$,
		we see from \eqref{E:SPLITFIRSTORDEREULERINTOEXTENDEDDATAPLUSREMAINDER}
		that $\vvariables_{(\text{ext})}+\tilde{\vvariables}$ solves the compressible Euler equations \eqref{FOEuler}
		in $\region_{(\text{in})}$. 
		To summarize,
		the array $\vvariables := \vvariables_{(\text{ext})}+\tilde{\vvariables}$ solves the compressible
		Euler equations \eqref{FOEuler} in $\region_{(\text{in})}$ and takes on the given initial data
		on $(\inihypersur \bigcup \inihypersurbar) \cap \region_{(\text{in})}$.
		
		It remains for us to prove uniqueness. Suppose we have two solutions 
		to \eqref{FOEuler} in $\region_{(\text{in})}$, denoted 
		${^{(1)}\vvariables}$, ${^{(2)}\vvariables}$,
		where $\region_{(\text{in})}$ is the region for 
		${^{(1)}\vvariables}$ generated from the existence part of the proof.
		Assume that ${^{(1)}\vvariables}$ and ${^{(2)}\vvariables}$
		both take on the given initial data on 
		$(\inihypersur \bigcup \inihypersurbar) \cap \region_{(\text{in})}$. 
		Let $\vvariables_{(\text{ext})}$ be the same $C^{\infty}$ Whitney extension of the initial data 
		from the existence part of the proof.
		Then ${^{(1)} \tilde{\vvariables}}
		:= 
		{^{(1)}\vvariables} - \vvariables_{(\text{ext})}$
		and $
		{^{(2)} \tilde{\vvariables}}
		:= 
		{^{(2)}\vvariables} - \vvariables_{(\text{ext})}$ are $C^{\infty}$ functions on $\region_{(\text{in})}$
		that vanish to all orders on 
		$(\inihypersur \bigcup \inihypersurbar) \cap \region_{(\text{in})}$.
		Hence, they can be extended to $C^{\infty}$ functions on all of $\region$
		by defining them to be $0$ in $\region_{(\text{out})}$,
		i.e., we set:
		\begin{align} \label{E:EXTENDINGTHEDIFFERENCESTOVANISHINTHEEXTERIORREGION}
			{^{(i)} \tilde{\vvariables}} \restriction_{\region_{(\text{out})}} 
			& := 0,
			& & \mbox{for } i=1,2.
		\end{align}
		It follows that each of the arrays ${^{(i)} \tilde{\vvariables}}$ is a solution to 
		the initial value problem
		\eqref{E:SPLITFIRSTORDEREULERINTOEXTENDEDDATAPLUSREMAINDER}--\eqref{E:0SPACELIKEDATAFORSPLITFIRSTORDEREULERINTOEXTENDEDDATAPLUSREMAINDER} on $\region$,
		where the inhomogeneous term
		$F$ is again defined by \eqref{E:INHOMOGENEOUSTERMFORSPLITFIRSTORDEREULERINTOEXTENDEDDATAPLUSREMAINDER}.
		Appealing to standard uniqueness for a hyperbolic PDE solution 
		with data on the spacelike surface $\widetilde{\Sigma}_0^{\mathcal{N}}$,
		we find that ${^{(1)} \tilde{\vvariables}} = {^{(2)} \tilde{\vvariables}}$ in $\region$,
		which implies in particular 
		that ${^{(1)}\vvariables} \restriction_{\region_{(\text{in})}}  
		= {^{(2)}\vvariables} \restriction_{\region_{(\text{in})}}$, as is desired.
		
	\end{proof}
	
	\subsection{The Acoustic Eikonal Functions $u,\ubar$ and the Geometric Angular Coordinates $(\sphere^1,\sphere^2)$}
	\label{SS:CONSTRUCTIONOFACOUSTICEIKONALFUNCTIONSANDGEOMETRICANGULARCOORDINATESINRENDALLREGION}
	Recall that in Section \ref{SS:ADAPTEDCOORDINATES}, we constructed
	adapted coordinates $(\iniu,\iniubar,\inisphere^1,\inisphere^2)$.
	Recall also the following limitation: even though 
	$\inihypersur$ is equal to the level set $\lbrace \iniu = 0 \rbrace$ 
	and $\inihypersurbar$ is equal to the level set $\lbrace \iniubar = 0 \rbrace$,
	the positive-value level sets of $\iniu$ and $\iniubar$ are not generally $\gfour$-null hypersurfaces.
	Hence, the level sets of $\iniu$ and $\iniubar$ are not generally good surfaces for carrying out 
	PDE analysis of the fluid solution.
	In this section, to remedy this difficulty, we construct eikonal functions
	$u$ and $\ubar$, whose level sets \emph{are} $\gfour$-null hypersurfaces. The main result is provided by 
	Prop.\,\ref{P:CONSTRUCTIONOFEIKONALANDGEOMETRICANGULARCOORDINATESINRENDALLREGION}.
	In later sections, we will derive estimates
	for the fluid on the level sets of $u$ and $\ubar$; 
	we need these estimates to prove our main existence and uniqueness result.
	
	Specifically, in Prop.\,\ref{P:CONSTRUCTIONOFEIKONALANDGEOMETRICANGULARCOORDINATESINRENDALLREGION}, 
	in a future neighborhood of $S_{0,0}$,
	we construct solutions $u$ and $\ubar$
	to the eikonal equation, which is the following fully nonlinear hyperbolic PDE whose coefficients
	depend on the acoustical metric $\gfour$ (and therefore on the fluid):
	\begin{align} \label{E:EIKONALFORENDALLREGION}
		(\gfour^{-1})^{\alpha\beta}\p_\alpha u\p_\beta u=&0,&
		(\gfour^{-1})^{\alpha\beta}\p_\alpha \ubar\p_\beta \ubar=&0.
	\end{align} 
	We assume the following initial conditions:
	\begin{align}  \label{E:INITIALCONDITIONSEIKONALFORENDALLREGION}
		u \restriction_{\inihypersur} & = \iniu_{\inihypersur},
		&
		\ubar \restriction_{\inihypersurbar} & = \iniubar_{\inihypersurbar}.
	\end{align}
	In particular, because of \eqref{E:INITIALCONDITIONSEIKONALFORENDALLREGION},
	$\inihypersur$ will be equal to the level set $\lbrace u = 0 \rbrace$ 
	and $\inihypersurbar$ will be equal to the level set $\lbrace \iniubar = 0 \rbrace$.
	
	In addition, in the same future neighborhood of $S_{0,0}$,
	we construct the geometric angular coordinates $(\sphere^1,\sphere^2)$,
	i.e., we solve the transport equation initial value problem 
	\eqref{E:TRANSPORTINGTHEANGULARCOORDINATEFUNCTIONS}.

	\begin{proposition}[Construction of the eikonal functions and geometric angular coordinates in the Rendall region]
		\label{P:CONSTRUCTIONOFEIKONALANDGEOMETRICANGULARCOORDINATESINRENDALLREGION}
		Assume the hypotheses and conclusions of Theorem~\ref{T:Rendall'sTheorem}.
		Let $\region_{(\text{in})}$ be the future-neighborhood of $S_{0,0}$ defined in \eqref{E:RENDALLREGION}, 
		on which Theorem~\ref{T:Rendall'sTheorem}
		guarantees the fluid solution exists uniquely and is $C^{\infty}$.
		There exists an open subset $\Rendallregionwithgeometry$ of $\region_{(\text{in})}$ in Cartesian coordinate space
		such that $S_{0,0} \subset \Rendallregionwithgeometry$ and such that
		on $\Rendallregionwithgeometry$, the following hold:
		\begin{itemize}
			\item There exist unique $C^{\infty}$ solutions $u,\ubar$ to the initial value problems
			\eqref{E:EIKONALFORENDALLREGION}--\eqref{E:INITIALCONDITIONSEIKONALFORENDALLREGION}
			with non-vanishing, non-parallel gradients.
			\item $\Rendallregionwithgeometry$ is foliated by level-set portions of $u$, and it is also
			foliated by level-set portions of $\ubar$.
			\item  In particular, there exists a $\timecutoffRendallregionwithgeometry > 0$ such that:
			\begin{align} \label{E:DESCRIPTIONOFLOCALREGIONNEARS00INEIKONALCOORDINATES}
				\Rendallregionwithgeometry
				& =
				\left\lbrace 
				(t,x^1,x^2,x^3)
				\ \left| \
				\begin{aligned}
					&u(t,x^1,x^2,x^3) \geq 0, 
				\
				\ubar(t,x^1,x^2,x^3) \geq 0,
				\\
				&\frac{1}{2}\left(u(t,x^1,x^2,x^3) + \ubar(t,x^1,x^2,x^3) \right) \leq \timecutoffRendallregionwithgeometry
				\end{aligned}\right.
				\right\rbrace.
			\end{align}
			\item Within $\Rendallregionwithgeometry$, the level sets of $u$ and $\ubar$ intersect transversally 
			in two-dimensional $\gfour$-spacelike submanifolds 
			$S_{u',\ubar'} := \lbrace u = u' \rbrace \cap \lbrace \ubar = \ubar' \rbrace$,
			all of which are diffeomorphic to $\mathbb{S}^2$.
			\item There exist unique $C^{\infty}$ solutions $(\sphere^1,\sphere^2)$ 
			to the transport equation initial value problem \eqref{E:TRANSPORTINGTHEANGULARCOORDINATEFUNCTIONS}
			that are local coordinates on the $S_{u',\ubar'}$.
			\item The map $\embedgeo$ defined by:
			\begin{align} \label{E:CHOVFROMGEOTOCARTESIAN}
				\embedgeo(u,\ubar,\sphere^1,\sphere^2)
				& := (t,x^1,x^2,x^3)
			\end{align}
			(i.e., the map from geometric coordinates to Cartesian coordinates)
			is a $C^{\infty}$ diffeomorphism from the domain
			$\lbrace (u,\ubar,\sphere) \ | \ u \geq 0, \ubar \geq 0, 
			\frac{1}{2}(u + \ubar) \leq \timecutoffRendallregionwithgeometry \rbrace$
			onto the image set $\Rendallregionwithgeometry$.
		\end{itemize}
	\end{proposition}

	\begin{proof}
		Let $\region_{(\text{in})}$ be the future-neighborhood of $S_{0,0}$ 
		provided by Theorem~\ref{T:Rendall'sTheorem}.
		Throughout this proof, we will silently use the results of Lemma~\ref{L:EXISTENCEOFALOCALDIFFEOMORPHISM},
		which shows that the adapted coordinates and Cartesian coordinates are in diffeomorphic correspondence on
		$\region_{(\text{in})}$; we can therefore use both coordinate systems in our analysis.
		
		We start by constructing $\ubar$ using a standard construction-by-null geodesics.
		To proceed, we recall the definition \eqref{D:uLgeo} of $\uLgeo$, which is a $\gfour$-null vectorfield
		transversal to $\inihypersur$. Recall that along $S_{0,0}$, $\uLgeo$ is equal to the null generator 
		$\uiniuLunit^{\alpha}$ of $\inihypersurbar$ (see \eqref{E:INITIALULUNITONHBAR}).
		For each $\iniubar,\inisphere$ and each corresponding 
		point $\embedding(\iniubar,\inisphere) \in \inihypersur \cap \region_{(\text{in})}$, 
		we solve the following null geodesic ODE system for the curve
		$\upsigma \rightarrow \underline{\Upsilon}_{\ubar,\theta}^\alpha(\upsigma)$:
		\begin{subequations}\label{E:LBARGEODESICINRENDALLREGION}
			\begin{align}
				\ddot{\underline{\Upsilon}}_{\iniubar,\inisphere}^\alpha(\upsigma)
				&
				= -\left.\Chfour^{\ \alpha}_{\kappa \ \lambda}\right|_{\underline{\Upsilon}_{\ubar,\theta}(\upsigma)}\dot{\underline{\Upsilon}}_{\ubar,\theta}^\kappa(\upsigma)\dot{\underline{\Upsilon}}_{\iniubar,\inisphere}^\lambda(\upsigma),
				\\
				\underline{\Upsilon}_{\ubar,\theta}^\alpha(0)
				& 
				= 
				\embedding(\iniubar,\inisphere),
				&
				\dot{\underline{\Upsilon}}_{\iniubar,\inisphere}^\alpha(0)
				& = 
				\left.\uLgeo^\alpha\right|_{\embedding(\iniubar,\inisphere)},
				\label{E:INITAILCONDITIONSFORLBARGEODESICINRENDALLREGION}
			\end{align} 
		\end{subequations}
		where 
		$\ddot{\underline{\Upsilon}}_{\iniubar,\inisphere}^\alpha
		: =\frac{\diff^2}{\diff \upsigma^2}\underline{\Upsilon}_{\iniubar,\inisphere}^\alpha$, 
		$\dot{\underline{\Upsilon}}_{\iniubar,\inisphere}^{\alpha}
		:=\frac{\diff}{\diff \upsigma}\underline{\Upsilon}_{\iniubar,\inisphere}^\alpha$, 
		$\Chfour^{\ \alpha}_{\kappa \ \lambda}
		=
		\frac{1}{2}
		(\gfour^{-1})^{\alpha \delta}
		\left\lbrace
		\partial_{\kappa} \gfour_{\delta \lambda}
		+
		\partial_{\lambda} \gfour_{\kappa \delta}
		-
		\partial_{\delta} \gfour_{\kappa \lambda}
		\right\rbrace
		$ 
		denotes the Cartesian Christoffel symbols of $\gfour$.
		By definition, $\upsigma$ is an affine parameter for the null geodesics.
		Theorem~\ref{T:Rendall'sTheorem} yields that the fluid is $C^{\infty}$ in
		$\region_{(\text{in})}$, and thus the terms 
		$\Chfour^{\ \alpha}_{\kappa \ \lambda}$ on RHS~\eqref{E:LBARGEODESICINRENDALLREGION} 
		are also $C^{\infty}$.
		Therefore, the standard existence + uniqueness + smooth dependence on initial conditions theorem for ODEs
		implies that 
		\eqref{E:LBARGEODESICINRENDALLREGION}--\eqref{E:INITAILCONDITIONSFORLBARGEODESICINRENDALLREGION}
		has a unique local $C^{\infty}$ solution that is also a $C^{\infty}$ function of 
		$\iniubar,\inisphere$.
		Since \eqref{E:LBARGEODESICINRENDALLREGION}--\eqref{E:INITAILCONDITIONSFORLBARGEODESICINRENDALLREGION}
		is the null geodesic equation, 
		straightforward computations yield that:
		\begin{align} \label{E:GEODESICSOLUTIONISNULLHASVANISHINGDERIVATIVE}
			\frac{d}{d \upsigma}
			\left(
			\gfour_{\alpha \beta}|_{\underline{\Upsilon}_{\iniubar,\inisphere}(\upsigma)} 
			\dot{\underline{\Upsilon}}_{\iniubar,\inisphere}^{\alpha}(\upsigma)
			\dot{\underline{\Upsilon}}_{\iniubar,\inisphere}^{\beta}(\upsigma)
			\right)
			& = 0
		\end{align}
		for all $\upsigma$ in the domain of maximal existence of the solution.
		In view of the second initial condition in \eqref{E:INITAILCONDITIONSFORLBARGEODESICINRENDALLREGION},
		we deduce from \eqref{E:GEODESICSOLUTIONISNULLHASVANISHINGDERIVATIVE} that:
		\begin{align} \label{E:GEODESICSOLUTIONISEVERYWHERENULL}
			\gfour_{\alpha \beta}|_{\underline{\Upsilon}_{\ubar,\theta}(\upsigma)} 
			\dot{\underline{\Upsilon}}_{\iniubar,\inisphere}^{\alpha}(\upsigma)
			\dot{\underline{\Upsilon}}_{\iniubar,\inisphere}^{\beta}(\upsigma)
			& = 0,
		\end{align}
		i.e.,
		the vector $\dot{\underline{\Upsilon}}_{\iniubar,\inisphere}^{\alpha}(\upsigma)$
		is everywhere $\gfour$-null. Similarly, 
		again using the geodesic equation \eqref{E:LBARGEODESICINRENDALLREGION}
		as well as the identity that results from differentiating \eqref{E:GEODESICSOLUTIONISEVERYWHERENULL}
		with $\frac{\partial}{\partial \inisphere^A}$, 
		we compute that:
		\begin{align} \label{E:EVOLUTIONEQUATIONFORGEODESICSOLUTIONISEVERYWHEREORTHOGONALTOANGULARDIRECTIONS}
			\frac{d}{d \upsigma}
			\left(
			\gfour_{\alpha \beta}|_{\underline{\Upsilon}_{\iniubar,\inisphere}(\upsigma)} 
			\dot{\underline{\Upsilon}}_{\iniubar,\inisphere}^{\alpha}(\upsigma)
			\frac{\partial}{\partial \inisphere^A}
			\underline{\Upsilon}_{\iniubar,\inisphere}^{\beta}(\upsigma)
			\right)
			& = 0,
			&  A & = 1,2.
		\end{align}
		From \eqref{E:EVOLUTIONEQUATIONFORGEODESICSOLUTIONISEVERYWHEREORTHOGONALTOANGULARDIRECTIONS}
		and the initial condition
		$\gfour_{\alpha \beta}|_{\underline{\Upsilon}_{\iniubar,\inisphere}(0)} 
		\dot{\underline{\Upsilon}}_{\iniubar,\inisphere}^{\alpha}(0)
		\frac{\partial}{\partial \inisphere^A}
		\underline{\Upsilon}_{\iniubar,\inisphere}^{\beta}(0)$
		(which holds because $\uLgeo$ is $\gfour$-orthogonal to $S_{0,0})$,
		we deduce that:
		\begin{align} \label{E:GEODESICSOLUTIONISEVERYWHEREORTHOGONALTOANGULARDIRECATIONS}
			\gfour_{\alpha \beta}|_{\underline{\Upsilon}_{\iniubar,\inisphere}(\upsigma)} 
			\dot{\underline{\Upsilon}}_{\iniubar,\inisphere}^{\alpha}(\upsigma)
			\frac{\partial}{\partial \inisphere^A}
			\underline{\Upsilon}_{\iniubar,\inisphere}^{\beta}(\upsigma)
			& = 0,
			&  A & = 1,2
		\end{align}
		for all $\upsigma$ in the domain of maximal existence of the solution.
		\eqref{E:GEODESICSOLUTIONISEVERYWHERENULL} and \eqref{E:GEODESICSOLUTIONISEVERYWHEREORTHOGONALTOANGULARDIRECATIONS}
		imply that the vectorfield
		$\dot{\underline{\Upsilon}}_{\iniubar,\inisphere}^{\alpha}(\upsigma)$
		is the null generator of the null hypersurface portions 
		parametrized by varying $\upsigma, \inisphere$ at fixed $\iniubar$,
		i.e., the image of the maps
		$(\upsigma,\inisphere) \rightarrow \underline{\Upsilon}_{\iniubar,\inisphere}^{\alpha}(\upsigma)$.
		In particular, $\dot{\underline{\Upsilon}}_{\iniubar,\inisphere}^{\alpha}(\upsigma)$ 
		is everywhere $\gfour$-orthogonal to these null hypersurface portions.
		Moreover, for solutions emanating from a point $\embedding(0,\theta) \in S_{0,0}$,
		$\dot{\underline{\Upsilon}}_{0,\inisphere}^{\alpha}(\upsigma)$
		must everywhere be parallel to the corresponding null generator $\uiniuLunit^{\alpha}$ of
		$\inihypersurbar$ (see \eqref{E:INITIALULUNITONHBAR}),
		and the solution trajectories 
		$\upsigma \rightarrow \underline{\Upsilon}_{0,\inisphere}^{\alpha}(\upsigma)$ 
		therefore exist long enough such that as $\inisphere$ varies over $\mathbb{S}^2$,
		they foliate $\inihypersurbar \cap \region_{(\text{in})}$.
		Appealing to smooth dependence of ODE solutions on the initial conditions,
		using compactness, 
		perhaps further shrinking the number $\Rendallexistencetime$ in \eqref{E:RENDALLREGION}
		(i.e., further shrinking $\inihypersurbar \cap \region_{(\text{in})}$),
		and using that
		$\uLgeo$ is transversal to $\inihypersur$,
		we further deduce that
		as $\iniubar$ varies over a small interval and $\inisphere$ varies over $\mathbb{S}^2$,
		the solution trajectories 
		$\upsigma \rightarrow \underline{\Upsilon}_{\iniubar,\inisphere}^{\alpha}(\upsigma)$
		emanating from the points in $\inihypersur \cap \region_{(\text{in})}$
		exist long enough to foliate all of $\region_{(\text{in})}$
		such that the map 
		$(\upsigma,\iniubar,\inisphere) \rightarrow \underline{\Upsilon}_{\iniubar,\inisphere}^{\alpha}(\upsigma)$
		is a $C^{\infty}$ diffeomorphism onto $\region_{(\text{in})}$.
		
		We then 
		define $\ubar$ to be the function of $(t,x^1,x^2,x^3)$ that is constant along 
		the curves $\upsigma \rightarrow \underline{\Upsilon}_{\iniubar,\inisphere}^{\alpha}(\upsigma)$
		and such that $\ubar \restriction_{\embedding(\iniubar,\inisphere)} = \iniubar$.
		$\ubar$ is uniquely defined by this construction and satisfies the desired second
		initial condition stated in \eqref{E:INITIALCONDITIONSEIKONALFORENDALLREGION}.
		By construction, the gradient one-form $\partial_{\alpha} \ubar$ is 
		non-vanishing along  $\inihypersur \cap \region_{(\text{in})}$, so by continuity,
		we can assume, after perhaps further shrinking the number $\Rendallexistencetime$ in \eqref{E:RENDALLREGION},
		that $\partial_{\alpha} \ubar$ is non-vanishing in all of $\region_{(\text{in})}$.
		Moreover, $\partial_{\alpha} \ubar$ is
		orthogonal to the level sets of $\ubar$,
		and its $\gfour$-dual must therefore be parallel to the 
		null generators of these surfaces, which we have shown to be
		$\dot{\underline{\Upsilon}}_{\iniubar,\inisphere}^{\alpha}$.
		From this fact and \eqref{E:GEODESICSOLUTIONISEVERYWHERENULL}, 
		we conclude that $\ubar$ satisfies the eikonal equation \eqref{E:EIKONALFORENDALLREGION}.
		
		We analogously construct the solution $u$ 
		to \eqref{E:EIKONALFORENDALLREGION}--\eqref{E:INITIALCONDITIONSEIKONALFORENDALLREGION}
		by solving the null geodesic equation
		emanating from points $\uembedding(\iniu,\inisphere) \in \inihypersurbar$: 
		\begin{subequations}\label{E:geodesiceq2}
			\begin{align}
				\ddot{\Upsilon}_{\iniu,\inisphere}^\alpha(\upsigma)
				& = -\left.\Chfour^{ \ \alpha}_{\kappa \ \lambda}\right|_{{\Upsilon}_{\iniu,\inisphere}(\upsigma)}\dot{\Upsilon}_{\iniu,\inisphere}^\kappa(\upsigma)\dot{\Upsilon}_{\iniu,\inisphere}^\lambda(\upsigma),\\
				\Upsilon_{\iniu,\inisphere}^\alpha(0)=&\uembedding(\iniu,\inisphere),&
				\dot{\Upsilon}_{\iniu,\inisphere}^\alpha(0) 
				& 
				= 
				\left.\Lgeo^\alpha\right|_{\uembedding(\iniu,\inisphere)},
			\end{align} 
		\end{subequations}
		where $\Lgeo^{\alpha}$ is as in \eqref{D:Lgeo}, i.e., 
		it is a $\gfour$-null vectorfield transversal to $\inihypersurbar$ 
		such that along $S_{0,0}$, it is equal
		to the null generator 
		$\iniLunit^{\alpha}$ of $\inihypersur$ (see \eqref{E:INITIALLUNITONH}).
		
		Given numbers $u',\ubar' \geq 0$, we define the sets
		$S_{u',\ubar'} := \lbrace u = u' \rbrace \cap \lbrace \ubar = \ubar' \rbrace$,
		$\mathcal{H}_{u'}([0,\ubar']) := \lbrace u = u' \rbrace \cap \lbrace 0 \leq \ubar \leq \ubar' \rbrace$,
		and
		$\underline{\mathcal{H}}_{\ubar'}([0,u']) := \lbrace \ubar = \ubar' \rbrace \cap \lbrace 0 \leq u \leq u' \rbrace$,
		viewed as subsets of Cartesian coordinate space. 
		Moreover, given a small number $\timecutoffRendallregionwithgeometry > 0$,
		$\Rendallregionwithgeometry$ denotes the set defined in
		\eqref{E:DESCRIPTIONOFLOCALREGIONNEARS00INEIKONALCOORDINATES},
		also viewed as a subset of Cartesian coordinate space.
		Our initial data are such that if $u'$ and $\ubar'$ are sufficiently small,
		the one-forms $\partial_{\alpha} u$ and $\partial_{\alpha} \ubar$ are non-parallel on
		$\mathcal{H}_0([0,\ubar']) \cup \underline{\mathcal{H}}_0([0,u'])$,
		and the sets $S_{0,\ubar'}$ and $S_{u',0}$ are $\gfour$-spacelike and diffeomorphic to $\mathbb{S}^2$.
		Hence, by continuity and the implicit function theorem, 
		if $\timecutoffRendallregionwithgeometry$ is sufficiently small,
		then $\partial_{\alpha} u$ and $\partial_{\alpha} \ubar$ 
		are everywhere non-parallel in $\Rendallregionwithgeometry$,
		and all the sets $S_{u',\ubar'}$ contained in
		$\timecutoffRendallregionwithgeometry$ are
		diffeomorphic to $\mathbb{S}^2$ and $\gfour$-spacelike.
		
		\medskip
		
		\noindent \textbf{Construction of $(\sphere^1,\sphere^2)$}:
		We now construct the geometric angular coordinates by 
		solving the transport equation initial value problem \eqref{E:TRANSPORTINGTHEANGULARCOORDINATEFUNCTIONS}.
		Part of our proof relies on the identity \eqref{EQ:determinant0}, which we independently derive later in the article.
		We start by noting that definition \eqref{D:breveL} implies that
		the Cartesian components of 
		the transport operator $\Lrenormal$
		from equation \eqref{E:TRANSPORTINGTHEANGULARCOORDINATEFUNCTIONS}
		are smooth functions
		of $\inip u, \inip \ubar,$ and the Cartesian components of $\gfour^{-1}$, 
		and that the Cartesian components of $\gfour^{-1}$ are
		smooth functions of $\vvariables$. In particular, by Theorem~\ref{T:Rendall'sTheorem} and the above results on the existence and smoothness of
		$u,\ubar$, it follows that
		the Cartesian components $\Lrenormal^{\alpha}$
		are $C^{\infty}$ functions on
		$\Rendallregionwithgeometry$.
		Moreover, \eqref{E:COORDINATENULLVECTORFIELDSAPPLIEDTOEIKONALFUNCTIONIS2}
		implies that $\Lrenormal$ is transversal to the hypersurface portion
		$\inihypersurbar \cap \Rendallregionwithgeometry$ where the initial data
		from \eqref{E:TRANSPORTINGTHEANGULARCOORDINATEFUNCTIONS}
		are posed, and it also implies that the $\ubar'$-parameterized flow map of $\Lrenormal$
		maps data spheres $S_{u',0}$ onto spheres $S_{u',2 \ubar'}$. Therefore, from the standard existence,
		uniqueness, and smooth dependence on initial conditions results from ODE theory,
		it follows that the solutions
		$(\sphere^1,\sphere^2)$ to \eqref{E:TRANSPORTINGTHEANGULARCOORDINATEFUNCTIONS}
		are $C^{\infty}$ functions of $(t,x^1,x^2,x^3)$ 
		in $\timecutoffRendallregionwithgeometry$.
		In addition, our assumptions on the initial embeddings from Sect.\ref{SS:initialembedding2} imply that
		the change of variables map
		$(t,x^1,x^2,x^3) \rightarrow (u,\ubar,\sphere^1,\sphere^2)$
		(i.e., the inverse of the map $\embedgeo$ from \eqref{E:CHOVFROMGEOTOCARTESIAN})
		is injective on the data surface portions
		$(\inihypersur \bigcup \inihypersurbar) \cap \Rendallregionwithgeometry$,
		and the aforementioned equation \eqref{EQ:determinant0} implies that the determinant of this change of variables map
		is non-vanishing on $(\inihypersur \bigcup \inihypersurbar) \cap \Rendallregionwithgeometry$.
		By compactness, continuity, and the inverse function theorem, it follows that
		if $\timecutoffRendallregionwithgeometry$ is sufficiently small,
		then this change of variables map is a $C^{\infty}$ diffeomorphism on
		$\Rendallregionwithgeometry$.
		
	\end{proof}

	\section{Acoustic Geometry}\label{S:geometry}
	In Prop.\,\ref{P:CONSTRUCTIONOFEIKONALANDGEOMETRICANGULARCOORDINATESINRENDALLREGION}, we extended 
	the acoustical functions $u,\ubar$ to a future-neighborhood 
	$\Rendallregionwithgeometry$ of $S_{0,0}$. Hence, on $\Rendallregionwithgeometry$,
	we are allowed to use the geometric framework from Section \ref{S:ACOUSTICGEOMETRY}.
	In this section, we introduce the acoustic geometry on an arbitrary spacetime region $\region$.
	Later on, we will use the geometry on the specific region $\Rendallregionwithgeometry$
	as well as other regions.
	
	We denote the Levi-Civita connection on $\suubar$ with respect to $\gsphere$ by $\angnabla$.
	We note that the acoustical quantities $\ubar,u$ that we defined in Section \ref{SS:InitialCharData} coincide with the $\iniubar,\iniu$ on the hypersurfaces $\inihypersur$ and $\inihypersurbar$ respectively. As a result, on $\inihypersur\bigcup\inihypersurbar$, we have: 
	\begin{subequations}\label{EQ:LapseConsistent}
		\begin{align}
			\iniu=&u,&\iniubar=&\ubar,\\
			\iniinvfolia=&\invfolia,&\iniuinvfolia=&\uinvfolia,\\
			\inilapse=&\lapse	
		\end{align}
	\end{subequations}
	
		
	\begin{lemma}[Properties related to $\gtau$]
		Let $K$ be the vectorfield defined in \eqref{Def:K}.
		Then the following identities hold:
		\begin{align}
			K\uptau=&0,& \gtauabs{K}^2:=& \gtau(K,K)=1-\left(\frac{2}{\lapse}\cdot\left(\sqrt{\frac{\uinvfolia}{\invfolia}}+\sqrt{\frac{\invfolia}{\uinvfolia}}\right)\right)^{2}.
		\end{align}
	\end{lemma}
	
	\begin{proof}
		The lemma follows from direct computation.
	\end{proof}

	\subsection{Subsets of Spacetime}\label{SS:SUBSETS}
	In this section, we define various subsets of spacetime that will play a role in our a priori estimates; see Fig.\,\ref{Figure4}
	for an illustration.
	
	In the rest of this article, we consider the subset of the spacetime region $I_{\ufinite}\times I_{\ubarfinite}\times\Stwo$, where $I_{\ufinite}:=[0,\ufinite]$ and$ I_{\ubarfinite}:=[0,\ubarfinite]$. That is, for every $u,\ubar$ in the rest of the article, $u\in[0,\ufinite]$ and $\ubar\in[0,\ubarfinite]$, where $\ufinite$ and $\ubarfinite$ are determined in Proposition \ref{P:DETERMININGTRANSVERALDERIVATIVESANDSATISFYINGEQUATIONSONINITIALNULLHYPERSURFACES}.
	
	Let $u_0\geq0,\ubar_0\geq0$. In this article, for $\uptau$ defined in Definition \ref{DE:TIMEFUNCTIONANDRADIALCOORDINATE}, we define the geometric hypersurface $\Stau$ as follows:
	\begin{align}
		\label{D:Sigmatau}\Stau:=\{u+\ubar=\uptau,u\geq u_0,\ubar\geq \ubar_0\}.
	\end{align}
	Moreover, for $u_1\geq u\geq u_0\geq0,\ubar_1\geq\ubar\geq\ubar_0\geq0$, we define the truncated version $\Stau^{u^\prime,\ubar^\prime}$ of $\Stau$ as follows:
	\begin{align}
		\label{DE:Stau}\Stau^{u^\prime,\ubar^\prime}:=\{u+\ubar=\uptau,u_1\geq u\geq u^\prime,\ubar_1\geq\ubar\geq\ubar^\prime\}.
	\end{align}
	
	$\Stau^{u,\ubar}$ is defined in (\ref{DE:Stau}). Let $\updelta_u:=u_1-u_0$, and $\updelta_{\ubar}:=\ubar_1-\ubar_0$. Let $\updelta:=\min\{\updelta_u,\updelta_{\ubar}\}$. For convenience, $\updelta_u,\updelta_{\ubar}\lesssim1$. Let $u+\ubar\leq\uptau\leq \min\{u+\ubar_1,u_1+\ubar\}$.  Then, we define the ``triangle'' space-time region $\region_1$ as follows:
	\begin{align}\label{DE:triangleM}
		\region_1^{u,\ubar,\uptau}:=
		[u+\ubar,\uptau]\times\Stau^{u,\ubar}.
	\end{align}
	
	For $u_0+\ubar_0+\min\{\updelta_u,\updelta_{\ubar}\}\leq\uptau_1\leq\uptau_2\leq\max\{u_0+\ubar_1,u_1+\ubar_0\}$, we define the ``parallelogram'' space-time region as follows:
	\begin{align}\label{DE:parallelogramM}
		\region_2^{\uptau_1,\uptau_2}:=
		[\uptau_1,\uptau_2]\times\Stau^{u_0,\ubar_0}.
	\end{align}
	
	For $\max\{u_0+\ubar_1,u_1+\ubar_0\}\leq\uptau_1\leq\uptau_2\leq u_1+\ubar_1$, we define the ``trapezoid'' space-time region as follows:
	\begin{align}\label{DE:trapezoidM}
		\region_3^{\uptau_1,\uptau_2}:=[\uptau_1,\uptau_2]\times\Stau^{u_0,\ubar_0}.
	\end{align}
	
	See Figure \ref{Figure4} for the regions $\region_1,\region_{2}$ and $\region_{3}$.
	In addition, we define $\hypersur_{u},\hypersurbar_{\ubar}$ to be the null hypersurfaces with constant $u,\ubar$ respectively. Then we define a portion of null hypersurfaces as follows:
	\begin{align}
		\hypersur_{u}^{\ubar_1,\ubar_2}:=&\hypersur_{u}\bigcap\{\ubar_1\leq\ubar\leq\ubar_2\},&
		\hypersur_{u}^{\ubar_2}:=&\hypersur_{u}\bigcap\{0\leq\ubar\leq\ubar_2\},\\
		\hypersurbar_{\ubar}^{u_1,u_2}:=&\hypersurbar_{\ubar}\bigcap\{u_1\leq u\leq u_2\},&
		\hypersurbar_{\ubar}^{u_2}:=&\hypersurbar_{\ubar}\bigcap\{0\leq u\leq u_2\}.
	\end{align} 
	
	For fixed $u,\ubar$, we define the 2-sphere $\suubar$ as follows:
	\begin{align}
		\suubar:=\hypersur_u\bigcap\hypersurbar_{\ubar}.
	\end{align}
	
	Also, for $0\leq u_1<u_2$ and $0\leq \ubar_1<\ubar_2$, we define $\region_{u_1,u_2}^{\ubar_1,\ubar_2}$ to be the ``diamond" region such that
	\begin{align}
		\region_{u_1,u_2}^{\ubar_1,\ubar_2}:=\{u_1\leq u\leq u_2\}\bigcap\{\ubar_1\leq \ubar\leq\ubar_2\}.
	\end{align}
	In particular, we define the ``diamond'' space-time region $\region_0$ as follows:
	\begin{align}\label{DE:region0}
		\region_0:=\{u_0\leq u\leq u_1,\ubar_0\leq\ubar\leq\ubar_1\}.
	\end{align}
	We denote $\region$ the finite union of $\region_{u_1,u_2}^{\ubar_1,\ubar_2}$. We note that such definition is consistent with the qualitative description in Section \ref{SS:Defregion}.
	\begin{figure}
		\centering
		\includegraphics[width=12cm]{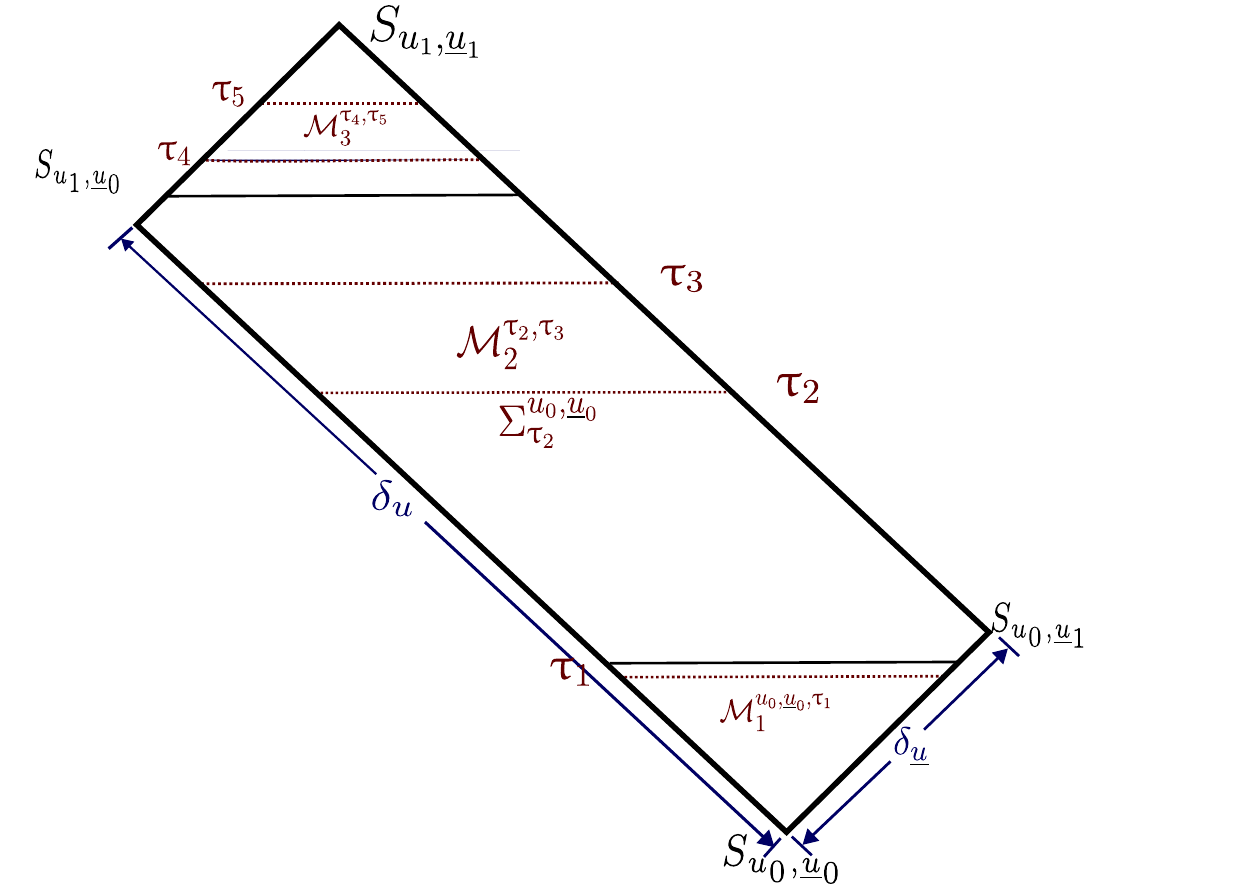}
		\caption{Regions in the case $\updelta_u>\updelta_{\ubar}$}
		\label{Figure4}
	\end{figure}

	\subsection{Metric under Geometric Coordinates}
	In this subsection, we present the components of $\gfour$ relative to the geometric coordinates in Lemma \ref{L:acousticmetric}. In addition, we define the volume forms on $\region,\St,\Stau,\suubar$ in Definition \ref{D:volumeforms}. Meanwhile, we list several useful identities that will be used later in our analysis.
	
	In our analysis, we aim to integrate along null geodesics. However, by (\ref{D:breveL}),(\ref{E:breveL}), Lie derivatives $\Lie_{\Lrenormal}$ and $\Lie_{\uLrenormal}$ differ slightly from $\Dfour_{4},\Dfour_{3}$. In the next lemma, we provide the distinctions between $\Lie_{\Lrenormal}$, $\Lie_{\uLrenormal}$ and $\Dfour_4,\Dfour_{ 3}$.
	\begin{lemma}\cite[(3.1.53), Difference between $\Lie_{\Lrenormal}$ and $\Lie_{\uLrenormal}$ and $\Dfour_{4},\Dfour_{3}$]{TheEvolutionProbleminGeneralRelativity}.\label{L:LieandCovariant}
		Let $\Lie_{V}$ be the standard Lie derivative with respect to vectorfield $V$. For scalar function $f$, $\suubar$-tensorfield $\xi_1\in\{\upeta,\ueta,\upzeta\},\xi_2\in\{\hchi,\huchi\}$, the following identities hold:
		\begin{subequations}
			\begin{align}
				\Lie_{\Lrenormal}f=&\uinvfolia\lapse\Dfour_{ 4}f,&
				\Lie_{\uLrenormal}f=&\invfolia\lapse\Dfour_{ 3}f,\\
				\Lie_{\Lrenormal}\xi_1=&\uinvfolia\lapse\left(\angD_{ 4}\xi_1+\upchi\cdot\xi_1\right),&
				\Lie_{\uLrenormal}\xi_1=&\invfolia\lapse\left(\angD_{ 3}\xi_1+\uchi\cdot\xi_1\right),\\
				\Lie_{\Lrenormal}\xi_2=&\uinvfolia\lapse\left(\angD_{ 4}\xi_2+2\upchi\cdot\xi_2\right),&
				\Lie_{\uLrenormal}\xi_2=&\invfolia\lapse\left(\angD_{ 3}\xi_2+2\uchi\cdot\xi_2\right).
			\end{align}
		\end{subequations}
	\end{lemma}
	\begin{lemma}[Acoustic metric under geometric coordinates]\label{L:acousticmetric}
		Let $\gfour$ be the acoustic metric defined in Definition \ref{D:acousticmetric}, and let $\gsphere$ be the first fundamental form of $\suubar$. We can express the spacetime metric $\gfour$ relative to the geometric coordinates $(u,\ubar,\sphere^1,\sphere^2)$:
		\begin{align}\label{E:gfouringeometriccoordinates}
			\gfour=-2\lapse^2\invfolia\uinvfolia\left(\diff u\otimes\diff\ubar+\diff \ubar\otimes\diff u\right)+\gsphere_{AB}\left(\diff\sphere^A+\Lambda^A\diff\ubar-\Lambda^A\diff u\right)\otimes\left(\diff\sphere^B+\Lambda^B\diff\ubar-\Lambda^B\diff u\right).
		\end{align}
		Moreover, for $\uptau,\rgeo$ defined in Definition \ref{DE:TIMEFUNCTIONANDRADIALCOORDINATE}, we have the following:
		\begin{align}\label{E:gfouringeometriccoordinates2}
			\gfour=-\lapse^2\invfolia\uinvfolia\diff \uptau\otimes\diff\uptau+\lapse^2\invfolia\uinvfolia\diff \rgeo\otimes\diff\rgeo+\gsphere_{AB}\left(\diff\sphere^A+\Lambda^A\diff\rgeo\right)\otimes\left(\diff\sphere^B+\Lambda^B\diff\rgeo\right).
		\end{align}
	\end{lemma}
	\begin{proof}[Discussion of the proof of Lemma \ref{L:acousticmetric}]
		The proof of Lemma \ref{L:acousticmetric} is omitted as it can be derived through straightforward computation using (\ref{D:breveL}), (\ref{E:breveL}), and (\ref{EQ:innerproduct}). Interested readers can refer to \cite[Lemma 3.45]{Smalldatashocks} for insights into the proof of (\ref{E:gfouringeometriccoordinates}).
	\end{proof}
	
	\begin{definition}[Geometric volume forms] \label{D:volumeforms} Let $\gfour$ be the acoustic metric defined in Definition \ref{D:acousticmetric}. Let $\esphere$ be the standard Euclidean metric on $\Stwo$. Let $g$ be the first fundamental form of $\St$. Let $\gtau$ be the first fundamental form of $\Stau$. Let $\gsphere$ be the first fundamental form of $\suubar$. We define the following volume forms:
		\begin{itemize}
			\item $\diff \flatspherevol$ denotes the volume form induced by $\esphere$ on $\Stwo$.
			\item $\diff \gvol$ denotes the volume form induced by $\gfour$ on $\region$.
			\item $\diff \tvol$ denotes the volume form induced by $\gt$ on $\St$.
			\item $\diff \tauvol$ denotes the volume form induced by $\gtau$ on $\Stau$.
			\item $\diff \spherevol$ denotes the volume form induced by $\gsphere$ on $\sTauu$.
		\end{itemize}
		We define the sphere-volume ratio as follows:
		\begin{align}\label{D:spherevolume}
			\volume:=&\frac{\sqrt{\det\gsphere}}{\sqrt{\det\esphere}}.
		\end{align}
		For scalar functions $f$, we define the average value of $f$, denoted by $\average{f}$, as follows:
		\begin{align}
			\average{f}&:=\frac{1}{\sgabs{S_{t,u}}}\int_{S_{t,u}}fd\spherevol,&
			\sgabs{S_{t,u}}&:=\int_{S_{t,u}}1d\spherevol
			.\end{align}
	\end{definition}
	\begin{lemma}\label{L:volumeform}
		Let $(t,x^1,x^2,x^3)$ be the Cartesian coordinates. For volume forms defined in Definition \ref{D:volumeforms}, the following identities hold:
		\begin{subequations}
			\begin{align}
				\label{E:tvol}\diff\tvol=&\speed^{-3}\diff x^1\diff x^2\diff x^3,\\
				\label{E:spherevol}\diff \spherevol=&\sqrt{\det\gsphere}\diff\flatspherevol,\\
				\label{E:gvol}\diff\gvol=&\diff\tvol\diff t=\lapse\sqrt{\invfolia\uinvfolia}\diff\uptau\diff\tauvol=2\lapse^2\invfolia\uinvfolia\diff \ubar\diff u\diff \spherevol.
			\end{align}
		\end{subequations}
	\end{lemma}
	\begin{proof}[Discussion of the proof of Lemma \ref{L:volumeform}]
		(\ref{E:tvol}), (\ref{E:gvol}) are direct results of (\ref{E:ACOUSTICALMETRIC}), (\ref{E:gfouringeometriccoordinates}) respectively. (\ref{E:spherevol}) follows from the definition of $\diff\flatspherevol$ and $\diff\spherevol$.
	\end{proof}
	\subsection{Norms}\label{SS:norms}
	In this article, for functions $f,g$ on a normed space $(X,\norm{\cdot}_X)$, we use the notation $\norm{f,g}_X:=\norm{f}_X+\norm{g}_X$. Similarly, for an array of functions $\vec{U}=(U^1,U^2,\dots,U^k)$, we have $\norm{\vec{U}}_X:=\sum\limits_{a=1}^k\norm{U^a}_X$. In particular, we use $\abs{\vec{U}}:=\sqrt{\sum\limits_{i=1}^{k}(U^i)^2}$. For functions $f$ and arrays $\vec{g}$, we also use $\norm{f,g}_X:=\norm{f}_X+\norm{\vec{g}}_X$.
	
	For $\suubar$-tangent tensorfield $\xi$, we define the norm on $\suubar$ as follows:
	\begin{align}
		\onenorm{\xi}{\gsphere}{2}{\suubar}:=&\left(\int_{\suubar}\sgabs{\xi}^2\diff\spherevol\right)^{1/2},&
		\onenorm{\xi}{\sangle}{2}{\suubar}:=&\left(\int_{\suubar}\sgabs{\xi}^2\diff\flatspherevol\right)^{1/2},
	\end{align}
	where $\sgabs{\xi}$ denotes the norm of $\suubar$-tangent tensorfield with respect to first fundamental form $\gsphere$. For example, for type $\binom{0}{2}$ tensorfield, $\sgabs{\xi}:=\sqrt{\gsphere^{-1}\cdot\gsphere^{-1}\cdot\xi\cdot\xi}$.

	We define the Sobolev norm on $\St$ and characteristic hypersurfaces $\hypersur_{u}$, $\hypersurbar_{\ubar}$ for $k\in\mathbb{N}$: 
	\begin{subequations}\label{D:Sobolevnorms}
		\begin{align}
			\sobolevnorm{F}{k}{\St}:=&\norm{\p^kF}_{L^2(\St)}+\norm{F}_{L^2(\St)},\\
			\sobolevnorm{F}{k}{\hypersur_{u}}:=&\norm{\left(\frac{\p }{\p \ubar},\angp\right)^kF}_{L^2(\hypersur_{u})}+\norm{F}_{L^2(\hypersur_{u})},\\
			\sobolevnorm{F}{k}{\hypersurbar_{\ubar}}:=&\norm{\left(\frac{\p }{\p u},\angp\right)^kF}_{L^2(\hypersurbar_{\ubar})}+\norm{F}_{L^2(\hypersurbar_{\ubar})}
		\end{align}
	\end{subequations}
	
	We denote the standard H\"older semi-norm $\dot{C}_x^{0,\beta}$ and H\"older norm $C_x^{0,\beta}$, where $0<\beta<1$, of a function $F$ with respect to the flat metric on the constant-time hypersurface $\Stau$ and characteristic hypersurfaces $\hypersur_{u}$, $\hypersurbar_{\ubar}$ by 
	\begin{subequations}
		\begin{align}
			\semiholdernorm{F}{x}{0}{\beta}{\Stau}:=&\sup\limits_{x\neq y\in\Stau}\frac{\abs{F(x)-F(y)}}{\abs{x-y}^\beta},&
			\holdernorm{F}{x}{0}{\beta}{\Stau}:=&\sup\limits_{x\in\Stau}\abs{F(x)}+\sup\limits_{x\neq y\in\Stau}\frac{\abs{F(x)-F(y)}}{\abs{x-y}^\beta},\\
			\semiholdernorm{F}{u}{0}{\beta}{\hypersur_{u}}:=&\sup\limits_{x\neq y\in\hypersur_{u}}\frac{\abs{F(x)-F(y)}}{\abs{x-y}^\beta},&
			\holdernorm{F}{u}{0}{\beta}{\hypersur_{u}}:=&\sup\limits_{x\in\hypersur_{u}}\abs{F(x)}+\sup\limits_{x\neq y\in\hypersur_{u}}\frac{\abs{F(x)-F(y)}}{\abs{x-y}^\beta},\\
			\semiholdernorm{F}{\ubar}{0}{\beta}{\hypersurbar_{\ubar}}:=&\sup\limits_{x\neq y\in\hypersurbar_{\ubar}}\frac{\abs{F(x)-F(y)}}{\abs{x-y}^\beta},&
			\holdernorm{F}{\ubar}{0}{\beta}{\hypersurbar_{\ubar}}:=&\sup\limits_{x\in\hypersurbar_{\ubar}}\abs{F(x)}+\sup\limits_{x\neq y\in\hypersurbar_{\ubar}}\frac{\abs{F(x)-F(y)}}{\abs{x-y}^\beta}.
		\end{align}
	\end{subequations}
	
	We also use the following mixed norms for function $F: \mathbb{R}^{1+3}\rightarrow\mathbb{R}$, where $1\leq q_1<\infty$, $1\leq q_2\leq\infty$, $I:=[0,u]$ is an interval of $u$ and $\underline{I}:=[0,\ubar]$ is an interval of $\ubar$:
	\begin{subequations}\label{DE:mixednorms}
		\begin{align}
			\twonorms{F}{u}{q_1}{\sangle}{q_2}{I\times\suubar}&:=\left\{\int_I\onenorm{F}{\sangle}{q_2}{\suprimeubar}^{q_1}\diff u^\prime\right\}^{1/q_1},&\twonorms{F}{u}{\infty}{\sangle}{q_2}{I\times\suubar}&:=\text{ess}\sup\limits_{u\in I}\onenorm{F}{\sangle}{q_2}{\suubar},\\
			\twonorms{F}{\ubar}{q_1}{\sangle}{q_2}{I\times\suubar}&:=\left\{\int_{\underline{I}}\onenorm{F}{\sangle}{q_2}{\suubarprime}^{q_1}\diff \ubar^\prime\right\}^{1/q_1},&\twonorms{F}{\ubar}{\infty}{\sangle}{q_2}{\underline{I}\times\suubar}&:=\text{ess}\sup\limits_{\ubar\in \underline{I}}\onenorm{F}{\sangle}{q_2}{\suubar},\\
			\holdertwonorms{F}{u}{q_1}{u}{0}{\beta}{I\times\hypersur_{u}}&:=\left\{\int_{I}\holdernorm{F}{u}{0}{q_2}{\hypersur_{u}}^{q_1}\diff\tau\right\}^{1/q_1},&\holdertwonorms{F}{u}{\infty}{u}{0}{q_2}{I\times\hypersur_{u}}&:=\text{ess}\sup\limits_{u\in I}\holdernorm{F}{u}{0}{q_2}{\hypersur_{u}},\\
			\holdertwonorms{F}{\ubar}{q_1}{\ubar}{0}{\beta}{\underline{I}\times\hypersurbar_{\ubar}}&:=\left\{\int_{I}\holdernorm{F}{u}{0}{q_2}{\hypersurbar_{\ubar}}^{q_1}\diff\tau\right\}^{1/q_1},&\holdertwonorms{F}{\ubar}{\infty}{\ubar}{0}{q_2}{\underline{I}\times\hypersurbar_{\ubar}}&:=\text{ess}\sup\limits_{\ubar\in \underline{I}}\holdernorm{F}{\ubar}{0}{q_2}{\hypersurbar_{\ubar}}.
		\end{align}
	\end{subequations}
	\subsection{Change of Variables Map, Elliptic Estimates on the Sphere, and Tensorial Identities}
	In this section, we define the change of variables map $\embedgeo$ in Definition \ref{DE:changofvariables}. We provide tensorial identities in Lemma \ref{TensorialId} and Corollary \ref{C:Idantisymmetric}.
	\begin{definition}[Change of variables map]\label{DE:changofvariables}
		We define $\embedgeo:[u_0,u_1]\times[\ubar_0,\ubar_1]\times \Stwo\rightarrow \region_0$ (defined in \eqref{DE:region0}) to be the change of variables map from geometric coordinates $(u,\ubar,\theta^1,\theta^2)$ to rectangular coordinates $(t,x^1,x^2,x^3)$.
	\end{definition}
	\begin{lemma}[The Jacobian of the change of variables map $\embedgeo$]\label{LE:Jacobian}
		$\embedgeo$ is the change of variables map defined in Definition \ref{DE:changofvariables}. Then the Jacobian determinant of $\embedgeo$ is as follows:
		\begin{align}\label{EQ:determinant0}
			\det\frac{\p(t,x^1,x^2,x^3)}{\p(u,\ubar,\theta^1,\theta^2)}=4\invfolia\uinvfolia\lapse^3(\det\gt)^{-\frac{1}{2}}(\det\gsphere)^{\frac{1}{2}}.
		\end{align}
	\end{lemma}
	\begin{proof}[Proof of Lemma \ref{LE:Jacobian}]
		By the construction of the geometry in Section \ref{SS:SUBSETSOFSPACETIME}, we have:
		\begin{subequations}\label{EQ:dPhi}
			\begin{align}
				\frac{\p t}{\p u}=&\uLrenormal^0=2\invfolia\lapse^2,&\frac{\p t}{\p \ubar}=&\Lrenormal^0=2\uinvfolia\lapse^2,\\
				\frac{\p x^i}{\p u}=&\uLrenormal^i-\Lambda^i,&\frac{\p x^i}{\p \ubar}=&\Lrenormal^i+\Lambda^i,\\
				\frac{\p t}{\p \theta^A}=&0,& \frac{\p x^i}{\p \theta^A}=&e_A^i.
			\end{align}
		\end{subequations}
		Then, we compute the determinant of the Jacobian matrix:
		\begin{align}\label{EQ:determinant1}
			\det\frac{\p(t,x^1,x^2,x^3)}{\p(u,\ubar,\theta^1,\theta^2)}=\det
			\begin{pmatrix}
				2\invfolia\lapse^2 & 2\uinvfolia\lapse^2 & 0 & 0 \\
				\uLrenormal^1-\Lambda^1 & \Lrenormal^1+\Lambda^1 & e_1^1& e_2^1 \\
				\uLrenormal^2-\Lambda^2  & \Lrenormal^2+\Lambda^2  & e_1^2& e_2^2  \\
				\uLrenormal^3-\Lambda^3 & \Lrenormal^3+\Lambda^3 & e_1^3& e_2^3 
			\end{pmatrix}
		\end{align}
		Noticing that $\Lambda\in \text{span}\{e_1,e_2\}$ and $\Lrenormal=\uinvfolia\lapse e_4,\uLrenormal=\invfolia\lapse e_3$, we compute \eqref{EQ:determinant1} as follows:
		\begin{align}\label{EQ:determinant2}
			\det\frac{\p(t,x^1,x^2,x^3)}{\p(u,\ubar,\theta^1,\theta^2)}=2\invfolia\uinvfolia\lapse^3\cdot\det
			\begin{pmatrix}
				e_4^1 & e_1^1& e_2^1 \\
				e_4^2 & e_1^2& e_2^2  \\
				e_4^3 & e_1^3& e_2^3 
			\end{pmatrix}
			-2\invfolia\uinvfolia\lapse^3\cdot\det
			\begin{pmatrix}
				e_3^1 & e_1^1& e_2^1 \\
				e_3^2 & e_1^2& e_2^2  \\
				e_3^3 & e_1^3& e_2^3 
			\end{pmatrix}
			=4\invfolia\uinvfolia\lapse^3\cdot\det
			\begin{pmatrix}
				\spherenormal^1 & e_1^1& e_2^1 \\
				\spherenormal^2 & e_1^2& e_2^2  \\
				\spherenormal^3 & e_1^3& e_2^3 
			\end{pmatrix}
		\end{align}
		The determinant of matrix on the RHS of (\ref{EQ:determinant2}) can be viewed as $\flattvol(\spherenormal,e_1,e_2)$, where $\flattvol$ is the volume form of the standard Euclidean metric $e$ on $\St$. By Lemma \ref{L:volumeform} and the fact that $\spherenormal$ is the $\gt$-unit normal to $\suubar$, we have:
		\begin{align}\label{EQ:determinant3}
			\flattvol(\spherenormal,e_1,e_2)=(\det\gt)^{-\frac{1}{2}}\tvol(\spherenormal,e_1,e_2)=(\det\gt)^{-\frac{1}{2}}\spherevol(e_1,e_2)=(\det\gt)^{-\frac{1}{2}}(\det\gsphere)^{\frac{1}{2}}.
		\end{align}
		Combining (\ref{EQ:determinant2}) and (\ref{EQ:determinant3}), we obtain the desired result (\ref{EQ:determinant0}).
	\end{proof}
	\begin{remark}
		$\embedgeo$ is an important object in deriving local existence and uniqueness for compressible Euler equations in Cartesian spacetime, as we obtain local existence and uniqueness only in geometric coordinates $(u,\ubar,\theta^1,\theta^2)$ via a priori estimates. In Section \ref{S:Lastslice}, we will show that under initial data assumptions in Section \ref{S:SectionAssumptions}, $\embedgeo$ is a $C^1$-diffeomorphism that maps onto a neighborhood of $\inihypersur\bigcup\inihypersurbar$.
	\end{remark}
	
	In below, we provide the transport lemma for integrals on the sphere $\suubar$.
	\begin{lemma}\cite[Lemma 3.1.3]{TheEvolutionProbleminGeneralRelativity}\label{LE:e4SuubarInt}.
		For $\suubar$ and scalar function $f$, we have the following equations:
		\begin{subequations}
			\begin{align}
				\label{EQ:LieL}\Lie_{\Lrenormal}\left(\int_{\suubar}f\diff\spherevol\right)=&\int_{\suubar}\left(\Lie_{\Lrenormal}f+\uinvfolia\lapse\gtr\upchi f\right)\diff\spherevol,\\
				\label{EQ:LieuL}\Lie_{\uLrenormal}\left(\int_{\suubar}f\diff\spherevol\right)=&\int_{\suubar}\left(\Lie_{\uLrenormal}f+\invfolia\lapse\gtr\uchi f\right)\diff\spherevol,
			\end{align}
		\end{subequations}
		where $\Lie_{V}$ is the standard Lie derivative with respect to vectorfield $V$.
	\end{lemma}
	\begin{lemma}\cite[Proposition 10, Sobolev embedding on $\suubar$]{GRCharIVP}.\label{LE:SobolevEmbedding}
		We assume that the metric $\gsphere$ on $\suubar$ satisfies the following:
		\begin{align}
			\tilde{a}^{-1}\seabs{\xi}^2\leq\gsphere_{AB}\xi^A\xi^B\leq& \tilde{a}\seabs{\xi}^2,&
			\text{ for any }& \suubar \text{-tangent vectorfield }\xi,
		\end{align}
		where $\esphere$ is the Euclidean metric on $\Stwo$, and $1\leq\tilde{a}\leq c<\infty$ for some constant $c$.
		Then, for scalar function $f$, the following Sobolev embedding holds on $\suubar$:
		\begin{align}
			\onenorm{f}{}{\infty}{\suubar}\leq C(\tilde{a})\sum\limits_{i=0}^2\onenorm{\angnabla^if}{}{2}{\suubar},
		\end{align}
		where $C(\tilde{a})$ is a smooth function of $\tilde{a}$.
	\end{lemma}

	We now provide several useful tensorial identities. We will use them when deriving spacetime elliptic estimates for the specific vorticity $\vortrenormalized$ and entropy gradient $\gradEnt$ in Section \ref{SS:elliptic}.
	\begin{lemma}\label{TensorialId}\cite[Lemma 5.6, Identities for $\St$-tangent vectorfields]{IntegralIdentities}. Let $V$ be a $\St$-tangent vectorfield. Then the following identities hold:
		\begin{subequations}
			\begin{align}
				V_0=&-g_{ab}v^aV^b=-v^aV_a,&
				V_i=&\speed^{-2}V^i,\\
				\materialderivative^\beta\p_\alpha V_\beta=&-V_{\beta}(\p_\alpha\materialderivative^\beta)=-V_a\p_\alpha v^a,\\
				\p_\alpha V_\beta-\p_\beta V_\alpha=&2(\p_\beta\ln\speed)V_{\alpha}-2(\p_\alpha\ln\speed)V_{\beta}+\updelta^0_\alpha V_a\p_\beta v^a-\delta^0_\beta V_a\p_\alpha v^a\\
				\notag&+\left(\delta^0_\alpha\gfour_{\beta\gamma}-\delta^0_\beta\gfour_{\alpha\gamma}\right)\materialderivative V^\gamma+\speed^{-2}\antisymmetic_{\alpha\beta\gamma\delta}\materialderivative^\gamma(\curl V)^\delta.
			\end{align}
		\end{subequations}
	\end{lemma}
	\begin{corollary}[Identities for anti-symmetric terms]\label{C:Idantisymmetric} Let $V^i:=\speed^2 \delta^{ia}\p_a \variables$ and $V^0:=0$, then $V_{i}=\p_i \variables$ and:
		\begin{align}\label{E:TensorialIdvariable}
			(\diff\p_i\variables)_{\alpha\beta}:=&\p_\alpha V_{\beta}-\p_\beta V_{\alpha}\\
			\notag=&\delta_\alpha^0\delta_\beta^i(\materialderivative\p_i\variables+\p_iv^a\p_a\variables)-\delta_\alpha^i\delta_\beta^0(\materialderivative\p_i\variables+\p_iv^a\p_a\variables)
		\end{align}
		Moreover, for $\vortrenormalized$ and $\gradEnt$, we have the following (see \cite[Corollary 5.7]{IntegralIdentities}):
		\begin{subequations}\label{E:vortofvortandgradEnt}
			\begin{align}
				\p_\alpha\vortrenormalized_\beta-\p_\beta\vortrenormalized_\alpha=&2(\p_\beta\ln\speed)\vortrenormalized_\alpha-2(\p_\alpha\ln\speed)\vortrenormalized_\beta+2\delta^0_\alpha\vortrenormalized_\alpha\p_\beta v^a-2\delta^0_\beta\vortrenormalized_a\p_\alpha v^a\\
				\notag&-\speed^{-4}e^{-2\logDensity}\frac{\pressure_{;\Ent}}{\backgroundDensity}\antisymmetic_{\alpha\beta\gamma\delta}(\materialderivative v^\gamma)\gradEnt^\delta+\speed^{-4}e^{-2\logDensity}\frac{\pressure_{;\Ent}}{\backgroundDensity}\antisymmetic_{\alpha\beta\gamma\delta}\materialderivative^\gamma\left[\gradEnt^\delta(\p_av^a)-\gradEnt^a\p_av^i\right]\\
				\notag&+\speed^{-2}e^{\logDensity}\antisymmetic_{\alpha\beta\gamma\delta}\materialderivative^\gamma\modivort^\delta,\\
				\p_\alpha\gradEnt_\beta-\p_\beta\gradEnt_\alpha=&2(\p_\beta\ln\speed)\gradEnt_\alpha-2(\p_\alpha\ln\speed)\gradEnt_\beta.
			\end{align}
		\end{subequations}
	\end{corollary}
	
	\section{Sobolev Regularity of the Initial Data}\label{S:REGULARITYOFDATA}
		In this section, we derive estimates for various norms of the fluid and acoustic geometry along the initial null hypersurfaces.
		We have already shown that the solution jet is $C^{\infty}$ along these initial hypersurfaces (see Sect.\,\ref{S:FLUIDCHARDATAFURTHERASPECTS}),
		so our analysis here can be viewed as consisting of a priori estimates for the data.

	\subsection{Regularities of the Lapses}\label{SS:lapseregularity}
	The following proposition reveals the reason why we (are allowed to) assume that the initial embedding $\embedding$ has the same regularity as the initial data $(\mathring{v}^i,\inilogDensity,\inient)$ (see (\ref{A:regularity})). Note that Proposition \ref{P:equationforlapses} holds true in the region $\region$. Moreover, on the initial hypersurfaces $\inihypersur,\inihypersurbar$, recalling (\ref{D:invfoliauinvfolia}), we have:
	\begin{subequations}
		\begin{align}
			\iniinvfolia\inilapse^2=&1,& \text{on }\inihypersurbar,\\
			\iniuinvfolia\inilapse^2=&1,& \text{on }\inihypersur.
		\end{align}
	\end{subequations}
	\begin{proposition}[Geodesic properties of the characteristic embedding]\label{P:equationforlapses}
		In this proposition, we consider the spacetime region $\region$ (defined in Section \ref{SS:SUBSETSOFSPACETIME}). For the first derivatives of $\invfolia$ and $\uinvfolia$ (defined in Section \ref{SS:SUBSETSOFSPACETIME}),
		
		\begin{subequations}
			\begin{align}
				\label{E:e_4e_43}e_4\ln\invfolia=&-\frac{1}{4}\dfour{4}{4}e_4\variables-\frac{1}{4}\left(2\dfour{4}{3}e_4\variables-\dfour{4}{4}e_3\variables\right)+\frac{1}{4\lapse}\left(2\angdfour^0_4e_4\variables-\dfour{4}{4}\angp^0\variables\right),\\
				\label{E:e_3e_3}e_3\ln\uinvfolia=&-\frac{1}{4}\dfour{3}{3}e_3\variables-\frac{1}{4}\left(2\dfour{4}{3}e_3\variables-\dfour{3}{3}e_4\variables\right)+\frac{1}{4\lapse}\left(2\angdfour^0_3e_3\variables-\dfour{3}{3}\angp^0\variables\right),
			\end{align}	
		\end{subequations}
		where $\angdfour$ and $\angp$ are defined in Definition \ref{DE:Projections}.
		
		Moreover, 
		for the first derivatives of $\Lunit^i$ and $\uLunit^i$ (defined in Section \ref{SS:SUBSETSOFSPACETIME}),
		\begin{subequations}
			\begin{align}
				\label{E:e_4e_432}
				e_4\Lunit^i=&\left(\frac{1}{4}\dfour{4}{4}e_4\variables\right)\uLunit^i-\left\{\frac{1}{4}\dfour{4}{4}e_4\variables-\frac{1}{4}\lapse^{-1}\left(2\angdfour^0_4e_4\variables-\dfour{4}{4}\angp^0\variables\right)\right\}\Lunit^i-\frac{1}{2}\lapse^{-1}\left(\angdfour^i_{4}e_4\variables-\frac{1}{2}\dfour{4}{4}\angp^i\variables\right),\\
				\label{E:e_3e_32}
				e_3\uLunit^i=&\left(\frac{1}{4}\dfour{3}{3}e_3\variables\right){\Lunit}^i-\left\{\frac{1}{4}\dfour{3}{3}e_3\variables-\frac{1}{4}\lapse^{-1}\left(2\angdfour^0_3e_3\variables-\dfour{3}{3}\angp^0\variables\right)\right\}\uLunit^i-\frac{1}{2}\lapse^{-1}\left(\angdfour^i_{3}e_3\variables-\frac{1}{2}\dfour{3}{3}\angp^i\variables\right).
			\end{align}
		\end{subequations}
		
		In addition, for the first derivatives of $\lapse$, we have:
		\begin{subequations}
			\begin{align}
				\label{E:e_3lapse}e_3\ln\lapse=&-\frac{1}{4}\dfour{3}{4}e_3\variables-\frac{1}{8}\dfour{4}{4}e_3\variables-\frac{1}{4}\lapse^{-1}\upeta^0+\frac{1}{8}\lapse^{-1}\left(\angdfour_{3}^0e_4\variables+\angdfour_{4}^{0}e_3\variables-\dfour{3}{4}\angp^0\variables\right)\\
				\notag&-\frac{1}{8}\dfour{3}{3}e_3\variables+\frac{1}{8}\lapse^{-1}\left(2\angdfour^0_3e_3\variables-\dfour{3}{3}\angp^0\variables\right),\\
				\label{E:e_4lapse}e_4\ln\lapse=&-\frac{1}{4}\dfour{3}{4}e_4\variables-\frac{1}{8}\dfour{3}{3}e_4\variables-\frac{1}{4}\lapse^{-1}\ueta^0+\frac{1}{8}\lapse^{-1}\left(\angdfour_{4}^0e_3\variables+\angdfour_{3}^{0}e_4\variables-\dfour{3}{4}\angp^0\variables\right)\\
				\notag&-\frac{1}{8}\dfour{4}{4}e_4\variables+\frac{1}{8}\lapse^{-1}\left(2\angdfour^0_4e_4\variables-\dfour{4}{4}\angp^0\variables\right).
			\end{align}
		\end{subequations}
	\end{proposition}
	\begin{proof}[Proof of Proposition \ref{P:equationforlapses}]
		First of all, recall that:
		\begin{align}
			\Dfour_{\Lgeo}\Lgeo=0,
		\end{align}
		where $\Lgeo=\invfolia^{-1}\Lunit$.
		
		We consider the equation: 
		\begin{align}\label{E:CLunit1}
			\Lunit\ln\invfolia=-\invfolia^{2}\Lgeo\Lgeo^0=(\gfour^{-1})^{0\gamma}\Chfour_{\gamma\alpha\beta}\Lunit^\alpha\Lunit^\beta.
		\end{align}
		We decompose $\gfour^{-1}$ as follows:
		\begin{align}\label{E:CLunit2}
			(\gfour^{-1})^{0\gamma}=-\lapse\left(e_4^\gamma+e_3^\gamma\right)+(\gsphere^{-1})^{0\gamma}=-2\lapse^{2}\left(\Lunit^\gamma+\uLunit^\gamma\right)+(\gsphere^{-1})^{0\gamma}.
		\end{align}
		By (\ref{E:CLunit1}) and (\ref{E:CLunit2}), we obtained (\ref{E:e_4e_43}). \eqref{E:e_3e_3} can be proved in a similar fashion.
		
		Now, to prove \eqref{E:e_4e_432}, we consider the equation:
		\begin{align}\label{E:e_4e_40}
			\Lunit\Lunit^\nu=\Dfour_\Lunit\Lunit^\nu-\frac{1}{2}(\gfour^{-1})^{\nu\kappa}\Lunit^\alpha \Lunit^\beta\left(\dfour{\kappa}{\beta}\p_\alpha\variables+\dfour{\alpha}{\kappa}\p_\beta\variables-\dfour{\alpha}{\beta}\p_\kappa\variables\right).
		\end{align}
		Contracting (\ref{E:e_4e_40}) with $e_\alpha$, we compute:
		\begin{align}\label{E:e_4e_420}
			\Lunit\Lunit^\nu=&\lapse^{2}\left(\dfour{\Lunit}{\Lunit}\Lunit\variables-2\gfour(\Dfour_\Lunit\Lunit,\Lunit)\right) \uLunit^\nu+\lapse^{2}\left\{-2\gfour(\Dfour_\Lunit\Lunit,\uLunit)+(2\dfour{\uLunit}{\Lunit}\Lunit\variables-\dfour{\Lunit}{\Lunit}\uLunit\variables)\right\}\Lunit^\nu\\
			\notag&-(\gsphere^{-1})^{AB}\left(\dfour{A}{\Lunit}\Lunit\variables-\frac{1}{2}\dfour{\Lunit}{\Lunit}e_A\variables+\gfour(\Dfour_\Lunit\Lunit,e_A)\right)e_B^\nu,
		\end{align}	
		By (\ref{D:geo}), we have:
		\begin{align}\label{E:Dfour_4e_4,e_40}
			\gfour(\Dfour_\Lunit\Lunit,\Lunit)=\gfour(\Dfour_\Lunit\Lunit,e_A)=0.
		\end{align}
		By the fact that $\Lgeo=\invfolia^{-1}\Lunit$, we also have:
		\begin{align}\label{E:DfourLL}
			\Dfour_\Lunit\Lunit=\Lunit\ln\invfolia\cdot \Lunit+\invfolia\Dfour_{\Lgeo}\Lgeo.
		\end{align}
		Contracting (\ref{E:DfourLL}) with $\uLunit$, we derive:
		\begin{align}
			\label{E:omegainG1}\gfour(\Dfour_\Lunit\Lunit,\uLunit)=&-\frac{1}{2}\lapse^{-2}\Lunit\ln\invfolia.
		\end{align}
		Then by (\ref{E:omegainG1}), (\ref{E:Dfour_4e_4,e_40}), (\ref{E:e_4e_420}) and (\ref{E:omegainG}), we have obtained (\ref{E:e_4e_432}). \eqref{E:e_3e_32} can be proved in a similar fashion.
		
		To prove (\ref{E:e_3lapse}), we consider
		\begin{align}\label{E:lapse1.1}
			e_3\left\{(2\lapse^2)^{-1}\right\}=e_3\left\{\gfour(\Lunit,\uLunit)\right\}=\dfour{\Lunit}{\uLunit}e_3\variables+\gfour(e_3\Lunit,\uLunit)+\gfour(\Lunit,e_3\uLunit).
		\end{align}
		To compute $\gfour(e_3\Lunit,\uLunit)$, we consider
		\begin{align}\label{E:lapse1.2}
			0=e_3\Lunit^0=-2\lapse^2\gfour(e_3\Lunit,\uLunit)\Lunit^0-2\lapse^2\gfour(e_3\Lunit,\Lunit)\uLunit^0+(\gsphere^{-1})^{AB}\gfour(e_3\Lunit,e_A)e_B^0.
		\end{align}
		For the second term on the RHS of \eqref{E:lapse1.2}, notice that
		\begin{align}\label{E:lapse1.25}
			0=e_3\left\{\gfour(\Lunit,\Lunit)\right\}=\dfour{\Lunit}{\uLunit}e_3\variables+2\gfour(e_3\Lunit,\Lunit).
		\end{align}
		For the last term on the RHS of \eqref{E:lapse1.2}, by \eqref{DE:ConnectionCoefficients}, we have:
		\begin{align}\label{E:lapse1.26}
			(\gsphere^{-1})^{AB}\gfour(e_3\Lunit,e_A)e_B^0=&\frac{1}{2\lapse}(\gsphere^{-1})^{AB}e_B^0\left[\upeta_A-\frac{1}{2}\left(\dfour{3}{A}e_4\variables+\dfour{4}{A}e_3\variables-\dfour{3}{4}e_A\variables\right)\right]\\
			\notag=&\frac{1}{2}\lapse^{-1}\upeta^0-\frac{1}{4}\lapse^{-1}\left(\angdfour_{3}^0e_4\variables+\angdfour_{4}^{0}e_3\variables-\dfour{3}{4}\angp^0\variables\right).
		\end{align}
		Combining \eqref{E:lapse1.2}, (\ref{E:lapse1.25}) and (\ref{E:lapse1.26}), we obtain:
		\begin{subequations}\label{E:lapse1.3}
			\begin{align}
				\gfour(e_3\Lunit,\Lunit)=&-\frac{1}{8}\lapse^{-2}\dfour{4}{4}e_3\variables,\\
				\gfour(e_3\Lunit,\uLunit)=&\frac{1}{8}\lapse^{-2}\dfour{4}{4}e_3\variables+\frac{1}{4}\lapse^{-3}\upeta^0-\frac{1}{8}\lapse^{-3}\left(\angdfour_{3}^0e_4\variables+\angdfour_{4}^{0}e_3\variables-\dfour{3}{4}\angp^0\variables\right),
			\end{align}
		\end{subequations}
		where $\upeta^0=(\gsphere^{-1})^{0\alpha}\upeta_\alpha$.
		
		For the last term on the RHS of \eqref{E:lapse1.1}, by (\ref{E:e_3e_32}), we have:
		\begin{align}\label{E:lapse1.4}
			\gfour(\Lunit,e_3\uLunit)=\frac{1}{8}\lapse^{-2}\dfour{3}{3}e_3\variables-\frac{1}{8}\lapse^{-3}\left(2\angdfour^0_3e_3\variables-\dfour{3}{3}\angp^0\variables\right).
		\end{align}
		Combining (\ref{E:lapse1.1}), (\ref{E:lapse1.3}) and (\ref{E:lapse1.4}), we get the desired result (\ref{E:e_3lapse}). (\ref{E:e_4lapse}) can be proved in a similar fashion.
	\end{proof}
	
	\begin{corollary}\label{CO:omegauomega}
		Using Proposition \ref{P:equationforlapses} and (\ref{omegauomega}), by direct computation, we have the following expression for $\omega$ and $\uomega$:
		\begin{subequations}\label{E:omegapde1}
			\begin{align}
				\label{omega}\omega=&-\frac{3}{8}\dfour{4}{4}e_4\variables+\frac{1}{4}\dfour{4}{4}e_3\variables\\
				\notag&-\frac{3}{4}\dfour{3}{4}e_4\variables-\frac{1}{4}\lapse^{-1}\ueta^0+\frac{1}{8}\lapse^{-1}\left(\angdfour_{4}^0e_3\variables+\angdfour_{3}^{0}e_4\variables-\dfour{3}{4}\angp^0\variables\right)\\
				\notag&+\frac{3}{8}\lapse^{-1}\left(2\angdfour^0_4e_4\variables-\dfour{4}{4}\angp^0\variables\right),\\
				\label{uomega}\uomega=&-\frac{3}{8}\dfour{3}{3}e_3\variables+\frac{1}{4}\dfour{3}{3}e_4\variables\\
				\notag&-\frac{3}{4}\dfour{3}{4}e_3\variables-\frac{1}{4}\lapse^{-1}\upeta^0+\frac{1}{8}\lapse^{-1}\left(\angdfour_{3}^0e_4\variables+\angdfour_{4}^{0}e_3\variables-\dfour{3}{4}\angp^0\variables\right)\\
				\notag&+\frac{3}{8}\lapse^{-1}\left(2\angdfour^0_3e_3\variables-\dfour{3}{3}\angp^0\variables\right).
			\end{align}
		\end{subequations}
	\end{corollary}
	
	\subsection{Regularities of the Initial Data for Fluid Variables}\label{SS:regularityoffluid}
	In Section \ref{S:CID}-Section \ref{S:FLUIDCHARDATAFURTHERASPECTS}, we  provided a description of the initial data for the fluid variables $(v^i,\logDensity,\Ent)$ and their first derivatives, as well as $\inimodivort,\iniDivGradEnt$ on the characteristic initial hypersurfaces. In this section, we develop a regularity theory for the characteristic initial data in terms of Sobolev norms. Specifically, we derive estimates for fluid variables and null frames.
	
	Let $\ufinite,\ubarfinite$ be as in Proposition \ref{P:DETERMININGTRANSVERALDERIVATIVESANDSATISFYINGEQUATIONSONINITIALNULLHYPERSURFACES}. Let the embedding $\embed$ be defined as in Definition \ref{D:LDIFFEOCONSTRUCTEDOUTOFEMBEDDINGS}. Let $I_{\ufinite}:=[0,\ufinite]$ and $ I_{\ubarfinite}:=[0,\ubarfinite]$. We focus solely on the initial hypersurface $\inihypersur_{[0,\ubarfinite]}=\embed\left(\{0\}\times I_{\ubarfinite}\times\Stwo\right)$; the regularity arguments on $\inihypersurbar_{[0,\ufinite]}=\embed\left(I_{\ufinite}\times\{0\}\times \Stwo\right)$ are similar.
	\begin{lemma}[Standard estimates for transport equations on initial hypersurfaces]\label{LE:transport}
		Assuming\footnote{This assumption is a direct result of (\ref{A:regularity}), (\ref{A:regularity2}) and (\ref{E:chiinG})} $\onenorm{\gtr\upchi}{ }{\infty}{\inihypersur}\lesssim\iniBAinfty$, we have the following standard $L^2$ estimate for transport equations:
		\begin{align}\label{ES:L2transport}
			\norm{\varphi}_{L^2(S_{0,\ubar})}\lesssim&\norm{\varphi}_{L^2(S_{0,0})}+\int_{0}^{\ubar}\norm{e_4\varphi}_{L^2(S_{0,\ubarprime})}\diff\ubarprime,& \text{ on }& \inihypersur.
		\end{align}
		
		Moreover, we have the following $L^\infty$ estimate:
		\begin{align}\label{ES:Linftytransport}
			\norm{\varphi}_{L^\infty(S_{0,\ubar})}\lesssim&\norm{\varphi}_{L^\infty(S_{0,0})}+\int_{0}^{\ubar}\norm{e_4\varphi}_{L^\infty(S_{0,\ubarprime})}\diff\ubarprime,& \text{ on }& \inihypersur.
		\end{align}
	\end{lemma}
	\begin{proof}[Proof of Lemma \ref{LE:transport}]
		(\ref{ES:L2transport})-(\ref{ES:Linftytransport}) are direct results of Lemma \ref{L:LieandCovariant}, Lemma \ref{LE:e4SuubarInt} and the fact that $\iniuinvfolia\inilapse=1$ on $\inihypersur$. We refer readers to \cite[Proposition 6-7]{GRCharIVP} for the proofs.
	\end{proof}
	\subsubsection{Regularities of the mixed derivatives of the initial data}\label{SSS:inisecondderivative}
	
	Now we focus on mixed derivatives of fluid variables. In this subsection, we use a schematic form to illustrate ideas instead of presenting detailed expressions, which does not affect our analysis. 
	
	\begin{proposition}[$L^\infty_{\ubar}L^\infty_{\gsphere}$-bound for the fluid variables of the characteristic initial data]\label{PR:InftyBoundfortheData}
		Let $\embed$ be the characteristic embedding defined in Definition \ref{D:LDIFFEOCONSTRUCTEDOUTOFEMBEDDINGS}. Let $(\inivel^i,\inilogDensity,\inient)$ be the smooth characteristic initial data that we prescribed/constructed in Section \ref{S:CID}-\ref{S:localexistence} so that their derivatives are well-defined at all orders on $\inihypersur_{[0,\ubarfinite]}\bigcup\inihypersurbar_{[0,\ufinite]}$, where $\ufinite,\ubarfinite$ are determined in Proposition \ref{P:DETERMININGTRANSVERALDERIVATIVESANDSATISFYINGEQUATIONSONINITIALNULLHYPERSURFACES}. We impose the regularity assumptions on the initial data in Sobolev norms \eqref{A:regularity}-\eqref{A:regularity2}. Then, the following estimates are satisfied by the mixed derivatives of the fluid variables on $\inihypersur_{[0,\ubarfinite]}$ with non-negative indices $j,k$:
		\begin{subequations}\label{ES:DataInfinity}
			\begin{align}
				\label{ES:DataInfinity1}&\twonorms{\mathring{e_3}^{j}\left(\mathring{e_4},\iniangp\right)^k(\inivel^i,\inilogDensity)}{\ubar}{\infty}{\gsphere}{\infty}{\inihypersur_{[0,\ubarfinite]}}\leq\tilde{\gensmoothfunction}_{j,k}(\iniLtwoS),&
				\text{ for }& 0\leq 2j+k\leq  \frac{N}{2}+1,\\
				\label{ES:DataInfinity2}&\twonorms{\mathring{e_3}^{j}\left(\mathring{e_4},\iniangp\right)^k(\mathring{e_3}^2\inient,\iniDivGradEnt)}{\ubar}{\infty}{\gsphere}{\infty}{\inihypersur_{[0,\ubarfinite]}}\leq\tilde{\gensmoothfunction}_{j,k}(\iniLtwoS),&
				\text{ for }& 0\leq 2j+k\leq  \frac{N}{2}-1.
			\end{align}
		\end{subequations}
		where $\tilde{\gensmoothfunction}_{j,k}(\iniLtwoS)<\infty$ are smooth increasing functions of $\iniLtwoS$ such that $\tilde{\gensmoothfunction}_{j,k}(0)=0$.
	\end{proposition}
	\begin{proof}[Proof of Proposition \ref{PR:InftyBoundfortheData}]
		We prove \eqref{ES:DataInfinity1}-\eqref{ES:DataInfinity2} by induction with respect to $j$. \eqref{ES:DataInfinity1} with $j=0$, \eqref{ES:DataInfinity1} with $j=1,k=0$,
		and \eqref{ES:DataInfinity2} with $j,k=0$,
		are ensured by the initial data assumptions \eqref{A:regularity} and \eqref{A:regularity2}. We first prove the case \eqref{ES:DataInfinity1} with $j=1$ by applying \eqref{ES:Linftytransport} to equation \eqref{EQ:e4eAlrho1.1} with $\varphi:=\mathring{e_3}\left(\mathring{e_4},\iniangp\right)^k\inilogDensity$:
		\begin{align}\label{ES:DataInfinity1.1}
			\norm{\mathring{e_3}\left(\mathring{e_4},\iniangp\right)^k\inilogDensity}_{L^\infty(S_{0,\ubar})}\lesssim&\norm{\mathring{e_3}\left(\mathring{e_4},\iniangp\right)^k\inilogDensity}_{L^\infty(S_{0,0})}+\int_{0}^{\ubar}\left(\norm{\mathring{e_3}\left(\mathring{e_4},\iniangp\right)^k\inilogDensity}_{L^\infty(S_{0,\ubarprime})}+1\right)\gensmoothfunction_{0,k}(\iniLtwoS)\diff\ubarprime,
		\end{align}
		where $\gensmoothfunction_{0,k}(\iniLtwoS)$ is a smooth increasing function of $\iniLtwoS$ such that $\gensmoothfunction_{0,k}(0)=0$, which represents the size of all the terms whose bound has been obtained directly from the initial data or the induction hypotheses. In particular, we emphasize that the term $\left(\mathring{e_4},\iniangp\right)^k\iniangp^2\inilogDensity$ (corresponding to the first term on the RHS of \eqref{EQ:e4eAlrho1.1}) causes the loss of derivative in the direction of the transversal derivatives compared to the tangential derivatives, since this term takes $k+2$ tangential derivatives on fluid variables while the LHS of \eqref{ES:DataInfinity1.1} takes one transversal derivative and $k$ tangential derivatives; this phenomenon appears even in the linear wave equation. Nevertheless, this term can be controlled by the initial data assumptions \eqref{A:regularity}.
		The desired bound  \eqref{ES:DataInfinity1} for $\inilogDensity$ with $j=1$ is obtained by applying Gr\"onwall's inequality to \eqref{ES:DataInfinity1.1}.
		Then, \eqref{ES:DataInfinity1} for $\inivel$ with $j=1$ follows from commuting \eqref{E:e3v} with $\left(\mathring{e_4},\iniangp\right)^k$, \eqref{ES:DataInfinity1} for $\inilogDensity$, initial data assumptions \eqref{A:regularity}, and the induction hypothesis. Then, \eqref{ES:DataInfinity2} with $j=0$ follows from commuting \eqref{E:ENTROPYEULEREQUATIONRELATIVETONULLFRAME} with $\mathring{e_3}\left(\mathring{e_4},\iniangp\right)^k$, commuting \eqref{EQ:SchematicDivGradEnt} with $\left(\mathring{e_4},\iniangp\right)^k$, initial data assumptions \eqref{A:regularity},  and previously proved results \eqref{ES:DataInfinity1}.
		
		Now, assuming \eqref{ES:DataInfinity1} holds for all $j_1\leq J_1$ with $J_1\leq \frac{N}{4}-\frac{1}{2}$, and \eqref{ES:DataInfinity2} holds for all $j_2\leq J_2$ with $J_2\leq \frac{N}{4}-\frac{3}{2}$, we first prove \eqref{ES:DataInfinity1} for $\inilogDensity$ with $j=J_1+1$. We apply \eqref{ES:Linftytransport} to equation \eqref{EQ:e4eAlrho1.1} with $\varphi:=\mathring{e_3}^{j}\left(\mathring{e_4},\iniangp\right)^k\inilogDensity$:
		\begin{align}\label{ES:DataInfinity1.2}
			\norm{\mathring{e_3}^{j}\left(\mathring{e_4},\iniangp\right)^k\inilogDensity}_{L^\infty(S_{0,\ubar})}\lesssim&\norm{\mathring{e_3}^{j}\left(\mathring{e_4},\iniangp\right)^k\inilogDensity}_{L^\infty(S_{0,0})}+\int_{0}^{\ubar}\left(\norm{\mathring{e_3}^{j}\left(\mathring{e_4},\iniangp\right)^k\inilogDensity}_{L^\infty(S_{0,\ubarprime})}+1\right)\gensmoothfunction_{j,k}(\iniLtwoS)\diff\ubarprime,
		\end{align}
		where $\gensmoothfunction_{j,k}(\iniLtwoS)$ is a smooth increasing function of $\iniLtwoS$ such that $\gensmoothfunction_{j,k}(0)=0$, which represents the size of all the terms whose bounds have been obtained directly by the initial data or the induction hypotheses. We again emphasize that the term $\mathring{e_3}^{j-1}\left(\mathring{e_4},\iniangp\right)^k\angp^2\inilogDensity$ (corresponding to the first term on the RHS of \eqref{EQ:e4eAlrho1.1}) is the term that causes the loss of derivatives in the transversal direction, which is controlled by the induction hypothesis. The desired bound  \eqref{ES:DataInfinity1} for $\inilogDensity$ is obtained by applying Gr\"onwall's inequality on \eqref{ES:DataInfinity1.2}. Then, \eqref{ES:DataInfinity1} for $\inivel$ follows from commuting \eqref{E:e3v} with $\mathring{e_3}^{j-1}\left(\mathring{e_4},\iniangp\right)^k$, \eqref{ES:DataInfinity1} for $\inilogDensity$, and the induction hypothesis. Then, \eqref{ES:DataInfinity2} follows from commuting \eqref{E:ENTROPYEULEREQUATIONRELATIVETONULLFRAME} with $\mathring{e_3}^{j}\left(\mathring{e_4},\iniangp\right)^k$, commuting \eqref{EQ:SchematicDivGradEnt} with $\mathring{e_3}^{j-1}\left(\mathring{e_4},\iniangp\right)^k$ and previously proved results \eqref{ES:DataInfinity1}.
	\end{proof}
	
	\begin{corollary}[$L^\infty_{\ubar}L^\infty_{\gsphere}$-bound for the acoustical quantities on the characteristic initial hypersurfaces]\label{CO:InftyBoundfortheData}
		We impose the assumptions in Proposition \ref{PR:InftyBoundfortheData}. Moreover, let $u,\ubar$ be the acoustical functions and $(\sphere^1,\sphere^2)$ be the geometric angular coordinates constructed in Section \ref{SS:InitialCharData} and Section \ref{SS:CONSTRUCTIONOFACOUSTICEIKONALFUNCTIONSANDGEOMETRICANGULARCOORDINATESINRENDALLREGION}. Let $\iniACC$ be the collection of connection coefficients (as in Definition \ref{D:Defconnectioncoeff}) on the initial hypersurfaces. Let $\iniinvfolia,\iniuinvfolia,\inilapse$ be the lapses (as in Section \ref{SS:SUBSETSOFSPACETIME}) on the initial hypersurfaces. Then, the following estimates are satisfied by the acoustical quantities on $\inihypersur_{[0,\ubarfinite]}$:
		\begin{subequations}\label{ES:GeometryInfinity}
			\begin{align}
				\label{ES:GeometryInfinity1}\norm{\iniACC}_{L^\infty(\inihypersur_{[0,\ubarfinite]})}\leq& \tilde{\gensmoothfunction}(\iniLtwoS),\\
				\label{ES:GeometryInfinity2}\norm{\ln\iniinvfolia,\ln\iniuinvfolia,\ln\inilapse}_{L^\infty(\inihypersur_{[0,\ubarfinite]})}\leq& \tilde{\gensmoothfunction}(\iniLtwoS).
			\end{align}
		\end{subequations}
		where $\tilde{\gensmoothfunction}(\iniLtwoS)<\infty$ is a smooth increasing function of $\iniLtwoS$ such that $\tilde{\gensmoothfunction}(0)=0$.
	\end{corollary}
	\begin{proof}[Proof of Corollary \ref{CO:InftyBoundfortheData}]
		By assumptions \eqref{A:regularity}, the embedding $\embed$ is smooth and $\norm{\pfour^k\embed^{\alpha}}_{L^\infty(\inihypersur_{[0,\ubarfinite]})}\leq \iniLtwoS$ for every integer $k\geq0$. By Definition \ref{D:LDIFFEOCONSTRUCTEDOUTOFEMBEDDINGS} and \ref{D:NULLFRAMEDATA}, we have:
		\begin{align}\label{ES:GeometryInfinity1.1}
			\mathring{e_I}=\gensmoothfunction(\pfour\embed,(\inivel^i,\inilogDensity,\inient)).
		\end{align}
		Moreover, connection coefficients can be written schematically as follows:
		\begin{align}\label{ES:GeometryInfinity1.2}
			\iniACC=\Dfour_{\mathring{e_I}}\mathring{e_J}=\mathring{e_I}(\mathring{e_J})+\Chfour_{IJ}^K{\mathring{e_I}}=\gensmoothfunction(\pfour^{\leq 2}\embed,\pfour^{\leq 1}(\inivel^i,\inilogDensity,\inient)).
		\end{align}
		Then, \eqref{ES:GeometryInfinity1} follows from the algebraic equation \eqref{ES:GeometryInfinity1.2} and Proposition \ref{PR:InftyBoundfortheData}. 
		
		Now we prove \eqref{ES:GeometryInfinity2}. First we recall the assumptions that $\frac{\p\embed}{\p\iniu},\frac{\p\embed}{\p\iniubar}$ are future directed vectorfields. In particular, by assumption \eqref{E:INITIALLAPSEBOUNDS} and Definition \ref{D:lapsequantities}, we have the following bounds as in \eqref{E:E40ANDE30INITIALBOUNDS} for the lapse $\inilapse$ on $\inihypersur_{[0,\ubarfinite]}$:
		\begin{align}
			\frac{1}{\boundinilapse}
			&
			\leq
			\inilapse
			\leq 
			\boundinilapse,
		\end{align} 
		where $\boundinilapse$ is a universal constant.
		Then, by \eqref{D:iniuinvfolia}, we have:
		\begin{align}
			\frac{1}{\boundinilapse^2}\leq\iniuinvfolia\leq&\boundinilapse^2,&
			\text{on } &\inihypersur_{[0,\ubarfinite]}.
		\end{align}
		Lastly, we prove \eqref{ES:GeometryInfinity2} for $\ln\iniinvfolia$. We apply \eqref{ES:Linftytransport} on \eqref{E:e_4e_43} with $\varphi:=\ln\iniinvfolia$. Then, \eqref{ES:GeometryInfinity2} for $\ln\iniinvfolia$ is obtained by previously proved estimates \eqref{ES:DataInfinity}, \eqref{ES:GeometryInfinity} and the fact that $\abs{\ln\iniinvfolia}\leq2\ln\boundinilapse$ on $S_{0,0}$.
	\end{proof}
	
	\begin{proposition}[$L^\infty_{\ubar}L^2_{\gsphere}$-bound for the fluid variables of the characteristic initial data]\label{PR:L2BoundfortheData}
		Let $\embed$ be the characteristic embedding defined in Definition \ref{D:LDIFFEOCONSTRUCTEDOUTOFEMBEDDINGS}. Let $(\inivel^i,\inilogDensity,\inient)$ be the smooth characteristic initial data that we prescribed/constructed in Section \ref{S:CID}-\ref{S:localexistence} so that their derivatives are well-defined at all orders on $\inihypersur_{[0,\ubarfinite]}\bigcup\inihypersurbar_{[0,\ufinite]}$, where $\ufinite,\ubarfinite$ are determined in Proposition \ref{P:DETERMININGTRANSVERALDERIVATIVESANDSATISFYINGEQUATIONSONINITIALNULLHYPERSURFACES}. We impose the regularity assumptions on the initial data in Sobolev norms \eqref{A:regularity}-\eqref{A:regularity2}. Then, the following estimates are satisfied by the mixed derivatives of the fluid variables on $\inihypersur_{[0,\ubarfinite]}$ with non-negative indices $j,k$:
		\begin{subequations}\label{ES:DataL2}
			\begin{align}
				\label{ES:DataL21}&\twonorms{\mathring{e_3}^{j}\left(\mathring{e_4},\iniangp\right)^k(\inivel^i,\inilogDensity)}{\ubar}{\infty}{\gsphere}{2}{\inihypersur_{[0,\ubarfinite]}}\leq\tilde{\gensmoothfunction}_{j,k}(\iniLtwoS),&
				\text{ for }& 0\leq 2j+k\leq  N,\\
				\label{ES:DataL22}&\twonorms{\mathring{e_3}^{j}\left(\mathring{e_4},\iniangp\right)^k(\mathring{e_3}^2\inient,\iniDivGradEnt)}{\ubar}{\infty}{\gsphere}{2}{\inihypersur_{[0,\ubarfinite]}}\leq\tilde{\gensmoothfunction}_{j,k}(\iniLtwoS),&
				\text{ for }& 0\leq 2j+k\leq  N-2.
			\end{align}
		\end{subequations}
		where $\tilde{\gensmoothfunction}_{j,k}(\iniLtwoS)<\infty$ are smooth increasing functions of $\iniLtwoS$ such that $\tilde{\gensmoothfunction}_{j,k}(0)=0$.
	\end{proposition}
	\begin{proof}[Proof of Proposition \ref{PR:InftyBoundfortheData}]
		We prove \eqref{ES:DataL21}-\eqref{ES:DataL22} by induction with respect to $j$ in a similar fashion as in Proposition \ref{PR:InftyBoundfortheData}. To start with, we only consider the case where \eqref{ES:DataL21} with $\frac{N}{2}+1\leq 2j+k\leq  N$ and \eqref{ES:DataL22} with $\frac{N}{2}-1\leq 2j+k\leq  N-2$ since the lower-order cases can be covered by Proposition \ref{PR:InftyBoundfortheData}.	
		
		Assuming \eqref{ES:DataL21} holds for all $j_1\leq J_1$ with $J_1\leq N-1$, and \eqref{ES:DataL22} holds for all $j_2\leq J_2$ with $J_2\leq N-3$, we first prove \eqref{ES:DataL21} for $\inilogDensity$ with $j=J_1+1$. We apply \eqref{ES:Linftytransport} to equation \eqref{EQ:e4eAlrho1.1} with $\varphi:=\mathring{e_3}^{j}\left(\mathring{e_4},\iniangp\right)^k\inilogDensity$:
		\begin{align}\label{ES:DataL21.1}
			\norm{\mathring{e_3}^{j}\left(\mathring{e_4},\iniangp\right)^k\inilogDensity}_{L^2(S_{0,\ubar})}\lesssim&\norm{\mathring{e_3}^{j}\left(\mathring{e_4},\iniangp\right)^k\inilogDensity}_{L^2(S_{0,0})}+\int_{0}^{\ubar}\left(\norm{\mathring{e_3}^{j}\left(\mathring{e_4},\iniangp\right)^k\inilogDensity}_{L^2(S_{0,\ubarprime})}+1\right)\gensmoothfunction_{j,k}(\iniLtwoS)\diff\ubarprime,
		\end{align}
		where $\gensmoothfunction_{j,k}(\iniLtwoS)$ is a smooth increasing function of $\iniLtwoS$ such that $\gensmoothfunction_{j,k}(0)=0$, which represents the size of all the terms whose bounds have been obtained directly by the initial data or the induction hypothesis. Similar to Proposition \ref{PR:InftyBoundfortheData}, we note that the term $\mathring{e_3}^{j-1}\left(\mathring{e_4},\iniangp\right)^k\angp^2\inilogDensity$ (corresponding to the first term on the RHS of \eqref{EQ:e4eAlrho1.1}) is the term that causes the loss of derivatives in the transversal direction, which is controlled by the induction hypothesis. The desired bound  \eqref{ES:DataL21} for $\inilogDensity$ is obtained by applying Gr\"onwall's inequality on \eqref{ES:DataL21.1}. Then, \eqref{ES:DataL21} for $\inivel$ follows from commuting \eqref{E:e3v} with $\mathring{e_3}^{j-1}\left(\mathring{e_4},\iniangp\right)^k$, \eqref{ES:DataL21} for $\inilogDensity$, and the induction hypothesis. Then, \eqref{ES:DataL22} follows from commuting \eqref{E:ENTROPYEULEREQUATIONRELATIVETONULLFRAME} with $\mathring{e_3}^{j}\left(\mathring{e_4},\iniangp\right)^k$, commuting \eqref{EQ:SchematicDivGradEnt} with $\mathring{e_3}^{j-1}\left(\mathring{e_4},\iniangp\right)^k$ and previously proved results \eqref{ES:DataL21}.
	\end{proof}
	
	\begin{corollary}[Bounds for the modified vorticity on the characteristic initial hypersurfaces]\label{CO:ModivortInfinity}
		We impose the assumptions in Proposition \ref{PR:InftyBoundfortheData} and Proposition \ref{PR:L2BoundfortheData}. Let $\modivort$ be the modified vorticity defined as in Definition \ref{D:modifiedfluidvariables} and $\inimodivort:=\modivort\restriction_{\inihypersur\cup\inihypersurbar}$. Then, the following estimates are satisfied by the modified vorticity data $\inimodivort$ on $\inihypersur_{[0,\ubarfinite]}$:
		\begin{subequations}\label{ES:ModivortInfinity}
			\begin{align}
				\label{ES:ModivortInfinity1}&\twonorms{\left(\mathring{e_4},\iniangp\right)^k\inimodivort}{\ubar}{\infty}{\gsphere}{\infty}{\inihypersur_{[0,\ubarfinite]}}\leq\tilde{\gensmoothfunction}_{j,k}(\iniLtwoS),&
				\text{ for }& 0\leq k\leq  \frac{N}{2}-2,\\
				\label{ES:ModivortInfinity2}&\twonorms{\mathring{e_3}^{j}\left(\mathring{e_4},\iniangp\right)^k\mathring{e_3}\inimodivort}{\ubar}{\infty}{\gsphere}{\infty}{\inihypersur_{[0,\ubarfinite]}}\leq\tilde{\gensmoothfunction}_{j,k}(\iniLtwoS),&
				\text{ for }& 0\leq 2j+k\leq  \frac{N}{2}-3,\\
				\label{ES:ModivortInfinity3}&\twonorms{\left(\mathring{e_4},\iniangp\right)^k\inimodivort}{\ubar}{\infty}{\gsphere}{2}{\inihypersur_{[0,\ubarfinite]}}\leq\tilde{\gensmoothfunction}_{j,k}(\iniLtwoS),&
				\text{ for }& 0\leq k\leq  N-2,\\
				\label{ES:ModivortInfinity4}&\twonorms{\mathring{e_3}^{j}\left(\mathring{e_4},\iniangp\right)^k\mathring{e_3}\inimodivort}{\ubar}{\infty}{\gsphere}{2}{\inihypersur_{[0,\ubarfinite]}}\leq\tilde{\gensmoothfunction}_{j,k}(\iniLtwoS),&
				\text{ for }& 0\leq 2j+k\leq  N-3.
			\end{align}
		\end{subequations}
	\end{corollary}
	\begin{proof}[Proof of Corollary \ref{CO:ModivortInfinity}]
		Expressing $\left(\mathring{e_4},\iniangp\right)^k\inimodivort$ by commuting \eqref{E:reforboxgv} with $\left(\mathring{e_4},\iniangp\right)^k$, \eqref{ES:ModivortInfinity1} and \eqref{ES:ModivortInfinity3} follow directly from the previously proved results \eqref{ES:DataInfinity} and \eqref{ES:DataL2}. Finally, we write \eqref{E:BC} schematically as follows:
		\begin{align}
			\mathring{	e_3} \inimodivort
			& 
			= \mathring{	e_4} \inimodivort
			+
			2 \inilapse^{-1} \iniangB \inimodivort+\gensmoothfunction\left(\inimodivort,\pfour^{\leq 1}(\inivel^i,\inilogDensity),\pfour^{\leq2}\inient)\right).
		\end{align}
		Then, \eqref{ES:ModivortInfinity2} and \eqref{ES:ModivortInfinity4} follow from commuting \eqref{E:BC} with $\mathring{e_3}^{j}\left(\mathring{e_4},\iniangp\right)^k$ and previously proved results \eqref{ES:DataInfinity}, \eqref{ES:DataL2}, \eqref{ES:ModivortInfinity1} and \eqref{ES:ModivortInfinity3}.
	\end{proof}
	
	\subsection{Table of Regularities of Quantities on Initial Characteristic Hypersurfaces}\label{SS:iniregularity}
	In this subsection, we gather the regularity of quantities on $\inihypersur$ that we obtained in Section \ref{SS:regularityoffluid}, where the integer $N\geq17$.
	\begin{center}

		\begin{tabular}{|c|c|c|} 
			\hline
			Symbol & Bounded $L^\infty_{\ubar}L^2_{\gsphere}$-norm&  Bounded $L^{\infty}_{\ubar}L^{\infty}_{\gsphere}$-norm  
			\\
			\hline\hline
			$\left[\mathring{e_3}^j(\mathring{e_4},\iniangp)^{k}\right](\inilogDensity,\mathring{v}^i,\inient)$ & $2j+k\leq N$ & $2j+k\leq\frac{N}{2}+1$\\
			\hline
			$[e_3^j(\mathring{e_4},\iniangp)^k]\left(\pfour\inient\right)$& $2j+k\leq N-1$&$2j+k\leq\frac{N}{2}$\\ 
			\hline
			$\left[\mathring{e_3}^j(\mathring{e_4},\iniangp)^k\right](\iniDivGradEnt,\pfour^2\inient)$ & $2j+k\leq N-2$& $2j+k\leq\frac{N}{2}-1$\\
			\hline
			$[e_3^j(\mathring{e_4},\iniangp)^k]\inimodivort$ & $2j+k\leq N-3$&$2j+k\leq \frac{N}{2}-2$ \\
			\hline
			$[e_3^j(\mathring{e_4},\iniangp)^k]\pfour\inimodivort$ & $2j+k\leq N-4$&$2j+k\leq \frac{N}{2}-3$\\
			\hline
		\end{tabular}
	\end{center}
	
	\section{Concise Statement of Sobolev Regularity of the Initial Data and the Bootstrap Assumptions}\label{S:SectionAssumptions}
	In this section, we state the Sobolev norms on the initial data that we will use in our main a priori estimates. We also state the bootstrap assumptions that we will
	use in deriving our a priori estimates.

	Let $\ufinite,\ubarfinite$ be as in Proposition \ref{P:DETERMININGTRANSVERALDERIVATIVESANDSATISFYINGEQUATIONSONINITIALNULLHYPERSURFACES}. Let the embedding $\embed$ be defined as in Definition \ref{D:LDIFFEOCONSTRUCTEDOUTOFEMBEDDINGS}. In the rest of this article, we consider $u\in I_{\ufinite}:=[0,\ufinite]$ and $\ubar\in I_{\ubarfinite}:=[0,\ubarfinite]$, $\inihypersur_{[0,\ubarfinite]}=\embed\left(\{0\}\times I_{\ubarfinite}\times\Stwo\right)$ and $\inihypersurbar_{[0,\ufinite]}=\embed\left(I_{\ufinite}\times\{0\}\times \Stwo\right)$.
	
	\subsection{Restatement of Sobolev Regularity of Data}\label{SS:iniregularity2}
	Recalling the regularities of the data in Section \ref{SS:iniregularity}, in this subsection, we provide precise assumptions for initial data that align with the computed regularities. That is, assumptions (\ref{A:regularity}) imply the assumptions in this Section\footnote{In Section \ref{SS:regularityoffluid}, we derive the estimates for the initial data in $L^\infty_{\ubar} L^2_{\gsphere}$-norm and $L^\infty_{\ubar} L^\infty_{\gsphere}$-norm provided the fact that $u_*,\ubar_*$ are finite constants.}. We again emphasize that although we assume that the initial data are $C^{\infty}$, all of our estimates are in terms of Sobolev norms. Let $\iniepsilon>0$ be the size of the initial data. 
	For $j,k\geq0$, $A=1,2$, integer $N\geq17$, on\footnote{For convenience in proving a priori estimates, in the rest of the article, we use the notation $\hypersur_{0}=\inihypersur$ and $\hypersurbar_{0}=\inihypersurbar$.} $\hypersur_{0}$, we assume that:
	\begin{subequations}\label{A:H}
		\begin{align}
			\twonorms{\left[\mathring{e_3}^j(\mathring{e_4},\iniangp)^{k}\right](\inilogDensity,\mathring{v}^i,\inient)}{\ubar}{2}{\gsphere}{2}{\hypersur_{0}}\leq& \iniepsilon,&
			\textnormal{for } & 0\leq 2j+k\leq N,\\
			\twonorms{\left[\mathring{e_3}^j(\mathring{e_4},\iniangp)^k\right](\iniDivGradEnt,\pfour\inigradEnt,\pfour^2\inient,\inivortrenormalized)}{\ubar}{2}{\gsphere}{2}{\hypersur_{0}}\leq& \iniepsilon,&
			\textnormal{for } & 0\leq 2j+k\leq N-2,\\
			\twonorms{(\mathring{e_4},\iniangp)^k\inimodivort}{\ubar}{2}{\gsphere}{2}{\hypersur_{0}}\leq& \iniepsilon,&
			\textnormal{for } & 0\leq k\leq N-3,\\
			\twonorms{[e_3^j(\mathring{e_4},\iniangp)^k]\pfour\inimodivort}{\ubar}{2}{\gsphere}{2}{\hypersur_{0}}\leq& \iniepsilon,&
			\textnormal{for } & 0\leq 2j+k\leq N-4.
		\end{align}
	\end{subequations}
	Similarly, on $\hypersurbar_{0}$, we assume that:
	\begin{subequations}\label{A:Hbar}
		\begin{align}
			\twonorms{\left[\mathring{e_4}^j(\mathring{e_3},\iniangp)^{k}\right](\inilogDensity,\mathring{v}^i,\inient)}{u}{2}{\gsphere}{2}{\hypersurbar_{0}}\leq& \iniepsilon,&
			\textnormal{for } & 0\leq 2j+k\leq N,\\
			\twonorms{\left[\mathring{e_4}^j(\mathring{e_3},\iniangp)^k\right](\iniDivGradEnt,\pfour\inigradEnt,\pfour^2\inient)}{u}{2}{\gsphere}{2}{\hypersurbar_{0}}\leq& \iniepsilon,&
			\textnormal{for } & 0\leq 2j+k\leq N-2,\\
			\twonorms{(\mathring{e_3},\iniangp)^k\inimodivort}{u}{2}{\gsphere}{2}{\hypersurbar_{0}}\leq& \iniepsilon,&
			\textnormal{for } & 0\leq k\leq N-3,\\
			\twonorms{[e_4^j(\mathring{e_3},\iniangp)^k]\pfour\inimodivort}{u}{2}{\gsphere}{2}{\hypersurbar_{0}}\leq& \iniepsilon,&
			\textnormal{for } & 0\leq 2j+k\leq N-4.
		\end{align}
	\end{subequations}
	Moreover, we assume that:
	\begin{subequations}\label{A:Boundness}
		\begin{align}
			\label{Datageometry1}\norm{\pfour^k(\mathring{v}^i,\inilogDensity,\inient),\pfour^{k-1}(\inimodivort,\iniDivGradEnt)}_{L^\infty(\hypersur_{0}\bigcup\hypersurbar_{0})}\leq& \iniepsilon,& &\textnormal{for }0\leq k\leq\frac{N-1}{4},\\
			\label{Datageometry2}\norm{\pfour^k(\iniu,\iniubar),\pfour^j\ln\iniinvfolia,\pfour^j\ln\iniuinvfolia,\pfour^j\ln\inilapse,\pfour^j\iniACC}_{L^\infty(\hypersur_{0}\bigcup\hypersurbar_{0})}\leq& \iniepsilon,&&\textnormal{for }0\leq k\leq 3, 0\leq j\leq 1.
		\end{align}
	\end{subequations}
	In addition, we assume that:
	\begin{align}
		\label{Datageometry3}\mathring{a}^{-1}\seabs{\xi}^2\leq\inigsphere_{AB}\xi^A\xi^B\leq& \mathring{a}\seabs{\xi}^2,&
		\text{ for any }& \suubar \text{-tangent vectorfield }\xi,
	\end{align}
	where $\esphere$ is the Euclidean metric on $\Stwo$, and $1\leq\mathring{a}\leq \iniepsilon$.
	\subsection{Bootstrap Assumptions}\label{SS:BA}
	Let\footnote{$\BAinfty>4\iniepsilon$ is enough to close the estimates in our argument, while it is not necessarily the sharp bound.} $\BAinfty:=\gensmoothfunction(\iniepsilon)>\max\{4\iniepsilon,1\}$. $(v^i,\logDensity,\Ent)$ is the solution to the compressible Euler equations \eqref{E:firstorderformulation}. $\modivort,\DivGradEnt$ are the modified fluid quantities defined in Definition \ref{D:modifiedfluidvariables}. $u,\ubar$ are the solutions to the acoustical eikonal equations \eqref{EQ:eikonalequations}. Let the Ricci coefficients $\ACC$ be as defined in Definition \ref{D:Defconnectioncoeff}. We make the following bootstrap assumptions:
	
	\noindent\textbf{Bootstrap assumptions of fluid variables}
	\begin{align}
		\label{BAfluid}\norm{\pfour^k(v^i,\logDensity,\Ent),\pfour^{k-1}(\modivort,\DivGradEnt)}_{L^\infty(\region)}\leq& \BAinfty,& 0\leq k\leq&\frac{N-1}{4}.
	\end{align}
	
	\noindent\textbf{Bootstrap assumptions on geometric structure}
	\begin{align}
		\label{BAacoustic}\norm{\pfour^k(u,\ubar),\pfour^j\ln\invfolia,\pfour^j\ln\uinvfolia,\pfour^j\ln\lapse,\pfour^j\ACC}_{L^\infty(\region)}\leq& \exp(\BAinfty),&
		0&\leq k\leq 3, 0\leq j\leq 1\\
		\label{BAgeometry3}\tilde{a}^{-1}\seabs{\xi}^2\leq\gsphere_{AB}\xi^A\xi^B\leq& \tilde{a}\seabs{\xi}^2,&
		\text{ for any }& \suubar \text{-tangent vectorfield }\xi,
	\end{align}
	where $\esphere$ is the Euclidean metric on $\Stwo$, and $\mathring{a}\leq\tilde{a}\leq\BAinfty$.
	
	\section{Energies, Energy Identities, and Energy Inequalities} \label{S:ENERGYANDENERGYIDENTITIES}
	In this section, we define the energies that we will use in our a priori estimates. 
	We also provide various energy identities and inequalities for these quantities.

	\subsection{Various Energies}\label{SS:energies}
	In this section, we construct and define various energies that we will use in deriving a priori energy estimates.
	
	\begin{definition}[Energy-momentum tensor, energy current, and deformation tensor]\label{D:energymomentum} We define the energy-momentum tensor $Q_{\mu\nu}[\varphi]$ associated to a scalar function $\varphi$ to be the following tensorfield:
		\begin{align}
			Q_{\mu\nu}[\varphi]:=\p_\mu\varphi\p_\nu\varphi-\frac{1}{2}\gfour_{\mu\nu}(\gfour^{-1})^{\alpha\beta}\p_\alpha\varphi\p_\beta\varphi.
		\end{align} 
		
		Given $\varphi$ and any multiplier vectorfield $\mathbf{X}$, we define the corresponding energy current $\Jenarg{\mathbf{X}}{\alpha}[\varphi]$ vectorfield as follows:
		\begin{align}
			\Jenarg{\mathbf{X}}{\alpha}[\varphi]:=Q^{\alpha\beta}[\varphi]\mathbf{X}_{\beta}-\varphi^2\mathbf{X}^\alpha.
		\end{align}
		
		We define the deformation tensor of $\mathbf{X}$ as follows:
		\begin{align}\label{D:deformtensor}
			\deform{\mathbf{X}}_{\alpha\beta}:=\Dfour_{\alpha}\mathbf{X}_\beta+\Dfour_{\beta}\mathbf{X}_\alpha,
		\end{align}
		where $\Dfour$ is the Levi-Civita connection of the acoustical metric $\gfour$.
		
		We have the following well-known divergence identity:
		\begin{align}\label{E:DalphaJalpha}
			\Dfour_\alpha\Jenarg{\mathbf{X}}{\alpha}[\varphi]=\boxg\varphi(\mathbf{X}\varphi)+\frac{1}{2}Q^{\mu\nu}[\varphi]\deform{\mathbf{X}}_{\mu\nu}-2\varphi\mathbf{X}\varphi-\frac{1}{2}\varphi^2(\gfour^{-1})^{\mu\nu}\deform{\mathbf{X}}_{\mu\nu}.
		\end{align}
	\end{definition}
	\begin{definition}[Energies for the ``wave-part'']\label{DE:waveenergy1}
		Recall (\ref{DE:Stau}) the definition of $\Stau^{u,\ubar}$. For $0\leq \ubar_0\leq\ubar\leq\ubar_1$ and $0\leq u_0\leq u\leq u_1$, $u+\ubar\leq\uptau\leq\max\{u_0+\ubar_1,u_1+\ubar_0\}$, we define the energy $\bEnergy^{(1)}[\varphi](\uptau)$ on $\Stau^{u,\ubar}$ as follows, where $\Timelike$ is the future-directed $\gfour$-timelike vectorfield defined in (\ref{D:Timelike}):
		\begin{align}\label{D:Energy1}
			\bEnergy_{u,\ubar}^{(1)}[\varphi](\uptau):=\int_{\Stau^{u,\ubar}}\Jenarg{\mathbf{\Timelike}}{\alpha}[\varphi](\Timelike)_\alpha\diff\tauvol=\int_{\Stau^{u,\ubar}}\left\{\frac{1}{2}(e_4\varphi)^2+\frac{1}{2}(e_3\varphi)^2+\frac{1}{2}\gsphere_{AB}e_A\varphi e_B\varphi+\varphi^2\right\}\diff\tauvol,
		\end{align}
		where $\gtau$ is the induced metric on $\Stau$ (see (\ref{D:defg})) and $\diff\tauvol$ is the volume form on $\Stau$ with respect to $\gtau$ (see Definition \ref{D:volumeforms}). In particular, for $\uptau=u_0+\ubar_0$ or $\uptau=u_1+\ubar_1$
		\begin{subequations}\label{D:Energy1.5}
			\begin{align}
				\bEnergy_{u_0,\ubar_0}^{(1)}[\varphi](u_0+\ubar_0):=&\int_{S_{u_0,\ubar_0}}\Jenarg{\mathbf{\Timelike}}{\alpha}[\varphi](\Timelike)_\alpha\diff\tauvol=\int_{S_{u_0,\ubar_0}}\left\{\frac{1}{2}(e_4\varphi)^2+\frac{1}{2}(e_3\varphi)^2+\frac{1}{2}\gsphere_{AB}e_A\varphi e_B\varphi+\varphi^2\right\}\diff\spherevol,\\
				\bEnergy_{u_0,\ubar_0}^{(1)}[\varphi](u_1+\ubar_1):=&\int_{S_{u_1,\ubar_1}}\Jenarg{\mathbf{\Timelike}}{\alpha}[\varphi](\Timelike)_\alpha\diff\tauvol=\int_{S_{u_1,\ubar_1}}\left\{\frac{1}{2}(e_4\varphi)^2+\frac{1}{2}(e_3\varphi)^2+\frac{1}{2}\gsphere_{AB}e_A\varphi e_B\varphi+\varphi^2\right\}\diff\spherevol,
			\end{align}
		\end{subequations}
		where $\spherevol$ is the volume form on $\suubar$ with respect to $\gsphere$ (see (\ref{D:defg})).
		
		For $0\leq \ubar_0<\ubar\leq\ubar_1$ and $0\leq u_0<u\leq u_1$, we define the energy fluxes $\Energy^{(1)}_{u}[\varphi](\ubar_0,\ubar_1)$ and $\uEnergy^{(1)}_{\ubar}[\varphi](u_0,u_1)$ on null hypersurfaces $\hypersur_{u}^{\ubar_0,\ubar_1}$ and $\hypersurbar_{\ubar}^{u_0,u_1}$ (defined in Section \ref{SS:SUBSETSOFSPACETIME}) as follows:
		\begin{subequations}\label{D:Energy2}
			\begin{align}
				\Energy^{(1)}_u[\varphi](\ubar_0,\ubar_1):=&\int_{\hypersur_u\bigcap\{\ubar_0\leq\ubar\leq\ubar_1\}}\Jenarg{\mathbf{\Timelike}}{\alpha}[\varphi](e_4)_\alpha\diff\spherevol\diff \ubar\\
				\notag=&\int_{\hypersur_{u}\bigcap\{\ubar_0\leq\ubar\leq\ubar_1\}}\left\{\frac{1}{2}(e_4\varphi)^2+\frac{1}{2}\gsphere_{AB}e_A\varphi e_B\varphi+\varphi^2\right\}\diff\spherevol\diff \ubar,\\
				\uEnergy^{(1)}_{\ubar}[\varphi](u_0,u_1):=&\int_{\hypersurbar_{\ubar}\bigcap\{u_0\leq u\leq u_1\}}\Jenarg{\mathbf{\Timelike}}{\alpha}[\varphi](e_3)_\alpha \diff\spherevol\diff u\\
				\notag=&\int_{\hypersurbar_{\ubar}\bigcap\{u_0\leq u\leq u_1\}}\left\{\frac{1}{2}(e_3\varphi)^2+\frac{1}{2}\gsphere_{AB}e_A\varphi e_B\varphi+\varphi^2\right\}\diff\spherevol\diff u.
			\end{align}
		\end{subequations}
	\end{definition}
	\begin{definition}[Energies for the ``transport-part'']\label{d:Energy34}
		For $0\leq \ubar_0\leq\ubar\leq\ubar_1$ and $0\leq u_0\leq u\leq u_1$, $u+\ubar\leq\uptau\leq\max\{u_0+\ubar_1,u_1+\ubar_0\}$, we define the energy $\bEnergy^{(2)}[\varphi]$ on $\Stau^{u,\ubar}$ (defined in (\ref{DE:Stau})) as follows:
		\begin{align}\label{D:Energy3}
			\bEnergy^{(2)}_{u,\ubar}[\varphi](\uptau):=\int_{\Stau^{u,\ubar}}\varphi^2\diff\tauvol,
		\end{align}
		where $\gtau$ is the induced metric on $\Stau$ (see (\ref{D:defg})) and $\diff\tauvol$ is the volume form on $\Stau$ with respect to $\gtau$ (see Definition \ref{D:volumeforms}). In particular, for $\uptau=u_0+\ubar_0$ or $\uptau=u_1+\ubar_1$,
		\begin{subequations}\label{D:Energy3.5}
			\begin{align}
				\bEnergy^{(2)}_{u_0,\ubar_0}[\varphi](\uptau):=&\int_{S_{u_0,\ubar_0}}\varphi^2\diff\spherevol,\\
				\bEnergy^{(2)}_{u_1,\ubar_1}[\varphi](\uptau):=&\int_{S_{u_1,\ubar_1}}\varphi^2\diff\spherevol,
			\end{align}
		\end{subequations}
		where $\spherevol$ is the volume form on $\suubar$ with respect to $\gsphere$ (see (\ref{D:defg})).
		
		For $0\leq \ubar_0\leq\ubar\leq\ubar_1$ and $0\leq u_0\leq u\leq u_1$, we define the energy fluxes $\Energy^{(2)}_{u}[\varphi](\ubar_0,\ubar_1)$ and $\uEnergy^{(2)}_{\ubar}[\varphi](u_0,u_1)$ on null hypersurfaces $\hypersur_{u}^{\ubar_0,\ubar_1}$ and $\hypersurbar_{\ubar}^{u_0,u_1}$ (defined in Section \ref{SS:SUBSETSOFSPACETIME}) as follows:
		\begin{subequations}\label{D:Energy4}
			\begin{align}
				\Energy^{(2)}_u[\varphi](\ubar_0,\ubar_1):=&\int_{\hypersur_u\bigcap\{\ubar_0\leq\ubar\leq\ubar_1\}}\varphi^2\diff\spherevol\diff \ubar,\\
				\uEnergy^{(2)}_{\ubar}[\varphi](u_0,u_1):=&\int_{\hypersurbar_{\ubar}\bigcap\{u_0\leq u\leq u_1\}}\varphi^2\diff\spherevol\diff u.
			\end{align}
		\end{subequations}
	\end{definition}
	\begin{definition}[Energies in the space-time region]\label{DE:Energyinspacetime}
		Let $0\leq \ubar_0\leq\ubar_1$ and $0\leq u_0\leq u_1$. Let $\updelta_u:=u_1-u_0$, and $\updelta_{\ubar}:=\ubar_1-\ubar_0$. For vectorfield $V$, we define the energy $\KEnergy^{(1)}[V]$ on the ``triangle" space-time region $\region_1^{u,\ubar,\uptau}$ (defined in (\ref{DE:triangleM})) as follows:
		\begin{align}\label{D:Energy5}
			\KEnergy^{(1)}[V](u,\ubar,\uptau):=\int_{\region_1}\abs{\pfour V}^2\diff\gvol+\int_{S_{u,\uptau-u}}\sgabs{\angV}^2\diff\spherevol+\int_{S_{\uptau-\ubar,\ubar}}\sgabs{\angV}^2\diff\spherevol.
		\end{align}
		
		Let $V$ be a vectorfield. Then, for $0\leq \ubar_0\leq\ubar_1$ and $0\leq u_0\leq u_1$, $u_0+\ubar_0+\min\{\updelta_u,\updelta_{\ubar}\}\leq\uptau_1\leq\uptau_2\leq\max\{u_0+\ubar_1,u_1+\ubar_0\}$, we define the energy $\KEnergy^{(2)}[V]$ on $\region^{\uptau_1,\uptau_2}_2$ ($\region^{\uptau_1,\uptau_2}_2$ is defined in Section \ref{SS:SUBSETSOFSPACETIME}) as follows:\\
		
		\textbf{When $\updelta_u\leq\updelta_{\ubar}$:}
		\begin{align}\label{D:KEnergy2}
			\KEnergy^{(2)}[V](\uptau_1,\uptau_2):=&\int_{\region_2}\abs{\pfour V}^2\diff\gvol+\int_{S_{u_1,\uptau_2-u_1}}\sgabs{\angV}^2\diff\spherevol.
		\end{align}
		
		\textbf{When $\updelta_u>\updelta_{\ubar}$:}
		\begin{align}\label{D:KEnergy3}
			\KEnergy^{(2)}[V](\uptau_1,\uptau_2):=&\int_{\region_2}\abs{\pfour V}^2\diff\gvol+\int_{S_{\uptau_2-\ubar_1,\ubar_1}}\sgabs{\angV}^2\diff\spherevol.
		\end{align}
		
		Moreover, for $\max\{u_0+\ubar_1,u_1+\ubar_0\}\leq\uptau_1\leq\uptau_2\leq u_1+\ubar_1$, we define $\KEnergy^{(3)}[V]$ on the ``trapezoid" region $\region_3$ (defined in Section \ref{SS:SUBSETSOFSPACETIME}) as follows:
		\begin{align}
			\KEnergy^{(3)}[V](\uptau_1,\uptau_2):=\int_{\region_3}\abs{\pfour V}^2\diff\gvol+\int_{S_{\uptau_2-\ubar_1,\ubar_1}}\sgabs{\angV}^2\diff\spherevol+\int_{S_{u_1,\uptau_2-u_1}}\sgabs{\angV}^2\diff\spherevol.
		\end{align}
	\end{definition}
	\begin{definition}[The controlling energy quantities]\label{D:ControllingEnergyQuantities}
		Let $0\leq \ubar_0\leq\ubar\leq\ubar_1$ and $0\leq u_0\leq u\leq u_1$.  For $1\leq k\leq \frac{N-1}{2}$, $1\leq j\leq N-2$, we define the fluid energy quantities on $\Stau$ for $u_0+\ubar_0<\uptau< u_1+\ubar_1$ as follows:
		\begin{subequations}\label{DE:Energy1}
			\begin{align}
				\bEnergy_{(i)}(u,\ubar,\uptau):=&\sum\limits_{\varphi\in\{\pfour^{i-1}(\logDensity,v^i,\Ent)\}}\bEnergy_{u,\ubar}^{(1)}[\varphi](\uptau)\\
				\notag&+\sum\limits_{\varphi\in\{\pfour^{i-1}(\vortrenormalized^i,\gradEnt^i,\modivort^i,\DivGradEnt)\}}\bEnergy_{u,\ubar}^{(2)}[\varphi](\uptau),\\
				\tilde{\bEnergy}_{(k)}(u,\ubar,\uptau):=&\sum\limits_{0\leq i\leq k}\bEnergy_{(i)}(u,\ubar,\uptau).
			\end{align}
		\end{subequations}
		
		Also, we define the fluid energy quantities on null hypersurfaces $\hypersur_u$ and $\hypersurbar_{\ubar}$:
		\begin{subequations}\label{DE:Energy2}
			\begin{align}
				\Energy_{(i)}(u,\ubar_0,\ubar_1):=&\sum\limits_{\varphi\in\{\pfour^{i-1}(\logDensity,v^i,\Ent)\}}\Energy_{u}^{(1)}[\varphi](\ubar_0,\ubar_1)\\
				\notag&+\sum\limits_{\varphi\in\{\pfour^{i-1}(\vortrenormalized^i,\gradEnt^i,\modivort^i,\DivGradEnt)\}}\Energy_{u}^{(2)}[\varphi](\ubar_0,\ubar_1),\\
				\uEnergy_{(i)}(\ubar,u_0,u_1):=&\sum\limits_{\varphi\in\{\pfour^{i-1}(\logDensity,v^i,\Ent)\}}\uEnergy_{\ubar}^{(1)}[\varphi](u_0,u_1)\\
				\notag&+\sum\limits_{\varphi\in\{\pfour^{i-1}(\vortrenormalized^i,\gradEnt^i,\modivort^i,\DivGradEnt)\}}\uEnergy_{\ubar}^{(2)}[\varphi](u_0,u_1),\\
				\tilde{\Energy}_{(k)}(u,\ubar_0,\ubar_1):=&\sum\limits_{0\leq i\leq k}\Energy_{(i)}(u,\ubar_0,\ubar_1),\\
				\tilde{\uEnergy}_{(k)}(\ubar,u_0,u_1):=&\sum\limits_{0\leq i\leq k}\uEnergy_{(i)}(\ubar,u_0,u_1),
			\end{align}
		\end{subequations}
		where, for convenience, we define $\Energy_{(i)}(u,\ubar^\prime_0,\ubar^\prime_1)=0$ if $\ubar^\prime_0\geq\ubar^\prime_1$, and $\uEnergy_{(i)}(\ubar,u^\prime_0,u^\prime_1)=0$ if $u^\prime_0\geq u^\prime_1$.
		
		Furthermore, we define fluid energy quantities on the space-time regions $\region_1$, $\region_2$, and $\region_3$ (defined in Section \ref{SS:SUBSETSOFSPACETIME}) as follows:
		\begin{subequations}\label{DE:Energy3}
			\begin{align}
				\KEnergy^{(1)}_{(i)}(u,\ubar,\uptau):=&\sum\limits_{V\in\{\pfour^{i-1}(\vVort,\vgradientEnt)\}}\KEnergy^{(1)}[V](u,\ubar,\uptau),\\
				\tilde{\KEnergy}^{(1)}_{(k)}(u,\ubar,\uptau):=&\sum\limits_{0\leq i\leq k}\KEnergy^{(1)}_{(i)}(u,\ubar,\uptau),\\
				\KEnergy^{(2)}_{(i)}(\uptau_1,\uptau_2):=&\sum\limits_{V\in\{\pfour^{i-1}(\vVort,\vgradientEnt)\}}\KEnergy^{(2)}[V](\uptau_1,\uptau_2),\\
				\tilde{\KEnergy}^{(3)}_{(k)}(\uptau_1,\uptau_2):=&\sum\limits_{0\leq i\leq k}\KEnergy^{(3)}_{(i)}(\uptau_1,\uptau_2),\\
				\KEnergy^{(3)}_{(i)}(\uptau_1,\uptau_2):=&\sum\limits_{V\in\{\pfour^{i-1}(\vVort,\vgradientEnt)\}}\KEnergy^{(3)}[V](\uptau_1,\uptau_2),\\
				\tilde{\KEnergy}^{(3)}_{(k)}(\uptau_1,\uptau_2):=&\sum\limits_{0\leq i\leq k}\KEnergy^{(3)}_{(i)}(\uptau_1,\uptau_2).
			\end{align}
		\end{subequations}
		
	\end{definition}
	\begin{remark}
		In Definition \ref{D:ControllingEnergyQuantities}, \eqref{DE:Energy1}-\eqref{DE:Energy3} will be used in the energy estimates in the general case, where one suffers a loss of derivatives in the transversal direction for the initial data. In this case, we only commute the compressible Euler equations \eqref{E:WaveEQ}-\eqref{E:divcurlsystem} with the standard Cartesian spacetime derivative $\pfour$. 
	\end{remark}

	\subsection{Standard Energy Inequalities}\label{SS:energyinequalities}
	\begin{lemma}[Basic energy inequality for the wave equations in ``triangle'' region]\label{basicenergy}
		Let $u_1\geq u\geq u_0\geq 0$ and $\ubar_1\geq\ubar\geq\ubar_0\geq 0$. Let $\updelta_u:=u_1-u_0$, and $\updelta_{\ubar}:=\ubar_1-\ubar_0$. Let $u+\ubar\leq\uptau\leq u+\ubar+\min\{\updelta_u,\updelta_{\ubar}\}$. Let $\varphi$ be smooth on $\region^{u,\ubar,\uptau}_1$ (defined in (\ref{DE:triangleM})). Then the following inequality holds: 
		\begin{align}
			\label{E:energyineq2}\bEnergy_{u,\ubar}^{(1)}[\varphi](\uptau)\leq &\Energy^{(1)}_{u}[\varphi](\ubar,\uptau-u)+\uEnergy^{(1)}_{\ubar}[\varphi](u,\uptau-\ubar)\\
			\notag&
			+\int_{\region^{u,\ubar,\uptau_1}_1}\left((\pfour\varphi)^2+\varphi^2\right)\cdot\deform{\frac{1}{2}(e_3+e_4)}+\varphi\cdot\pfour\varphi\diff\gvol\\
			\notag&+\int_{\region^{u,\ubar,\uptau}_1}\left(\boxg\varphi\right)\left(\Timelike\varphi\right)\diff\gvol.
		\end{align}
	\end{lemma}
	
	\begin{proof}[Proof of Lemma \ref{basicenergy}]
		We apply the divergence theorem for $\Jen{\Timelike}[\varphi]$ on the space-time region $\region_1^{u,\ubar,\uptau}$ relative to the volume form $\diff\gvol$. By (\ref{D:deformtensor}), (\ref{E:DalphaJalpha}), (\ref{D:Energy1}), with $\mathbf{X}:=\Timelike$, we have:
		\begin{align}
			\bEnergy_{u,\ubar}^{(1)}[\varphi](\uptau)\leq &\Energy^{(1)}_{u}[\varphi](\ubar,\uptau-u)+\uEnergy^{(1)}_{\ubar}[\varphi](u,\uptau-\ubar)\\
			\notag&
			+\int_{\region^{u,\ubar,\uptau}_1}\left(\boxg\varphi(\Timelike\varphi)+\frac{1}{2}Q^{\mu\nu}[\varphi]\deform{\Timelike}_{\mu\nu}-2\varphi\Timelike\varphi-\frac{1}{2}\varphi^2(\gfour^{-1})^{\mu\nu}\deform{\Timelike}_{\mu\nu}\right)\diff\gvol,
		\end{align}
		which, by rearranging terms, is exactly (\ref{E:energyineq2}).
	\end{proof}
	
	\begin{lemma}[Basic energy inequality for the wave equations in the ``parallelogram'' region]\label{LE:basicenergy2}
		Let $u_1\geq u\geq u_0\geq 0$ and $\ubar_1\geq\ubar\geq\ubar_0\geq 0$. Let $\updelta_u:=u_1-u_0$, and $\updelta_{\ubar}:=\ubar_1-\ubar_0$. For $u_0+\ubar_0+\min\{\updelta_u,\updelta_{\ubar}\}\leq\uptau_1\leq\uptau_2\leq\max\{u_0+\ubar_1,u_1+\ubar_0\}$, Let $\varphi$ be smooth on $\region^{\uptau_1,\uptau_2}_2$ (defined in (\ref{DE:parallelogramM})). Then the following inequality holds: \\
		
		\textbf{When $\updelta_u\leq\updelta_{\ubar}$}:
		\begin{align}
			\label{EQ:energyineq3}\Energy^{(1)}_{u_1}[\varphi](\uptau_1-u_1,\uptau_2-u_2)+\bEnergy_{u_0,\ubar_0}^{(1)}[\varphi](\uptau_2)\leq &\Energy^{(1)}_{u_0}[\varphi](\uptau_1-u_0,\uptau_2-u_0)+\bEnergy_{u_0,\ubar_0}^{(1)}[\varphi](\uptau_1)\\
			\notag&
			+\int_{\region^{\uptau_1,\uptau_2}_2}\left((\pfour\varphi)^2+\varphi^2\right)\cdot\deform{\frac{1}{2}(e_3+e_4)}+\varphi\cdot\pfour\varphi\diff\gvol\\
			\notag&+\int_{\region^{\uptau_1,\uptau_2}_2}\left(\boxg\varphi\right)\left(\Timelike\varphi\right)\diff\gvol.
		\end{align}
		
		\textbf{When $\updelta_u>\updelta_{\ubar}$}:
		\begin{align}
			\label{EQ:energyineq4}\uEnergy^{(1)}_{\ubar_1}[\varphi](\uptau_1-\ubar_1,\uptau_2-\ubar_2)+\bEnergy_{u_0,\ubar_0}^{(1)}[\varphi](\uptau_2)\leq &\uEnergy^{(1)}_{\ubar_0}[\varphi](\uptau_1-\ubar_1,\uptau_2-\ubar_2)+\bEnergy_{u_0,\ubar_0}^{(1)}[\varphi](\uptau_1)\\
			\notag&
			+\int_{\region^{\uptau_1,\uptau_2}_2}\left((\pfour\varphi)^2+\varphi^2\right)\cdot\deform{\frac{1}{2}(e_3+e_4)}+\varphi\cdot\pfour\varphi\diff\gvol\\
			\notag&+\int_{\region^{\uptau_1,\uptau_2}_2}\left(\boxg\varphi\right)\left(\Timelike\varphi\right)\diff\gvol.
		\end{align}
	\end{lemma}
	
	\begin{lemma}[Basic energy inequality for the wave equations in the ``trapezoid'' region]\label{LE:basicenergy3.5}
		Let $u_1\geq u\geq u_0\geq 0$ and $\ubar_1\geq\ubar\geq\ubar_0\geq 0$. Let $\updelta_u:=u_1-u_0$, and $\updelta_{\ubar}:=\ubar_1-\ubar_0$. Let $\max\{\updelta_u,\updelta_{\ubar}\}\leq\uptau_1\leq\uptau_2\leq u_1+\ubar_1$. Let $\varphi$ be smooth on $\region^{\uptau_1,\uptau_2}_3$ (defined in (\ref{DE:trapezoidM})). Then the following inequality holds: 
		\begin{align}
			\label{EQ:energyineq1}&\bEnergy_{u,\ubar}^{(1)}[\varphi](\uptau_2)+\Energy^{(1)}_{u_1}[\varphi](\uptau_1-u_1,\uptau_2-u_1)+\uEnergy^{(1)}_{\ubar_1}[\varphi](\uptau_1-\ubar_1,\uptau_2-\ubar_1)\\
			\notag\leq &\bEnergy_{u,\ubar}^{(1)}[\varphi](\uptau_1)+
			\int_{\region^{\uptau_1,\uptau_2}_3}\left((\pfour\varphi)^2+\varphi^2\right)\cdot\deform{\frac{1}{2}(e_3+e_4)}+\varphi\cdot\pfour\varphi\diff\gvol\\
			\notag&+\int_{\region^{\uptau_1,\uptau_2}_3}\left(\boxg\varphi\right)\left(\Timelike\varphi\right)\diff\gvol.
		\end{align}
	\end{lemma}
	\begin{proof}[Proof of Lemma \ref{LE:basicenergy2} and Lemma \ref{LE:basicenergy3.5}]
		The proofs of Lemma \ref{LE:basicenergy2} and Lemma \ref{LE:basicenergy3.5} are similar to the proof of Lemma \ref{basicenergy} except that we apply the divergence theorem for $\Jen{\Timelike}[\varphi]$ on the space-time regions $\region_2^{\uptau_1,\uptau_2}$ and $\region^{\uptau_1,\uptau_2}_3$ respectively.
	\end{proof}
	\begin{lemma}[Basic energy inequality for the transport equations in the ``triangle'' region]\label{LE:basicenergy4}
		Let $u_1\geq u\geq u_0\geq 0$ and $\ubar_1\geq\ubar\geq\ubar_0\geq 0$. Let $\updelta_u:=u_1-u_0$, and $\updelta_{\ubar}:=\ubar_1-\ubar_0$. Let $u+\ubar\leq\uptau\leq u+\ubar+\min\{\updelta_u,\updelta_{\ubar}\}$. Let $\varphi$ be smooth on $\region^{u,\ubar,\uptau}_1$ (defined in (\ref{DE:triangleM})). Then, under the bootstrap assumptions (\ref{BAfluid})-\eqref{BAgeometry3}, the following inequality holds: 
		\begin{align}
			\label{EQ:energyineq4.1}\bEnergy_{u,\ubar}^{(2)}[\varphi](\uptau)\leq &\Energy^{(2)}_{u}[\varphi](\ubar,\uptau-u)+\uEnergy^{(2)}_{\ubar}[\varphi](u,\uptau-\ubar)\\
			\notag&
			+\int_{u+\ubar}^{\uptau}\norm{\pfour\vvariables}_{L^\infty(\region^{u,\ubar,\uptau}_{1})}\bEnergy_{u,\ubar}^{(2)}[\varphi](\uptau^\prime)\diff\uptau^\prime\\
			\notag&+\int_{u+\ubar}^{\uptau}\left(\bEnergy_{u,\ubar}^{(2)}[\varphi](\uptau^\prime)\cdot\bEnergy_{u,\ubar}^{(2)}[\materialderivative\varphi](\uptau^\prime)\right)^{\frac{1}{2}}\diff\uptau^\prime.
		\end{align}
	\end{lemma}
	\begin{proof}[Proof of Lemma \ref{LE:basicenergy4}]
		Let $\mathbf{J}^\alpha:=\varphi^2\materialderivative^\alpha$; then $\p_\alpha \mathbf{J}^\alpha=2\varphi\materialderivative\varphi+\left(\p_\alpha\materialderivative^\alpha\right)\varphi^2$. We apply the divergence theorem for $\mathbf{J}^\alpha$ on the space-time region $\region_1^{u,\ubar,\uptau}$ relative to the Cartesian coordinates. Note that $\gfour(\mathbf{J},e_4)=2\varphi^2\lapse$. By the Cauchy-Schwarz inequality and bootstrap assumptions (\ref{BAacoustic}), we obtain the desired estimates.
	\end{proof}
	\begin{lemma}[Basic energy inequality for the transport equations in the ``parallelogram'' region]\label{LE:basicenergy5}
		Let $u_1\geq u\geq u_0\geq 0$ and $\ubar_1\geq\ubar\geq\ubar_0\geq 0$. Let $\updelta_u:=u_1-u_0$, and $\updelta_{\ubar}:=\ubar_1-\ubar_0$. For $u_0+\ubar_0+\min\{\updelta_u,\updelta_{\ubar}\}\leq\uptau_1\leq\uptau_2\leq\max\{u_0+\ubar_1,u_1+\ubar_0\}$, let $\varphi$ be smooth on $\region^{\uptau_1,\uptau_2}_2$ (defined in (\ref{DE:parallelogramM})). $\Energy,\bEnergy$ are defined in Section \ref{SS:energies}. Then the following inequality holds: \\
		
		\textbf{When $\updelta_u\leq\updelta_{\ubar}$}:
		\begin{align}
			\label{EQ:energyineq5}\Energy^{(2)}_{u_1}[\varphi](\uptau_1-u_1,\uptau_2-u_2)+\bEnergy_{u_0,\ubar_0}^{(2)}[\varphi](\uptau_2)\leq &\Energy^{(2)}_{u_0}[\varphi](\uptau_1-u_0,\uptau_2-u_0)+\bEnergy_{u_0,\ubar_0}^{(2)}[\varphi](\uptau_1)\\
			\notag&+\int_{\uptau_1}^{\uptau_2}\norm{\pfour\vvariables}_{L^\infty(\region^{\uptau_1,\uptau_2}_2)}\bEnergy_{u_0,\ubar_0}^{(2)}[\varphi](\uptau^\prime)\diff\uptau^\prime\\
			\notag&+\int_{\uptau_1}^{\uptau_2}\left(\bEnergy_{u_0,\ubar_0}^{(2)}[\varphi](\uptau^\prime)\cdot\bEnergy_{u_0,\ubar_0}^{(2)}[\materialderivative\varphi](\uptau^\prime)\right)^{\frac{1}{2}}\diff\uptau^\prime.
		\end{align}
		
		\textbf{When $\updelta_u>\updelta_{\ubar}$}:
		\begin{align}
			\label{EQ:energyineq6}\uEnergy^{(2)}_{\ubar_1}[\varphi](\uptau_1-\ubar_1,\uptau_2-\ubar_2)+\bEnergy_{u_0,\ubar_0}^{(2)}[\varphi](\uptau_2)\leq &\uEnergy^{(2)}_{\ubar_0}[\varphi](\uptau_1-\ubar_1,\uptau_2-\ubar_2)+\bEnergy_{u_0,\ubar_0}^{(2)}[\varphi](\uptau_1)\\
			\notag&+\int_{\uptau_1}^{\uptau_2}\norm{\pfour\vvariables}_{L^\infty(\region^{\uptau_1,\uptau_2}_2)}\bEnergy_{u,\ubar}^{(2)}[\varphi](\uptau^\prime)\diff\uptau^\prime\\
			\notag&+\int_{\uptau_1}^{\uptau_2}\left(\bEnergy_{u,\ubar}^{(2)}[\varphi](\uptau^\prime)\cdot\bEnergy_{u,\ubar}^{(2)}[\materialderivative\varphi](\uptau^\prime)\right)^{\frac{1}{2}}\diff\uptau^\prime.
		\end{align}
	\end{lemma}
	\begin{lemma}[Basic energy inequality for the transport equations in the ``trapezoid'' region]\label{LE:basicenergy5.5}
		Let $u_1\geq u\geq u_0\geq 0$ and $\ubar_1\geq\ubar\geq\ubar_0\geq 0$. Let $\updelta_u:=u_1-u_0$, and $\updelta_{\ubar}:=\ubar_1-\ubar_0$. Let $\max\{u_0+\ubar_1,u_1+\ubar_0\}\leq\uptau_1\leq\uptau_2\leq u_1+\ubar_1$. Let $\varphi$ be smooth on $\region^{\uptau_1,\uptau_2}_3$ (defined in (\ref{DE:trapezoidM})). Then the following inequality holds: 
		\begin{align}
			\label{ES:energyineq1}&\bEnergy_{u_0,\ubar_0}^{(2)}[\varphi](\uptau_2)+\Energy^{(2)}_{u_1}[\varphi](\uptau_1-u_1,\uptau_2-u_1)+\uEnergy^{(2)}_{\ubar_1}[\varphi](\uptau_1-\ubar_1,\uptau_2-\ubar_1)\\
			\notag\leq&\bEnergy_{u_0,\ubar_0}^{(2)}[\varphi](\uptau_1)+\int_{\uptau_1}^{\uptau_2}\norm{\pfour\vvariables}_{L^\infty(\region^{\uptau_1,\uptau_2}_3)}\bEnergy_{u_0,\ubar_0}^{(2)}[\varphi](\uptau^\prime)\diff\uptau^\prime\\
			\notag&+\int_{\uptau_1}^{\uptau_2}\left(\bEnergy_{u_0,\ubar_0}^{(2)}[\varphi](\uptau^\prime)\cdot\bEnergy_{u_0,\ubar_0}^{(2)}[\materialderivative\varphi](\uptau^\prime)\right)^{\frac{1}{2}}\diff\uptau^\prime.
		\end{align}
	\end{lemma}
	\begin{proof}[Proof of Lemma \ref{LE:basicenergy5} and Lemma \ref{LE:basicenergy5.5}] 
		The proofs of Lemma \ref{LE:basicenergy5} and Lemma \ref{LE:basicenergy5.5} are similar to the proof of Lemma \ref{LE:basicenergy4} except that we apply the divergence theorem for $\mathbf{J}^\alpha$ (defined in the proof of Lemma \ref{LE:basicenergy4}) on the space-time regions $\region_2^{\uptau_1,\uptau_2}$ and $\region^{\uptau_1,\uptau_2}_3$ relative to the Cartesian coordinates respectively.
	\end{proof}
	
	\subsection{Space-time Elliptic Estimates for Vorticity and Entropy Gradient}\label{SS:elliptic}
	Inspired by \cite{IntegralIdentities}, in this section, we introduce a space-time version of elliptic estimates for vorticity $\vortrenormalized$ and entropy gradient $\gradEnt$. 
	\begin{definition}\cite[Definition 4.5. Elliptic-hyperbolic current]{IntegralIdentities}.\label{D:elliptichyperboliccurrent}
		Let $V$ be a $\St$-tangent vectorfield; we define $\tilde{\mathbf{J}}[V]$ to be the following $\Stau$-tangent tensorfield\footnote{We note that $\tilde{\mathbf{J}}[V]$ is a $\Stau$-tangent tensorfield, while $V$ is a $\St$-tangent vectorfield.}, where $\Sigmatauproject$ is the projection defined in Definition \ref{D:projectionoperators}:
		\begin{align}
			\tilde{\mathbf{J}}^\alpha[V]:=V^\gamma\Sigmatauproject^\lambda_\gamma\Sigmatauproject^\alpha_\kappa\p_\lambda V^\kappa-V^\gamma\Sigmatauproject^\alpha_\gamma\Sigmatauproject^\kappa_\lambda\p_\kappa V^\lambda=V^\gamma\Sigmatauproject^\alpha_\kappa\taup_\gamma V^\kappa-V^\gamma\Sigmatauproject^\alpha_\gamma\taup_\lambda V^\lambda.
		\end{align}
	\end{definition}
	
	\begin{lemma}\cite[Lemma 4.7. Divergence identity along $\Stau$]{IntegralIdentities}.\label{L:Divergenceidentity}
		Let $V$ be a $\St$-tangent vectorfield, and $\tilde{\mathbf{J}}[V]$ be as defined in Definition \ref{D:elliptichyperboliccurrent}. $K$ is defined in (\ref{Def:K}). $\Timelike$ is defined in (\ref{D:Timelike}). $\gtau$ is defined in (\ref{D:defg}). $\taup$ is defined in (\ref{D:taupangp}). $\tilde{\nabla}_\alpha\tilde{\mathbf{J}}^\alpha$ is the $\Stau$-projected divergence of $\Stau$-tangent vectorfields (see \cite[Lemma 3.19]{IntegralIdentities} for multiple expressions of $\tilde{\nabla}_\alpha\tilde{\mathbf{J}}^\alpha$). Then the following divergence identity holds:
		\begin{align}
			\label{E:IdentityinSigmatau1}
			\gtauabs{\pfour V}^2-(K_\alpha K V)^2=&\tilde{\nabla}_\alpha(\tilde{\mathbf{J}}^\alpha[V])\\
			\notag&+\mathfrak{J}_{(\widetilde{\text{Antisym}})}[\pfour V,\pfour V]+\mathfrak{J}_{(\widetilde{\text{Div}})}[\pfour V,\pfour V]\\
			\notag&+\mathfrak{J}_{(\widetilde{\text{Coeff}})}[\pfour V,\pfour V],
		\end{align}
		where 
		\begin{subequations}
			\begin{align}
				\label{E:JAnti}\mathfrak{J}_{(\widetilde{\text{Antisym}})}[\pfour V,\pfour V]:=&\frac{1}{2}\gtau^{\alpha\gamma}\gtau^{\beta\delta}\left(\p_\alpha V_\beta-\p_\beta V_\alpha\right)\left(\p_\gamma V_\delta-\p_\delta V_\gamma\right),\\
				\label{E:JDiv}\mathfrak{J}_{(\text{Div})}[\pfour V,\pfour V]:=&(\p_a V^a)^2+(2\lapse)^{-2}(K_\alpha\materialderivative V^\alpha)^2\\
				\notag&+2(K_\alpha KV^\alpha)\p_aV^a-\lapse^{-1}(K_\alpha\materialderivative V^\alpha)(K_\alpha K V^\alpha)-\lapse^{-1}(\p_aV^a)K_\alpha\materialderivative V^\alpha,\\
				\mathfrak{J}_{(\text{Coeff})}[V,\pfour V]:=&V^\alpha\gfour_{\beta\gamma}(\taup_\alpha\Timelike^\beta)\Timelike V^\gamma-V_\alpha(\taup_\beta\Timelike^\alpha)\Timelike V^\beta\\
				\notag&+V^\alpha\Timelike_\alpha(\taup_\beta\Timelike^\beta)\taup_\gamma V^\gamma-V^\alpha\Timelike_\alpha(\taup_\beta\Timelike^\gamma)\taup_\gamma V^\beta\\
				\notag&+V^\alpha\Timelike_\beta(\taup_\alpha\Timelike^\gamma)\taup_\gamma V^\beta-V^\alpha\Timelike_\beta(\taup_\gamma\Timelike^\gamma)\taup_\alpha V^\beta\\
				\notag&+V^\alpha\Timelike^\beta(\taup_\alpha\gfour_{\beta\gamma})\Timelike V^\gamma-V^\alpha\Timelike^\beta(\taup_\gamma\gfour_{\alpha\beta})\Timelike V^\gamma\\
				\notag&+\frac{1}{2}V^\alpha\Timelike_\beta\Timelike^\gamma\Timelike^\delta(\taup_\alpha\gfour_{\gamma\delta})\Timelike V^\beta-\frac{1}{2}V^\alpha\Timelike_\alpha\Timelike^\beta\Timelike^\gamma(\taup_\delta\gfour_{\beta\gamma})\Timelike V^\delta\\
				\notag&+2V^\alpha(\taup_\beta\gfour_{\alpha\beta})\taup^\gamma V^\beta-2\Sigmatauproject^\alpha_\beta V^\gamma(\taup^\gamma\gfour_{\alpha\gamma})\taup_\delta V^\beta\\
				\notag&+\frac{1}{2}\gtau^{\alpha\beta}V^\gamma(\taup_\gamma\gfour_{\alpha\beta})\taup_\delta V^\delta-\frac{1}{2}\gtau^{\alpha\beta}V^\gamma(\taup_\delta\gfour_{\alpha\beta})\taup_\gamma V^\delta\\
				\notag&+V^\alpha V^\beta(\taup^\gamma\gfour_{\alpha\delta})(\taup^\delta\gfour_{\beta\gamma})-V^\alpha V^\beta(\gfour^{-1})^{\gamma\delta}(\taup^\kappa\gfour_{\alpha\gamma})(\taup_\kappa\gfour_{\beta\delta}).
			\end{align}
		\end{subequations}
	\end{lemma}
	\begin{lemma}\cite[Lemma 4.8. Preliminaries on boundary terms]{IntegralIdentities}.\label{L:preboundaryterms}
		Let $V$ be a $\St$-tangent vectorfield, and $\tilde{\mathbf{J}}[V]$ be as defined in Definition \ref{D:elliptichyperboliccurrent}. $\angV:=\sphereproject V$ is a $\suubar$-projected vectorfield, where $\sphereproject$ is defined in Definition \ref{D:projectionoperators}. $\lapse,\uinvfolia,\Timelike,\spherenormal,\Lrenormal,\Theta$ are defined in Section \ref{SS:SUBSETSOFSPACETIME}. Then the following divergence identity holds for $\spherenormal_\alpha\tilde{\mathbf{J}}^\alpha[V]$:
		\begin{subequations}
			\begin{align}
				\Timelike_\alpha\tilde{\mathbf{J}}^\alpha[V]=&0,\\
				\label{E:NJ1}\spherenormal_\alpha\tilde{\mathbf{J}}^\alpha[V]=&-\Lrenormal\left(\frac{1}{\lapse\uinvfolia}\sgabs{V}^2\right)+\frac{1}{\lapse}(\Lunit^\alpha+\Theta^\alpha)\angV^\beta(\p_\alpha V_\beta-\p_\beta V_\alpha)+\Lrenormal\left\{\frac{1}{\lapse\uinvfolia}\gsphere^{\alpha\beta}\right\}V_\alpha V_\beta\\
				\notag&+\left\{\Theta\left(\frac{1}{2\lapse}\gsphere^{\alpha\beta}\right)\right\}V_\alpha V_\beta-2V_\alpha(\angV\materialderivative^\alpha)+\frac{1}{2\lapse}\sgabs{\angV}^2\angdiv\Theta+\spherenormal_\alpha V^\alpha V^\beta\angdiv\angp_\beta+V_\alpha\angV \spherenormal^\alpha\\
				\notag&+V^\alpha\spherenormal^\beta\angV\gfour_{\alpha\beta}-\angdiv\left(\frac{1}{2\lapse}\sgabs{\angV}^2\Theta\right)-\angdiv\left(\spherenormal_\alpha V^\alpha\angV\right),\\
				\label{E:NJ2}\spherenormal_\alpha\tilde{\mathbf{J}}^\alpha[V]=&-\uLrenormal\left(\frac{1}{\lapse\invfolia}\sgabs{V}^2\right)+\frac{1}{\lapse}(\uLunit^\alpha+\Theta^\alpha)\angV^\beta(\p_\alpha V_\beta-\p_\beta V_\alpha)+\uLrenormal\left\{\frac{1}{\lapse\invfolia}\gsphere^{\alpha\beta}\right\}V_\alpha V_\beta\\
				\notag&+\left\{\Theta\left(\frac{1}{2\lapse}\gsphere^{\alpha\beta}\right)\right\}V_\alpha V_\beta-2V_\alpha(\angV\materialderivative^\alpha)+\frac{1}{2\lapse}\sgabs{\angV}^2\angdiv\Theta+\spherenormal_\alpha V^\alpha V^\beta\angdiv\angp_\beta+V_\alpha\angV \spherenormal^\alpha\\
				\notag&+V^\alpha\spherenormal^\beta\angV\gfour_{\alpha\beta}-\angdiv\left(\frac{1}{2\lapse}\sgabs{\angV}^2\Theta\right)-\angdiv\left(\spherenormal_\alpha V^\alpha\angV\right)
			\end{align}
		\end{subequations}
	\end{lemma}
	\begin{lemma}\cite[Theorem 10.6. Hodge-type integral identities for $\vortrenormalized$ and $\gradEnt$]{IntegralIdentities}.\label{L:integralidentitiesforS}
		Let $u_1\geq u\geq u_0\geq 0$ and $\ubar_1\geq\ubar\geq\ubar_0\geq 0$. Let $\updelta_u:=u_1-u_0$, and $\updelta_{\ubar}:=\ubar_1-\ubar_0$. Let $u+\ubar\leq\uptau\leq u+\ubar+\min\{\updelta_u,\updelta_{\ubar}\}$. $\min\{u_0+\ubar_1,u_1+\ubar_0\}\leq\uptau_1\leq\uptau_2\leq u_1+\ubar_1$. $\region^{\uptau_1,\uptau_2}_2,\region^{\uptau_1,\uptau_2}_3$ are defined in (\ref{DE:parallelogramM}) and (\ref{DE:trapezoidM}). $\region^{u,\ubar,\uptau}_1$ is defined in (\ref{DE:triangleM}).  $\lapse,\uinvfolia,\invfolia,K,\Lrenormal,\uLrenormal,\Theta$ are defined in Section \ref{SS:SUBSETSOFSPACETIME}. $\angV:=\sphereproject V$ is $\suubar$-projected vectorfield, where $\sphereproject$ is defined in Definition \ref{D:projectionoperators}. For $v^i,\logDensity,\Ent$ the solution of the compressible Euler equations (\ref{E:Brho})-(\ref{E:Bs}). $V\in\{\vortrenormalized,\gradEnt\}$ satisfies the following integral identities in $\region^{u,\ubar,\uptau}_1$ , $\region^{\uptau_1,\uptau_2}_2$, and $\region^{\uptau_1,\uptau_2}_3$:\\
		
		\noindent\textbf{On $\region^{u,\ubar,\uptau}_1$:}
		\begin{align}\label{E:integralidentities1}
			&\int_{\region^{u,\ubar,\uptau}_1}(\lapse^2\invfolia\uinvfolia)^{-1/2}\quadsmoothfunction(\pfour V,\pfour V)\diff\gvol+\int_{S_{u,\uptau-u}}\frac{1}{\lapse(\invfolia+\uinvfolia)}\sgabs{\angV}^2\diff\spherevol+\int_{S_{\uptau-\ubar,\ubar}}\frac{1}{\lapse(\invfolia+\uinvfolia)}\sgabs{\angV}^2\diff\spherevol\\
			\notag=&\int_{S_{u_0,\ubar_0}}\frac{2}{\lapse(\invfolia+\uinvfolia)}\sgabs{\angV}^2\diff\spherevol\\
			\notag&-\int_{\hypersur_{u}^{\ubar,\uptau-u}}\mathfrak{H}[ V]+\mathfrak{H}_{(1)}[ V]\diff\gvol-\int_{\hypersur_{\ubar}^{u,\uptau-\ubar}}\underline{\mathfrak{H}}[ V]+\underline{\mathfrak{H}}_{(1)}[V]\diff\gvol\\
			\notag&+\int_{\region^{u,\ubar,\uptau}_1}(\lapse^2\invfolia\uinvfolia)^{-1/2}\left\{\frac{1}{2}\gtauabs{\mathfrak{A}^{( V)}}^2+\gtauabs{\mathfrak{B}_{( V)}}^2+\mathfrak{C}^{( V)}+\mathfrak{D}^{( V)}+\mathfrak{J}_{(\text{Coeff})}[ V,\pfour V]\right\}\diff\gvol.
		\end{align}
		
		\noindent\textbf{On $\region^{\uptau_1,\uptau_2}_2$:}
		
		\textbf{When $\updelta_u\leq\updelta_{\ubar}$:}
		\begin{align}
			\label{E:integralidentities2}
			&\int_{\region^{\uptau_1,\uptau_2}_2}(\lapse^2\invfolia\uinvfolia)^{-1/2}\quadsmoothfunction(\pfour V,\pfour V)\diff\gvol+\int_{S_{u_1,\uptau_2-u_1}}\frac{1}{\lapse(\invfolia+\uinvfolia)}\sgabs{\angV}^2\diff\spherevol\\
			\notag=&\int_{S_{u_0,\uptau_2-u_0}}\frac{1}{\lapse(\invfolia+\uinvfolia)}\sgabs{\angV}^2\diff\spherevol-\int_{S_{u_0,\uptau_1-u_0}}\frac{2}{\lapse(\invfolia+\uinvfolia)}\sgabs{\angV}^2\diff\spherevol+\int_{S_{u_1,\uptau_1-u_1}}\frac{2}{\lapse(\invfolia+\uinvfolia)}\sgabs{\angV}^2\diff\spherevol\\
			\notag&-\int_{\hypersur_{u_0}^{\uptau_1-u_0,\uptau_2-u_0}}\mathfrak{H}[ V]+\mathfrak{H}_{(1)}[ V]\diff\gvol+\int_{\hypersur_{u_1}^{\uptau_1-u_1,\uptau_2-u_1}}\mathfrak{H}[ V]+\mathfrak{H}_{(1)}[ V]\diff\gvol\\
			\notag&+\int_{\region^{\uptau_1,\uptau_2}_2}(\lapse^2\invfolia\uinvfolia)^{-1/2}\left\{\frac{1}{2}\gtauabs{\mathfrak{A}^{( V)}}^2+\gtauabs{\mathfrak{B}_{( V)}}^2+\mathfrak{C}^{( V)}+\mathfrak{D}^{( V)}+\mathfrak{J}_{(\text{Coeff})}[ V,\pfour V]\right\}\diff\gvol.
		\end{align}
		
		\textbf{When $\updelta_u>\updelta_{\ubar}$:}
		\begin{align}
			\label{E:integralidentities3}
			&\int_{\region^{\uptau_1,\uptau_2}_2}(\lapse^2\invfolia\uinvfolia)^{-1/2}\quadsmoothfunction(\pfour V,\pfour V)\diff\gvol+\int_{S_{\uptau_2-\ubar_1,\ubar_1}}\frac{1}{\lapse(\invfolia+\uinvfolia)}\sgabs{\angV}^2\diff\spherevol\\
			\notag=&\int_{S_{\uptau_2-\ubar_0,\ubar_0}}\frac{1}{\lapse(\invfolia+\uinvfolia)}\sgabs{\angV}^2\diff\spherevol-\int_{S_{\uptau_1-\ubar_0,\ubar_0}}\frac{2}{\lapse(\invfolia+\uinvfolia)}\sgabs{\angV}^2\diff\spherevol+\int_{S_{\uptau_1-\ubar_1,\ubar_1}}\frac{2}{\lapse(\invfolia+\uinvfolia)}\sgabs{\angV}^2\diff\spherevol\\
			\notag&-\int_{\hypersurbar_{\ubar_0}^{\uptau_1-\ubar_0,\uptau_2-\ubar_0}}\mathfrak{H}[ V]+\mathfrak{H}_{(1)}[ V]\diff\gvol+\int_{\hypersurbar_{\ubar_1}^{\uptau_1-\ubar_1,\uptau_2-\ubar_1}}\mathfrak{H}[ V]+\mathfrak{H}_{(1)}[ V]\diff\gvol\\
			\notag&+\int_{\region^{\uptau_1,\uptau_2}_2}(\lapse^2\invfolia\uinvfolia)^{-1/2}\left\{\frac{1}{2}\gtauabs{\mathfrak{A}^{( V)}}^2+\gtauabs{\mathfrak{B}_{( V)}}^2+\mathfrak{C}^{( V)}+\mathfrak{D}^{( V)}+\mathfrak{J}_{(\text{Coeff})}[ V,\pfour V]\right\}\diff\gvol,
		\end{align}
		
		\noindent\textbf{On $\region^{\uptau_1,\uptau_2}_3$:}
		\begin{align}\label{E:integralidentities4}
			&\int_{\region^{\uptau_1,\uptau_2}_3}(\lapse^2\invfolia\uinvfolia)^{-1/2}\quadsmoothfunction(\pfour V,\pfour V)\diff\gvol+\int_{S_{u_1,\uptau_2-u_1}}\frac{1}{\lapse(\invfolia+\uinvfolia)}\sgabs{\angV}^2\diff\spherevol+\int_{S_{\uptau_2-\ubar_1,\ubar_1}}\frac{1}{\lapse(\invfolia+\uinvfolia)}\sgabs{\angV}^2\diff\spherevol\\
			\notag=&\int_{S_{u_1,\uptau_1-u_1}}\frac{2}{\lapse(\invfolia+\uinvfolia)}\sgabs{\angV}^2\diff\spherevol+\int_{S_{\uptau_2-\ubar_1,\ubar_1}}\frac{2}{\lapse(\invfolia+\uinvfolia)}\sgabs{\angV}^2\diff\spherevol\\
			\notag&+\int_{\hypersur_{u_1}^{\uptau_1-u_1,\uptau_2-u_1}}\mathfrak{H}[ V]+\mathfrak{H}_{(1)}[ V]\diff\gvol+\int_{\hypersur_{\ubar_1}^{\uptau_1-\ubar_1,\uptau_2-\ubar_1}}\underline{\mathfrak{H}}[ V]+\underline{\mathfrak{H}}_{(1)}[V]\diff\gvol\\
			\notag&+\int_{\region^{\uptau_1,\uptau_2}_3}(\lapse^2\invfolia\uinvfolia)^{-1/2}\left\{\frac{1}{2}\gtauabs{\mathfrak{A}^{( V)}}^2+\gtauabs{\mathfrak{B}_{( V)}}^2+\mathfrak{C}^{( V)}+\mathfrak{D}^{( V)}+\mathfrak{J}_{(\text{Coeff})}[ V,\pfour V]\right\}\diff\gvol.
		\end{align}
		
		where 
		\begin{align}
			\quadsmoothfunction(\pfour V,\pfour V):=\gtauabs{\Sigmatauproject(\pfour V)}^2-(K_\alpha KV^\alpha)^2+\gtau_{\alpha\beta}(\materialderivative V^\alpha)(\materialderivative V^\beta).
		\end{align}
		\begin{subequations}
			\begin{align}
				\mathfrak{A}_{\alpha\beta}^{(\vortrenormalized)}:=&2(\p_\beta\ln \speed)\vortrenormalized_\alpha-2(\p_\alpha\ln \speed)\vortrenormalized_\beta+2\delta^0_\alpha\vortrenormalized_a\p_\beta v^a-2\delta^0_\beta\vortrenormalized_a\p_\alpha v^a\\
				\notag&-\speed^{-4}e^{-2\logDensity}\frac{\pressure_{;\Ent}}{\backgroundDensity}\antisymmetic_{\alpha\beta\gamma\delta}(\materialderivative v^\gamma)\gradEnt^\delta+\speed^{-4}e^{-2\logDensity}\frac{\pressure_{;\Ent}}{\backgroundDensity}\antisymmetic_{\alpha\beta\gamma\delta}\materialderivative^\gamma\left\{\gradEnt^\delta(\p_a v^a)-\gradEnt^a\p_a v^\delta\right\}\\
				\notag&+\speed^{-2}e^{\logDensity}\antisymmetic_{\alpha\beta\gamma\delta}\materialderivative^\gamma\modivort^{\delta},\\
				\mathfrak{A}_{\alpha\beta}^{(\gradEnt)}:=&2(\p_\beta\ln\speed)\gradEnt_\alpha-2(\p_\alpha\ln\speed)\gradEnt_\beta,\\
				\mathfrak{B}_{(\vortrenormalized)}^i:=&\vortrenormalized^a\p_av^i-e^{-2\logDensity}\speed^{-2}\frac{\pressure_{;\Ent}}{\backgroundDensity}\antisymmetic_{iab}(\materialderivative v^a)\gradEnt^b,\\
				\mathfrak{B}_{(\gradEnt)}^i:=&-\delta_{ab}(\p_iv^a)\gradEnt^b,\\
				\mathfrak{C}^{(\vortrenormalized)}:=&-2(K_aK\vortrenormalized^a)\vortrenormalized^b\p_b\logDensity-\lapse^{-1}(K_aK\vortrenormalized^a)K_b\mathfrak{B}_{(\vortrenormalized)}^b+\dive\FV,\\
				\mathfrak{C}^{(\gradEnt)}:=&2(K_aK\gradEnt^a)\left(e^{2\logDensity}\DivGradEnt+\gradEnt^b\p_b\logDensity\right)-\lapse^{-1}(K_aK\gradEnt^a)K_b\mathfrak{B}_{(\gradEnt)}^b,\\
				\label{E:mathfrakD}\mathfrak{D}^{(\vortrenormalized)}:=&\left\{(\vortrenormalized^a\p_a\logDensity+\dive\FV)+(2\lapse)^{-1}K_a\mathfrak{B}^a_{(\vortrenormalized)}\right\}^2,\\
				\mathfrak{D}^{(\gradEnt)}:=&\left\{e^{2\logDensity}\DivGradEnt+\gradEnt^a\p_a\logDensity+(2\lapse)^{-1}K_a\mathfrak{B}^a_{(\gradEnt)}\right\}^2.
			\end{align}
		\end{subequations}
		\begin{subequations}\label{D:integralquantities2}
			\begin{align}	
				\mathfrak{H}[V]:=&\frac{1}{2\lapse\uinvfolia}\sgabs{\angV}^2\gsphere^{\alpha\beta}\Lie_{\Lrenormal}\gsphere_{\alpha\beta}+\Lrenormal\left\{\frac{1}{\lapse\uinvfolia}\gsphere^{\alpha\beta}\right\}V_\alpha V_\beta+\left\{\Theta\left(\frac{1}{2\lapse}\gsphere^{\alpha\beta}\right)\right\}V_\alpha V_\beta\\
				\notag&-2V_\alpha(\angV\materialderivative^\alpha)+\frac{1}{2\lapse}\sgabs{\angV}^2\angdiv\Theta+\spherenormal_\alpha V^\alpha V^\beta\angdiv\angp_\beta+V_\alpha\angV \spherenormal^\alpha+V^\alpha\spherenormal^\beta\angV\gfour_{\alpha\beta},\\
				\mathfrak{H}_{(1)}[\vortrenormalized]:=&\frac{1}{\lapse}(\Lunit^\alpha+\Theta^\alpha)\angvort^\beta\left\{2(\p_\beta\ln\speed)\vortrenormalized_\alpha-2(\p_\alpha\ln\speed)\vortrenormalized_\beta+2\delta^0_\alpha\vortrenormalized_\alpha\p_\beta v^a-2\delta^0_\beta\vortrenormalized_a\p_\alpha v^a\right.\\
				\notag&\left.-\speed^{-4}e^{-2\logDensity}\frac{\pressure_{;\Ent}}{\backgroundDensity}\antisymmetic_{\alpha\beta\gamma\delta}(\materialderivative v^\gamma)\gradEnt^\delta+\speed^{-2}\antisymmetic_{\alpha\beta\gamma\delta}\materialderivative^\gamma(\curl\vortrenormalized)^\delta\right\},\\
				\mathfrak{H}_{(1)}[\gradEnt]:=&\frac{1}{\lapse}(\Lunit^\alpha+\Theta^\alpha)\anggradEnt^\beta\left\{2(\p_\beta\ln\speed)\gradEnt_\alpha-2(\p_\alpha\ln\speed)\gradEnt_\beta\right\},\\
				\underline{\mathfrak{H}}[V]:=&\frac{1}{2\lapse\uinvfolia}\sgabs{\angV}^2\gsphere^{\alpha\beta}\Lie_{\uLrenormal}\gsphere_{\alpha\beta}+\uLrenormal\left\{\frac{1}{\lapse\invfolia}\gsphere^{\alpha\beta}\right\}V_\alpha V_\beta+\left\{\Theta\left(\frac{1}{2\lapse}\gsphere^{\alpha\beta}\right)\right\}V_\alpha V_\beta\\
				\notag&-2V_\alpha(\angV\materialderivative^\alpha)+\frac{1}{2\lapse}\sgabs{\angV}^2\angdiv\Theta-\spherenormal_\alpha V^\alpha V^\beta\angdiv\angp_\beta-V_\alpha\angV \spherenormal^\alpha-V^\alpha\spherenormal^\beta\angV\gfour_{\alpha\beta},\\
				\underline{\mathfrak{H}}_{(1)}[\vortrenormalized]:=&\frac{1}{\lapse}(\uLunit^\alpha+\Theta^\alpha)\angvort^\beta\left\{2(\p_\beta\ln\speed)\vortrenormalized_\alpha-2(\p_\alpha\ln\speed)\vortrenormalized_\beta+2\delta^0_\alpha\vortrenormalized_\alpha\p_\beta v^a-2\delta^0_\beta\vortrenormalized_a\p_\alpha v^a\right.\\
				\notag&\left.-\speed^{-4}e^{-2\logDensity}\frac{\pressure_{;\Ent}}{\backgroundDensity}\antisymmetic_{\alpha\beta\gamma\delta}(\materialderivative v^\gamma)\gradEnt^\delta+\speed^{-2}\antisymmetic_{\alpha\beta\gamma\delta}\materialderivative^\gamma(\curl\vortrenormalized)^\delta\right\},\\
				\underline{\mathfrak{H}}_{(1)}[\gradEnt]:=&\frac{1}{\lapse}(\uLunit^\alpha+\Theta^\alpha)\anggradEnt^\beta\left\{2(\p_\beta\ln\speed)\gradEnt_\alpha-2(\p_\alpha\ln\speed)\gradEnt_\beta\right\}.
			\end{align}
		\end{subequations}
	\end{lemma}
	\begin{proof}[Discussion of the proof of Lemma \ref{L:integralidentitiesforS}]
		To prove (\ref{E:integralidentities1}), we substitute $V$ by $\vortrenormalized,\gradEnt$ in Lemma \ref{L:Divergenceidentity} respectively, integrate (\ref{E:IdentityinSigmatau1}) along $\Stau^{u,\ubar}$ (defined in (\ref{DE:Stau})), apply the divergence theorem, use Lemma \ref{L:preboundaryterms} for boundary terms, then integrate over $u+\ubar\leq\uptau\leq u+\ubar+\min\{\updelta_u,\updelta_{\ubar}\}$. We refer\footnote{In \cite{IntegralIdentities}, the authors integrate in a different region. However, the techniques and quantities are the same as here. That is, compared to integral identities in \cite{IntegralIdentities}, regions/boundaries are different in (\ref{E:integralidentities1})-(\ref{E:integralidentities3}) while the integrands are the same.} readers to \cite[Theorem 10.6]{IntegralIdentities} for the derivation of quantities\footnote{For $\mathfrak{H}_{(1)}[V]$, for convenience, we provide its expression $\mathfrak{H}_{(1)}[V]=\frac{1}{\lapse}(\uLunit^\alpha+\Theta^\alpha)\angvort^\beta\mathfrak{A}_{\alpha\beta}^{(V)}$ as in the second term on the RHS of \cite[(115)]{IntegralIdentities}, while its geometric property is important, as in Lemma \ref{L:Geometricproperty}} in (\ref{D:integralquantities2}), where the related quantities in \cite{IntegralIdentities} are $\mathcal{W}=1,\underline{\iota}=2\lapse^2\uinvfolia$, $\underline{z}=\frac{1}{2\lapse},\underline{h}=\frac{1}{4\lapse^2}$.
		
		The proofs of (\ref{E:integralidentities2})-(\ref{E:integralidentities4}) are similar to the proof of (\ref{E:integralidentities1}) except that we integrate (\ref{E:IdentityinSigmatau1}) along $\Stau^{u_0,\ubar_0}$ (defined in (\ref{DE:Stau})), apply the divergence theorem, then integrate over $\uptau_1\leq\uptau\leq\uptau_2$.
	\end{proof}

	
	\begin{lemma}\cite[Lemma 9.9. Coerciveness of the quadratic form $\quadsmoothfunction$]{IntegralIdentities}.\label{LE:integralID}
		Let $0\leq \ubar_0\leq\ubar_1$ and $0\leq u_0\leq u_1$. On the space-time regions $\region^{u,\ubar,\uptau}_1$ (defined in (\ref{DE:triangleM})), $\region^{\uptau_1,\uptau_2}_2$ (defined in (\ref{DE:parallelogramM})), and $\region^{\uptau_1,\uptau_2}_3$ (defined in (\ref{DE:trapezoidM})), the following estimates hold for the quadratic form $\quadsmoothfunction(\pfour V,\pfour V)$ and energy $\KEnergy$ defined in Definition \ref{DE:Energyinspacetime}:\\
		
		\noindent\textbf{On $\region^{u,\ubar,\uptau}_1$:}
		\begin{align}
			\int_{\region^{u,\ubar,\uptau}_1}\quadsmoothfunction(\pfour V,\pfour V)\diff\gvol+\int_{S_{u,\uptau-u}}\frac{1}{\lapse(\invfolia+\uinvfolia)}\sgabs{\angV}^2\diff\spherevol+\int_{S_{\uptau-\ubar,\ubar}}\frac{1}{\lapse(\invfolia+\uinvfolia)}\sgabs{\angV}^2\diff\spherevol\approx\KEnergy^{(1)}[V](u,\ubar,\uptau).
		\end{align}
		
		\noindent\textbf{On $\region^{\uptau_1,\uptau_2}_2$:}
		
		\textbf{When $\updelta_u\leq\updelta_{\ubar}$:}
		\begin{align}
			\int_{\region^{\uptau_1,\uptau_2}_2}\quadsmoothfunction(\pfour V,\pfour V)\diff\gvol+\int_{S_{u_1,\uptau_2-u_1}}\frac{1}{\lapse(\invfolia+\uinvfolia)}\sgabs{\angV}^2\diff\spherevol\approx\KEnergy^{(2)}[V](\uptau_1,\uptau_2).
		\end{align}
		
		\textbf{When $\updelta_u>\updelta_{\ubar}$:}
		\begin{align}
			\int_{\region^{\uptau_1,\uptau_2}_2}\quadsmoothfunction(\pfour V,\pfour V)\diff\gvol+\int_{S_{\uptau_2-\ubar_1,\ubar_1}}\frac{1}{\lapse(\invfolia+\uinvfolia)}\sgabs{\angV}^2\diff\spherevol\approx\KEnergy^{(2)}[V](\uptau_1,\uptau_2).
		\end{align}
		
		\noindent\textbf{On $\region^{\uptau_1,\uptau_2}_3$:}
		\begin{align}
			\int_{\region^{\uptau_1,\uptau_2}_3}\quadsmoothfunction(\pfour V,\pfour V)\diff\gvol+\int_{S_{\uptau_2-\ubar_1,\ubar_1}}\sgabs{\angV}^2\diff\spherevol+\int_{S_{u_1,\uptau_2-u_1}}\sgabs{\angV}^2\diff\spherevol\approx\KEnergy^{(3)}[V](\uptau_1,\uptau_2).
		\end{align}
	\end{lemma}
	\begin{lemma}\cite[Theorem 7.2. Geometric structure of the boundary terms]{IntegralIdentities}.\label{L:Geometricproperty}
		Let $v^i,\logDensity,\Ent$ be the solutions of the compressible Euler equations (\ref{E:Brho})-(\ref{E:Bs}). Let $\vec{e}=\{e_1,e_2,e_3,e_4\}$ be the array of the null frame defined in Section \ref{SS:SUBSETSOFSPACETIME}. Then, for $V\in\{\vortrenormalized,\gradEnt\}$, we have the following structure for $\mathfrak{H}[V]$, $\mathfrak{H}_{(1)}[V]$, $\underline{\mathfrak{H}}[V]$, and $\underline{\mathfrak{H}}_{(1)}[V]$, which are defined in (\ref{D:integralquantities2}):
		\begin{subequations}
			\begin{align}
				\mathfrak{H}[V],\mathfrak{H}_{(1)}[V]=\linsmoothfunction[\vvariables,\vVort,\vgradientEnt,\vec{e}]\left((e_4,\angp)(\vvariables,\vec{e}),\vortrenormalized,\gradEnt,\vectvort\right),\\
				\underline{\mathfrak{H}}[V],\underline{\mathfrak{H}}_{(1)}[V]=\linsmoothfunction[\vvariables,\vVort,\vgradientEnt,\vec{e}]\left((e_3,\angp)(\vvariables,\vec{e}),\vortrenormalized,\gradEnt,\vectvort\right).
			\end{align}
		\end{subequations}
	\end{lemma}

	\section{Base-Level Energy Estimates}\label{S:BASELEVELENERGYESTIMATES}
	In this section, we consider the geometric formulation of compressible Euler flow from Theorem~\ref{T:geometricformulation}. 
	We use it to derive energy estimates for solutions to these equations, but not yet for their higher derivatives. We refer to these as ``base-level'' energy estimates.
	Our estimates apply on regions of type $\region_{1},\region_{2},\region_{3}$,
	as defined in Section~\ref{SS:SUBSETS}.
	In Sect.\,\ref{S:APRIORIESTIMATESUPTOTOPORDER}, we will derive similar energy estimates for the solution's higher derivatives.

	\subsection{Base-Level Energy Estimates in the ``Triangle" Region $\region_{1}^{u_0,\ubar_0,\uptau}$}
	In this section, we derive base-level energy estimates for fluid variables in the ``triangular" region $\region_{1}$ defined as in Section \ref{SS:SUBSETS}.
	
	\begin{theorem}[Base-level energy estimates in a neighborhood of $S_{u_0,\ubar_0}$]\label{T:EnergyEstimates0}
		Let $\RBAinfty$ be as in Definition \ref{DE:LOT}. Let $\updelta_u:=u_1-u_0, \updelta_{\ubar}:=\ubar_1-\ubar_0$. $v^i,\logDensity,\Ent$ are the solutions of the compressible Euler equations (\ref{E:Brho})-(\ref{E:Bs}). Then, for various energies $\bEnergy,\KEnergy,\Energy,\uEnergy$ defined in Definition \ref{D:ControllingEnergyQuantities}, where $u_0+\ubar_0\leq\uptau\leq\min\{u_0+\ubar_1,u_1+\ubar_0\}$, the following energy estimates hold:
		\begin{subequations}\label{ES:CKEnergies1}
			\begin{align}
				\label{ES:KEnergy1}\KEnergy^{(1)}_{(1)}(u_0,\ubar_0,\uptau)\lesssim&\KEnergy^{(1)}_{(1)}(u_0,\ubar_0,\uptau)+\Energy_{(1)}(u_0,\ubar_0,\uptau-u_0)+\uEnergy_{(1)}(\ubar_0,u_0,\uptau-\ubar_0)\\
				\notag&+\left\{\KEnergy^{(1)}_{(1)}(u_0,\ubar_0,\uptau)\int_{u_0+\ubar_0}^{\uptau}\bEnergy_{(1)}(u_0,\ubar_0,\uptau^\prime)\diff \uptau^\prime\right\}^{1/2}\\
				\notag&+\int_{u_0+\ubar_0}^{\uptau}\left[\bEnergy_{(1)}(u_0,\ubar_0,\uptau^\prime)\right]\RBAinfty\diff \uptau^\prime,\\
				\label{ES:Energy1}\bEnergy_{(1)}(u_0,\ubar_0,\uptau)\lesssim&\Energy_{(1)}(u_0,\ubar_0,\uptau-u_0)+\uEnergy_{(1)}(\ubar_0,u_0,\uptau-\ubar_0)\\
				\notag&+\left\{\left[\KEnergy^{(1)}_{(1)}(u_0,\ubar_0,\uptau)\right]\int_{u_0+\ubar_0}^{\uptau}\bEnergy_{(1)}(u_0,\ubar_0,\uptau^\prime)\diff \uptau^\prime\right\}^{1/2}\\
				\notag&+\int_{u_0+\ubar_0}^{\uptau}\left[\bEnergy_{(1)}(u_0,\ubar_0,\uptau^\prime)\right]\RBAinfty\diff \uptau^\prime.
			\end{align} 
		\end{subequations}
		Moreover, by bootstrap assumptions (\ref{BAfluid})-(\ref{BAgeometry3}), we have:
		\begin{align}\label{ES:CKEnergies2}	&\KEnergy^{(1)}_{(1)}(u_0,\ubar_0,\uptau),\bEnergy_{(1)}(u_0,\ubar_0,\uptau)\\
			\notag\lesssim&\left\{\KEnergy^{(1)}_{(1)}(u_0,\ubar_0,u_0+\ubar_0)+\Energy_{(1)}(u_0,\ubar_0,\uptau-u_0)+\uEnergy_{(1)}(\ubar_0,u_0,\uptau-\ubar_0)\right\}\exp(\RBAinfty\min\{\updelta_u,\updelta_{\ubar}\}).
		\end{align}
	\end{theorem}
	
	\begin{proof}[Discussion of the proof of Theorem~\ref{T:EnergyEstimates0}]
		We follow the proof strategy used in \cite[Theorem 9.10]{IntegralIdentities}.
		
		To prove (\ref{ES:KEnergy1}), we use (\ref{E:integralidentities1}) and Lemma \ref{LE:integralID} for $\quadsmoothfunction(\pfour V,\pfour V)$, Lemma \ref{L:Geometricproperty} for boundary terms. 
		
		To prove (\ref{ES:Energy1}), we use Lemma \ref{basicenergy} with $\varphi\in\{\logDensity,v^i\}$, Lemma \ref{LE:basicenergy4} with $\varphi\in\{\Ent,\gradEnt,\modivort^i,\DivGradEnt\}$, (\ref{E:WaveEQ}) for $v,\logDensity$, (\ref{E:Bs}) for $\Ent$,  (\ref{E:BC})-(\ref{E:BD}) for $\modivort,\DivGradEnt$. 
		
		(\ref{ES:CKEnergies2}) follows by combining (\ref{ES:CKEnergies1}), bootstrap assumptions (\ref{BAfluid}), (\ref{BAacoustic}), Gr\"onwall's inequality (with respect to $\uptau$), and Young's inequality.
	\end{proof}
	
	\begin{corollary}[Base-level energy estimates in a neighborhood of $S_{0,0}$]\label{C:Energybound1}
		Let $\updelta_u:=u_1-u_0, \updelta_{\ubar}:=\ubar_1-\ubar_0$. Assume that $u_0=\ubar_0=0$. $\region_{1}^{0,0,\uptau}$ is defined in (\ref{DE:triangleM}). $v^i,\logDensity,\Ent$ are the solutions of the compressible Euler equations (\ref{E:Brho})-(\ref{E:Bs}). Then, under the initial data assumptions (\ref{A:H})-(\ref{A:Boundness}) and bootstrap assumptions (\ref{BAfluid})-(\ref{BAgeometry3}), for various energies $\KEnergy^{(1)},\bEnergy,\Energy,\uEnergy$ defined in Definition \ref{D:ControllingEnergyQuantities}, where $0\leq\uptau\leq\min\{\updelta_u,\updelta_{\ubar}\}$, the following energy estimates hold:
		\begin{align}\label{ES:Energybound1}
			\KEnergy^{(1)}_{(1)}(0,0,\uptau),\bEnergy_{(1)}(0,0,\uptau)\lesssim \iniepsilon\cdot\min\{\updelta_u,\updelta_{\ubar}\}\cdot\exp(\RBAinfty\min\{\updelta_u,\updelta_{\ubar}\}),
		\end{align}
		where $\iniepsilon$ is the size of the initial data; see Section \ref{SS:iniregularity2}.
	\end{corollary}
	
	\begin{corollary}[Base-level energy estimates in a neighborhood of $S_{u_0,\ubar_0}$]\label{C:Energybound1.2}
		Let $\updelta_u:=u_1-u_0, \updelta_{\ubar}:=\ubar_1-\ubar_0$. $\region_{1}^{u_0,\ubar_0,\uptau}$ is defined in (\ref{DE:triangleM}). $v^i,\logDensity,\Ent$ are the solutions of the compressible Euler equations (\ref{E:Brho})-(\ref{E:Bs}). Then, under the initial data assumptions (\ref{A:H})-(\ref{A:Boundness}) and bootstrap assumptions (\ref{BAfluid})-(\ref{BAgeometry3}), for various energies $\KEnergy^{(1)},\bEnergy,\Energy,\uEnergy$ defined in Definition \ref{D:ControllingEnergyQuantities}, where $0\leq\uptau\leq\min\{u_0+\ubar_1,u_1+\ubar_0\}$, and\footnote{Such $\updelta_u$ (or $\updelta_{\ubar}$) exists because of the absolute continuity of the integral such that, initial energy on $\hypersurbar_{\ubar_0}^{u_0,u_1}$ (or $\hypersur_{u_0}^{\ubar_0,\ubar_1}$) is small (which can be controlled by $\iniepsilon\cdot\min\{\updelta_u,\updelta_{\ubar}\}$).} 
		\begin{align}\label{AS:energy1}
			\KEnergy^{(1)}_{(1)}(u_0,\ubar_0,u_0+\ubar_0),\Energy_{(1)}(u_0,\ubar_0,\ubar_0+\min\{\updelta_u,\updelta_{\ubar}\}),\uEnergy_{(1)}(\ubar_0,u_0,u_0+\min\{\updelta_u,\updelta_{\ubar}\})\lesssim\iniLtwo\cdot\min\{\updelta_u,\updelta_{\ubar}\},
		\end{align}
		for some constant\footnote{In the rest of this article, $\iniLtwo$ denotes the size of the (various) energies of the fluid solutions and geometric quantities on $\hypersur_{u_0}\bigcup\hypersurbar_{\ubar_0}$.} $\iniLtwo$. Then, the following energy estimates hold:
		\begin{align}\label{ES:Energybound1.2}
			\KEnergy^{(1)}_{(1)}(u_0,\ubar_0,\uptau),\bEnergy_{(1)}(u_0,\ubar_0,\uptau)\lesssim \iniLtwo \cdot\min\{\updelta_u,\updelta_{\ubar}\}\cdot\exp(\RBAinfty\min\{\updelta_u,\updelta_{\ubar}\}).
		\end{align}
	\end{corollary}
	
	\begin{proof}[Proof of Corollary \ref{C:Energybound1} and Corollary \ref{C:Energybound1.2}]	
		By Section \ref{SS:iniregularity2} and the fact that $\uptau\leq\min\{\updelta_u,\updelta_{\ubar}\}$, we have 
		\begin{align}
			\KEnergy^{(1)}_{(1)}(0,0,0),\Energy_{(1)}(0,0,\uptau),\uEnergy_{(1)}(0,0,\uptau)\leq\iniepsilon\cdot\min\{\updelta_u,\updelta_{\ubar}\}.
		\end{align}
		Using (\ref{ES:CKEnergies2}), we obtain the desired result. Corollary \ref{C:Energybound1.2} can be proved in a similar fashion.
	\end{proof}

	\subsection{Base-Level Energy Estimates in the ``Parallelogram" Region $\region_2^{\uptau_1,\uptau_2}$}
	In this section, we derive energy estimates for fluid variables in the ``parallelogram" region $\region_2^{\uptau_1,\uptau_2}$ defined as in Section \ref{SS:SUBSETS}.
	
	\begin{theorem}[Base-level energy estimates for the fluid variables in the ``parallelogram" region $\region_2^{\uptau_1,\uptau_2}$]\label{TH:EnergyEstimates1}
		Let $\updelta_u:=u_1-u_0, \updelta_{\ubar}:=\ubar_1-\ubar_0$. $v^i,\logDensity,\Ent$ are the solutions of the compressible Euler equations (\ref{E:Brho})-(\ref{E:Bs}). Then, for various energies $\bEnergy,\KEnergy,\Energy,\uEnergy$ defined in Definition \ref{D:ControllingEnergyQuantities}, where $\min\{u_0+\ubar_1,u_1+\ubar_0\}\leq\uptau_1\leq\uptau_2\leq\max\{u_0+\ubar_1,u_1+\ubar_0\}$, we have the following energy estimates:\\
		
		\textbf{When $\updelta_u\leq\updelta_{\ubar}$:}
		\begin{subequations}\label{ES:CKEnergies3}
			\begin{align}
				\label{ES:KEnergy3}&\KEnergy^{(2)}_{(1)}(\uptau_1,\uptau_2)\lesssim\KEnergy^{(2)}_{(1)}(\uptau_1,\uptau_1)+\Energy_{(1)}(u_1,\uptau_1-u_1,\uptau_2-u_1)-\Energy_{(1)}(u_0,\uptau_1-u_0,\uptau_2-u_0)\\
				\notag&+\left\{\KEnergy^{(2)}_{(1)}(\uptau_1,\uptau_2)\int_{\uptau_1}^{\uptau_2}\bEnergy_{(1)}(u_0,\ubar_0,\uptau^\prime)\diff \uptau^\prime\right\}^{1/2}\\
				\notag&+\int_{\uptau_1}^{\uptau_2}\left[\bEnergy_{(1)}(u_0,\ubar_0,\uptau^\prime)\right]\RBAinfty\diff \uptau^\prime,\\
				\label{ES:Energy3}&\bEnergy_{(1)}(u_0,\ubar_0,\uptau_2)+\Energy_{(1)}(u_1,\uptau_1-u_1,\uptau_2-u_1)\\
				\notag\lesssim&\Energy_{(1)}(u_0,\uptau_1-u_0,\uptau_2-u_0)+\bEnergy_{(1)}(u_0,\ubar_0,\uptau_1)\\
				\notag&+\left\{\left[\KEnergy^{(2)}_{(1)}(\uptau_1,\uptau_2)\right]\int_{\uptau_1}^{\uptau_2}\bEnergy_{(1)}(u_0,\ubar_0,\uptau^\prime)\diff \uptau^\prime\right\}^{1/2}\\
				\notag&+\int_{\uptau_1}^{\uptau_2}\left[\bEnergy_{(1)}(u_0,\ubar_0,\uptau^\prime)\right]\RBAinfty\diff \uptau^\prime.
			\end{align} 
		\end{subequations}
		
		\textbf{When $\updelta_u>\updelta_{\ubar}$:}
		\begin{subequations}\label{ES:CKEnergies4}
			\begin{align}
				\label{ES:KEnergy4}&\KEnergy^{(2)}_{(1)}(\uptau_1,\uptau_2)\lesssim\KEnergy^{(2)}_{(1)}(\uptau_1,\uptau_1)+\uEnergy_{(1)}(\ubar_1,\uptau_1-\ubar_1,\uptau_2-\ubar_1)-\uEnergy_{(1)}(\ubar_0,\uptau_1-\ubar_0,\uptau_2-\ubar_0)\\
				\notag&+\left\{\KEnergy^{(2)}_{(1)}(\uptau_1,\uptau_2)\int_{\uptau_1}^{\uptau_2}\bEnergy_{(1)}(u_0,\ubar_0,\uptau^\prime)\diff \uptau^\prime\right\}^{1/2}\\
				\notag&+\int_{\uptau_1}^{\uptau_2}\left[\bEnergy_{(1)}(u_0,\ubar_0,\uptau^\prime)\right]\RBAinfty\diff \uptau^\prime,\\
				\label{ES:Energy4}&\bEnergy_{(1)}(u_0,\ubar_0,\uptau_2)+\uEnergy_{(1)}(\ubar_1,\uptau_1-\ubar_1,\uptau_2-\ubar_1)\\
				\notag\lesssim&\uEnergy_{(1)}(\ubar_0,\uptau_1-\ubar_0,\uptau_2-\ubar_0)+\bEnergy_{(1)}(u_0,\ubar_0,\uptau_1)\\
				\notag&+\left\{\left[\KEnergy^{(2)}_{(1)}(\uptau_1,\uptau_2)\right]\int_{\uptau_1}^{\uptau_2}\bEnergy_{(1)}(u_0,\ubar_0,\uptau^\prime)\diff \uptau^\prime\right\}^{1/2}\\
				\notag&+\int_{\uptau_1}^{\uptau_2}\left[\bEnergy_{(1)}(u_0,\ubar_0,\uptau^\prime)\right]\RBAinfty\diff \uptau^\prime.
			\end{align} 
		\end{subequations}
		Moreover, by bootstrap assumptions (\ref{BAfluid})-(\ref{BAacoustic}), we have:\\
		
		\textbf{When $\updelta_u\leq\updelta_{\ubar}$:}
		\begin{align}\label{ES:CKEnergies41}	&\KEnergy^{(2)}_{(1)}(\uptau_1,\uptau_2),\bEnergy_{(1)}(u_0,\ubar_0,\uptau_2),\\
			\notag\lesssim&\left\{\KEnergy^{(2)}_{(1)}(\uptau_1,\uptau_1)+\bEnergy_{(1)}(u_0,\ubar_0,\uptau_1)+\Energy_{(1)}(u_0,\uptau_1-u_0,\uptau_2-u_0)\right\}\exp\left(\RBAinfty\left[\max\{\updelta_u,\updelta_{\ubar}\}-\min\{\updelta_u,\updelta_{\ubar}\}\right]\right).
		\end{align}
		
		\textbf{When $\updelta_u>\updelta_{\ubar}$:}
		\begin{align}\label{ES:CKEnergies42}	&\KEnergy^{(2)}_{(1)}(\uptau_1,\uptau_2),\bEnergy_{(1)}(u_0,\ubar_0,\uptau_2)\\
			\notag\lesssim&\left\{\KEnergy^{(2)}_{(1)}(\uptau_1,\uptau_1)+\bEnergy_{(1)}(u_0,\ubar_0,\uptau_1)+\uEnergy_{(1)}(u_0,\uptau_1-\ubar_0,\uptau_2-\ubar_0)\right\}\exp\left(\RBAinfty\left[\max\{\updelta_u,\updelta_{\ubar}\}-\min\{\updelta_u,\updelta_{\ubar}\}\right]\right).
		\end{align}
	\end{theorem}
	
	\begin{proof}[Discussion of the proof of Theorem~\ref{TH:EnergyEstimates1}]
		We refer readers to \cite[Theorem 9.10]{IntegralIdentities} for details of the proof. 
		
		To prove (\ref{ES:KEnergy3}), (\ref{ES:KEnergy4}), we use (\ref{E:integralidentities2}), (\ref{E:integralidentities3}) respectively. We also use Lemma \ref{LE:integralID} for $\quadsmoothfunction(\pfour V,\pfour V)$, Lemma \ref{L:Geometricproperty} for boundary terms. 
		
		To prove (\ref{ES:Energy3}) and (\ref{ES:Energy4}), we use Lemma \ref{LE:basicenergy2} with $\varphi\in\{\logDensity,v^i\}$, Lemma \ref{LE:basicenergy5} with $\varphi\in\{\Ent,\gradEnt,\modivort^i,\DivGradEnt\}$, (\ref{E:WaveEQ}) for $v,\logDensity$, (\ref{E:Bs}) for $\Ent$,  (\ref{E:BC})-(\ref{E:BD}) for $\modivort,\DivGradEnt$.
		
		(\ref{ES:CKEnergies41})/(\ref{ES:CKEnergies42}) follows by combining (\ref{ES:CKEnergies3})/(\ref{ES:CKEnergies4}), bootstrap assumptions (\ref{BAfluid}), (\ref{BAacoustic}), Gr\"onwall's inequality (with respect to $\uptau_2$), and Young's inequality.
	\end{proof}
	\begin{corollary}[A priori estimates in the ``parallelogram" region near $\hypersur_{0}$ and $\hypersurbar_{0}$]\label{C:Energybound2.1}
		Let $\updelta_u:=u_1-u_0, \updelta_{\ubar}:=\ubar_1-\ubar_0$. Assume $u_0=\ubar_0=0$. $\region_{2}^{\min\{\updelta_u,\updelta_{\ubar}\},\uptau}$ is defined in (\ref{DE:parallelogramM}). $v^i,\logDensity,\Ent$ are the solutions of the compressible Euler equations (\ref{E:Brho})-(\ref{E:Bs}). Then, under the initial data assumptions (\ref{A:H})-(\ref{A:Boundness}) and bootstrap assumptions (\ref{BAfluid})-(\ref{BAacoustic}), for various energies $\KEnergy^{(1)},\bEnergy,\Energy,\uEnergy$ defined in Definition \ref{D:ControllingEnergyQuantities}, where $\min\{\updelta_u,\updelta_{\ubar}\}\leq\uptau\leq\max\{\updelta_u,\updelta_{\ubar}\}$, the following energy estimates hold:
		\begin{align}\label{ES:Energybound2.1}
			\KEnergy^{(2)}_{(1)}(\min\{\updelta_u,\updelta_{\ubar}\},\uptau),\bEnergy_{(1)}(0,0,\uptau)\lesssim \exp\left(\RBAinfty\max\{\updelta_u,\updelta_{\ubar}\}\right)\iniepsilon,
		\end{align}
		where $\iniepsilon$ is the size of the initial data, see Section \ref{SS:iniregularity2}.
	\end{corollary}
	\begin{corollary}[A priori estimates in the ``parallelogram" region near $\hypersur_{u_0}$ and $\hypersurbar_{\ubar_0}$]\label{C:Energybound2.2}
		Let $\updelta_u:=u_1-u_0, \updelta_{\ubar}:=\ubar_1-\ubar_0$. $\region_{2}^{\min\{u_0+\ubar_1,u_1+\ubar_0\},\uptau}$ is defined in (\ref{DE:parallelogramM}). $v^i,\logDensity,\Ent$ are the solutions of the compressible Euler equations (\ref{E:Brho})-(\ref{E:Bs}). Then, under the initial data assumptions (\ref{A:H})-(\ref{A:Boundness}) and bootstrap assumptions (\ref{BAfluid})-(\ref{BAacoustic}), for various energies $\KEnergy,\Energy,\uEnergy$ defined in Definition \ref{D:ControllingEnergyQuantities}, where $\min\{u_0+\ubar_1,u_1+\ubar_0\}\leq\uptau\leq\max\{u_0+\ubar_1,u_1+\ubar_0\}$, and:\\
		
		\textbf{When $\updelta_u\leq\updelta_{\ubar}$:}
		\begin{subequations}
			\begin{align}\label{ES:Assumption2.1}
				\KEnergy^{(1)}_{(1)}(u_0,\ubar_0,\ubar_0),\uEnergy_{(1)}(\ubar_0,u_0,u_1)\lesssim&\iniLtwo\cdot\min\{\updelta_u,\updelta_{\ubar}\},\\
				\Energy_{(1)}(u_0,\ubar_0,\ubar_1)\lesssim&\iniLtwo.
			\end{align}
		\end{subequations}
		
		\textbf{When $\updelta_u>\updelta_{\ubar}$:}
		\begin{subequations}
			\begin{align}\label{ES:Assumption2.2}
				\KEnergy^{(1)}_{(1)}(u_0,\ubar_0,\ubar_0),\Energy_{(1)}(u_0,\ubar_0,\ubar_1)\lesssim&\iniLtwo\cdot\min\{\updelta_u,\updelta_{\ubar}\},\\
				\uEnergy_{(1)}(\ubar_0,u_0,u_1)\lesssim&\iniLtwo,
			\end{align}
		\end{subequations}
		for some constant $\iniLtwo$.
		Then, we have the following energy estimates:
		\begin{align}\label{ES:Energybound2.2}
			\KEnergy^{(2)}_{(1)}(\min\{u_0+\ubar_1,u_1+\ubar_0\},\uptau),\bEnergy_{(1)}(u_0,\ubar_0,\uptau)\lesssim \iniLtwo\exp\left(\RBAinfty\max\{\updelta_u,\updelta_{\ubar}\}\right).
		\end{align}
		
	\end{corollary}
	\begin{proof}[Proof of Corollary \ref{C:Energybound2.1}-\ref{C:Energybound2.2}]
		We only prove Corollary \ref{C:Energybound2.2}; Corollary \ref{C:Energybound2.1} can be proved by the same approach. 
		
		By Corollary \ref{C:Energybound1.2} with $\uptau=\min\{u_0+\ubar_1,u_1+\ubar_0\}$, we have
		\begin{align}\label{ES:Energybound2.3}
			\KEnergy^{(2)}_{(1)}(\min\{u_0+\ubar_1,u_1+\ubar_0\},\min\{u_0+\ubar_1,u_1+\ubar_0\}),\bEnergy_{(1)}(u_0,\ubar_0,\min\{u_0+\ubar_1,u_1+\ubar_0\})\lesssim \iniLtwo\cdot\min\{\updelta_u,\updelta_{\ubar}\}.
		\end{align}
		Using (\ref{ES:CKEnergies41})-(\ref{ES:CKEnergies42}) with $\uptau_1=\min\{u_0+\ubar_1,u_1+\ubar_0\}$, by (\ref{ES:Assumption2.1}), (\ref{ES:Assumption2.2}) and (\ref{ES:Energybound2.3}), we obtain the desired result (\ref{ES:Energybound2.2}).
	\end{proof}
	\subsection{Energy Estimates for the Fluid Variables in the ``Trapezoid'' Region $\region_{3}^{\uptau_1,\uptau_2}$}
	In this section, we derive first-order energy estimates for fluid variables in the ``trapezoid'' region $\region_{3}$ (defined in (\ref{DE:trapezoidM})).
	
	\begin{theorem}[Base-level energy estimates for the fluid variables in the ``Trapezoid" Region $\region_{3}^{\uptau_1,\uptau_2}$]\label{TH:EnergyEstimates04}
		Let $\updelta_u:=u_1-u_0, \updelta_{\ubar}:=\ubar_1-\ubar_0$. $v^i,\logDensity,\Ent$ are the solutions of the compressible Euler equations (\ref{E:Brho})-(\ref{E:Bs}). Then, for various energies $\bEnergy,\KEnergy,\Energy,\uEnergy$ defined in Definition \ref{D:ControllingEnergyQuantities}, where $\max\{u_0+\ubar_1,u_1+\ubar_0\}\leq\uptau_1\leq\uptau_2\leq u_1+\ubar_1$, the following energy estimates hold:
		\begin{subequations}\label{ES:CKEnergies40}
			\begin{align}
				\label{ES:KEnergy41}&\KEnergy^{(3)}_{(1)}(\uptau_1,\uptau_2)+\Energy_{(1)}(u_1,\uptau_1-u_1,\uptau_2-u_1)+\uEnergy_{(1)}(\ubar_1,\uptau_1-\ubar_1,\uptau_2-\ubar_1)\\
				\notag\lesssim&\KEnergy^{(3)}_{(1)}(\uptau_1,\uptau_1)+\left\{\KEnergy^{(3)}_{(1)}(\uptau_1,\uptau_2)\int_{\uptau_1}^{\uptau_2}\bEnergy_{(1)}(u_0,\ubar_0,\uptau^\prime)\diff \uptau^\prime\right\}^{1/2}\\
				\notag&+\int_{\uptau_1}^{\uptau_2}\left[\bEnergy_{(1)}(u_0,\ubar_0,\uptau^\prime)\right]\RBAinfty\diff \uptau^\prime,\\
				\label{ES:Energy41}&\bEnergy_{(1)}(u_0,\ubar_0,\uptau_2)+\Energy_{(1)}(u_1,\uptau_1-u_1,\uptau_2-u_1)+\uEnergy_{(1)}(\ubar_1,\uptau_1-\ubar_1,\uptau_2-\ubar_1)\\
				\notag\lesssim&\bEnergy_{(1)}(u_0,\ubar_0,\uptau_1)+\left\{\left[\KEnergy^{(3)}_{(1)}(\uptau_1,\uptau_2)\right]\int_{\uptau_1}^{\uptau_2}\bEnergy_{(1)}(u_0,\ubar_0,\uptau^\prime)\diff \uptau^\prime\right\}^{1/2}\\
				\notag&+\int_{\uptau_1}^{\uptau_2}\left[\bEnergy_{(1)}(u_0,\ubar_0,\uptau^\prime)\right]\RBAinfty\diff \uptau^\prime.
			\end{align} 
		\end{subequations}
		Moreover, by bootstrap assumptions (\ref{BAfluid})-(\ref{BAacoustic}), we have:
		\begin{align}\label{ES:CKEnergies43}	&\bEnergy_{(1)}(u_0,\ubar_0,\uptau_2),\KEnergy^{(3)}_{(1)}(\uptau_1,\uptau_2),\Energy_{(1)}(u_1,\uptau_1-u_1,\uptau_2-u_1),\uEnergy_{(1)}(\ubar_1,\uptau_1-\ubar_1,\uptau_2-\ubar_1)\\
			\notag\lesssim&\left\{\KEnergy^{(3)}_{(1)}(\uptau_1,\uptau_1)+\bEnergy_{(1)}(u_0,\ubar_0,\uptau_1)\right\}\exp\left(\RBAinfty\min\{\updelta_u,\updelta_{\ubar}\}\right).
		\end{align}
	\end{theorem}
	
	\begin{proof}[Discussion of the proof of Theorem~\ref{TH:EnergyEstimates04}]
		We use the approach used in the proof of \cite[Theorem 9.10]{IntegralIdentities}. 
		
		To prove (\ref{ES:KEnergy41}), we use (\ref{E:integralidentities4}) and Lemma \ref{LE:integralID} for $\quadsmoothfunction(\pfour V,\pfour V)$, Lemma \ref{L:Geometricproperty} for boundary terms.
		
		To prove (\ref{ES:Energy41}), we use Lemma \ref{LE:basicenergy3.5} with $\varphi\in\{\logDensity,v^i\}$, Lemma \ref{LE:basicenergy5.5} with $\varphi\in\{\Ent,\gradEnt,\modivort^i,\DivGradEnt\}$, (\ref{E:WaveEQ}) for $v,\logDensity$, (\ref{E:Bs}) for $\Ent$,  (\ref{E:BC})-(\ref{E:BD}) for $\modivort,\DivGradEnt$.
		
		(\ref{ES:CKEnergies43}) follows by combining (\ref{ES:CKEnergies40}), bootstrap assumptions (\ref{BAfluid}), (\ref{BAacoustic}), Gr\"onwall's inequality (with respect to $\uptau_2$), and Young's inequality.
	\end{proof}
	\begin{corollary}[A priori estimates in the ``trapezoid" region $\region_{3}^{\uptau_1,\uptau_2}$]\label{C:Energybound3.1}
		Let $\updelta_u:=u_1-u_0, \updelta_{\ubar}:=\ubar_1-\ubar_0$. Assume $u_0=\ubar_0=0$. $\region_{3}^{\max\{u_0+\ubar_1,u_1+\ubar_0\},\uptau}$ is defined in (\ref{DE:trapezoidM}). $v^i,\logDensity,\Ent$ are the solutions of the compressible Euler equations (\ref{E:Brho})-(\ref{E:Bs}). Then, under the initial data assumptions (\ref{A:H})-(\ref{A:Boundness}) and bootstrap assumptions (\ref{BAfluid})-(\ref{BAacoustic}), for various energies $\KEnergy^{(1)},\bEnergy,\Energy,\uEnergy$ defined in Definition \ref{D:ControllingEnergyQuantities}, where $\max\{\updelta_u,\updelta_{\ubar}\}\leq\uptau\leq\updelta_u+\updelta_{\ubar}$, the following energy estimates hold:
		\begin{align}\label{ES:Energybound3.1}
			&\KEnergy^{(3)}_{(1)}(\max\{\updelta_u,\updelta_{\ubar}\},\uptau),\bEnergy_{(1)}(0,0,\uptau),\Energy_{(1)}(u_1,\max\{\updelta_u,\updelta_{\ubar}\}-u_1,\uptau-u_1),\uEnergy_{(1)}(\ubar_1,\max\{\updelta_u,\updelta_{\ubar}\}-\ubar_1,\uptau-\ubar_1)\\
			\notag&\lesssim\exp\left(\RBAinfty\left[\updelta_u+\updelta_{\ubar}\right]\right)\iniepsilon,
		\end{align}
		where $\iniepsilon$ is the size of the initial data, see Section \ref{SS:iniregularity2}.
	\end{corollary}
	\begin{corollary}[A priori estimates in the ``trapezoid" region $\region_{3}^{\uptau_1,\uptau_2}$ near $\hypersur_{u_0}$ and $\hypersurbar_{\ubar_0}$]\label{C:Energybound3.2}
		Let $\updelta_u:=u_1-u_0, \updelta_{\ubar}:=\ubar_1-\ubar_0$. $\region_{3}^{\max\{u_0+\ubar_1,u_1+\ubar_0\},\uptau}$ is defined in (\ref{DE:trapezoidM}). $v^i,\logDensity,\Ent$ are the solutions of the compressible Euler equations (\ref{E:Brho})-(\ref{E:Bs}). Then, under the initial data assumptions (\ref{A:H})-(\ref{A:Boundness}) and bootstrap assumptions (\ref{BAfluid})-(\ref{BAacoustic}), for various energies $\KEnergy,\Energy,\uEnergy$ defined in Definition \ref{D:ControllingEnergyQuantities}, where $\max\{u_0+\ubar_1,u_1+\ubar_0\}\leq\uptau\leq u_1+\ubar_1$, and:\\
		
		\textbf{When $\updelta_u\leq\updelta_{\ubar}$:}
		\begin{subequations}
			\begin{align}\label{ES:Assumption3.1}
				\KEnergy^{(1)}_{(1)}(u_0,\ubar_0,\ubar_0),\uEnergy_{(1)}(\ubar_0,u_0,u_1)\lesssim&\iniLtwo\cdot\min\{\updelta_u,\updelta_{\ubar}\},\\
				\Energy_{(1)}(u_0,\ubar_0,\ubar_1)\lesssim&\iniLtwo.
			\end{align}
		\end{subequations}
		
		\textbf{When $\updelta_u>\updelta_{\ubar}$:}
		\begin{subequations}
			\begin{align}\label{ES:Assumption3.2}
				\KEnergy^{(1)}_{(1)}(u_0,\ubar_0,\ubar_0),\Energy_{(1)}(u_0,\ubar_0,\ubar_1)\lesssim&\iniLtwo\cdot\min\{\updelta_u,\updelta_{\ubar}\},\\
				\uEnergy_{(1)}(\ubar_0,u_0,u_1)\lesssim&\iniLtwo.,
			\end{align}
		\end{subequations}
		for some constant $\iniLtwo$.
		Then, we have the following energy estimates:
		\begin{align}\label{ES:Energybound3.2}
			&\KEnergy^{(3)}_{(1)}(\max\{u_0+\ubar_1,u_1+\ubar_0\},\uptau),\bEnergy_{(1)}(u_0,\ubar_0,\uptau),\\
			\notag&\Energy_{(1)}(u_1,\max\{u_0+\ubar_1,u_1+\ubar_0\}-u_1,\uptau-u_1),\uEnergy_{(1)}(\ubar_1,\max\{u_0+\ubar_1,u_1+\ubar_0\}-\ubar_1,\uptau-\ubar_1)\\
			\notag&\lesssim \iniLtwo\exp\left(\RBAinfty\left[\updelta_u+\updelta_{\ubar}\right]\right).
		\end{align}
		
	\end{corollary}
	\begin{proof}[Proof of Corollary \ref{C:Energybound3.1}-\ref{C:Energybound3.2}]
		We only prove Corollary \ref{C:Energybound3.2}, Corollary \ref{C:Energybound3.1} can be proved by the same approach. 
		
		By Corollary \ref{C:Energybound1.2} with $\uptau=\min\{u_0+\ubar_1,u_1+\ubar_0\}$, and Corollary \ref{C:Energybound2.2} with $\uptau=\max\{u_0+\ubar_1,u_1+\ubar_0\}$, we have
		\begin{align}\label{ES:Energybound3.3}
			&\KEnergy^{(3)}_{(1)}(\max\{u_0+\ubar_1,u_1+\ubar_0\},\max\{u_0+\ubar_1,u_1+\ubar_0\}),\bEnergy_{(1)}(u_0,\ubar_0,\max\{u_0+\ubar_1,u_1+\ubar_0\})\\
			\notag&\lesssim \iniLtwo\exp(\RBAinfty\max\{\updelta_u,\updelta_{\ubar}\}).
		\end{align}
		Using (\ref{ES:CKEnergies43}) with $\uptau_1=\max\{u_0+\ubar_1,u_1+\ubar_0\}$, by  (\ref{ES:Energybound3.3}), we obtain the desired result (\ref{ES:Energybound3.2}).
	\end{proof}

	\section{A Priori Estimates up to Top-Order and Recovering the Bootstrap Assumptions}\label{S:APRIORIESTIMATESUPTOTOPORDER}
	In this section, we use the bootstrap assumptions to derive energy estimates up to top-order. 
	We emphasize that our estimates are obtained by commuting the equations \underline{only} with Cartesian coordinate partial derivatives. 
	We use the acoustic geometry (i.e., the eikonal functions $u$ and $\ubar$) only to control the ``shape'' of the regions on which the estimates hold.
	As in Section~\ref{S:BASELEVELENERGYESTIMATES}, 
	our estimates apply on regions of type $\region_{1},\region_{2},\region_{3}$, as defined in Sect.\,\ref{SS:SUBSETS}.
	All three types of regions will play a role in the proof of our main existence result, Theorem~\ref{TH:localwellposedness}.
	We then derive strict improvements of the bootstrap assumptions.
	
	\subsection{Pointwise Estimates Satisfied by Derivatives of the Solution Variables}\label{SS:SectionCommutators}
	To start, we use the bootstrap assumptions to derive simple pointwise estimates for the wave equations after they are commuted with
	Cartesian partial derivatives $\pfour$.
	
	\begin{definition}[Low-order terms' notation]\label{DE:LOT}
		Let $\LOTsmoothfunction:\mathbb{R}\rightarrow\mathbb{R}$ be an increasing function such that $\LOTsmoothfunction(\BAinfty)\leq \max\{(\BAinfty)^{2N},e^{\BAinfty}\}=:\widetilde{\BAinfty}$.
	\end{definition}
	\begin{lemma}[Pointwise estimates satisfied by derivatives of the solution variables]\label{L:lemmacommute}
		Let $\frac{N-1}{4}\leq k\leq\frac{N-1}{2}$. We consider the solutions to the equations of Theorem~\ref{T:geometricformulation}, that is, if  $\variables\in\{v^0,v^1,v^2,v^3,\logDensity,\Ent\}$ solves the compressible Euler equations (\ref{E:firstorderformulation}), by bootstrap assumptions \eqref{BAfluid}, the following estimates hold:
		\begin{subequations}\label{E:Commutedequations}
			\begin{align}
				\label{boxg2}&\abs{\boxg\pfour^{k-1}\variables}\lesssim\sum\limits_{\frac{N-1}{4}\leq i\leq k}\left(\abs{\pfour^{i-1}\modivort,\pfour^{i-1}\DivGradEnt,\pfour^{i}\vvariables}+1\right)\cdot\RBAinfty,\\
				\label{bc2}&\abs{\materialderivative\pfour^{k-1}\modivort^a}\lesssim\sum\limits_{\frac{N-1}{4}\leq i\leq k}\left(\abs{\pfour^{i}\vvariables,\pfour^{i}\vvort,\pfour^{i}\vgradientEnt}+1\right)\cdot\RBAinfty,\\
				\label{bd2}&\abs{\materialderivative\pfour^{k-1}\DivGradEnt}\lesssim \sum\limits_{\frac{N-1}{4}\leq i\leq k}\left(\abs{\pfour^i\vvort,\pfour^i\vvariables,\pfour^i\vgradientEnt}+1\right)\cdot\RBAinfty.
			\end{align}
		\end{subequations}
	\end{lemma}
	\begin{proof}[Sketch of the proof of Lemma \ref{L:lemmacommute}]
		By commuting (\ref{E:WaveEQ}), (\ref{E:BC}) and (\ref{E:BD}) with $\pfour^{k-1}$, using relations (by Definition \ref{D:modifiedfluidvariables}) $\modivort=\linsmoothfunction(\vvariables,\vvort,\vgradientEnt)[\pfour\vvort,\pfour\vvariables]$ and $\DivGradEnt=\linsmoothfunction[\pfour\vgradientEnt]+\linsmoothfunction(\vgradientEnt)[\pfour\vvariables]$, (\ref{E:Commutedequations}) are derived by bootstrap assumptions \eqref{BAfluid} and straightforward computations.
	\end{proof}
	
	
	\subsection{A Priori Energy Estimates up to Top-order}
	\label{SS:APRIORIENERGYESTIMATESWITHISOTROPICREGULARITY}
	The next theorem yields our main a priori energy estimates.
	
	\begin{theorem}[A priori energy estimates up to top-order] \label{TH:UPTOTOPORDERENERGYESTIMATES}
		Let $\region_0:=\{u_0\leq u \leq u_1, \ubar_0 \leq\ubar\leq\ubar_1\}$. Let $\updelta_u:=u_1-u_0, \updelta_{\ubar}:=\ubar_1-\ubar_0$. Assume that 
		$v^i,\logDensity,\Ent$ are $C^{\infty}$ solutions to the compressible Euler equations \eqref{E:firstorderformulation}
		(and hence also solutions to the equations of Theorem~\ref{T:geometricformulation}). 
		Then under the initial data assumptions (\ref{A:H})-(\ref{A:Boundness}) and the bootstrap assumptions (\ref{BAfluid})-(\ref{BAacoustic}), for the energies $\tilde{\KEnergy},\tilde{\bEnergy},\tilde{\Energy},\tilde{\uEnergy}$ defined in Definition~\ref{D:ControllingEnergyQuantities} where $u_0+\ubar_0<\uptau\leq u_1+\ubar_1$, $1\leq k\leq \frac{N-1}{2}$ (see Section \ref{SS:SOBOLEVREGULARITYOFDATA} for $N$), assuming $\tilde{\KEnergy}^{(1)}_{(j)},\tilde{\KEnergy}^{(2)}_{(j)},\tilde{\KEnergy}^{(3)}_{(j)},\tilde{\bEnergy}_{(j)}\lesssim \iniLtwo\exp(\updelta_u+\updelta_{\ubar})$ for $j<k$, and\footnote{This assumption is the same as $\tilde{\KEnergy}_{(k)}(u_0,\ubar_0,u_0+\ubar_0)\lesssim\min\{\updelta_u,\updelta_{\ubar}\}$, and $	\left\{\tilde{\Energy}_{(k)}(u_0,\ubar_0,\ubar_1)+\tilde{\uEnergy}_{(k)}(\ubar_0,u_0,u_1)\right\}\lesssim 1$, which are satisfied on the initial hypersurfaces $\hypersur_{0},\hypersurbar_0$. Such $\updelta_u$ (or $\updelta_{\ubar}$) exists because of the absolute continuity of the integral such that, the energies on $\hypersurbar_{\ubar_0}^{u_0,u_0+\min\{\updelta_u,\updelta_{\ubar}\}}$ and  $\hypersur_{u_0}^{\ubar_0,\ubar_0+\min\{\updelta_u,\updelta_{\ubar}\}}$ are small (which can be controlled by $1\cdot\min\{\updelta_u,\updelta_{\ubar}\}$).}:\\
		
		\textbf{When $\updelta_u\leq\updelta_{\ubar}$:}
		\begin{subequations}
			\begin{align}\label{ES:Assumption4.1}
				\KEnergy_{(k)}^{(1)}(u_0,\ubar_0,\ubar_0),\uEnergy_{(k)}(\ubar_0,u_0,u_1)\lesssim&\iniLtwo\cdot\min\{\updelta_u,\updelta_{\ubar}\},\\
				\Energy_{(k)}(u_0,\ubar_0,\ubar_1)\lesssim&\iniLtwo.
			\end{align}
		\end{subequations}
		
		\textbf{When $\updelta_u>\updelta_{\ubar}$:}
		\begin{subequations}
			\begin{align}\label{ES:Assumption4.2}
				\KEnergy^{(1)}_{(k)}(u_0,\ubar_0,\ubar_0),\Energy_{(k)}(u_0,\ubar_0,\ubar_1)\lesssim&\iniLtwo\cdot\min\{\updelta_u,\updelta_{\ubar}\},\\
				\uEnergy_{(k)}(\ubar_0,u_0,u_1)\lesssim&\iniLtwo,
			\end{align}
		\end{subequations}
		for some constant $\iniLtwo$.
		Then, letting $\RBAinfty$ be as in Definition \ref{DE:LOT}, we have the following energy estimates:
		\begin{align}\label{ES:EnergyRegularity1}
			&\tilde{\KEnergy}^{(1)}_{(k)}(u_0,\ubar_0,\min\{u_0+\ubar_1,u_1+\ubar_0\}),\tilde{\KEnergy}^{(2)}_{(k)}(\min\{u_0+\ubar_1,u_1+\ubar_0\},\max\{u_0+\ubar_1,u_1+\ubar_0\}),\\
			\notag&\tilde{\KEnergy}^{(3)}_{(k)}(\max\{u_0+\ubar_1,u_1+\ubar_0\},u_1+\ubar_1),\tilde{\bEnergy}_{(k)}(u_0,\ubar_0,\uptau),\\
			\notag&\tilde{\Energy}_{(k)}(u_1,\ubar_0,\ubar_1),\tilde{\uEnergy}_{(k)}(\ubar_1,u_0,u_1)\\
			\notag&\lesssim \iniLtwo\exp\left(\RBAinfty\left[\updelta_u+\updelta_{\ubar}\right]\right).
		\end{align}
		
		In particular, if $u_0=\ubar_0=0$, then $\iniLtwo=\iniepsilon$.
	\end{theorem}
	\begin{proof}[Proof of Theorem~\ref{TH:UPTOTOPORDERENERGYESTIMATES}]
		We take the same approach as in the proof in Section \ref{S:BASELEVELENERGYESTIMATES}. It suffices to prove the following:\\ 
		
		\noindent\textbf{On $\region^{u_0,\ubar_0,\min\{u_0+\ubar_1,u_1+\ubar_0\}}_1$, where $u_0+\ubar_0\leq\uptau\leq\min\{u_0+\ubar_1,u_1+\ubar_0\}$:}
		\begin{subequations}\label{ES:REGEnergy1}
			\begin{align}
				\label{ES:REGEnergy1.1}\KEnergy^{(1)}_{(k)}(u_0,\ubar_0,\uptau)\lesssim&\KEnergy^{(1)}_{(k)}(u_0,\ubar_0,u_0+\ubar_0)+\tilde{\Energy}_{(k)}(u_0,\ubar_0,\uptau-u_0)+\tilde{\uEnergy}_{(k)}(\ubar_0,u_0,\uptau-\ubar_0)\\
				\notag&+\left\{\tilde{\KEnergy}^{(1)}_{(k)}(u_0,\ubar_0,\uptau)\int_{u_0+\ubar_0}^{\uptau}\tilde{\bEnergy}_{(k)}(u_0,\ubar_0,\uptau^\prime)\diff \uptau^\prime\right\}^{1/2}\\
				\notag&+\int_{u_0+\ubar_0}^{\uptau}\tilde{\bEnergy}_{(k)}(u_0,\ubar_0,\uptau^\prime)\RBAinfty\diff \uptau^\prime,\\
				\label{ES:REGEnergy1.2}\bEnergy_{(k)}(u_0,\ubar_0,\uptau)\lesssim&\tilde{\Energy}_{(k)}(u_0,\ubar_0,\uptau-u_0)+\tilde{\uEnergy}_{(k)}(\ubar_0,u_0,\uptau-\ubar_0)\\
				\notag&+\left\{\tilde{\KEnergy}^{(1)}_{(k)}(u_0,\ubar_0,\uptau)\int_{u_0+\ubar_0}^{\uptau}\tilde{\bEnergy}_{(k)}(u_0,\ubar_0,\uptau^\prime)\diff \uptau^\prime\right\}^{1/2}\\
				\notag&+\int_{u_0+\ubar_0}^{\uptau}\tilde{\bEnergy}_{(k)}(u_0,\ubar_0,\uptau^\prime)\RBAinfty\diff \uptau^\prime.
			\end{align} 
		\end{subequations}
		
		\noindent\textbf{On $\region_{2}^{u_1+\ubar_0,\uptau}$, when $\updelta_u\leq\updelta_{\ubar}$, $u_1+\ubar_0\leq\uptau \leq u_0+\ubar_1$:}
		\begin{subequations}\label{ES:REGEnergy2}
			\begin{align}
				\label{ES:REGEnergy2.1}\KEnergy^{(2)}_{(k)}(u_1+\ubar_0,\uptau)\lesssim&\tilde{\KEnergy}^{(2)}_{(k)}(u_1+\ubar_0,u_1+\ubar_0)+\tilde{\Energy}_{(k)}(u_1,\ubar_0,\uptau-u_1)-\tilde{\Energy}_{(k)}(u_0,\ubar_0+\updelta_u,\uptau-u_0)\\
				\notag&+\left\{\tilde{\KEnergy}^{(2)}_{(k)}(u_1+\ubar_0,\uptau)\int_{u_1+\ubar_0}^{\uptau}\tilde{\bEnergy}_{(k)}(u_0,\ubar_0,\uptau^\prime)\diff \uptau^\prime\right\}^{1/2}\\
				\notag&+\int_{u_1+\ubar_0}^{\uptau}\tilde{\bEnergy}_{(k)}(u_0,\ubar_0,\uptau^\prime)\RBAinfty\diff \uptau^\prime,\\
				\label{ES:REGEnergy2.2}\bEnergy_{(k)}(u_0,\ubar_0,\uptau)+\tilde{\Energy}_{(k)}(u_1,\ubar_0,\uptau-u_1)\lesssim&\tilde{\Energy}_{(k)}(u_0,\ubar_0+\updelta_u,\uptau-u_0)+\tilde{\bEnergy}_{(k)}(u_0,\ubar_0,u_1+\ubar_0)\\
				\notag&+\left\{\tilde{\KEnergy}^{(2)}_{(k)}(u_1+\ubar_0,\uptau)\int_{u_1+\ubar_0}^{\uptau}\tilde{\bEnergy}_{(k)}(u_0,\ubar_0,\uptau^\prime)\diff \uptau^\prime\right\}^{1/2}\\
				\notag&+\int_{u_1+\ubar_0}^{\uptau}\tilde{\bEnergy}_{(k)}(u_0,\ubar_0,\uptau^\prime)\RBAinfty\diff \uptau^\prime.
			\end{align} 
		\end{subequations}
		
		\noindent\textbf{On $\region_{2}^{u_0+\ubar_1,\uptau}$, when $\updelta_u>\updelta_{\ubar}$, $u_0+\ubar_1\leq\uptau \leq u_1+\ubar_0$:}
		\begin{subequations}\label{ES:REGEnergy3}
			\begin{align}
				\label{ES:REGEnergy3.1}\KEnergy^{(2)}_{(k)}(u_0+\ubar_1,\uptau)\lesssim&\tilde{\KEnergy}^{(2)}_{(k)}(u_0+\ubar_1,u_0+\ubar_1)+\tilde{\Energy}_{(k)}(\ubar_1,u_0,\uptau-\ubar_1)-\tilde{\Energy}_{(k)}(\ubar_0,u_0+\updelta_{\ubar},\uptau-\ubar_0)\\
				\notag&+\left\{\tilde{\KEnergy}^{(2)}_{(k)}(u_0+\ubar_1,\uptau)\int_{u_0+\ubar_1}^{\uptau}\tilde{\bEnergy}_{(k)}(u_0,\ubar_0,\uptau^\prime)\diff \uptau^\prime\right\}^{1/2}\\
				\notag&+\int_{u_0+\ubar_1}^{\uptau}\tilde{\bEnergy}_{(k)}(u_0,\ubar_0,\uptau^\prime)\RBAinfty\diff \uptau^\prime,\\
				\label{ES:REGEnergy3.2}\bEnergy_{(k)}(u_0,\ubar_0,\uptau)+\tilde{\Energy}_{(k)}(u_1,\ubar_0,\uptau-u_1)\lesssim&\tilde{\Energy}_{(k)}(u_0,\ubar_0+\updelta_u,\uptau-u_0)+\tilde{\bEnergy}_{(k)}(u_0,\ubar_0,u_0+\ubar_1)\\
				\notag&+\left\{\tilde{\KEnergy}^{(2)}_{(k)}(u_0+\ubar_1,\uptau)\int_{u_0+\ubar_1}^{\uptau}\tilde{\bEnergy}_{(k)}(u_0,\ubar_0,\uptau^\prime)\diff \uptau^\prime\right\}^{1/2}\\
				\notag&+\int_{u_0+\ubar_1}^{\uptau}\tilde{\bEnergy}_{(k)}(u_0,\ubar_0,\uptau^\prime)\RBAinfty\diff \uptau^\prime.
			\end{align} 
		\end{subequations}
		
		\noindent\textbf{On $\region_{3}^{\max\{u_0+\ubar_1,u_1+\ubar_0\},\uptau}$, where $\max\{u_0+\ubar_1,u_1+\ubar_0\}\leq\uptau \leq u_1+\ubar_1$:}
		\begin{subequations}\label{ES:REGEnergy4}
			\begin{align}
				\label{ES:REGEnergy4.1}&\KEnergy^{(3)}_{(k)}(\max\{u_0+\ubar_1,u_1+\ubar_0\},\uptau)+\tilde{\Energy}_{(k)}(u_1,\max\{u_0+\ubar_1,u_1+\ubar_0\}-u_1,\uptau-u_1)\\
				\notag&+\tilde{\uEnergy}_{(k)}(\ubar_1,\max\{u_0+\ubar_1,u_1+\ubar_0\}-\ubar_1,\uptau-\ubar_1)\\
				\notag\lesssim&\KEnergy^{(3)}_{(k)}(\max\{u_0+\ubar_1,u_1+\ubar_0\},\uptau)+\left\{\tilde{\KEnergy}^{(3)}_{(k)}(\max\{u_0+\ubar_1,u_1+\ubar_0\},\uptau)\int_{\max\{u_0+\ubar_1,u_1+\ubar_0\}}^{\uptau}\tilde{\bEnergy}_{(k)}(u_0,\ubar_0,\uptau^\prime)\diff \uptau^\prime\right\}^{1/2}\\
				\notag&+\int_{\max\{u_0+\ubar_1,u_1+\ubar_0\}}^{\uptau}\tilde{\bEnergy}_{(k)}(u_0,\ubar_0,\uptau^\prime)\RBAinfty\diff \uptau^\prime,\\
				\label{ES:REGEnergy4.2}&\bEnergy_{(k)}(u_0,\ubar_0,\uptau)+\tilde{\Energy}_{(k)}(u_1,\max\{u_0+\ubar_1,u_1+\ubar_0\}-u_1,\uptau-u_1)\\
				\notag&+\tilde{\uEnergy}_{(k)}(\ubar_1,\max\{u_0+\ubar_1,u_1+\ubar_0\}-\ubar_1,\uptau-\ubar_1)\\
				\notag\lesssim&\tilde{\bEnergy}_{(k)}(u_0,\ubar_0,\uptau_1)+\left\{\tilde{\KEnergy}^{(3)}_{(k)}(\max\{u_0+\ubar_1,u_1+\ubar_0\},\uptau)\int_{\max\{u_0+\ubar_1,u_1+\ubar_0\}}^{\uptau}\tilde{\bEnergy}_{(k)}(u_0,\ubar_0,\uptau^\prime)\diff \uptau^\prime\right\}^{1/2}\\
				\notag&+\int_{\max\{u_0+\ubar_1,u_1+\ubar_0\}}^{\uptau}\tilde{\bEnergy}_{(k)}(u_0,\ubar_0,\uptau^\prime)\RBAinfty\diff \uptau^\prime.
			\end{align} 
		\end{subequations}
		The proofs of (\ref{ES:REGEnergy1})-(\ref{ES:REGEnergy4}) are the same as in Theorem~\ref{T:EnergyEstimates0}, Theorem~\ref{TH:EnergyEstimates1}, and Theorem~\ref{TH:EnergyEstimates04}. We omit the details of the proof, while we point out the differences between (\ref{ES:REGEnergy1})-(\ref{ES:REGEnergy4}) and the proofs in Section \ref{S:BASELEVELENERGYESTIMATES}.
		
		To prove (\ref{ES:REGEnergy1.1}), (\ref{ES:REGEnergy2.1}), (\ref{ES:REGEnergy3.1}), (\ref{ES:REGEnergy4.1}), we use Lemma \ref{L:integralidentitiesforS} and substitute $V\in\{\vortrenormalized,\gradEnt\}$ by $\pfour^jV\in\{\pfour^{j}\vortrenormalized,\pfour^{j}\gradEnt\}$ with $0\leq j\leq k-1$. We note that for $\pfour^jV\in\{\pfour^j\vortrenormalized,\pfour^j\gradEnt\}$, we have:
		\begin{subequations}
			\begin{align}
				\label{E:A1}\mathfrak{A}_{\alpha\beta}^{(\pfour^j V)}:=&\pfour^j\mathfrak{A}_{\alpha\beta}^{(V)},\\
				\label{E:B1}\mathfrak{B}^i_{(\pfour^j V)}:=&\sum\limits_{m+n=j}\pfour^m\mathfrak{B}^i_{(V)}\cdot\pfour^n\vvariables,\\
				\mathfrak{C}^{(\pfour^j V)}:=&2(K_aK\pfour^j V)(\pfour^j\dive V)-\lapse^{-1}(K_aK(\pfour^j V)^a)K_b\mathfrak{B}^b_{(\pfour^j V)},\\
				\label{E:D1}\mathfrak{D}^{(\pfour^j V)}:=&\left\{\pfour^j\dive V+(2\lapse)^{-1}K_a\mathfrak{B}^a_{(\pfour^j V)}\right\},\\
				\label{E:H1}\mathfrak{H}_{(1)}[\pfour^j V]:=&\frac{1}{\lapse}(\Lunit^\alpha+\Theta^\alpha)(\sphereproject(\pfour^{j}V))^\beta\mathfrak{A}_{\alpha\beta}^{(\pfour^jV)},\\
				\label{E:uH1}\underline{\mathfrak{H}}_{(1)}[\pfour^j V]:=&\frac{1}{\lapse}(\uLunit^\alpha+\Theta^\alpha)(\sphereproject(\pfour^{j}V))^\beta\mathfrak{A}_{\alpha\beta}^{(\pfour^jV)}.
			\end{align}
		\end{subequations}
		where (\ref{E:A1}) is obtained by substituting $V$ by $\pfour^jV$ in (\ref{E:JAnti}), (\ref{E:B1})-(\ref{E:D1}) are obtained by Lemma \ref{L:lemmacommute} and substituting $V$ by $\pfour^jV$ in (\ref{E:JDiv}), (\ref{E:H1}) and (\ref{E:uH1}) are obtained by substituting $V$ by $\pfour^jV$ for the second term on the RHS of (\ref{E:NJ1}) and (\ref{E:NJ2}) respectively.
		
		To prove (\ref{ES:REGEnergy1.2}), (\ref{ES:REGEnergy2.2}), (\ref{ES:REGEnergy3.2}), (\ref{ES:REGEnergy4.2}), we use Lemma \ref{L:lemmacommute} in the following way: 1) for the ``wave-part'',  we use Lemma \ref{basicenergy}, Lemma \ref{LE:basicenergy2}, Lemma \ref{LE:basicenergy3.5}, respectively, with $\varphi\in\{\pfour^{j-1}\logDensity,\pfour^{j-1}v^i,\pfour^{j-1}\Ent\}$; 2) for the ``transport-part'', we use Lemma \ref{LE:basicenergy4}, Lemma \ref{LE:basicenergy5}, Lemma \ref{LE:basicenergy5.5}, respectively, with $\varphi\in\{\pfour^{j-1}\modivort^i,\pfour^{j-1}\DivGradEnt\}$. We note that the lower order terms can be absorbed by the highest order term and $\RBAinfty$, by bootstrap assumptions \eqref{BAfluid}-\eqref{BAacoustic}, and the induction assumptions $\tilde{\KEnergy}^{(1)}_{(j)},\tilde{\KEnergy}^{(2)}_{(j)},\tilde{\KEnergy}^{(3)}_{(j)}, \tilde{\bEnergy}_{(j)}\lesssim\iniLtwo\exp(\updelta_u+\updelta_{\ubar})$ for $j\leq k$.
		
		Then, combining estimates (\ref{ES:REGEnergy1})-(\ref{ES:REGEnergy4}), by Gr\"onwall's inequality, Young's inequality, and bootstrap assumptions (\ref{BAfluid}) and (\ref{BAacoustic}), we obtained the desired results (\ref{ES:EnergyRegularity1}).
	\end{proof}

		\subsection{A Priori $L^{\infty}$ Estimates}
	\label{SS:APRIORILINFINITYESTIMATESWITHISOTROPICREGULARITY}
	
	\begin{proposition}[$L^\infty$-bound of the solution and improvements of the bootstrap assumptions]\label{PR:Linftybound}
		Let $\RBAinfty$ be as in Definition \ref{DE:LOT}. Let $u_0\leq u\leq u_1$, $\ubar_0\leq u\leq \ubar_1$, $u_0+\ubar_0\leq\uptau\leq u_1+\ubar_1$. $\Stau$ is defined in (\ref{D:Sigmatau}). Let $ d_u:=u_1-u_0,  d_{\ubar}:=\ubar_1-\ubar_0$. Let $u_0\leq u^\prime\leq u_1$ and $\ubar_0\leq \ubar^\prime\leq \ubar_1$. $v^i,\logDensity,\Ent$ are solutions of the original geometric formulation in Theorem~\ref{T:geometricformulation} under the initial data assumptions (\ref{A:H})-(\ref{A:Boundness}) and bootstrap assumptions (\ref{BAfluid})-(\ref{BAacoustic}). Assuming $\min\{ d_u, d_{\ubar}\}\leq \delta$, where $\delta$ is sufficiently small such that (\ref{ES:Assumption4.1})-(\ref{ES:Assumption4.2}) hold and 
		\begin{align}\label{ES:smalldelta}
			\iniLtwo\cdot\exp\left\{\RBAinfty+\iniLtwo\exp\left(\RBAinfty\left[d_u+ d_{\ubar}\right]\right)\right\}\cdot\delta\lesssim\frac{1}{4}\BAinfty,
		\end{align} 
		where $\iniLtwo$ is the size of energy of the solution at $\hypersur_{u_0}$ and $\hypersurbar_{\ubar_0}$ as in Theorem~\ref{TH:UPTOTOPORDERENERGYESTIMATES}.
		Also, for $0\leq n\leq \frac{N-9}{2}$, assuming:
		\begin{align}\label{ES:Linftyassumption}
			\onenorm{\pfour^n(v^i,\logDensity,\Ent),\pfour^{n-1}(\modivort,\DivGradEnt)}{}{\infty}{\hypersur_{u_0}},\onenorm{\pfour^n(v^i,\logDensity,\Ent),\pfour^{n-1}(\modivort,\DivGradEnt)}{}{\infty}{\hypersurbar_{\ubar_0}}\leq\frac{1}{4}\BAinfty.
		\end{align}
		Then, we have the following $L^\infty$ bound for the solutions:
		\begin{align}
			\label{C:Linftyfluid2}&\onenorm{\pfour^n(v^i,\logDensity,\Ent),\pfour^{n-1}(\modivort,\DivGradEnt)}{}{\infty}{\region}\leq \frac{1}{2}\BAinfty,&
			\textnormal{for } n\leq&\frac{N-9}{2}.
		\end{align}
		Since $N\geq17$, we have $\frac{N-9}{2}\geq\frac{N-1}{4}$. Hence, (\ref{C:Linftyfluid2}) is a strict improvement of (\ref{BAfluid}).
	\end{proposition}
	\begin{proof}[Proof of Proposition \ref{PR:Linftybound}]
		We prove the case where $d_u\leq\delta$, while the proof of the case $d_{\ubar}\leq\delta$ is similar. For $k\leq\frac{N-1}{2}$, substituting $f$ by $\abs{\pfour^{k-1}(v^i,\logDensity,\Ent),\pfour^{k-2}(\modivort,\DivGradEnt)}^2$ in equation \eqref{EQ:LieuL}, by Gronwall's lemma, H\"older's inequality and bootstrap assumptions \eqref{BAfluid}-\eqref{BAgeometry3}, we have:
		\begin{align}
			\onenorm{\pfour^{k-1}(v^i,\logDensity,\Ent),\pfour^{k-2}(\modivort,\DivGradEnt)}{}{2}{S_{u^\prime,\ubar^\prime}}^2\leq\onenorm{\pfour^{k-1}(v^i,\logDensity,\Ent),\pfour^{k-2}(\modivort,\DivGradEnt)}{}{2}{S_{u_0,\ubar^\prime}}^2\cdot\exp\left([\BAinfty+\tilde{\uEnergy}_{(k)}(\ubar^\prime,u_0,u_1)]\cdot\delta^{1/2}\right).
		\end{align}
		By Sobolev embedding on spheres in Lemma \ref{LE:SobolevEmbedding}, Theorem~\ref{TH:UPTOTOPORDERENERGYESTIMATES}, and the consequence of the initial data regularity assumption in Section \ref{SSS:inisecondderivative}, we obtain: 
		\begin{align}
			\onenorm{\pfour^{k-3}(v^i,\logDensity,\Ent),\pfour^{k-4}(\modivort,\DivGradEnt)}{}{\infty}{S_{u^\prime,\ubar^\prime}}\lesssim \Gamma(\iniLtwo,\BAinfty,d_u,d_{\ubar}),
		\end{align}
		where 
		\begin{align}
			\Gamma(\iniLtwo,\BAinfty,d_u,d_{\ubar})=\iniLtwo\cdot\exp\left\{\RBAinfty+\iniLtwo\exp\left(\RBAinfty\left[d_u+ d_{\ubar}\right]\right)\right\}.
		\end{align}
		By the fundamental theorem of calculus, the assumption on the data on $\hypersur_{u_0}$ \eqref{ES:Linftyassumption} and smallness of $\delta$ as in \eqref{ES:smalldelta}, we have
		\begin{align}
			\onenorm{\pfour^{k-4}(v^i,\logDensity,\Ent),\pfour^{k-5}(\modivort,\DivGradEnt)}{}{\infty}{S_{u^\prime,\ubar^\prime}}\leq\frac{1}{4}\BAinfty+\delta C\cdot\Gamma(\iniLtwo,\BAinfty,d_u,d_{\ubar})\leq\frac{1}{2}\BAinfty,
		\end{align}
		where $C$ is a universal constant.
	\end{proof}
	\begin{proposition}[$L^\infty$-bound for geometric quantities]\label{PR:BoundforGeometricQuantities}
		The acoustical functions $u,\ubar$ are defined in Section \ref{SS:SUBSETSOFSPACETIME}. Let $u_0\leq u\leq u_1$, $\ubar_0\leq u\leq \ubar_1$, $u_0+\ubar_0\leq\uptau\leq u_1+\ubar_1$. $\Stau$ is defined in (\ref{D:Sigmatau}). Let $u_0\leq u^\prime\leq u_1$ and $\ubar_0\leq \ubar^\prime\leq \ubar_1$. Let $ d_u:=u_1-u_0,  d_{\ubar}:=\ubar_1-\ubar_0$. We note that $d_u,d_{\ubar}\lesssim 1$ and $\min\{d_u,d_{\ubar}\}\leq\delta$. The lapse variables $\invfolia,\uinvfolia,\lapse$ are defined in Section \ref{SS:SUBSETSOFSPACETIME}. Ricci coefficients $\ACC$ are defined in Definition \ref{D:Defconnectioncoeff}. Assuming bootstrap assumptions \eqref{BAfluid}-\eqref{BAgeometry3}, and\footnote{ $\pfour^3u,\pfour^3\ubar,\pfour^{\leq1}\ACC,\pfour^{\leq1}\ln\lapse,\pfour^{\leq1}\ln\invfolia,\pfour^{\leq1}\ln\uinvfolia$ satisfy \eqref{ES:Linftyassumption2} on $\hypersur_{0}\bigcup\hypersurbar_{0}$, see (\ref{Datageometry2}).}
		\begin{align}\label{ES:Linftyassumption2}
			\onenorm{\pfour^{\leq3}u,\pfour^{\leq3}\ubar,\pfour^{\leq1}\ACC,\pfour^{\leq1}\ln\lapse,\pfour^{\leq1}\ln\invfolia,\pfour^{\leq1}\ln\uinvfolia}{}{\infty}{\hypersur_{u_0}}&\leq\frac{1}{4}\BAinfty,\\
			\notag\onenorm{\pfour^{\leq3}u,\pfour^{\leq3}\ubar,\pfour^{\leq1}\ACC,\pfour^{\leq1}\ln\lapse,\pfour^{\leq1}\ln\invfolia,\pfour^{\leq1}\ln\uinvfolia}{}{\infty}{\hypersurbar_{\ubar_0}}&\leq\frac{1}{4}\BAinfty,
		\end{align}
		where $\pfour^{\leq k}\phi$ means the collection of all $\pfour^j\phi$ with $0\leq j\leq k$.
		Then, the following bounds hold:
		\begin{subequations}\label{ES:GeometricQuantities}
			\begin{align}
				\label{ES:uubar}&\onenorm{\pfour^{\leq3}u,\pfour^{\leq3}\ubar}{}{\infty}{S_{u^\prime,\ubar^\prime}}\leq\frac{1}{2}\exp(\BAinfty),\\	
				\label{ES:ACCGaussian}&\onenorm{\pfour^{\leq1}\ACC}{}{\infty}{S_{u^\prime,\ubar^\prime}}\leq\frac{1}{2}\exp(\BAinfty),\\
				\label{ES:Lapses}&\onenorm{\pfour^{\leq1}\ln\lapse,\pfour^{\leq1}\ln\invfolia,\pfour^{\leq1}\ln\uinvfolia}{}{\infty}{S_{u^\prime,\ubar^\prime}}\leq\frac{1}{2}\exp(\BAinfty).
			\end{align}
		\end{subequations}
		We note that \eqref{ES:GeometricQuantities} strictly improves \eqref{BAacoustic}.
	\end{proposition}
	\begin{proof}[Proof of Proposition \ref{PR:BoundforGeometricQuantities}]
		We first prove \eqref{ES:uubar}. Commuting \eqref{EQ:eikonalequations} with $\pfour^{k}$ with $0\leq k\leq 3$,
		we deduce the following estimates:
		\begin{subequations}\label{E:Commutedequations2}
			\begin{align}
				\label{ES:Lu1}&\abs{e_4\pfour^{k} u}\lesssim\sum\limits_{i\leq k}\abs{\gensmoothfunction_3\left(\pfour^{\leq i}\vvariables,\pfour^{\leq i} u\right)},\\
				\label{ES:Lbarubar1}&\abs{e_3\pfour^{k} \ubar}\lesssim\sum\limits_{i\leq k}\abs{\gensmoothfunction_3\left(\pfour^{\leq i}\vvariables,\pfour^{\leq i} \ubar\right)},
			\end{align}
		\end{subequations}
		where $\gensmoothfunction_3$ is an up-to-cubic polynomial of its arguments. Integrating \eqref{ES:Lu1} over $e_4$, applying Gr\"onwall's inequality, using the fact that $d_{\ubar}\lesssim 1$, we obtain:
		\begin{align}
			\onenorm{\pfour^{\leq3}u}{}{\infty}{S_{u^\prime,\ubar^\prime}}\lesssim\onenorm{\pfour^{\leq3}u}{}{\infty}{S_{u^\prime,\ubar_0}}\exp\left(\int_{\ubar_0}^{\ubar^\prime}\onenorm{\pfour^{\leq3}\vvariables}{}{\infty}{S_{u^\prime,\ubar}}\diff\ubar\right)\leq \frac{1}{2}\exp(\BAinfty).
		\end{align}
		Similarly, integrating \eqref{ES:Lbarubar1} over $e_3$, we have:
		\begin{align}
			\onenorm{\pfour^{\leq3}\ubar}{}{\infty}{S_{u^\prime,\ubar^\prime}}\lesssim\onenorm{\pfour^{\leq3}\ubar}{}{\infty}{S_{u_0,\ubar^\prime}}\exp\left(\int_{u_0}^{u^\prime}\onenorm{\pfour^{\leq3}\vvariables}{}{\infty}{S_{u,\ubar^\prime}}\diff u\right)\leq \frac{1}{2}\exp(\BAinfty).
		\end{align}
		Recall the  definition of $\ACC:=\Dfour e_I=\pfour[ \gensmoothfunction(\variables)(\pfour u,\pfour\ubar)]+\Chfour\cdot(\pfour u,\pfour\ubar)$, where $\gensmoothfunction$ is a smooth function (allowed to vary) of its argument. Then, schematically, we have:
		\begin{align}\label{EQ:ACCSchematic}
			\ACC=\gensmoothfunction(\pfour\vvariables)(\pfour^2u,\pfour^2{\ubar})+\gensmoothfunction(\pfour\vvariables,\pfour u,\pfour\ubar).
		\end{align}
		Then the desired result \eqref{ES:ACCGaussian} for $\ACC$ follows from \eqref{EQ:ACCSchematic} and \eqref{ES:uubar}.
		Finally, we prove \eqref{ES:Lapses}.
		We first consider transport equations (\ref{E:e_3e_3}) and (\ref{E:e_3lapse}). By the fundamental theorem of calculus, 
		\begin{align}
			\norm{\pfour^{\leq1}\ln\lapse,\pfour^{\leq1}\ln\uinvfolia}_{L^\infty(\hypersur_{u})}\lesssim\norm{\pfour^{\leq1}\ln\lapse,\pfour^{\leq1}\ln\uinvfolia}_{L^\infty(\hypersur_{u_0})}+\int_{u_0}^{u}\norm{\pfour^{\leq2}(v^i,\logDensity,\Ent),\pfour^{\leq1}\upeta}_{L^\infty(\hypersur_{u^\prime})}\diff u^\prime.
		\end{align}
		Similarly, considering (\ref{E:e_4e_43}) and (\ref{E:e_4lapse}), 
		\begin{align}
			\norm{\pfour^{\leq1}\ln\lapse,\pfour^{\leq1}\ln\invfolia}_{L^\infty(\hypersurbar_{\ubar})}\lesssim\norm{\pfour^{\leq1}\ln\lapse,\pfour^{\leq1}\ln\invfolia}_{L^\infty(\hypersurbar_{\ubar_0})}+\int_{\ubar_0}^{\ubar}\norm{\pfour^{\leq2}(v^i,\logDensity,\Ent),\pfour^{\leq1}\ueta}_{L^\infty(\hypersurbar_{\ubar^\prime})}\diff \ubar^\prime.
		\end{align}
		The desired estimates (\ref{ES:Lapses}) follow from \eqref{ES:uubar}-(\ref{ES:ACCGaussian}).
	\end{proof}
	\begin{corollary}[Ellipticity of the induced metric $\gsphere$]\label{CO:Ellipticity}
		The acoustical functions $u,\ubar$ are defined in Section \ref{SS:SUBSETSOFSPACETIME}. Let $u_0\leq u\leq u_1$, $\ubar_0\leq u\leq \ubar_1$, $u_0+\ubar_0\leq\uptau\leq u_1+\ubar_1$. $\Stau$ is defined in (\ref{D:Sigmatau}). Let $u_0\leq u^\prime\leq u_1$ and $\ubar_0\leq \ubar^\prime\leq \ubar_1$. Let $ d_u:=u_1-u_0,  d_{\ubar}:=\ubar_1-\ubar_0$. We note that $d_u,d_{\ubar}\lesssim 1$ and $\min\{d_u,d_{\ubar}\}\leq\delta$. Ricci coefficients $\ACC$ are defined in Definition \ref{D:Defconnectioncoeff}. Assume the bootstrap assumptions \eqref{BAfluid}-\eqref{BAgeometry3}, and that there exists $\mathring{a}\leq a_1\leq \frac{1}{4}\BAinfty$ such that the following bounds hold on $\hypersur_{u_0}$ and $\hypersurbar_{\ubar_0}$:
		\begin{align}\label{AS:Ellipticity1}
			a_1^{-1}\seabs{\xi}^2\leq\gsphere_{AB}\xi^A\xi^B\leq& a_1\seabs{\xi}^2,&
			\text{ for any }& \suubar \text{-tangent vectorfield }\xi.
		\end{align}
		Then, there exists $\mathring{a}\leq a\leq \BAinfty$ such that the following bounds hold:
		\begin{align}\label{ES:Ellipticity1}
			a^{-1}\seabs{\xi}^2\leq\gsphere_{AB}\xi^A\xi^B\leq& a\seabs{\xi}^2,&
			\text{ for any }& \suubar \text{-tangent vectorfield }\xi.
		\end{align}
	\end{corollary}
	\begin{proof}[Proof of Corollary \ref{CO:Ellipticity}]
		We only prove \eqref{ES:Ellipticity1} in the neighborhood of $\hypersur_{u_0}$ where $\updelta_u\leq \delta$. \eqref{ES:Ellipticity1} in the neighborhood of $\hypersurbar_{\ubar_0}$ can be proved in a similar fashion. Recall the definition of $\gsphere$ in \eqref{DE:gsphere}:
		\begin{align}\label{EQ:gsphere}
			\gsphere_{\alpha\beta}:=\gfour_{\alpha\beta}+\frac{1}{2}(e_4)_\alpha (e_3)_\beta+\frac{1}{2}(e_3)_\alpha (e_4)_\beta.
		\end{align}
		Using the fact that $e_3,e_4$ depend on $\vvariables,u,\ubar$ (as in Section \ref{SS:SUBSETSOFSPACETIME}) and $\gfour$ depends on $\vvariables$ (as in Definition \ref{D:acousticmetric}), we apply the fundamental theorem of calculus on $e_3$-derivative of \eqref{EQ:gsphere}, use the previously proved Proposition \ref{PR:BoundforGeometricQuantities} and smallness of $\updelta_u$ as in \eqref{ES:smalldelta}:
		\begin{align}\label{ES:gsphere}
			\abs{\gsphere_{\alpha\beta}\restriction_{S_{u,\ubar}}-\gsphere_{\alpha\beta}\restriction_{S_{u_0,\ubar}}}\leq 	\int_{0}^{u}\onenorm{\gensmoothfunction(\pfour^{\leq2}u,\pfour^{\leq2}\ubar,\pfour\vvariables,\ln\lapse,\ln\invfolia)}{\gsphere}{\infty}{S_{u^\prime,\ubar}}\diff u^\prime\leq \frac{1}{4}\BAinfty,
		\end{align}
		where $\gensmoothfunction$ is a linear combination of polynomials of order $\leq 5$ with respect to its arguments.
		Then the desired result \eqref{ES:Ellipticity1} follows from \eqref{ES:gsphere} and assumptions on $\hypersur_{u_0}$ \eqref{AS:Ellipticity1}.
	\end{proof}
	
	\subsection{Propagation of $C^{\infty}$ Regularity}\label{SS:SmoothnessoftheSolutions}
	While our a priori energy estimates close at finite regularity, for ``soft reasons,'' our approach to continuation relies on having 
	$C^{\infty}$ solutions that extend as $C^{\infty}$ functions to the boundary of their domain of definition.
	For this reason, in this section, we derive estimates for all of the solution's higher derivatives to all orders.

	More precisely, we claim that for any $j,k\geq0$, there exists a constant $\Lambda_{j,k}$ such that:
	\begin{subequations}\label{A:SmoothBound1}
		\begin{align}
			\twonorms{\left[\mathring{e_3}^j(\mathring{e_4},\iniangp)^{k}\right](\inilogDensity,\mathring{v}^i,\inient)}{\ubar}{2}{\gsphere}{2}{\hypersur_{0}}\leq& \Lambda_{j,k},\\
			\twonorms{\left[\mathring{e_4}^j(\mathring{e_3},\iniangp)^{k}\right](\inilogDensity,\mathring{v}^i,\inient)}{u}{2}{\gsphere}{2}{\hypersurbar_{0}}\leq& \Lambda_{j,k}.
		\end{align}
	\end{subequations}
	Moreover, for any $i\geq 0$, there exists $\Lambda_{i}$ such that:
	\begin{subequations}\label{A:SmoothBound2}
		\begin{align}
			\label{SmoothDatageometry1}\norm{\pfour^i(\mathring{v}^i,\inilogDensity,\inient)}_{L^\infty(\hypersur_{0}\bigcup\hypersurbar_{0})}\leq& \Lambda_{i},\\
			\label{SmoothDatageometry2}\norm{\pfour^i(\iniu,\iniubar),\pfour^i\ln\iniinvfolia,\pfour^i\ln\iniuinvfolia,\pfour^i\ln\inilapse,\pfour^i\iniACC}_{L^\infty(\hypersur_{0}\bigcup\hypersurbar_{0})}\leq& \Lambda_{i}.
		\end{align}
	\end{subequations}
	We note that $\Lambda_{j,k}$ and $\Lambda_{i}$ can be large and divergent in terms of $i,j,k$. 
	
	By the same approach\footnote{The only difference is replacing $\RBAinfty$ by $\norm{\pfour^k(v^i,\logDensity,\Ent),\pfour^{k-1}(\modivort,\DivGradEnt)}_{L^\infty(\region)}$ with higher order $k$ in Theorem~\ref{TH:UPTOTOPORDERENERGYESTIMATES} and Proposition \ref{PR:Linftybound}, which can be bounded by the same approach as in Proposition \ref{PR:Linftybound} and \eqref{SmoothDatageometry1}.} as in Section \ref{SS:APRIORIENERGYESTIMATESWITHISOTROPICREGULARITY} and an induction argument, we preserve the smoothness of the fluid variables. That is, letting $\updelta_u:=u_1-u_0, \updelta_{\ubar}:=\ubar_1-\ubar_0$, for various energies $\tilde{\KEnergy},\tilde{\bEnergy},\tilde{\Energy},\tilde{\uEnergy}$ defined in Definition \ref{D:ControllingEnergyQuantities} where $u_0+\ubar_0<\uptau\leq u_1+\ubar_1$, $1\leq k\leq \frac{N-1}{2}$ (see Section \ref{SS:SOBOLEVREGULARITYOFDATA} for $N$), there exist $\widetilde{\Lambda_{k}}$ and $\Delta_{i}$ (similar to $\Lambda_{i},\Lambda_{j,k}$, $\widetilde{\Lambda_{k}}$ and $\Delta_{i}$ can be large and divergent in terms of $i,k$) such that 
	\begin{align}\label{ES:SmoothL2}
		&\tilde{\KEnergy}^{(1)}_{(k)}(u_0,\ubar_0,\min\{u_0+\ubar_1,u_1+\ubar_0\}),\tilde{\KEnergy}^{(2)}_{(k)}(\min\{u_0+\ubar_1,u_1+\ubar_0\},\max\{u_0+\ubar_1,u_1+\ubar_0\}),\\
		\notag&\tilde{\KEnergy}^{(3)}_{(k)}(\max\{u_0+\ubar_1,u_1+\ubar_0\},u_1+\ubar_1),\tilde{\bEnergy}_{(k)}(u_0,\ubar_0,\uptau),\\
		\notag&\tilde{\Energy}_{(k)}(u_1,\ubar_0,\ubar_1),\tilde{\uEnergy}_{(k)}(\ubar_1,u_0,u_1)\\
		\notag&\lesssim \widetilde{\Lambda_{k}}.
	\end{align}
	Moreover, for $i\geq 0$, there exists $\Delta_i$ such that the following $L^\infty$-bound holds:
	\begin{align}\label{ES:SmoothLinfty}
		\norm{\pfour^i(v^i,\logDensity,\Ent)}_{L^\infty(\region)}\leq& \Delta_{i}.
	\end{align}
	We emphasize that \eqref{ES:SmoothL2}-\eqref{ES:SmoothLinfty} are needed in the (domain) extension argument of the existence of the solution, since we are going to utilize Rendall's argument as in Section \ref{S:localexistence} on a neighborhood of $S_{0,\ubar'}$ and $S_{u',\ubar}$. To run this argument, one requires solutions to be smooth functions on a portion of $\hypersur_{u'},\hypersurbar_{\ubar'}$. Nevertheless, the uniform bounds for the energy or $L^\infty$-norm are not needed in Section \ref{S:localexistence}. Hence, we proved a uniform bound under the regularity $N\leq 17$ as well as controlled the geometry including the change of variables map in Section \ref{SS:APRIORIENERGYESTIMATESWITHISOTROPICREGULARITY}, and proved a non-uniform bound for higher regularities in this section. Then, combining both pieces of the regularity results, we prove the local existence and uniqueness of the characteristic initial value problem and continuation criterion in Section \ref{S:Lastslice}.

	\section{The Main Theorem on Local Existence and Uniqueness and Continuation Criteria}\label{S:Lastslice}
	Our first main goal in this section is to prove Theorem~\ref{TH:localwellposedness}, which is the main existence and uniqueness result in the paper.
	Our second goal is to prove Theorem~\ref{TH:continuationcriteria}, which yields continuation criteria.
	As a precursor, we have to control the change of variables map $\embedgeo$ between geometric and Cartesian coordinates
	i.e., we give criteria that ensure $\embedgeo$ is a diffeomorphism in the regions under consideration.

	\subsection{A Quantitative Version of the Inverse Function Theorem}
	\label{SS:QUANTITATIVEVERSIONOFTHEINVERSEFUNCTION}
	To prove existence in a region whose size can be controlled by the data, 
	we need to control the change of variables map between geometric coordinates and Cartesian coordinates,
	in particular to show that it is a diffeomorphism in a region whose size can be controlled.
	To achieve this, we rely on a quantitative version of the inverse function theorem,
	which we prove as Lemma~\ref{L:INVERSEFUNCTIONTHEOREMWITHCONTROLOVERTHESIZEOFTHEREGION}.
	Before proving the lemma, we first introduce a simple distance function that we will use in its 
	statement and proof.
	
	\begin{definition}[A distance function on geometric coordinate space]
		\label{D:DISTANCEFUNCTIONONGEOMETRICCOORDINATESPACE}
		For pairs of points 
		$p := (u_p,\ubar_p,\uptheta_p), q := (u_q,\ubar_q,\uptheta_q) \in [0,\infty) \times [0,\infty) \times \mathbb{S}^2$ in geometric coordinate space,
		let $d$ be the corresponding distance function defined by:
		\begin{align} \label{E:DISTANCEFUNCTIONONGEOMETRICCOORDINATESPACE}
			d(p,q)
			& := |u_p - u_q|
			+
			|\ubar_p - \ubar_q|
			+
			d_{\mathring{\upgamma}}(\uptheta_p,\uptheta_q),
		\end{align}
		where $d_{\mathring{\upgamma}}(\uptheta_p,\uptheta_q)$ is the distance function on $\mathbb{S}^2$ induced by the round metric
		$\mathring{\upgamma}$ on $\mathbb{S}^2$.
	\end{definition}

	\begin{lemma}[Quantitative version of the inverse function theorem]
		\label{L:INVERSEFUNCTIONTHEOREMWITHCONTROLOVERTHESIZEOFTHEREGION}
		Let $\mathcal{K} := \mathcal{D} \times \mathbb{S}^2$ be a compact 
		subset of geometric coordinate space, where $\mathcal{D} \subset [0,\infty) \times [0,\infty)$
		is a compact subset of ``$(u,\ubar)''$-space.
		Assume that there are constants $M_1 > 0, M_2 > 0$ such
		that every pair of points $(u_1,\ubar_1), (u_2,\ubar_2) \in \mathcal{D}$ can be joined by
		a piecewise smooth curve $\ell$ in $\mathcal{D}$ consisting of a finite number of constant-$\ubar$
		line segments and/or constant-$u$ line segments such that the total length of these segments, denoted by $|\ell|$,
		satisfies the following estimates: 
		\begin{align} \label{E:QUASICONVEXITYSTRAIGHTSEGMENTBOUNDFORDOMAIN}
			M_1 \left( |u_2 - u_1| + |\ubar_2 - \ubar_1| \right)
			& \leq |\ell|
			\leq M_2 \left( |u_2 - u_1| + |\ubar_2 - \ubar_1| \right).
		\end{align}
		Let $F: \mathcal{K} \rightarrow \mathbb{R}^{1+3}$
		be a $C^2$ function. Let $C > 0$ be a constant
		such that the following bounds hold for all $p \in \mathcal{K}$, where
		$\mathrm{d} F$ denotes the Jacobian matrix of $F$ and
		$\mathrm{d}^2 F$ denotes the array of all second order partial derivatives of $F$
		(relative to a local chart on $\mathbb{S}^2$ for angular derivatives):
		\begin{align} \label{E:UNIFORMLOWEROBUNDFORJACOBIANDETCHOV} 
			\frac{1}{|\mbox{\upshape det }  (\mathrm{d} F(p))|} 
			& \leq C,
			\\
			|F(p)|,
			\,
			|\mathrm{d} F(p)|,
			\,
			|\mathrm{d}^2 F(p)|
			& \leq C.
			\label{E:UNIFORMBOUNDSFORCHOV}
		\end{align} 
		Given $\varepsilon > 0$ and a point $\widetilde{p} \in \mathcal{D}$,
		let $B_{\widetilde{p}}(\varepsilon)$ denote the intersection of $\mathcal{D}$ with the
		$d$-ball of radius $\varepsilon$ centered at $\widetilde{p}$,
		where $d$ is the distance function from 
		Def.\,\ref{D:DISTANCEFUNCTIONONGEOMETRICCOORDINATESPACE}.
		
		Then there is a number $\varepsilon_0 > 0$
		and constants $K_1, K_2 > 0$, depending only on $M_1$, $M_2$, and $C$,
		such that for every $\widetilde{p} \in \mathcal{D}$,
		$F$ is a diffeomorphism on $B_{\widetilde{p}}(\varepsilon_0)$, and
		the following estimates hold:
		\begin{align} \label{E:CHOVMAPHOMEOMORPHISMESTIMATEINSMALLBALL}
			K_1 d(p,q)
			&
			\leq
			|F(p) - F(q)|
			\leq K_2 d(p,q),
			&& \mbox{for every } p,q \in B_{\widetilde{p}}(\varepsilon_0).
		\end{align}
		
	\end{lemma}
	
	\begin{proof}
		The inverse function theorem guarantees that $F$ is a local diffeomorphism, so it only remains to prove
		\eqref{E:CHOVMAPHOMEOMORPHISMESTIMATEINSMALLBALL}. 
		Let $q \in \mathcal{K}$. 
		Note that the function $p \rightarrow [\mathrm{d} F(q)]^{-1} F(p)$ has Jacobian matrix equal to the identity
		when $p=q$.
		By assumption, for pairs of points $p,q \in \mathcal{K}$,
		we can connect their ``$(u,\ubar)$'' coordinates with a curve $\ell$ of the type from the assumptions of the lemma.
		Hence, repeatedly 
		using the mean value theorem along the straight line segments of $\ell$,
		the assumption \eqref{E:QUASICONVEXITYSTRAIGHTSEGMENTBOUNDFORDOMAIN},
		and the bounds
		\eqref{E:UNIFORMLOWEROBUNDFORJACOBIANDETCHOV}--\eqref{E:UNIFORMBOUNDSFORCHOV}
		we deduce the following estimate for the function:
		\begin{align} \label{E:MEANVALUETHEOREMESTIMATE}
			\left|
			[\mathrm{d} F(q)]^{-1} F(q)
			-
			[\mathrm{d} F(q)]^{-1} F(p)
			-
			(p - q)
			\right|
			& \leq C (d(p,q))^2,
		\end{align}
		where we have expressed the terms on LHS~\eqref{E:MEANVALUETHEOREMESTIMATE} with respect to
		a coordinate chart on $\mathcal{D} \times \mathbb{S}^2$ (in particular, we are viewing $p - q$ as a point in $\mathbb{R}^4$, 
		and $|p - q| \approx d(p,q)$ in the chart). From \eqref{E:UNIFORMBOUNDSFORCHOV} and \eqref{E:MEANVALUETHEOREMESTIMATE},
		we further deduce that:
		\begin{align} \label{E:SECONDMEANVALUETHEOREMESTIMATE}
			\left|
			F(q)
			-
			F(p)
			-
			dF(q)\cdot(p - q)
			\right|
			& \leq C (d(p,q))^2.
		\end{align}
		In addition, the bounds
		\eqref{E:UNIFORMLOWEROBUNDFORJACOBIANDETCHOV}--\eqref{E:UNIFORMBOUNDSFORCHOV}
		imply that there are constants $C_1 > 0$, $C_2 > 0$, independent of $p$ and $q$, such that when $d(p,q)$ is sufficiently small, we have:
		\begin{align} \label{E:DFINVERSETIMESVECTORISAPPROXIMATELYSAMELENGTHASVECTOR}
			C_1 d(p,q)
			& \leq
			\left|
			[\mathrm{d} F(q)] \cdot (p - q)
			\right|
			\leq
			C_2 d(p,q).
		\end{align}
		Finally, from \eqref{E:SECONDMEANVALUETHEOREMESTIMATE},
		\eqref{E:DFINVERSETIMESVECTORISAPPROXIMATELYSAMELENGTHASVECTOR},
		and \eqref{E:UNIFORMBOUNDSFORCHOV},
		we conclude that if $d(p,q)$ is sufficiently small,
		then $F(q) \neq F(p)$, and there are constants $K_1 > 0$ and $K_2 > 0$
		such that \eqref{E:CHOVMAPHOMEOMORPHISMESTIMATEINSMALLBALL} holds.
		
	\end{proof}
	
	\begin{proposition}[Quantitative global diffeomorphism estimate for the change of variables map]
		\label{E:QUANTITATIVEDIFFEOMORPHISMESTIMATEFORCHOVMAP}
		Let $\embedgeo$ be the change of variables map from geometric coordinates to Cartesian coordinates defined in 
		\eqref{E:CHOVFROMGEOTOCARTESIAN}.
		Under the initial data and bootstrap assumptions, 
		there exist constants $C_1 > 0$ and $C_2 > 0$
		such that if $\varsigma > 0$ is sufficiently small, then
		$\embedgeo$ is a global diffeomorphism on the bootstrap region, 
		and for all pairs of points $p,q$ in the bootstrap region, 
		we have the following estimates, where $d(\cdot,\cdot)$ is the distance function defined in \eqref{E:DISTANCEFUNCTIONONGEOMETRICCOORDINATESPACE}:
		\begin{align} \label{E:CHOVMAPHOMEOMORPHISMESTIMATE}
			C_1 d(p,q)
			&
			\leq
			|\embedgeo(p) - \embedgeo(q)|
			\leq C_2 d(p,q).
		\end{align}
	\end{proposition}
	
	\begin{proof}
		For $F=\embedgeo$ as in Lemma \ref{L:INVERSEFUNCTIONTHEOREMWITHCONTROLOVERTHESIZEOFTHEREGION},  by bootstrap assumptions \eqref{BAfluid}-\eqref{BAgeometry3}, and $\det\gt=\speed^{-6}$ with $\speed\geq\speed_0>0$ as in Definition \ref{DE:SoundSpeed} and $\speed(\logDensity,\Ent)\leq \BAinfty$ by \eqref{BAfluid}, the assumption \eqref{E:UNIFORMLOWEROBUNDFORJACOBIANDETCHOV} holds since we have \eqref{EQ:determinant0}:
		\begin{align}
			\det(\diff \embedgeo)=4\invfolia\uinvfolia\lapse^3(\det\gt)^{-\frac{1}{2}}(\det\gsphere)^{\frac{1}{2}}>0.
		\end{align}
		Moreover, by \eqref{EQ:dPhi} and bootstrap assumptions \eqref{BAfluid}-\eqref{BAgeometry3}, \eqref{E:UNIFORMBOUNDSFORCHOV} also holds for $F=\embedgeo$ since 
		\begin{align}\label{ES:d2Phi}
			\abs{\diff^2 \embedgeo(p)}\leq \abs{\gensmoothfunction(\pfour^{\leq2}(u,\ubar),\pfour^{\leq1}(\ln\lapse,\ln\invfolia,\ln\uinvfolia),\pfour^{\leq 1}\vvariables)}\leq \gensmoothfunction(\BAinfty),
		\end{align}
		where $\gensmoothfunction$ are finite-order polynomials with respect to their arguments. $\abs{\diff F(p)},\abs{F(p)}$ are bounded by the fact that $u,\ubar$ are bounded and \eqref{ES:d2Phi}. 
		
		Therefore, the bootstrap assumptions imply that
		the bootstrap domain and
		$\embedgeo$ satisfy the assumptions and conclusions of
		Lemma~\ref{L:INVERSEFUNCTIONTHEOREMWITHCONTROLOVERTHESIZEOFTHEREGION}.
		In particular, $\embedgeo$ is a local diffeomorphism.
		It remains for us to prove the bound \eqref{E:CHOVMAPHOMEOMORPHISMESTIMATE}
		when $\varsigma > 0$ is small enough,
		which will in particular show that $\embedgeo$ is globally injective on the bootstrap region.
		By Lemma~\ref{L:INVERSEFUNCTIONTHEOREMWITHCONTROLOVERTHESIZEOFTHEREGION}, there is an $\varepsilon_0 > 0$
		such that \eqref{E:CHOVMAPHOMEOMORPHISMESTIMATE} holds for any pair of points $p,q$
		in the bootstrap region such that $d(p,q) < \varepsilon_0$.
		
		It remains for us to consider the case $d(p,q) > \varepsilon_0$.
		To proceed, we note that our assumptions on the initial embeddings stated in Sect.\,\ref{SS:initialembedding2}
		imply that $\embedgeo \restriction_{\inihypersur \bigcup \inihypersurbar}$ is injective. By compactness,
		there are constants $C_1', C_2' > 0$ such that
		for every 
		$\mathring{p}, \mathring{q} \in \inihypersur \bigcup \inihypersurbar$,
		we have:
		\begin{align} \label{E:CHOVMAPHOMEOMORPHISMESTIMATEFORDATA}
			C_1' d(\mathring{p},\mathring{q})
			&
			\leq
			|\embedgeo(\mathring{p}) - \embedgeo(\mathring{q})|
			\leq C_2' d(\mathring{p},\mathring{q}).
		\end{align}
		
		Next, we note that we can connect $p$ to a point
		$\mathring{p} \in \inihypersur \bigcup \inihypersurbar$
		by a line segment of constant $u$ or constant $\ubar$ whose geometric coordinate
		length is $\leq \varsigma$ such that the segment is contained in the bootstrap region.
		We can similarly connect $q$ to a point $\mathring{q} \in \inihypersur \bigcup \inihypersurbar$
		by the same kind of line segment.
		Since $\varsigma$ is the width of the bootstrap region, if $\varsigma$ 
		is sufficiently small compared to $\varepsilon_0$, 
		then the following inequalities hold:
		\begin{align} \label{E:POINTSNEARPOINTSONDATAHYPERSURFACE}
			d(p,\mathring{p})
			& \leq \varsigma,
			&
			d(q,\mathring{p})
			& \leq \varsigma,
			\\
			d(\mathring{p},\mathring{q})
			& \geq \frac{d(p,q)}{2}.	
			\label{E:TRIANGLEINEQUALITYFORPOINTSONINITIALHYPERSURFACEVSGENERALPOINTS}
		\end{align}
		From
		\eqref{E:CHOVMAPHOMEOMORPHISMESTIMATEINSMALLBALL} (with $\embedgeo$ in the role of $F$),
		\eqref{E:CHOVMAPHOMEOMORPHISMESTIMATEFORDATA},
		\eqref{E:POINTSNEARPOINTSONDATAHYPERSURFACE},
		\eqref{E:TRIANGLEINEQUALITYFORPOINTSONINITIALHYPERSURFACEVSGENERALPOINTS},
		and the triangle inequality, we deduce that:
		\begin{align} \label{E:KEYLOWERBOUNDFORDISTANCEBETWEENTWOIMAGEPOINTS}
			\begin{split}
				|\embedgeo(p) - \embedgeo(q)|
				& \geq |\embedgeo(\mathring{p}) - \embedgeo(\mathring{q})|
				-
				|\embedgeo(p) - \embedgeo(\mathring{p})|
				-
				|\embedgeo(q) - \embedgeo(\mathring{q})|
				\\
				& \geq C_1' d(\mathring{p},\mathring{q})
				- C \varsigma
				\geq \frac{C_1'}{2} 
				d(p,q)
				- 
				C \varsigma
				\\
				& \geq \frac{C_1'}{4} d(p,q),
			\end{split}
		\end{align}
		where to obtain the last inequality, we have assumed that $\varsigma$ is small enough such that
		$C \varsigma \leq \frac{C_1'}{4} \varepsilon_0$ and used our assumption $d(p,q) > \varepsilon_0$.
		From \eqref{E:CHOVMAPHOMEOMORPHISMESTIMATEINSMALLBALL} 
		and
		\eqref{E:KEYLOWERBOUNDFORDISTANCEBETWEENTWOIMAGEPOINTS},
		we conclude the first inequality in \eqref{E:CHOVMAPHOMEOMORPHISMESTIMATE}.
		The second inequality in \eqref{E:CHOVMAPHOMEOMORPHISMESTIMATE} can be proved using a similar argument.	
		
	\end{proof}
	
	\subsection{The Main Theorem on Local Existence and Uniqueness}
	\label{SS:MAINTHEOREMLOCALEXISTNCEANDUNIQUENESS}
	We are now ready to state and prove our main result on local existence and uniqueness.
	
	\begin{figure}
		\centering
		\includegraphics[width=12cm]{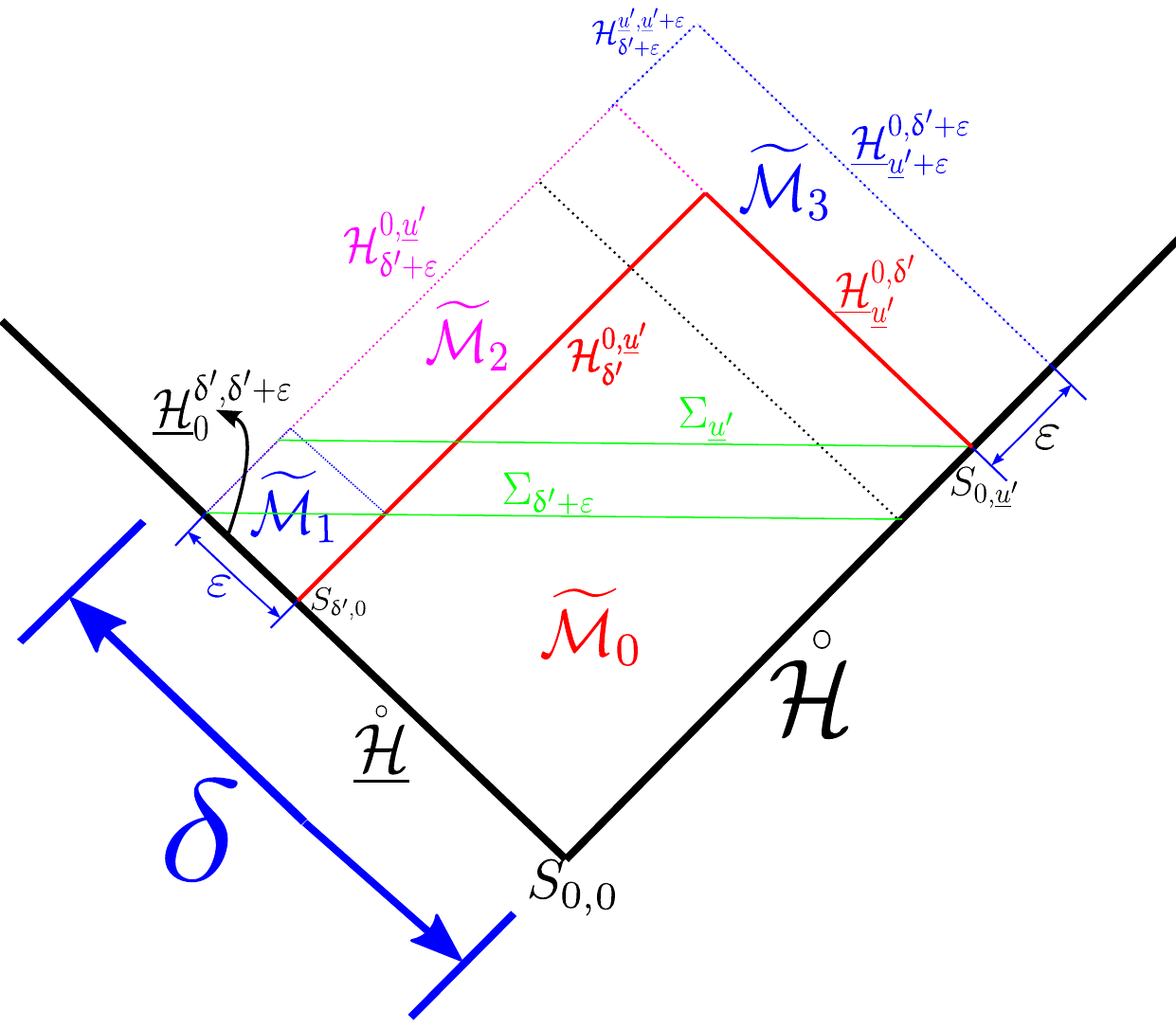}
		\caption{Continuation of the solutions in a neighborhood of $\hypersur_{0}$}
		\label{Figure3}
	\end{figure}

	\begin{theorem}[Main Theorem: Local existence and uniqueness in a neighborhood of the initial null hypersurfaces]\label{TH:localwellposedness}
		Consider $C^{\infty}$ fluid data for the compressible Euler equations \eqref{E:Brho}--\eqref{E:Bs}
		that are prescribed on the union 
		$\inihypersur \cup \inihypersurbar$,
		where $\inihypersur$ and $\inihypersurbar$ are $C^{\infty}$ characteristic hypersurface portions 
		that intersect transversally in a set $S_{0,0}$ that is diffeomorphic to $\mathbb{S}^2$,
		as is described in Sect.\,\ref{S:CID}.
		Assume that the Euler constraint equation \eqref{E:constraint1} is satisfied on
		$\inihypersur$, that the Euler constraint equation \eqref{E:constraint2} is satisfied on $\inihypersurbar$,
		and that the data satisfy the corner compatibility conditions stated in Def.\,\ref{D:CORNERCOMPATIBILITYCONDITIONSONDATA}.
		Recall that in Prop.\,\ref{P:DETERMININGTRANSVERALDERIVATIVESANDSATISFYINGEQUATIONSONINITIALNULLHYPERSURFACES}, 
		we proved that the solution's tangential and transversal derivatives to all orders are finite on
		the sub-portions $\inihypersur([0,\ubarfinite]) \cup \inihypersurbar([0,\ufinite])$.
		Let $N \geq 17$ be an integer, 
		and assume that 
		on $\inihypersur([0,\ubarfinite]) \cup \inihypersurbar([0,\ufinite])$
		the initial data satisfy the norm bounds stated 
		in \eqref{A:H}--\eqref{A:Boundness}.
		Let $\updelta > 0$ be the number from a priori $L^{\infty}$ estimates of 
		Props.\,\ref{PR:Linftybound} and~\ref{PR:BoundforGeometricQuantities}.
		In particular, these propositions and
		Theorem~\ref{TH:UPTOTOPORDERENERGYESTIMATES} show that any $C^{\infty}$ solution
		on the compact region
		$\ExistenceRegionthreearg{\updelta}{\ufinite}{\ubarfinite}$ defined by:
		\begin{align} 
			\begin{split} \label{E:EXISTENCENEIGHBORHOODOFNULLHYPERSURFACESINMAINTHEOREM}
				\ExistenceRegionthreearg{\updelta}{\ufinite}{\ubarfinite}
				& := 	
				\left\lbrace
				(t,x^1,x^2,x^3) \in \mathbb{R}^{1+3}
				\ | \
				0 \leq \ubar(t,x^1,x^2,x^3)  \leq \ubarfinite,
				0 \leq u(t,x^1,x^2,x^3)  \leq \updelta
				\right\rbrace
				\\
				& 
				\bigcup
				\left\lbrace
				(t,x^1,x^2,x^3) \in \mathbb{R}^{1+3}
				\ | \
				0 \leq \ubar(t,x^1,x^2,x^3) \leq \updelta,
				0 \leq u(t,x^1,x^2,x^3) \leq \ufinite
				\right\rbrace.
			\end{split}
		\end{align}
		must satisfy the a priori estimates \eqref{ES:EnergyRegularity1}
		on $\ExistenceRegionthreearg{\updelta}{\ufinite}{\ubarfinite}$,
		as long as the map $\embedgeo$ defined by:
		\begin{align} \label{E:CHOVFROMGEOTOCARTESIANMAINLWPTHEOREM}
			\embedgeo(u,\ubar,\sphere^1,\sphere^2)
			& := (t,x^1,x^2,x^3)
		\end{align}
		(i.e., the map from geometric coordinates to Cartesian coordinates)
		is a $C^{\infty}$ diffeomorphism from the domain
		$[0,\updelta] \times [0,\ubarfinite] \times \mathbb{S}^2 \bigcup [0,\ufinite] \times [0,\updelta] \times \mathbb{S}^2$
		onto the image set $\ExistenceRegionthreearg{\updelta}{\ufinite}{\ubarfinite}$.
		
		Then the following hold on the future-neighborhood
		$\ExistenceRegionthreearg{\updelta}{\ufinite}{\ubarfinite}$
		of $\inihypersur([0,\ubarfinite]) \cup \inihypersurbar([0,\ufinite])$:
		\begin{itemize}
			\item There exists a unique $C^{\infty}$ solution $\vvariables: =(\logDensity,v^1,v^2,v^3,\Ent)$
			to the compressible Euler equations \eqref{E:Brho}--\eqref{E:Bs} 
			whose restriction to $\inihypersur([0,\ubarfinite]) \cup \inihypersurbar([0,\ufinite])$
			coincides with the given initial data on $\inihypersur$ and $\inihypersurbar$.
			\item There exist unique $C^{\infty}$ solutions $u,\ubar$ to the initial value problems
			\eqref{E:EIKONALFORENDALLREGION}--\eqref{E:INITIALCONDITIONSEIKONALFORENDALLREGION}
			with non-vanishing, non-parallel gradients.
			\item Within $\ExistenceRegionthreearg{\updelta}{\ufinite}{\ubarfinite}$, the level sets of $u$ and $\ubar$ intersect transversally 
			in two-dimensional $\gfour$-spacelike submanifolds 
			$S_{u',\ubar'} := \lbrace u = u' \rbrace \cap \lbrace \ubar = \ubar' \rbrace$,
			all of which are diffeomorphic to $\mathbb{S}^2$.
			\item There exist unique $C^{\infty}$ solutions $(\sphere^1,\sphere^2)$ 
			to the transport equation initial value problem \eqref{E:TRANSPORTINGTHEANGULARCOORDINATEFUNCTIONS}
			that are local coordinates on the $S_{u',\ubar'}$.
			\item $\embedgeo$ is a $C^{\infty}$ diffeomorphism from the domain
			$[0,\updelta] \times [0,\ubarfinite] \times \mathbb{S}^2 \bigcup [0,\ufinite] \times [0,\updelta] \times \mathbb{S}^2$
			onto the image set $\ExistenceRegionthreearg{\updelta}{\ufinite}{\ubarfinite}$.
			\item The solution satisfies the energy estimates stated in Theorem~\ref{TH:UPTOTOPORDERENERGYESTIMATES}
			(where $u_0 = \ubar_0 := 0$ in the theorem).
		\end{itemize}
		
	\end{theorem}

	\begin{proof}
		We will prove existence in detail only for a region of type
		\begin{align*}
			\left\lbrace
		(t,x^1,x^2,x^3) \in \mathbb{R}^{1+3}
		\ | \
		0 \leq \ubar(t,x^1,x^2,x^3)  \leq \ubarfinite,
		0 \leq u(t,x^1,x^2,x^3)  \leq \updelta
		\right\rbrace;
		\end{align*} 
	existence in the region of type 
	\begin{align*}
		\left\lbrace
		(t,x^1,x^2,x^3) \in \mathbb{R}^{1+3}
		\ | \
		0 \leq \ubar(t,x^1,x^2,x^3) \leq \updelta,
		0 \leq u(t,x^1,x^2,x^3) \leq \ufinite
		\right\rbrace
	\end{align*} can be proved using similar arguments. 
		
		The proof of the theorem is based on continuing the fluid solution and acoustic geometry into several different regions in geometric coordinate space, depicted in Fig.\,\ref{Figure3}.
		In all of the steps, we can use Prop.\,\ref{E:QUANTITATIVEDIFFEOMORPHISMESTIMATEFORCHOVMAP} to ensure that the change of variables map
		also extends as a $C^{\infty}$ diffeomorphism; we will sometimes make silent use of this fact.
		
		First, by Theorem~\ref{T:Rendall'sTheorem}, there exist numbers 
		$u',\ubar',\updelta'$ satisfying
		$0 < u' \leq \ufinite$, $0 < \ubar' \leq \ubarfinite$,
		and $0 < \updelta' \leq \updelta$
		such that the conclusions of the theorem 
		hold on $\ExistenceRegionthreearg{\updelta'}{u'}{\ubar'}=:\widetilde{\mathcal{M}}_0$. 
		We will first prove that
		there is a number $\varepsilon > 0$ such that the conclusions of Theorem~\ref{TH:localwellposedness}
		hold on $\ExistenceRegionthreearg{\updelta'+ \varepsilon}{u' }{\ubar'}$.
		
		\medskip
		
		\noindent \textbf{(First continuation)}.
		 The results of Section~\ref{SS:SmoothnessoftheSolutions}
		guarantee that the solution is $C^{\infty}$ on the compact region
		$\ExistenceRegionthreearg{\updelta'}{u' }{\ubar'}$. We can therefore 
		use Theorem~\ref{T:Rendall'sTheorem} to deduce that there is an $\varepsilon > 0$ such that
		the solution can uniquely be extended as a $C^{\infty}$ solution to a future-neighborhood 
		$\widetilde{\mathcal{M}}_1:=[\updelta',\updelta'+\varepsilon]\times[0,\varepsilon]\times\Stwo$ of the sphere
		$S_{u',0}$. 
		\medskip
		
		\noindent \textbf{(Second continuation)}. 
		Let $\Sigma_{\timefunction'}$ (defined as in \eqref{D:Sigmatau}) denote the intersection of 
		$\ExistenceRegionthreearg{\updelta'+ \varepsilon}{u' }{\ubar'}$ with
		the $\timefunction'$-level set of
		the geometric time function from \eqref{E:GEOMETRICTIMEFUNCTION}.
		The estimates \eqref{ES:EnergyRegularity1} imply that the fluid and acoustic geometry extend as
		classical solutions to $\overline{\widetilde{\mathcal{M}}_1\cup\widetilde{\mathcal{M}}_0}$, the closure of $\widetilde{\mathcal{M}}_1\cup\widetilde{\mathcal{M}}_0$, which satisfies $\Sigma_{\timefunction'}\subset\widetilde{\mathcal{M}}_1\cup\widetilde{\mathcal{M}}_0$ for all $0\leq \timefunction'\leq \updelta'+\varepsilon$. Note that $\Sigma_{\updelta' + \varepsilon}$ is a spacelike hypersurface that
		is contained in the region $\overline{\widetilde{\mathcal{M}}_1\cup\widetilde{\mathcal{M}}_0}$ where the solution is known to be $C^{\infty}$.
		Hence, using standard well-posedness theory for the Cauchy problem with data
		on $\Sigma_{\updelta' + \varepsilon}$,
		the a priori estimates \eqref{ES:EnergyRegularity1} and those of Sect.\,\ref{SS:SmoothnessoftheSolutions}
		(note that our a priori estimates include estimates on spacelike surfaces $\Sigma_{\timefunction'}$),
		and standard continuation criteria for solutions on regions foliated by spacelike hypersurfaces,
		we deduce that
		the $C^{\infty}$ solution extends uniquely to the entire future domain of dependence of
		$\Sigma_{\updelta' + \varepsilon}$, which is a triangular-shaped region that includes the region
		$\widetilde{\mathcal{M}}_2:=[\updelta',\updelta'+\varepsilon]\times[\varepsilon,\ubar']\times\Stwo$
		in Fig.\,\ref{Figure3}.
\medskip		

\noindent \textbf{(Third continuation)}. Using arguments similar to the ones we used to extend the solution to 
$\widetilde{\mathcal{M}}_1\cup\widetilde{\mathcal{M}}_2$, we can, after perhaps shrinking $\varepsilon$, 
also uniquely extend the $C^{\infty}$ solution to the region
$\widetilde{\mathcal{M}}_3:=[0,\updelta'+\varepsilon]\times[\ubar',\ubar'+\varepsilon]\times\Stwo$ in Fig.\,\ref{Figure3}.
	
We have therefore extended the solution to the larger region $\ExistenceRegionthreearg{\updelta'+ \varepsilon}{u' }{\ubar'}$.
Let $M$ be the supremum of all numbers $\updelta' > 0$ such that a unique $C^{\infty}$ solution exists on 
$\ExistenceRegionthreearg{\updelta'}{u' }{\ubar'}$. We claim that $M \geq \updelta$, where
$\updelta > 0$ be the number from Props.\,\ref{PR:Linftybound} and~\ref{PR:BoundforGeometricQuantities}.
$\ExistenceRegionthreearg{\updelta'+ \varepsilon}{u' }{\ubar'}$. Assume for the sake of deriving a contradiction that
$M < \updelta$. The a priori estimates \eqref{ES:EnergyRegularity1} and those of Sect.\,\ref{SS:SmoothnessoftheSolutions},
as well as continuity,
imply that a unique $C^{\infty}$ solution exists and satisfies the estimates on the compact region
$\ExistenceRegionthreearg{M}{u' }{\ubar'}$. The arguments given above imply that we can uniquely extend the $C^{\infty}$ solution
to a larger region $\ExistenceRegionthreearg{M + \varepsilon}{u' }{\ubar'}$
on which the a priori estimates hold. This contradicts the definition of $M$.

\end{proof}


	\subsection{Continuation Criteria} \label{SS:continuationcriterion}
	In the next theorem, we show that under
	suitable assumptions on the fluid solution and the null hypersurfaces, 
	the region of existence can be extended by a small amount, i.e., we provide continuation criteria.
	
	\begin{theorem}[Continuation criteria] \label{TH:continuationcriteria}
		Assume the hypotheses on the initial data from the assumptions of Theorem~\ref{TH:localwellposedness},
		and let $\ufinite > 0$ and $\ubarfinite > 0$ be the (data-dependent) numbers featured in those assumptions.
		Let $u',\ubar',\updelta_1,\updelta_2$ be any real numbers satisfying
		$0 \leq u' \leq \ufinite$, $0 \leq \ubar' \leq \ubarfinite$,
		and $0 < \updelta_1, \updelta_2$.
		Assume that the following hold:
			\begin{enumerate}
			\item The change of variables map $\embedgeo$ 
			from Def.\,\ref{DE:changofvariables} extends to a global $C^{\infty}$-diffeomorphism from the compact subset 
			$([0,u'] \times [0,\updelta_1] \times \mathbb{S}^2) \bigcup  ([0,\updelta_2] \times [0,\ubar'] \times \mathbb{S}^2)$ 
			of geometric coordinate space
			onto 
			its compact image 
			\begin{align} \label{E:CLOSEDOFFREGION}
			\mathcal{R}_{u',\updelta_1,\ubar',\updelta_2}
			& := 
			\embedgeo\left(([0,u'] \times [0,\updelta_1] \times \mathbb{S}^2) \bigcup  
			([0,\updelta_2] \times [0,\ubar'] \times \mathbb{S}^2) \right)
		\end{align}
		in Cartesian coordinate space.
		\item The fluid solution $\vvariables:=(\logDensity,v^1,v^2,v^3,\Ent)$ exists and is $C^{\infty}$ 
			on the ``open-at-the-top'' subset 
			\begin{align} \label{E:OPENATTOPREGION}
			\widehat{\mathcal{R}}_{u',\updelta_1,\ubar',\updelta_2} 
			& := 
			\embedgeo\left(([0,u'] \times [0,\updelta_1) \times \mathbb{S}^2) \bigcup  
			([0,\updelta_2) \times [0,\ubarfinite_2] \times \mathbb{S}^2) \right)
		\end{align} 
		of $\mathcal{R}_{u',\updelta_1,\ubar',\updelta_2}$.
			\item The eikonal functions $u,\ubar$ are $C^{\infty}$ solutions to the eikonal equation initial value problems
			\eqref{E:EIKONALFORENDALLREGION}--\eqref{E:INITIALCONDITIONSEIKONALFORENDALLREGION}
			on $\widehat{\mathcal{R}}_{u',\updelta_1,\ubar',\updelta_2}$.
			\item  The geometric angular coordinates $(\sphere^1,\sphere^2)$ are $C^{\infty}$ solutions 
			to the transport equation initial value problem \eqref{E:TRANSPORTINGTHEANGULARCOORDINATEFUNCTIONS}
			on $\widehat{\mathcal{R}}_{u',\updelta_1,\ubar',\updelta_2}$.
			\item For all integers $k,i \geq 0$, there exist
				numbers $\widetilde{\Lambda}_{k}$ and $\Delta_i > 0$ such that the estimates \eqref{ES:SmoothL2}--\eqref{ES:SmoothLinfty} hold
				on $\widehat{\mathcal{R}}_{u',\updelta_1,\ubar',\updelta_2}$.
		\end{enumerate}
		Then there exists an $\varepsilon > 0$ such that the fluid, eikonal functions, and geometric angular coordinates
		uniquely extend as $C^{\infty}$ solutions to the larger compact region
		$\mathcal{R}_{u',\updelta_1 + \varepsilon,\ubar',\updelta_2 + \varepsilon} 
		$
		such that $\embedgeo$ is a
		$C^{\infty}$-diffeomorphism on 
		$([0,u'] \times [0,\updelta_1+\varepsilon] \times \mathbb{S}^2) \bigcup  ([0,\updelta_2+\varepsilon] \times [0,\ubar'] \times \mathbb{S}^2)$.
	\end{theorem}
	
	\begin{proof}
		The assumptions of the theorem
		imply that the fluid and acoustic geometry 
		have derivatives to all orders that are each uniformly bounded on the pre-compact region $\widehat{\mathcal{R}}_{u',\updelta_1,\ubar',\updelta_2}$.
		Hence, all these quantities extend as
		$C^{\infty}$ solutions to the compact region 
		$\mathcal{R}_{u',\updelta_1,\ubar',\updelta_2}$.
		The data induced on the ``boundary'' null hypersurface portions
		$\hypersur_{\updelta_1}^{0,\ubar'}$
		and $\hypersurbar_{\updelta_2}^{0,u'}$
		satisfy the hypotheses of Theorem~\ref{TH:localwellposedness}.
		We can therefore apply Theorem~\ref{TH:localwellposedness}
		to uniquely extend the solution
		to $\mathcal{R}_{u',\updelta_1 + \varepsilon,\ubar',\updelta_2 + \varepsilon} 
		$ for some $\varepsilon > 0$.
\end{proof}

	\appendix
	\section*{Appendices}
	\addcontentsline{toc}{section}{Appendices}
	\renewcommand{\thesubsection}{\Alph{subsection}}
	\subsection{Scaling of the Null Vectorfields}
	For convenience, we list our null vectorfields with different scaling as follows:
	\begin{subequations}
		\begin{align}
			\Lgeo^\alpha=&-(\gfour^{-1})^{\alpha\beta}\p_\beta u,&
			\uLgeo^\alpha=&-(\gfour^{-1})^{\alpha\beta}\p_\beta \ubar,\\
			\Lunit^\alpha=&-\invfolia(\gfour^{-1})^{\alpha\beta}\p_\beta u,&
			\uLunit^\alpha=&-\uinvfolia(\gfour^{-1})^{\alpha\beta}\p_\beta \ubar,\\
			e_4^\alpha=&-2\invfolia\lapse(\gfour^{-1})^{\alpha\beta}\p_\beta u,&
			e_3^\alpha=&-2\uinvfolia\lapse(\gfour^{-1})^{\alpha\beta}\p_\beta \ubar,\\
			\Lrenormal^\alpha=&-2\invfolia\uinvfolia\lapse^2(\gfour^{-1})^{\alpha\beta}\p_\beta u,&
			\uLrenormal^\alpha=&-2\invfolia\uinvfolia\lapse^2(\gfour^{-1})^{\alpha\beta}\p_\beta \ubar.
		\end{align}
	\end{subequations}
	\subsection{Notations}\label{AP:notations}
	In this appendix, we gather some of the notation that we use throughout the article.
	\begin{center}

		\begin{tabular}{|c|c|} 
			\hline
			Symbol & Reference  \\ 
			\hline\hline
			$\lesssim,\p,\pfour,\antisymmetic_{\alpha\beta\gamma\delta}$ & Section \ref{Notation} \\
			\hline
			$\St$ & Section \ref{Notation} \\
			\hline
			$v^i,\density,\Ent,\pressure$, $\speed$ & Section \ref{D:speedofsound} \\
			\hline
			$\vort,\logDensity,\backgroundDensity,\vortrenormalized,\gradEnt$ & Section \ref{SS:morefluid} \\
			\hline
			$\materialderivative,\delta^{ab},\vvariables$ & Section \ref{SS:firstorderequations} \\
			\hline
			$\modivort,\DivGradEnt$ & Definition \ref{DE:CD} \\
			\hline
			$\gfour,\gfour^{-1}$ & Definition \ref{D:acousticmetric} \\
			\hline
			$\boxg,\Dfour$ & Definition \ref{D:LeviCivita} \\
			\hline
			$\dfour{\alpha}{\beta}$ & Definition \ref{DE:dfour} \\
			\hline
			$\Chfour$ & Definition \ref{D:Christoffelsymbols} \\
			\hline
			$\mathfrak{Q}_{(v)}^i, \mathfrak{Q}_{(\logDensity)}, \mathfrak{Q}^i_{(\modivort)},\mathfrak{Q}_{(\DivGradEnt)}$ & Theorem~\ref{T:geometricformulation} \\
			\hline
			$\mathfrak{L}_{(v)}^i,\mathfrak{L}_{(\logDensity)},\mathfrak{L}_{(\Ent)},\mathfrak{L}_{(\vortrenormalized)}^i,\mathfrak{L}_{(\dive\vortrenormalized)},\mathfrak{L}^i_{(\modivort)}$ & Theorem~\ref{T:geometricformulation} \\
			\hline 
			$\gt,\gtau,\Sigmatproject,\Sigmatauproject,\sphereproject,\taup,\angp,\angdfour$ & Section \ref{SS:fundamentalforms}\\
			\hline
			$\angnabla,u,\ubar,\theta,\invfolia,\uinvfolia,\Lgeo,\uLgeo$ & Section \ref{SS:SUBSETSOFSPACETIME} \\
			\hline
			$\deform{},\upchi,\uchi,\upeta,\ueta,\omega,\uomega,\upzeta,\ACC$ & Definition \ref{D:Defconnectioncoeff}\\
			\hline
			$\uptau,\rgeo,\lapse,\gsphere$ & Section \ref{SS:SUBSETSOFSPACETIME} \\
			\hline
			$\Lunit,\uLunit,e_I,e_J,\Lrenormal,\uLrenormal,\Theta,\Timelike,\spherenormal,\Trenomal,K,\Lambda$ & Section \ref{SS:SUBSETSOFSPACETIME} \\
			\hline
			$\mathfrak{R},\mathfrak{V}^i,\mathfrak{S},\mathfrak{W}$ & Section \ref{SS:NOTAITONFOREQUATIONS} \\
			\hline
			$\inihypersur,\inihypersurbar,\inigfour,\inivvariables,\iniu,\iniubar,\inisphere,\embedding,\uembedding$ & Section \ref{SS:initialembedding2} \\
			\hline
			$\mathring{\mathcal{S}},\mathring{\mathcal{T}},\inimaterialderivative,\iniLunit,\iniuLunit,\uiniuLunit,\uiniLunit,\inimaterialderivative$ & Section \ref{SS:GEOMETRYONDATAHYPERSURFACES} \\
			\hline
			$\inip,\iniLbargeo,\iniLgeo,\embeddomain_{\embedkappa,\ubar_*,u_*},\embed$ & Section \ref{SS:ADAPTEDCOORDINATES} \\
			\hline
		\end{tabular}
		\begin{tabular}{|c|c|} 
			\hline
			Symbol & Reference  \\ 
			\hline\hline	
			$E_I,\mathring{\nullapse},\mathring{e}_I,K,\inigsphere,\inilapse$ & Section \ref{SS:ADAPTEDCOORDINATES} \\
			\hline	
			$\widetilde{\FW}$ & Section \ref{E:WAVEEQUATIONFORDENSITYINADPATEDCOORDINATESALONGINITIALNULL} \\
			\hline
			$\inilogDensity_{(K)},\inient_{(K)},\inivel_{(K)}$ & Section \ref{SS:TRANSVERALDERIVATIVESONINITIALNULLHYPERSURFACES} \\
			\hline
			$\inimodivort,\mathring{\omega},\mathring{\upchi},\mathring{\uchi},\mathring{\ueta},\mathring{\upzeta}$ & Section \ref{SSS:inimodivort} \\
			\hline
			$\region$ & Section \ref{SS:Defregion} \\	
			\hline
			$\Stau,\region_1,\region_2,\region_3,\hypersur_{u}^{\ubar_1,\ubar_2},\hypersurbar_{\ubar}^{u_1,u_2}$ & Section \ref{SS:SUBSETSOFSPACETIME}\\
			\hline
			$\updelta_u,\updelta_{\ubar},\suubar,\region_{u_1,u_2}^{\ubar_1,\ubar_2},\region_0$ & Section \ref{SS:SUBSETSOFSPACETIME}\\
			\hline
			$\Lie$ & Lemma \ref{L:LieandCovariant}\\
			\hline
			$\diff\flatspherevol,\diff\gvol,\diff\tvol,\diff\tauvol,\diff\spherevol,\volume$ & Definition \ref{D:volumeforms}\\
			\hline
			$\embedgeo$ & Lemma \ref{LE:Jacobian}\\
			\hline
			$\iniepsilon,\BAinfty,\mathring{a},\tilde{a}$ & Section \ref{S:SectionAssumptions}\\
			\hline
			$\RBAinfty$ & Definition \ref{DE:LOT}\\
			\hline
			$\bEnergy^{(1)}_{u,\ubar},\Energy^{(1)}_u,\uEnergy^{(1)}_{\ubar}$ & Definition \ref{DE:waveenergy1}\\
			\hline
			$\bEnergy^{(2)}_{u,\ubar},\Energy^{(2)}_u,\uEnergy^{(2)}_{\ubar}$ & Definition \ref{d:Energy34}\\
			\hline
			$\KEnergy^{(1)},\KEnergy^{(2)},\KEnergy^{(3)}$ & Definition \ref{DE:Energyinspacetime}\\
			\hline
			$\bEnergy,\Energy,\uEnergy$ & Definition \ref{D:ControllingEnergyQuantities}\\
			\hline
			$\tilde{\mathbf{J}}$ & Definition \ref{D:elliptichyperboliccurrent}\\
			\hline
			$\mathfrak{A}^{(V)},\mathfrak{B}_{(V)},\mathfrak{C}^{(V)},\mathfrak{D}^{(V)}$ & Lemma \ref{L:integralidentitiesforS}\\
			\hline
			$\mathfrak{H},\mathfrak{H}_{(1)},\underline{\mathfrak{H}},\underline{\mathfrak{H}}_{(1)}$& Lemma \ref{L:integralidentitiesforS}\\
			\hline
			$\iniLtwo$& Section \ref{S:BASELEVELENERGYESTIMATES}\\
			\hline
			$d(p,q),\updelta,\ufinite,\ubarfinite$& Section \ref{SS:QUANTITATIVEVERSIONOFTHEINVERSEFUNCTION} \\
			\hline
			$\ExistenceRegionthreearg{\updelta}{\ufinite}{\ubarfinite},\widetilde{\mathcal{M}}_1,\widetilde{\mathcal{M}}_2,\widetilde{\mathcal{M}}_3$& Section 
			\ref{SS:MAINTHEOREMLOCALEXISTNCEANDUNIQUENESS}\\
			\hline
		\end{tabular}
	\end{center}
	\bibliography{bib.bib}
	\bibliographystyle{plain} 
\end{document}